\documentclass[12pt, a4paper]{amsart}
\usepackage{amsmath}
\usepackage{amscd}
\usepackage{amssymb}
\usepackage{color}

\newtheorem{thm}{Theorem}
\newtheorem{cor}[thm]{Corollary}
\newtheorem{lem}[thm]{Lemma}
\newtheorem{pro}[thm]{Proposition}

\theoremstyle{definition}
\newtheorem{examp}[thm]{Example}
\newtheorem{defn}[thm]{Definition}
\newtheorem{rem}[thm]{Remark}
\newtheorem{ass}[thm]{Assumption}
\newtheorem{conv}[thm]{Convention}
\numberwithin{equation}{section}
\numberwithin{thm}{section}
\begin{document}
\title{Bloch's principle for holomorphic maps into subvarieties of semi-abelian varieties}
\author{Katsutoshi Yamanoi}
\address{Department of Mathematics, Graduate School of Science,
Osaka University, Toyonaka, Osaka 560-0043, Japan}
\email{yamanoi@math.sci.osaka-u.ac.jp}

\maketitle

\newcommand{\cal}[0]{\mathcal}
\newcommand{\supp}[0]{\operatorname{supp}}
\newcommand{\m}[2]{m\left( #1,#2\right)}
\newcommand{\N}[2]{N\left( #1,#2\right)}
\newcommand{\tn}[3]{\overline{n}\left( #1,#2,#3\right)}
\newcommand{\h}[2]{h_{#1}\left( #2\right)}
\newcommand{\Sm}[2]{S\left( r,#1\right)}
\newcommand{\laoo}[2]{O\left(\log ^{+}\T{#1}{#2}\right)+O(\log ^{+}r)+O(1) \quad \Vert}
\newcommand{\mul}[3]{\operatorname{mult}_{#1}{#2}\cdot {#3}}
\newcommand{\St}[1]{\operatorname{St}(#1)}
\newcommand{\Ind}[3]{\operatorname{Ind}_{#1,#2}(#3)}
\newcommand{\Der}[2]{\operatorname{Der}_{#1}(#2)}
\newcommand{\cod}[2]{\operatorname{codim}(#1,#2)}
\newcommand{\Lie}[1]{\operatorname{Lie}(#1)}
\newcommand{\Spec}[0]{\operatorname{Spec}}
\newcommand{\I}[2]{\mathcal I(#1,#2)}
\newcommand{\Pe}[2]{\mathcal P(#1,#2)}
\def\ord{\mathop{\mathrm{ord}}\nolimits}
\def\ve{\mathop{\mathrm{vert}}\nolimits}
\def\edge{\mathop{\mathrm{edge}}\nolimits}
\def\Pic{\mathop{\mathrm{Pic}}\nolimits}
\def\vol{\mathop{\mathrm{vol}}\nolimits}
\def\supp{\mathop{\mathrm{supp}}\nolimits}
\def\ram{\mathop{\mathrm{ram}}\nolimits}

\begin{abstract}
We generalize a fundamental theorem in higher dimensional value distribution theory about entire curves in subvarieties $X$ of semi-abelian varieties to the situation of the sequences of holomorphic maps from the unit disc into $X$.
This generalization implies, among other things, that subvarieties of log general type in semi-abelian varieties are pseudo-Kobayashi hyperbolic.
As another application, we improve a classical theorem due to Cartan in 1920's about the system of nowhere vanishing holomorphic functions on the unit disc satisfying Borel's identity. 
\end{abstract}

\setcounter{tocdepth}{1}
\tableofcontents

\section{Introduction}

Bloch's principle is a widely recognized guiding principle in the study of complex function theory.
The origin of this principle goes back to a famous Bloch's dictum ``{\it Nihil est in infinito quod non prius fuerit in finito}''\footnote{According to Schiff \cite[p. 101]{Sch}, this may be translated as: {\it Nothing exists in the infinite plane that has not been previously done in the finite disc}.} made in his several papers written in 1926 (eg. \cite[p. 311]{Blo1}).
In many contexts, this statement is interpreted more concretely as the following heuristic principle that a family of holomorphic functions in a domain, all of which have a property $P$, is likely to be normal if $P$ cannot be possessed by non-constant entire functions in the plane (eg. \cite[p. 101]{Sch}).
A typical example is the correspondence between Picard's little theorem and Montel's theorem that a family of holomorphic functions on a domain, all of which omit two values $0$ and $1$, is normal. 
There are several very good references for Bloch's principle including \cite{Ber}, \cite[Chapter 4]{Sch}, \cite{Zal}.

In the light of his principle, Bloch \cite{Blo1} investigated a theory over the finite disc that corresponds to Borel's generalization of Picard's little theorem (cf. Theorem \ref{thm:borelthm}).
This investigation was succeeded by Cartan \cite{Car} who in particular  generalized Montel's theorem above (cf. Theorem \ref{thm:cartanthm}).
After a half century, Kiernan and Kobayashi \cite{KK} interpreted the works of Bloch and Cartan in the context of Kobayashi hyperbolic geometry.
These developments are fully explained by Lang \cite[Ch. VIII]{Lan}.
We shall discuss an improvement of the theorem of Cartan later (cf. Theorem \ref{cor:bc}).

In this paper, we are interested in holomorphic maps into subvarieties of semi-abelian varieties.
In the situation over the complex plane $\mathbb C$, we have the following theorem due to Bloch \cite{Blo2}, Ochiai \cite{O}, Kawamata \cite{Ka} and Noguchi \cite{N81}.

\begin{thm}[Bloch, Ochiai, Kawamata, Noguchi]\label{thm:lbo}
Let $A$ be a semi-abelian variety and let $X\subsetneqq A$ be a proper closed algebraic subvariety of $A$. 
Let  $f:\mathbb C\to X$ be a holomorphic map.
Then there exists a proper semi-abelian subvariety $B\subsetneqq A$ with the following two properties:
\begin{enumerate}
\item
$\varpi\circ f:\mathbb C\to A/B$ is a constant map, where $\varpi:A\to A/B$ is the quotient map.
\item $f(\mathbb C)\subset \bigcap_{b\in B}(X+b)$.
\end{enumerate}
\end{thm}

The statement of this theorem is possibly unfamiliar, but convenient to discuss Bloch's principle.
An equivalent statement is that the Zariski closure of the image $f(\mathbb C)$ is a translate of a semi-abelian subvariety $B'\subset A$ (cf. \cite[Thm 3.9.19]{Kbook}).
We may take $B\subset A$ in Theorem \ref{thm:lbo} to be this $B'$. 
As noted by Noguchi (cf. \cite[p. 156]{NW}), this theorem includes Borel's generalization of Picard's little theorem (cf. Theorem \ref{thm:borelthm}); We apply the theorem to the case that $A$ is an algebraic torus $(\mathbb G_m)^n$ and $X\subset (\mathbb G_m)^n$ is a subvariety defined by the linear equation $x_1+\cdots+x_n+1=0$ using coordinates $x_1,\ldots,x_n$ of $(\mathbb G_m)^n\subset \mathbb C^n$.
For more discussion about Theorem \ref{thm:lbo}, we refer the readers to \cite[Sec. 3.9]{Kbook} and \cite[Sec. 4.8]{NW}.

Now we are going to discuss a corresponding generalization of Theorem \ref{thm:lbo} for holomorphic mappings from the unit disc $\mathbb D$.
To state our main theorem, we first introduce one terminology from \cite[p. 242]{Lan}.
Let $\gamma >0$ and let $W\subset \mathbb C$ be an open set.
An assertion concerning points $w\in W$ will be said to hold for $\gamma$-almost all $w\in W$ if it holds for all $w\in W$ possibly except for $w$ contained in at most countably many closed discs such that the sum of the radii is less than $\gamma$.

Let $V\subset \overline{A}$ be a Zariski closed set, where $\overline{A}$ is an equivariant compactification of a semi-abelian variety $A$.
See Appendix \ref{sec:a} for the necessary matters on semi-abelian varieties.
Let $B\subset A$ be a semi-abelian variety.
We set
$$
\mathrm{Sp}_BV=\bigcap_{b\in B}(V+b)\subset V.
$$
Then $\mathrm{Sp}_BV\subset \overline{A}$ is a Zariski closed subset.

The following is the main result of this paper.
In the following statement, we denote by $\mathrm{Hol}(\mathbb D,X)$ the set of all holomorphic mappings from $\mathbb D$ to $X$.
For $0<s<1$, we set $\mathbb D(s)=\{z\in\mathbb D; |z|<s\}$.

\begin{thm}\label{thm:mstrong}
Let $A$ be a semi-abelian variety.
Let $X\subsetneqq A$ be a proper closed algebraic subvariety.
Let $(f_n)_{n\in \mathbb N}$  be a sequence of holomorphic maps in $\mathrm{Hol}(\mathbb D,X)$.
Then there exist a proper semi-abelian subvariety $B\subsetneqq A$ and a subsequence $(f_{n_k})_{k\in\mathbb N}$ with the following two properties:
\begin{enumerate}
\item $(\varpi\circ f_{n_k})_{k\in \mathbb N}$ converges uniformly on compact subsets of $\mathbb D$ to a holomorphic map $g:\mathbb D\to A/B$, where $\varpi:A\to A/B$ is the quotient map.
\item Let $\overline{A}$ be an equivariant compactification and let $\overline{X}\subset \overline{A}$ be the Zariski closure of $X$ in $\overline{A}$.
Then for every $0<s<1$, $\gamma >0$, and open neighbourhood $U\subset \overline{A}$ of $\mathrm{Sp}_B\overline{X}$, there exists $k_0\in \mathbb N$ such that, for all $k\geq k_0$, we have $f_{n_k}(z)\in U$ for $\gamma$-almost all $z\in \mathbb D(s)$.
\end{enumerate}
\end{thm}

Before going to discuss the applications of the theorem, we derive Theorem \ref{thm:lbo} from Theorem \ref{thm:mstrong}. 
Given $f:\mathbb C\to X$, we define a sequence $(\varphi_n)_{n\in \mathbb N}$ in $\mathrm{Hol}(\mathbb D,X)$ by $\varphi_n(z)=f(nz)$.
Then by Theorem \ref{thm:mstrong}, there exist a subsequence $\{\varphi_{n_k}\}_{k=1}^{\infty}$ and a proper semi-abelian subvariety $B\subsetneqq A$ such that $\{\varpi\circ \varphi_{n_k}\}$ converges uniformly on compact subsets of $\mathbb D$ to $g:\mathbb D\to A/B$.
We claim that $\varpi\circ f:\mathbb C\to A/B$ is constant.
So assume contrary that $\varpi\circ f$ is non-constant.
Then there exists a global holomorphic one-form $\eta\in \Gamma(A/B,\Omega_{A/B}^1)$ such that $(\varpi\circ f)^*\eta=\zeta(z)dz$ is non-zero on $\mathbb C$.
Hence there exists $l\geq 0$ such that $\zeta^{(l)}(0)\not=0$.
We set $(\varpi\circ \varphi_n)^*\eta=\xi_n(z)dz$.
Then $\xi_n(z)=n\zeta(nz)$.
Hence $|\xi_{n}^{(l)}(0)|=|n^{l+1}\zeta^{(l)}(0)|\to \infty$ as $n\to\infty$.
On the other hand, $\xi_{n_k}$ converges uniformly on compact subsets of $\mathbb D$ to $\xi(z)$ on $\mathbb D$, where $g^*\eta=\xi(z)dz$.
Hence $\xi_{n_k}^{(l)}(0)\to \xi^{(l)}(0)$ as $k\to\infty$.
This is a contradiction.
Hence $\varpi\circ f$ is constant.
To ensure the assertion (2) of Theorem \ref{thm:lbo}, we take a minimum semi-abelian variety $B\subsetneqq A$ such that $\varpi\circ f$ is constant.
By considering the translations by $f(0)$, we may assume without loss of generality that $f(0)=0_A$, i.e., the identity element of $A$.
Then we have $f(\mathbb C)\subset B$.
Let $X'\subset B$ be the Zariski closure of $f(\mathbb C)$.
To show $X'=B$, we assume contrary that $X'\subsetneqq B$.
Then by the argument above applied for $f:\mathbb C\to X'\subsetneqq B$, we get a proper semi-abelian subvariety $B'\subsetneqq B$ such that $\varpi'\circ f:\mathbb C\to B/B'$ is constant, where $\varpi':B\to B/B'$ is the quotient map.
This contradicts to the choice of $B$.
Thus $X'=B$.
By $X'\subset X\cap B$, we have $X\cap B=B$.
Hence $f(\mathbb C)\subset X\cap B\subset\bigcap_{b\in B}(b+X)$ as desired.
This completes the implication of Theorem \ref{thm:lbo} from Theorem \ref{thm:mstrong}. 

Note that in the implication above, the assertion (2) of Theorem \ref{thm:mstrong} plays no role.
However, in the applications of Theorem \ref{thm:mstrong} below, we need the assertion (2).

Next we discuss applications of Theorem \ref{thm:mstrong} to Kobayashi hyperbolic geometry.
We first introduce some terminologies from \cite[p. 245]{Kbook}.
Let $W$ be a relatively compact open domain of $M$, and let $\Delta$ be a closed subset of $M$.
We say that $W$ is tautly imbedded modulo $\Delta$ in $M$ if for each sequence $(f_n)_{n\in\mathbb N}$ in $\mathrm{Hol}(\mathbb D,W)$, one of the following holds:
\begin{enumerate}
\item 
$(f_n)_{n\in\mathbb N}$ has a subsequence $(f_{n_k})_{k\in\mathbb N}$ which converges uniformly on compact subsets of $\mathbb D$ to some $f\in \mathrm{Hol}(\mathbb D,M)$;
\item
for each compact set $K\subset \mathbb D$ and each compact set $L\subset M- \Delta$ there exists an integer $n_0$ such that $f_n(K)\cap L=\emptyset$ for all $n\geq n_0$.
\end{enumerate}

Let $X\subset A$ be a closed algebraic subvariety of a semi-abelian variety $A$.
Let $\overline{X}\subset \overline{A}$ be the compactification, where $\overline{A}$ is an equivariant compactification.
We set
$$
Z=\{ x\in \overline{X}\ ;\  \exists B\subset A, \text{a semi-abelian variety s.t. $\mathrm{dim}(x+B)\geq 1$ and $x+B\subset \overline{X}$}
\}.
$$
Let $S\subset \overline{X}$ be the Zariski closure of $Z$.
Then $S\subset \overline{X}$ is a Zariski closed set.
By \cite[Lemma 4.1]{N81}, $S$ is a proper subset of $\overline{X}$, provided $X$ is of log-general type.
As an application of Theorem \ref{thm:mstrong}, we obtain the following theorem.

\begin{thm}\label{thm:pkob}
Let $X\subset A$ be a closed algebraic subvariety.
Let $\overline{A}$ be a smooth equivariant compactification and let $\overline{X}\subset \overline{A}$ be the compactification.
Then $X$ is tautly imbedded modulo $S$ in $\overline{X}$.
\end{thm}

A theorem of Kiernan and Kobayashi \cite{KK} claims that if $W$ is tautly imbedded modulo $\Delta$ in $M$, then $W$ is hyperbolically imbedded modulo $\Delta$ in $M$ (cf. \cite[Thm 5.1.13]{Kbook}, \cite[Thm 1.4]{Lan}).
Hence by Theorem \ref{thm:pkob}, we immediately get the following corollary.

\begin{cor}
Let $X\subset A$ be a closed algebraic subvariety.
Let $\overline{A}$ be a smooth equivariant compactification and let $\overline{X}\subset \overline{A}$ be the compactification.
Then $X$ is hyperbolically imbedded modulo $S$ in $\overline{X}$.
In particular, if $X$ is of log-general type, then $X$ is pseudo-Kobayashi hyperbolic.
\end{cor}

When $A$ is compact, these results are previously proved in \cite{Y}.
In the compact case, Theorem \ref{thm:mstrong} yields a result on infinitesimal Kobayashi-Royden pseudo-metric $F_X$ defined as follows.
Let $X$ be an algebraic variety.
For each $x\in X$, we call the set $\breve{T}_xX$ of all 1-jets the tangent cone of $X$ at $x$, and $\breve{T}X=\cup_{x\in X}\breve{T}_xX$ the tangent cone of $X$ (cf. \cite[p. 31]{Kbook}).
If $X$ is smooth, then $\breve{T}X$ coincides with the usual tangent bundle.
For $v\in \breve{T}X$, we set
$$
F_X(v)=\inf \left\{ r>0 ; \exists f:\mathbb D\to X\ \text{s.t. $f'(0)=\frac{1}{r}v$}\right\}.
$$
Suppose there exists a holomorphic map $f:\mathbb C\to X$ such that $f'(0)=v\in \breve{T}X$, then we have $F_X(v)=0$.
However the converse is not true in general (cf. Example \ref{examp:20220518}).
Theorem \ref{thm:mstrong} yields the following theorem.
We do not know whether this statement is true or not for subvarieties of non-compact semi-abelian varieties. 

\begin{thm}\label{thm:20220429}
Let $A$ be an abelian variety.
Let $X\subset A$ be a closed subvariety.
Suppose $v\in \breve{T}X$ satisfies $F_X(v)=0$.
Then there exists a holomorphic map $f:\mathbb C\to X$ such that $f'(0)=v$.
\end{thm}

Now we return to the classical topics concerning Bloch-Cartan theorem.
To simplify the complicated indices in the description, we employ the following conventions.

\begin{conv}\label{conv:20220519}
For a set $E$ and an infinite set $I$, an indexed family $(x_i)_{i\in I}$ is a function $I\to E$.
When $E$ is an infinite set, we consider $E$ as an indexed family indexed by $E$ itself by the identity map $E\to E$.
Let $S$ be a metric space with the distance function $d$.
Let $(f_i)_{i\in I}$ be an indexed family of functions defined in $\mathbb D$, and with values in $S$.
We say that $(f_i)_{i\in I}$ converges uniformly on compact subsets on $\mathbb D$ to some $f:\mathbb D\to S$ if for every compact subset $K\subset \mathbb D$ and every $\varepsilon>0$, there exists a finite subset $F\subset I$ such that $d(f_i(z),f(z))<\varepsilon$ for all $z\in K$ and $i\in I-F$.
Of course, when $I=\mathbb N$, this definition coincides with the usual definition of the uniform convergence on compact subsets on $\mathbb D$.
If each $f_i$ is continuous, then $f$ is continuous.
Indeed, we take a sequence $i_1,i_2,i_3,\ldots$ of distinct elements in $I$.
Then the sequence $(f_{i_k})_{k\in\mathbb N}$ converges uniformly on compact subsets on $\mathbb D$ to $f$.
Hence $f$ is continuous.
Similarly, if $S$ is a complex manifold and each $f_i$ is holomorphic, then $f$ is holomorphic.
\end{conv}

\begin{thm}\label{cor:bc}
Let $\mathcal{F}$ be an infinite set of $p$-tuples $f=(f_1,\ldots,f_p)$ of nowhere vanishing holomorphic functions on $\mathbb D$ satisfying the following identity
\begin{equation}\label{eqn:20220429}
f_1+f_2+\cdots +f_p=0.
\end{equation}
Then there exist disjoint non-empty subsets $I_1,\ldots,I_n\subset \{ 1,\ldots,p\}$ and an infinite subset $\mathcal{G}\subset \mathcal{F}$ with the following properties:
\begin{enumerate}
\item
$n\geq 1$.
\item 
Let $k\in \{ 1,\ldots,n\}$.
Then:
\begin{enumerate}
\item
For all $i,j\in I_k$, the indexed family $( f_i/f_j)_{f\in\mathcal{G}}$ converges uniformly on compact subsets of $\mathbb D$ to a nowhere vanishing holomorphic function.
\item 
For all $j\in I_k$, the indexed family $(\sum_{i\in I_k}f_i/f_j)_{f\in\mathcal{G}}$ converges uniformly on compact subsets of $\mathbb D$ to $0$.
In particular, $I_k$ contains at least two elements.
\end{enumerate}
\item
Set $I=I_1\sqcup\cdots\sqcup I_n$.
Let $i\in \{1,\ldots,p\}- I$.
Let $\varepsilon>0$, $s\in (0,1)$ and $\gamma>0$.
Then there exists a finite subset $\mathcal{E}\subset \mathcal{G}$ such that
for all $f\in \mathcal{G}-\mathcal{E}$, we have
$$\frac{|f_i(z)|}{\sqrt{\sum_{j\in I}|f_j(z)|^2}}<\varepsilon$$ 
for $\gamma$-almost all $z\in\mathbb D(s)$.
\end{enumerate}
\end{thm}

Some historical remarks are required.
The identity \eqref{eqn:20220429} was considered by Borel \cite{Borel} who generalized Picard's little theorem as follows.

\begin{thm}[Borel]\label{thm:borelthm}
Let $f_1,\ldots,f_p$ be nowhere vanishing holomorphic functions on $\mathbb C$ satisfying the identity \eqref{eqn:20220429}.
Then there exists a partition of indices $\{1,\ldots,p\}=I_1\sqcup\cdots\sqcup I_n$ such that each $I_k$ satisfies that
\begin{enumerate}
\item
for every $i,j\in I_k$, the quotient $f_i/f_j$ is constant, and 
\item
$\sum_{i\in I_k}f_i=0$.
\end{enumerate}
\end{thm}

By the second conclusion, each $I_k$ has at least two elements.
In particular, when $p=3$, we have $n=1$ and $I_1=\{ 1,2,3\}$.
This implies Picard's little theorem as follows.
Let $g:\mathbb C\to\mathbb C-\{0,1\}$.
We set $f_1(z)=g(z)$, $f_2(z)=1-g(z)$ and $f_3(z)=-1$.
Then $f_1$, $f_2$ and $f_3$ are nowhere vanishing holomorphic functions on $\mathbb C$ satisfying the identity $f_1+f_2+f_3=0$.
Hence by Borel's theorem, $f_1/f_3=-g$ is constant, as desired.

As we already mentioned above, the corresponding theory over the disc $\mathbb D$ was investigated by Bloch \cite{Blo1} and Cartan \cite[p. 312]{Car}.
See also \cite[Ch. VIII]{Lan}.
Let $\mathcal{F}$ be an infinite set of $p$-tuples $f=(f_1,\ldots,f_p)$ of nowhere vanishing holomorphic functions on $\mathbb D$ satisfying the identity \eqref{eqn:20220429}.
A subset $I\subset \{1,\ldots,p\}$ is called C-class if there exists $i\in I$ such that for all $j\in I$, the sequence $(f_j/f_i)_{f\in\mathcal{F}}$ is uniformly bounded on every compact set of $\mathbb D$ and $(\sum_{j\in I}f_j/f_i)_{f\in\mathcal{F}}$ converges uniformly on compact subsets of $\mathbb D$ to $0$.

\begin{thm}[Cartan]\label{thm:cartanthm}
There exists an infinite subset $\mathcal{G}\subset \mathcal{F}$ such that $\{1,\ldots,p\}$ itself is C-class, or there exist two disjoint subsets $I_1,I_2\subset \{1,\ldots,p\}$ such that both are C-classes.
\end{thm}

Cartan conjectured that there exists an infinite subset $\mathcal{G}\subset \mathcal{F}$ such that the set $\{1,\ldots,p\}$ can be partitioned into C-classes.
However this conjecture was disproved by Eremenko \cite{Ere1}.
In \cite{Ere2}, Eremenko proposed a modified version of Cartan's conjecture and proved it for the case $p=5$.

\medskip

{\it Theorem \ref{cor:bc} implies Theorem \ref{thm:cartanthm} as follows.}
We apply Theorem \ref{cor:bc} to get an infinite subset $\mathcal{G}\subset \mathcal{F}$ and disjoint subsets $I_1,\ldots,I_n\subset \{ 1,\ldots,p\}$.
If $n\geq 2$, there is nothing to do for each $I_k$ is a C-class by the second assertion of Theorem \ref{cor:bc}.
Thus we consider the case $n=1$.
Let $j\in I_1$.
The assertion (2a) of Theorem \ref{cor:bc} implies that $(\sqrt{\sum_{i\in I_1}|f_i|^2}/|f_j|)_{f\in \mathcal{G}}$ converges uniformly on compact subsets of $\mathbb D$.
Hence the third assertion of Theorem \ref{cor:bc} reads as follows:
\begin{quote}
Let $i\not\in I_1$, $\varepsilon>0$, $s\in (0,1)$ and $\gamma>0$.
Then there exists a finite subset $\mathcal{E}\subset \mathcal{G}$ such that for all $f\in \mathcal{G}-\mathcal{E}$, we have $|f_i|/|f_j|<\varepsilon$ for $\gamma$-almost all $z\in\mathbb D(s)$.
\end{quote}
Let $i\not\in I_1$.
Let $K\subset \mathbb D$ be a compact set and let $\varepsilon>0$.
We take $s\in (0,1)$ and $\gamma>0$ such that $K\subset \mathbb D(s-2\gamma)$.
Then by the third assertion of Theorem \ref{cor:bc}, there exists a finite set $\mathcal{E}\subset \mathcal{G}$ such that for all $f\in \mathcal{G}-\mathcal{E}$, we have $|f_i|/|f_j|<\varepsilon$ for $\gamma$-almost all $z\in\mathbb D(s)$.
Let $f\in \mathcal{G}-\mathcal{E}$.
We may take $s'\in (s-2\gamma,s)$ such that $|f_i|/|f_j|<\varepsilon$ holds over the circle $\partial\mathbb D(s')$.
Since $f_i/f_j$ is holomorphic, the maximal principle yields that $|f_i|/|f_j|<\varepsilon$ for all $z\in \mathbb D(s')$, hence for all $z\in K$.
Hence $(f_i/f_j)_{f\in\mathcal{G}}$ converges uniformly to $0$ on $K$.
Hence if $n=1$, then $\{1,\ldots,p\}$ is C-class.
This completes the derivation of Theorem \ref{thm:cartanthm} from Theorem \ref{cor:bc}.
\hspace{\fill} $\square$

\medskip

The contents of this paper is as follows:
The sections \ref{sec:dem}-\ref{sec:10} are devoted for the proof of Theorem \ref{thm:mstrong}.
Although the proof of Theorem \ref{thm:mstrong} is lengthy, the structure is rather simple.
The proof is divided into two parts.
In the first part, we shall establish a new normality criterion for families $\mathcal{F}\subset \mathrm{Hol}(\mathbb D,A)$.
This is Proposition \ref{pro:320}.
We recall notions of Demailly jet space in section \ref{sec:dem}, which is used in the statement of this normality criterion.
The main technical tool for the proof of Proposition \ref{pro:320} is Nevanlinna theory, which is the theme of section \ref{sec:4}.
After preparations, we prove Proposition \ref{pro:320} in section \ref{sec:8}.
This is the first part.
In the second part, we shall find $B\subset A$ such that $\{ \varpi\circ f\}_{f\in\mathcal{F}}$ satisfies the normality criterion.
This is stated in Proposition \ref{cor:20210306}.
The proof of this proposition is the main theme of sections \ref{sec:8.5}-\ref{sec:9}.
Then we prove Theorem \ref{thm:mstrong} in section \ref{sec:10}.

In each of the sections \ref{sec:11}, \ref{sec:20220429} and \ref{sec:12}, we prove Theorems \ref{thm:pkob}, \ref{thm:20220429} and \ref{cor:bc} in this order.
Some needed facts in this paper for semi-abelian varieties are treated in appendix \ref{sec:a}.
Appendix \ref{sec:b} is devoted for the proof of algebraic geometrical proposition needed in this paper, namely a flattening result using blow-ups.

\begin{conv}\label{rem:var}
In this paper, an algebraic variety (or simply a variety) is an integral, separated scheme of finite type over the complex number field $\mathbb C$ (cf., e.g., \cite[p. 105]{H}).
In particular, every variety is reduced, irreducible and non-empty (cf. \cite[Chap. II, Prop. 3.1]{H}).
Every variety has a canonically associated complex space structure (cf. \cite[p. 439]{H}).
\end{conv}

\medskip

\section{Demailly jet spaces}\label{sec:dem}
We introduce Demailly jet spaces (cf. \cite{Dem}).
Let $M$ be a positive dimensional smooth algebraic variety.
Let $V\subset TM$ be an algebraic vector subbundle, whose bundle rank is positive.
Set $\tilde{M}= P(V)$.
Let $\pi :\tilde{M}\to M$ be the projection.
We define a vector subbundle $\tilde{V}\subset T\tilde{M}$ by the following: for every point $(x,[v])\in \tilde{M}$ associated with a vector $v\in V_x\backslash \{ 0\}$, we set
\begin{equation}\label{eqn:20211104}
\tilde{V}_{(x,[v])}=\{ \xi \in T_{(x,[v])}\tilde{M}\ ;\ \pi_*(\xi )\in \mathbb Cv\} ,
\end{equation}
where $\pi_*:T\tilde{M}\to TM$ is the induced map.
Let $f:\mathbb D\to M$ be a non-constant holomorphic map.
We say that $f$ is tangent to $V$ if $f'(z)\in V_{f(z)}$ for all $z\in \mathbb D$.
If $f$ is tangent to $V$, we may define $f_{[1]}:\mathbb D\to \tilde{M}$ by $f_{[1]}(z)=(f(z),[f'(z)])$.
Then $f_{[1]}$ is tangent to $\tilde{V}$.

\par

We inductively define the Demailly jet space $M_k$ together with vector subbundle $V_k\subset TM_k$ by 
$$
(M_0,V_0)=(M,TM), \qquad (M_k,V_k)=(\widetilde{M_{k-1}},\widetilde{V_{k-1}}).
$$
For a non-constant holomorphic map $f:\mathbb D\to M$, we define $f_{[k]}:\mathbb D\to M_k$ inductively by $f_{[0]}=f$ and $f_{[k]}=(f_{[k-1]})_{[1]}$.

For $k\geq 2$, we define the singular locus $M_k^{\mathrm{sing}}\subset M_k$ as follows.
We note $M_k\subset PTM_{k-1}$.
We have a natural map $M_{k-1}\to M$, from which we get the relative tangent bundle $T_{M_{k-1}/M}\subset TM_{k-1}$.
We set
\begin{equation}\label{eqn:202206093}
M_k^{\mathrm{sing}}=M_k\cap PT_{M_{k-1}/M}.
\end{equation}
Then $M_k^{\mathrm{sing}}\subset M_k$ is a Zariski closed set.
We claim that this is a divisor.
To show this, we consider the maps $M_{k-1}\to M_{k-2}\to M$, which induces $TM_{k-1}\to TM_{k-2}\to TM$.
For $v\in TM_{k-1}$, we have $v\in T_{M_{k-1}/M}$ if and only if $(\pi_{k-1})_*(v)\in T_{M_{k-2}/M}$, where $\pi_{k-1}:M_{k-1}\to M_{k-2}$ is the natural projection.
The rational map $PTM_{k-1}\dashrightarrow PTM_{k-2}$ induces the holomorphic map $p:PTM_{k-1}-PT_{M_{k-1}/M_{k-2}}\to PTM_{k-2}$.
Then
\begin{equation}\label{eqn:20211119}
PT_{M_{k-1}/M}-PT_{M_{k-1}/M_{k-2}}= p^{-1}(PT_{M_{k-2}/M}).
\end{equation}
Note that the subbundle $T_{M_{k-1}/M_{k-2}}\subset TM_{k-1}$ satisfies 
\begin{equation}\label{eqn:202203281}
T_{M_{k-1}/M_{k-2}}\subset V_{k-1}.
\end{equation}
The rank of  $T_{M_{k-1}/M_{k-2}}$ is equal to $\dim M-1$.
Set $D_k=PT_{M_{k-1}/M_{k-2}}\subset PV_{k-1}=M_k$.
Then $D_k$ is a divisor on $M_k$.
By \eqref{eqn:20211119}, we have 
\begin{equation}\label{eqn:202203282}
M_k^{\mathrm{sing}}=D_k\cup \pi_k^{-1}(M_{k-1}^{\mathrm{sing}}),
\end{equation}
where $\pi_k:M_k\to M_{k-1}$.
Hence $M_k^{\mathrm{sing}}$ is a divisor on $M_k$, using the induction on $k$.

Next we consider the case of semi-abelian varieties.
Let $A$ be a semi-abelian variety.
Let $m:A\times A\to A$ be the natural action such that $(x,a)\mapsto x+a$.
This induces
\begin{equation*}
m_*:T(A\times A) \to TA.
\end{equation*}
We have a subbundle
$$
A\times TA\subset T(A\times A).
$$
Thus we get
$$
(A\times TA)|_{A\times \{0_A\}}\to TA.
$$
By $T_{0_A}A=\mathrm{Lie}A$, we get
\begin{equation}\label{eqn:20211108}
\psi:A\times\mathrm{Lie}A\to TA.
\end{equation}
Then $\psi$ is an isomorphism of vector bundles over $A$.
For each $a\in A$, we denote by $t_a:A\to A$ the translation defined by $a$.
This induces an isomorphism 
$$(t_a)_*:TA\to TA.$$
Then we have
\begin{equation}\label{eqn:202112031}
(\psi^{-1}\circ(t_a)_*\circ\psi)(x,v)=(x+a,v).
\end{equation}
Let $f\in \mathrm{Hol}(\mathbb D,A)$.
We define $f_{\mathrm{Lie}A}:\mathbb D\to \mathrm{Lie}A$ by the composite of
$$
\mathbb D\overset{f'}{\to}TA\overset{\psi^{-1}}{\to}A\times \mathrm{Lie}A\to \mathrm{Lie}A.
$$
For $a\in A$, we define $f_a:\mathbb D\to A$ by $f_a(z)=f(z)+a$.
By \eqref{eqn:202112031}, we have 
\begin{equation}\label{eqn:20211213}
(f_a)_{\mathrm{Lie}A}=f_{\mathrm{Lie}A}.
\end{equation}

We consider the Demailly jet space for the case $M=A\times S$, where $S$ is a smooth algebraic variety.
We construct a smooth algebraic variety $S_{k,A}$ and a vector subbundle 
\begin{equation}\label{eqn:20230312}
V_k^{\dagger}\subset TS_{k,A}\times \mathrm{Lie}(A)
\end{equation}
as follows.
Set $S_{0,A}=S$ and $V_0^{\dagger}=TS\times \mathrm{Lie}(A)$.
Suppose $S_{k-1,A}$ and $V_{k-1}^{\dagger}\subset TS_{k-1,A}\times \mathrm{Lie}(A)$ are given.
We set 
\begin{equation}\label{eqn:202112033}
S_{k,A}=P(V_{k-1}^{\dagger}).
\end{equation}
Then $S_{k,A}$ is a smooth algebraic variety.
Let $\tau :S_{k,A}\to S_{k-1,A}$ be the projection.
We have a vector bundle map $(\tau _*,\mathrm{id}_{\mathrm{Lie}(A)}):TS_{k,A}\times \mathrm{Lie}(A)\to TS_{k-1,A}\times \mathrm{Lie}(A)$.
We define $V_k^{\dagger}\subset TS_{k,A}\times \mathrm{Lie}(A)$ as follows.
For each $(x,[v])\in S_{k,A}$, where $x\in S_{k-1,A}$ and $v\in V_{k-1}^{\dagger}\backslash \{ 0\}$, we set
\begin{equation}\label{eqn:20220621}
(V_k^{\dagger})_{(x,[v])}=\{ \xi \in T_{(x,[v])}S_{k,A}\times \mathrm{Lie}(A);\ (\tau _*,\mathrm{id}_{\mathrm{Lie}(A)})(\xi )\in \mathbb C\cdot v\} .
\end{equation}
By the isomorphism \eqref{eqn:20211108}, we have an isomorphism $T(A\times S_{k,A})\simeq A\times TS_{k,A}\times \mathrm{Lie}(A)$.
By this isomorphism, we consider $A\times V_k^{\dagger}$ as a vector subbundle of $T(A\times S_{k,A})$.

Next we construct an isomorphism
\begin{equation}\label{eqn:20211203}
\varphi_k:A\times S_{k,A}\to (A\times S)_k
\end{equation}
as follows.
For $k=0$, we set $\varphi_0=\mathrm{id}_{A\times S}$.
Note that $(\varphi_0)_*(A\times V_0^{\dagger})=V_0$.
Suppose we are given an isomorphism $\varphi_{k-1}:A\times S_{k-1,A}\to (A\times S)_{k-1}$ such that $(\varphi_{k-1})_*(A\times V_{k-1}^{\dagger})=V_{k-1}$.
Then the projectivization of $(\varphi_{k-1})_*$ induces an isomorphism $\varphi_k:A\times S_{k,A}\to (A\times S)_k$.
Under this isomorphism, we have $(\varphi_k)_*(A\times V_{k}^{\dagger})=V_{k}$.
Thus inductively, we have constructed the isomorphism \eqref{eqn:20211203}.

In the following, we identify $(A\times S)_k$ with $A\times S_{k,A}$ by the isomorphism \eqref{eqn:20211203}.
When $S$ is a single point, we denote
$$
P_{k,A}=\{ \mathrm{pt}\}_{k,A}.
$$
Then under the isomorphism \eqref{eqn:20211203}, we have $A_k=A\times P_{k,A}$.

Let $f\in \mathrm{Hol}(\mathbb D,A\times S)$ be non-constant.
We denote by $f_S:\mathbb D\to S$ the composite of $f:\mathbb D\to A\times S$ and the second projection $A\times S\to S$.  
We define 
\begin{equation}\label{eqn:20220206}
f_{S_{k,A}}:\mathbb D\to S_{k,A}
\end{equation}
as follows.
We set $f_{S_{0,A}}=f_S$.
Note that $((f_{S_{0,A}})',f_{\mathrm{Lie}A})(z)\in V_0^{\dagger}$.
Suppose that $f_{S_{k-1,A}}:\mathbb D\to S_{k-1,A}$ is given such that $((f_{S_{k-1,A}})',f_{\mathrm{Lie}A})(z)\in V_{k-1}^{\dagger}$.
We define $f_{S_{k,A}}:\mathbb D\to S_{k,A}$ by the projectivization of $((f_{S_{k-1,A}})',f_{\mathrm{Lie}A})$.
Then we have  $((f_{S_{k,A}})',f_{\mathrm{Lie}A})(z)\in V_{k}^{\dagger}$.
Thus we have constructed $f_{S_{k,A}}:\mathbb D\to S_{k,A}$ inductively for all $k$.
Let $f_A:\mathbb D\to A$ be the composite of $f:\mathbb D\to A\times S$ and the first projection $A\times S\to A$.
We have $\varphi_k\circ (f_A,f_{S_{k,A}})=f_{[k]}$, which follows from the construction.
For $a\in A$, we set $f_a:\mathbb D\to A\times S$ by $f_a(z)=(f_A(z)+a,f_S(z))$.
Then by \eqref{eqn:20211213}, we have
\begin{equation}\label{eqn:202112035}
(f_a)_{S_{k,A}}=f_{S_{k,A}}
\end{equation}
for all $a\in A$.

Let $k\geq 2$.
We define $S_{k,A}^{\mathrm{sing}}\subset S_{k,A}$ by
\begin{equation}\label{eqn:20220129}
S_{k,A}^{\mathrm{sing}}=S_{k,A}\cap P(T_{S_{k-1,A}/S}\times \{ 0\}),
\end{equation}
where $T_{S_{k-1,A}/S}\times \{ 0\}\subset TS_{k-1,A}\times \mathrm{Lie}A$.
Then we have 
\begin{equation}\label{eqn:20220611}
A\times S_{k,A}^{\mathrm{sing}}=(A\times S)_k^{\mathrm{sing}}
\end{equation}
under the isomorphism of \eqref{eqn:20211203}.

The following definition plays an important role in this paper.

\begin{defn}\label{defn:20201225}
Let $B\subset A$ be a semi-abelian subvariety.
For $k\geq 1$, we define $E_{k,A,A/B}\subset P_{k,A}$ by $E_{k,A,A/B}=P_{k,A}\cap P(TP_{k-1,A}\times \mathrm{Lie}B)$, where $TP_{k-1,A}\times \mathrm{Lie}B\subset TP_{k-1,A}\times \mathrm{Lie}A$.
\end{defn}

Then we have $E_{k,A,A}\subset E_{k,A,A/B}$.
Moreover by $T_{P_{k-1,A}/\{\mathrm{pt}\}}=TP_{k-1,A}$, we have
\begin{equation}\label{eqn:202206235}
P_{k,A}^{\mathrm{sing}}= E_{k,A,A}.
\end{equation}

\begin{lem}\label{lem:202206211}
Let $k\geq 1$.
Let $\tau:P_{k+1,A}\to P_{k,A}$ be the projection.
Then $\tau^{-1}(E_{k,A,A/B})\subset E_{k+1,A,A/B}$.
\end{lem}

{\it Proof.}\
By the definition \eqref{eqn:202112033}, we have $P_{k+1,A}=P(V_{k}^{\dagger})$, where $V_{k}^{\dagger}\subset TP_{k,A}\times \mathrm{Lie}(A)$.
Let $(y,[\xi])\in P_{k+1,A}\backslash E_{k+1,A,A/B}$, where $y\in P_{k,A}$ and $\xi\in V_{k}^{\dagger}\backslash\{0\}$.
Then $\xi\not\in TP_{k,A}\times\mathrm{Lie}(B)$.
Let $y=(x,[v])\in P(V_{k-1}^{\dagger})$, where $x\in P_{k-1,A}$ and $v\in V_{k-1}^{\dagger}\backslash\{0\}$.
Then by the definition of $V_k^{\dagger}$ (cf. \eqref{eqn:20220621}), the image of $\xi$ under the map $TP_{k,A}\times\mathrm{Lie}(A)\to TP_{k-1,A}\times\mathrm{Lie}(A)$ is contained in the linear space $\mathbb C\cdot v$.
Hence $v\not\in TP_{k-1,A}\times\mathrm{Lie}(B)$.
Hence $y\in P_{k,A}\backslash E_{k,A,A/B}$.
Hence we have proved $\tau^{-1}(E_{k,A,A/B})\subset E_{k+1,A,A/B}$.
\hspace{\fill} $\square$

\begin{rem}\label{rem:20230312}
Let $k\geq 0$.
We have the subbundle $V^{\dagger}_{k}\subset TP_{k,A}\times \mathrm{Lie}(A)$ so that $P_{k+1,A}=P(V^{\dagger}_k)$.
Set $S=P_{k,A}$.
For each $l\geq 0$, we denote by $V^{\dagger}_{l,S}\subset TS_{l,A}\times \mathrm{Lie}(A)$ the object in \eqref{eqn:20230312} so that $S_{l+1,A}=P(V^{\dagger}_{l,S})$.
Then for each $l\geq 0$, there exists a natural embedding 
\begin{equation}\label{eqn:20230313}
P_{k+l,A}\subset  S_{l,A}
\end{equation}
such that $V^{\dagger}_{k+l}\subset V^{\dagger}_{l,S}\cap (TP_{k+l,A}\times\mathrm{Lie}(A))$.
This is constructed inductively as follows.
For $l=0$, we note $P_{k,A}=S=S_{0,A}$ and $V^{\dagger}_{k}\subset TP_{k,A}\times \mathrm{Lie}(A)=V^{\dagger}_{0,S}$.
We discuss the induction step from $l$ to $l+1$.
By $P_{k+l,A}\subset S_{l,A}$ and $V^{\dagger}_{k+l}\subset V^{\dagger}_{l,S}\cap (TP_{k+l}\times\mathrm{Lie}(A))$, we have $P_{k+l+1,A}=P(V^{\dagger}_{k+l})\subset P(V^{\dagger}_{l,S})=S_{l+1,A}$.
The constructions of $V^{\dagger}_{k+l+1}$ and $V^{\dagger}_{l+1,S}$ (cf. \eqref{eqn:20220621}) yield $V^{\dagger}_{k+l+1}\subset V^{\dagger}_{l+1,S}\cap (TP_{k+l+1,A}\times\mathrm{Lie}(A))$.
This completes the induction step.
\end{rem}

\section{Sufficient condition for normality: Statement of Proposition \ref{pro:320}}\label{sec:2}

The goal of this section is to introduce Proposition \ref{pro:320}.
This proposition gives a sufficient condition for a subset of $\mathrm{Hol}(\mathbb D,A)$ to be normal, where $A$ is a semi-abelian variety.
The proof of this proposition is rather lengthy, so we devote sections \ref{sec:4}-\ref{sec:8} for the proof.
To state our proposition, we need to prepare several terminologies, which we describe below.

\subsection{Family of holomorphic maps and Zariski closed sets}
We start from the following two definitions.

\begin{defn}\label{defn:20211007}
Let $S$ be a variety and let $Z\subset S$ be a Zariski closed set.
\begin{enumerate}
\item
By a $Z$-{\it admissible modification} $\varphi:S'\to S$, we assume that 
\begin{enumerate}
\item
$\varphi$ is projective and birational, and 
\item
there exists a Zariski open set $U\subset S$ such that $Z\cap U\not=\emptyset$ and $\varphi^{-1}(U)\to U$ is an isomorphism.
\end{enumerate}
\item
For a $Z$-admissible modification $\varphi:S'\to S$, we define the {\it minimal transform} $Z'\subset S'$ as follows.
Let $\mathcal{U}$ be the set of all Zariski open subsets $U\subset S$ with the property (b) above.
We set $Z'=\cap_{U\in\mathcal{U}}\overline{\varphi^{-1}(Z\cap U)}$, where $\overline{\varphi^{-1}(Z\cap U)}\subset S'$ is the Zariski closure.
\end{enumerate}
\end{defn}

\begin{defn}\label{defn:20211006}
Let $\mathcal{F}=(f_i)_{i\in I}$ be an infinite indexed family in $\mathrm{Hol}(\mathbb D,S)$, where $S$ is a variety.
Let $Z\subset S$ be a Zariski closed set.
\begin{enumerate}
\item
We write $\mathcal{F}\to Z$ if the following holds:
For every $0<s<1$, $\gamma >0$, and open neighbourhood $U\subset S$ of $Z$, there exists a finite subset $E\subset I$ such that, for all $i\in I-E$, we have $f_i(z)\in U$ for $\gamma$-almost all $z\in \mathbb D(s)$.
\item
We write $\mathcal{F}\Rightarrow Z$ if the followings hold:
\begin{enumerate}
\item
Let $V\subset S$ be a Zariski closed set such that $Z\not\subset V$.
Then $f_i(\mathbb D)\not\subset V$ for all $i\in I$ with finite exception.
\item
Let $S'\to S$ be a $Z$-admissible modification and let $Z'\subset S'$ be the minimal transform.
Then $\mathcal{F}\to Z'$.
\end{enumerate}
\end{enumerate}
\end{defn}

\begin{rem}\label{rem:202205291}
We supplement the condition $\mathcal{F}\to Z'$ in the assertion (2b).  
There exists a Zariski open set $U\subset S$ such that $\varphi:S'\to S$ satisfies the assertion (1b) in Definition \ref{defn:20211007}.
Then by $Z\not\subset (S-U)$, there exists a finite subset $E\subset I$ such that $f_i(\mathbb D)\not\subset (S-U)$  for all $i\in I-E$.
Then for each $i\in I-E$, there is a natural lift $f_i':\mathbb D\to S'$ of $f_i$.
By $\mathcal{F}\to Z'$ in the assertion (2b), we  mean $(f_i')_{i\in I-E}\to Z'$.
\end{rem}

\begin{rem}\label{rem:20220523}
If $\mathcal{F}\to Z$, then $\mathcal{G}\to Z$ for all infinite indexed subfamily $\mathcal{G}$ of $\mathcal{F}$.
Here we call an infinite indexed family $\mathcal{G}=(f_j)_{j\in J}$ a subfamily of $\mathcal{F}=(f_i)_{i\in I}$ if $J\subset I$ is an infinite subset.
If $\mathcal{F}\Rightarrow Z$, then $\mathcal{G} \Rightarrow Z$ for all infinite indexed subfamily $\mathcal{G}$ of $\mathcal{F}$.
\end{rem}

We prove several basic properties related to Definition \ref{defn:20211006}.

\begin{lem}\label{lem:20220626}
Let $\mathcal{F}=(f_i)_{i\in I}$ be an infinite indexed family in $\mathrm{Hol}(\mathbb D,S)$.
Let $Z_1$ and $Z_2$ be Zariski closed subsets of $S$.
If $\mathcal{F}\to Z_1$ and $\mathcal{F}\to Z_2$, then $\mathcal{F}\to Z_1\cap Z_2$.
\end{lem}

{\it Proof.}\
Set $S_o=S\backslash (Z_1\cap Z_2)$.
Then $(Z_1\cap S_o)\cap (Z_2\cap S_o)=\emptyset$.
Hence there exist open neighbourhoods $U_1\subset S_o$ of $Z_1\cap S_o$ and $U_2\subset S_o$ of $Z_2\cap S_o$ such that $U_1\cap U_2=\emptyset$.

 We take $0<s<1$, $\gamma >0$, and open neighbourhood $U\subset S$ of $Z_1\cap Z_2$.
 Then $U_1\cup U\subset S$ is an open neighbourhood of $Z_1$.
 Hence by $\mathcal{F}\to Z_1$, there exists a finite set $E_1\subset I$ such that for all $i\in I\backslash E_1$, we have $f_i(z)\in U_1\cup U$ for $\gamma/2$-almost all $z\in\mathbb D(s)$.
 Similary, by $\mathcal{F}\to Z_2$, there exists a finite set $E_2\subset I$ such that for all $i\in I\backslash E_2$, we have $f_i(z)\in U_2\cup U$ for $\gamma/2$-almost all $z\in\mathbb D(s)$. 
 Set $E=E_1\cup E_2$.
 Then $E$ is finite.
 Note that $(U_1\cup U)\cap (U_2\cup U)=U$.
 Hence for all $i\in I\backslash E$, we have $f_i(z)\in U$ for $\gamma$-almost all $z\in\mathbb D(s)$.  
 Thus $\mathcal{F}\to Z_1\cap Z_2$.
 \hspace{\fill} $\square$
 
\begin{lem}\label{lem:202206261}
Let $\mathcal{F}=(f_i)_{i\in I}$ be an infinite indexed family in $\mathrm{Hol}(\mathbb D,S)$.
Then there exists a Zariski closed subset $Y$ such that Zariski closed subsets $Z\subset S$ satisfy $\mathcal{F}\to Z$ if and only if $Y\subset Z$. 
In particular, $\mathcal{F}\to Y$.
\end{lem} 
 
 {\it Proof.}\
 We denote by $\mathcal{Z}$ the set of all Zariski closed subsets $Z\subset S$ such that $\mathcal{F}\to Z$.
  We have $\mathcal{F}\to S$, hence $S\in\mathcal{Z}$.
 Hence $\mathcal{Z}\not=\emptyset$.
By the Noetherian property, there exists a minimal element $Y\in \mathcal{Z}$.
Then $\mathcal{F}\to Y$.
We note that if $Z\in\mathcal{Z}$, then $Y\subset Z$.
Indeed, if not, then $Y\cap Z\subsetneqq Y$.
By $\mathcal{F}\to Y$ and $\mathcal{F}\to Z$, we have $\mathcal{F}\to Y\cap Z$ (cf. Lemma \ref{lem:20220626}), hence $Y\cap Z\in\mathcal{Z}$.
This contradicts to the choice of $Y$.
Hence $Y\subset Z$.
Conversely, if a Zariski closed set $Z\subset S$ satisfies $Y\subset Z$, then by $\mathcal{F}\to Y$, we have $\mathcal{F}\to Z$.
Hence $Z\in\mathcal{Z}$.
Thus $Z\in\mathcal{Z}$ if and only if $Y\subset Z$.
 \hspace{\fill} $\square$

\begin{rem}\label{rem:20220626}
Suppose $\mathcal{F}=(f_i)_{i\in I}$ converges uniformly on compact subsets of $\mathbb D$ to a holomorphic map $g:\mathbb D\to S$.
Let $V\subset S$ be the Zariski closure of $g(\mathbb D)$.
Then $Y$ in Lemma \ref{lem:202206261} coincides with $V$.
To check this, we note $\mathcal{F}\to V$, hence $Y\subset V$.
To check the converse, we assume contrary that $V\not\subset Y$.
Then $g(\mathbb D)\not\subset Y$.
Hence $g^{-1}(Y)\subset \mathbb D$ is a discrete subset.
Hence we may take a closed disc $K\subset \mathbb D\backslash g^{-1}(Y)$ of positive radius.
We take an open neighbourhood $U\subset S$ of $Y$ such that $g(K)\cap \overline{U}=\emptyset$. 
Then since $\mathcal{F}$ converges to $g$ uniformly on $K$, there exists a finite subset $E\subset I$ such that for all $i\in I\backslash E$, we have $f_i(K)\cap U=\emptyset$.
This contradicts to $\mathcal{F}\to Y$.
Hence $V\subset Y$.
Thus $V=Y$.
\end{rem}
  
\begin{lem}\label{lem:202206260}
Let $\varphi:S'\to S$ be a morphism and let $Z\subset S$ be a Zariski closed set.
Let $\mathcal{F}=(f_i)_{i\in I}$ be an infinite indexed family in $\mathrm{Hol}(\mathbb D,S')$ such that $(\varphi\circ f_i)_{i\in I}\to Z$.
Suppose $\varphi:S'\to S$ is proper.
Then $\mathcal{F}\to \varphi^{-1}(Z)$.
\end{lem}

{\it Proof.}\
Let $s\in (0,1)$ and $\gamma>0$.
Let $U\subset S'$ be an open neighbourhood of $\varphi^{-1}(Z)$.
Since $\varphi$ is proper, $\varphi(S'\backslash U)$ is a closed subset.
Set $W=S\backslash \varphi(S'\backslash U)$.
Then $W\subset S$ is an open set such that $Z\subset W$.
We have $\varphi^{-1}(W)\subset U$.
By $(\varphi\circ f_i)_{i\in I}\to Z$, there exists a finite subset $E\subset I$ such that for all $i\in I\backslash E$ we have $\varphi\circ f_i(z)\in W$ for $\gamma$-almost all $z\in\mathbb D(s)$.
Then for all $i\in I\backslash E$ we have $f_i(z)\in U$ for $\gamma$-almost all $z\in\mathbb D(s)$.
Hence $\mathcal{F}\to \varphi^{-1}(Z)$.
\hspace{\fill} $\square$

\begin{lem}\label{lem:202105271}
Let $\mathcal{F}$ be an infinite indexed family in $\mathrm{Hol}(\mathbb D,S)$.
Let $Z_1$ and $Z_2$ be Zariski closed subsets of $S$.
If $\mathcal{F}\Rightarrow Z_1$ and $\mathcal{F}\Rightarrow Z_2$, then either $Z_1\subset Z_2$ or $Z_2\subset Z_1$.
\end{lem}

{\it Proof.}\
Assume contrary.
Set $V=Z_1\cap Z_2$ in the sense of scheme theory, namely $\mathcal{I}_V=\mathcal{I}_{Z_1}+\mathcal{I}_{Z_2}$ for the defining ideal sheaves in $\mathcal{O}_S$.
Set $S'=\mathrm{Bl}_VS$ with $\varphi:S'\to S$.
Let $Z_1'$ and $Z_2'$ be the strict transforms of $Z_1$ and $Z_2$, respectively.
Then 
\begin{equation}\label{eqn:20220522}
Z_1'\cap Z_2'=\emptyset.
\end{equation}

We prove this.
Note that $\varphi^*V\subset S'$ is a Cartier divisor.
Hence we may take a closed subscheme $Z_1''\subset S'$ such that $\varphi^*Z_1=\varphi^*V+Z_1''$.
Indeed we have $\mathcal{I}_{\varphi^*Z_1}\subset \mathcal{I}_{\varphi^*V}\subset \mathcal{O}_{S'}$, where $\mathcal{I}_{\varphi^*V}$ is an invertible sheaf.
So we take $Z_1''\subset S'$ so that $\mathcal{I}_{Z_1''}=\mathcal{I}_{\varphi^*Z_1}\otimes (\mathcal{I}_{\varphi^*V})^{-1}\subset \mathcal{O}_{S'}$.
Then $\mathcal{I}_{Z_1''}\cdot \mathcal{I}_{\varphi^*V}=\mathcal{I}_{\varphi^*Z_1}$.
Similarly there exists $Z_2''\subset S'$ such that $\mathcal{I}_{Z_2''}\cdot \mathcal{I}_{\varphi^*V}=\mathcal{I}_{\varphi^*Z_2}$.
Then by $\mathcal{I}_{\varphi^*Z_1}+\mathcal{I}_{\varphi^*Z_2}=\mathcal{I}_{\varphi^*V}$, we have $(\mathcal{I}_{Z_1''}+\mathcal{I}_{Z_2''})\cdot \mathcal{I}_{\varphi^*V}=\mathcal{I}_{\varphi^*V}$.
Hence $\mathcal{I}_{Z_1''}+\mathcal{I}_{Z_2''}=\mathcal{O}_{S'}$.
This shows $Z_1''\cap Z_2''=\emptyset$.
By $Z_1'\subset Z_1''$ and $Z_2'\subset Z_2''$, we get \eqref{eqn:20220522}.

Now note that $\mathrm{Bl}_VS\to S$ is $Z_1$-admissible.
Hence by $\mathcal{F}\Rightarrow Z_1$, we have $\mathcal{F}\to Z_1'$.
Similarly, $\mathcal{F}\to Z_2'$.
This is a contradiction.
\hspace{\fill} $\square$

\begin{lem}\label{lem:202105272}
Let $\mathcal{F}$ be an infinite indexed family in $\mathrm{Hol}(\mathbb D,S)$.
Let $Z\subset S$ be a Zariski closed set such that $\mathcal{F}\Rightarrow Z$.
Then $Z$ is irreducible.
\end{lem}

{\it Proof.}\
We assume contrary that $Z=Z_1\cup Z_2$.
Set $V=Z_1\cap Z_2$ in the sense of scheme theory.
Then $Z_1'\cap Z_2'=\emptyset$ in $\mathrm{Bl}_VS$, where $Z_1'$ and $Z_2'$ are the strict transforms (cf. \eqref{eqn:20220522}).
Let $\varphi:\mathrm{Bl}_VS\to S$, which is $Z$-admissible.
Set $U_1=S-Z_2$.
Then $V\cap U_1=\emptyset$.
Hence $\varphi^{-1}(U_1)\to U_1$ is an isomorphism.
Moreover we have $Z\cap U_1=Z_1\cap U_1\not=\emptyset$.
Hence the minimal transform $Z'\subset \mathrm{Bl}_ZS$ satisfies $Z'\subset \overline{\varphi^{-1}(Z_1\cap U_1)}\subset Z_1'$.
Similarly $Z'\subset Z_2'$.
Hence $Z'=\emptyset$.
This contradicts to $\mathcal{F}\to Z'$.
Hence $Z$ is irreducible.
\hspace{\fill} $\square$

\begin{lem}\label{lem:20220623}
Let $\mathcal{F}=(f_i)_{i\in I}$ be an infinite indexed family in $\mathrm{Hol}(\mathbb D,S)$.
Let $Z\subset S$ be a Zariski closed set such that $\mathcal{F}\Rightarrow Z$.
Let $\varphi:S'\to S$ be a $Z$-admissible modification and let $Z'\subset S'$ be the minimal transform.
Then $\mathcal{F}\Rightarrow Z'$.
\end{lem}

{\it Proof.}\
If $V\subset S'$ is a Zariski closed set such that $Z'\not\subset V$, then $Z\not\subset \varphi(V)$.
For all $i\in I$ with finite exception, we have $f_i(\mathbb D)\not\subset \varphi(V)$, hence $f_i'(\mathbb D)\not\subset V$, where $f_i':\mathbb D\to S'$ is the lift of $f_i:\mathbb D\to S$.
Let $S''\to S'$ be a $Z'$-admissible modification and let $Z''\subset S''$ be the minimal transform of $Z'$.
Then $S''\to S$ is $Z$-admissible and the minimal transform of $Z$ coincides with $Z''$.
These easily follow from the irreducibility of $Z$ (cf. Lemma \ref{lem:202105272}).
Hence $\mathcal{F}\to Z''$.
Hence $\mathcal{F}\Rightarrow Z'$.
\hspace{\fill} $\square$

\begin{defn}\label{defn:20220205}
Let $\mathcal{F}$ be an infinite indexed family in $\mathrm{Hol}(\mathbb D,S)$ and let $\{ E_i\}_{i\in I}$ be a countable family of Zariski closed sets in $S$.
$\mathrm{LIM}(\mathcal{F};\{ E_i\})$ is a Zariski closed set of $S$ with the following properties:
\begin{enumerate}
\item $\mathcal{F}\Rightarrow \mathrm{LIM}(\mathcal{F};\{ E_i\})$ and $\mathrm{LIM}(\mathcal{F};\{ E_i\})\not\subset \cup_{i\in I}E_i$
\item for every proper Zariski closed subset $W\subsetneqq \mathrm{LIM}(\mathcal{F};\{ E_i\})$ and every infinite indexed subfamily $\mathcal{G}$ of $\mathcal{F}$, either $\mathcal{G}\nRightarrow W$ or $W\subset \cup_{i\in I}E_i$. 
\end{enumerate}
\end{defn}

\begin{rem}\label{rem:20210805}
(1) $\mathrm{LIM}(\mathcal{F};\{ E_i\})$ is unique if it exists.
This follows from Lemma \ref{lem:202105271}.
Indeed, suppose that both $Z_1,Z_2\subset S$ satisfy the two conditions of Definition \ref{defn:20220205}.
Then $\mathcal{F}\Rightarrow Z_1$ and $\mathcal{F}\Rightarrow Z_2$.
By Lemma \ref{lem:202105271}, we may assume that $Z_1\subset Z_2$.
Assume contrary that $Z_1\not=Z_2$.
Then by the second condition of Definition \ref{defn:20220205}, we have $Z_1\subset \cup_{i\in I}E_i$, for $\mathcal{F}\Rightarrow Z_1$.
This contradicts to the first condition of Definition \ref{defn:20220205} that $Z_1\not\subset \cup_{i\in I}E_i$.
Hence $Z_1=Z_2$.
This shows that $\mathrm{LIM}(\mathcal{F};\{ E_i\})$ is unique if it exists.

(2) $\mathrm{LIM}(\mathcal{F};\{ E_i\})$ is irreducible, if it exists.
This follows from Lemma \ref{lem:202105272}.

(3) If $\mathrm{LIM}(\mathcal{F};\{ E_i\})$ exists, then for all infinite subfamily $\mathcal{G}$ of $\mathcal{F}$, $\mathrm{LIM}(\mathcal{G};\{ E_i\})$ exists and $\mathrm{LIM}(\mathcal{G};\{ E_i\})=\mathrm{LIM}(\mathcal{F};\{ E_i\})$.
This follows directly from Definition \ref{defn:20220205}
(cf. Remark \ref{rem:20220523}).
\end{rem}

\begin{lem}\label{lem:202105284}
Let $\mathcal{F}=(f_i)_{i\in I}$ be an infinite indexed family in $\mathrm{Hol}(\mathbb D,S)$.
Let $\{ E_j\}_{j\in J}$ be a countable family of Zariski closed sets of $S$.
Assume that for each $j\in J$, all but fintely many $i\in I$ satisfy $f_i(\mathbb D)\not\subset E_j$.
Then there exists an infinite indexed subfamily $\mathcal{G}$ of $\mathcal{F}$ such that $\mathrm{LIM}(\mathcal{G};\{ E_j\})$ exists.
\end{lem}

{\it Proof.}\
We denote by $\mathcal{Z}$ the set of Zariski closed subsets $Z\subset S$ such that $f_i(\mathbb D)\subset Z$ for infinitely many $i\in I$. 
We have $S\in \mathcal{Z}$, hence $\mathcal{Z}\not=\emptyset$.
By the Noetherian property, we may take a minimal $Z\in\mathcal{Z}$.
We take an infinite subset $I'\subset I$ such that $f_i(\mathbb D)\subset Z$ for all $i\in I'$.
Set $\mathcal{F}'=(f_i)_{i\in I'}$.
Then $\mathcal{F}'\Rightarrow Z$.
By Lemma \ref{lem:202105272}, $Z$ is irreducible.
We have $Z\not\subset \cup E_j$.
Indeed suppose $Z\subset \cup E_j$.
Then since $J$ is countable, there exists $j\in J$ such that $Z\subset E_j$.
Hence $f_i(\mathbb D)\subset E_j$ for all $i\in I'$, a contradiction.
Thus $Z\not\subset \cup E_j$.

For an infinite subfamily $\mathcal{G}$ of $\mathcal{F}'$, we denote by $\mathcal{W}_{\mathcal{G}}$ the set of Zariski closed subsets $W\subset S$ such that $\mathcal{G}\Rightarrow W$ and $W\not\subset \cup_{j\in J}E_j$.
By the above argument, we have $Z\in \mathcal{W}_{\mathcal{G}}$. 
Hence $\mathcal{W}_{\mathcal{G}}\not=\emptyset$.
By Lemma \ref{lem:202105271} and the Noetherian property, we may take a unique minimal element $W_{\mathcal{G}}\subset \mathcal{W}_{\mathcal{G}}$.
If $\mathcal{G}'$ is an infinite subfamily of $\mathcal{G}$, then $W_{\mathcal{G}}\subset \mathcal{W}_{\mathcal{G}'}$.
Hence $W_{\mathcal{G
}'}\subset W_{\mathcal{G}}$.
Hence by the Noetherian property, we may take an infinite subfamily $\mathcal{G}$ of $\mathcal{F}'$ such that $W_{\mathcal{G}}=W_{\mathcal{G}'}$ for all $\mathcal{G}'\subset \mathcal{G}$.
Then $\mathrm{LIM}(\mathcal{G};\{ E_j\})=W_{\mathcal{G}}$.
\hspace{\fill} $\square$

\begin{lem}\label{lem:20220219}
Let $p:S_1\to S_2$ be a morphism of varieties.
Let $\mathcal{F}=(f_i)_{i\in I}$ be an infinite indexed family in $\mathrm{Hol}(\mathbb D,S_1)$.
Let $Z\subset S_1$ be a Zariski closed set such that $\mathcal{F}\Rightarrow Z$.
Then $\{ p\circ f_i\}_{i\in I}\Rightarrow \overline{p(Z)}$, where $\overline{p(Z)}\subset S_2$ is the Zariski closure of $p(Z)$.
\end{lem}

{\it Proof.}\
Let $V\subset S_2$ be a Zariski closed set such that $\overline{p(Z)}\not\subset V$.
Then we have $Z\not\subset p^{-1}(V)$.
We have $f_i(\mathbb D)\not\subset p^{-1}(V)$ for all $i\in I$ with finite exception.
Hence $p\circ f_i(\mathbb D)\not\subset V$ for all $i\in I$ with finite exception.

Now let $S_2'\to S_2$ be a $\overline{p(Z)}$-admissible modification.
We may take a $Z$-admissible modification $S_1'\to S_1$ such that $p':S_1'\to S_2'$ exists. 
Let $Z'\subset S_1'$ and $\overline{p(Z)}'\subset S_2'$ be the minimal transforms of $Z$ and $\overline{p(Z)}$, respectively.
Then $p'(Z')\subset \overline{p(Z)}'$.
We have $\mathcal{F}\to Z'$.
Hence $\{ p'\circ f_i\}_{i\in I}\to \overline{p(Z)}'$.
Thus $\{ p\circ f_i\}_{i\in I}\Rightarrow \overline{p(Z)}$.
\hspace{\fill} $\square$

\subsection{Horizontally integrable}\label{subsec:h}

Let $\iota:U\hookrightarrow A\times S$ be an immersion, i.e., open of closed immersion.
Assume that $q:U\to S$ is {\'e}tale, where $q$ is the composite of the immersion $U\hookrightarrow A\times S$ and the second projection $A\times S\to S$.
Then we get an immersion $D\iota:PTU\hookrightarrow PT(A\times S)$.
By \eqref{eqn:20211203}, we have $PT(A\times S)=A\times S_{1,A}$.
Hence by the composite of the immersion $D\iota$ and the second projection $PT(A\times S)\to S_{1,A}$, we get $\iota' :PTU\to S_{1,A}$.

\begin{defn}\label{defn:20201122}
Let $Z\subset S_{1,A}$ be an irreducible Zariski closed set.
We say that $Z$ is {\it horizontally integrable} if there exists an immersion $\iota:U\hookrightarrow A\times S$ such that $q:U\to S$ is {\'e}tale and $Z\cap \iota'(PTU)\subset Z$ is Zariski dense in $Z$.
\end{defn}

We describe a simple example.
Let $S\hookrightarrow A\times S$ be an immersion induced from a constant map $S\to A$, and let $S_{1,A}^*\subset S_{1,A}$ be the image of $PTS\to PT(A\times S)\to S_{1,A}$.
Namely, we set 
\begin{equation}\label{eqn:20220707}
S_{1,A}^*=P(TS\times \{0\}).
\end{equation}
Then every irreducible Zariski closed set $Z\subset S_{1,A}$ such that $Z\subset S_{1,A}^*$ is horizontally integrable.

In general, let $\iota:U\hookrightarrow A\times S$ be an immersion such that $q:U\to S$ is {\'e}tale.
Then the induced map $q':U_{1,A}\to S_{1,A}$ is also {\'e}tale.
Let $\varphi:U\to A$ be the composite of $\iota$ and the first projection $A\times S\to A$.
We note that the immersion $\iota:U\hookrightarrow A\times S$ implies an isomorphism
\begin{equation}\label{eqn:20220708}
\Phi:A\times U\to  A\times U
\end{equation}
over $U$ by 
\begin{equation}\label{eqn:202302248}
\Phi (a,u)=(a-\varphi (u),u).
\end{equation}
This defines an isomorphism $PT(A\times U )\to PT(A\times U )$.
This induces
$$D\Phi :U_{1,A}\to U_{1,A},$$
which is an isomorphism over $U$.

\begin{lem}\label{lem:20220707}
Suppose $Z\subset S_{1,A}$ is horizontally integrable.
Let $\iota:U\hookrightarrow A\times S$ be an immersion as in Definition \ref{defn:20201122}.
Then there exists an irreducible component $\Xi\subset U_{1,A}$ of $(q')^{-1}(Z)$ such that $D\Phi(\Xi)\subset U_{1,A}^*$.
\end{lem}

{\it Proof.}\
The immersion $\iota$ induces the section $\iota_U:U\hookrightarrow A\times U$ of the second projection $A\times U\to U$.
Then $\iota_U$ induces $\iota_U':PTU\to U_{1,A}$ as above.
Note that $\iota':PTU\to S_{1,A}$ is the composite of $\iota_U':PTU\to U_{1,A}$ and $q':U_{1,A}\to S_{1,A}$.
Let $\Xi_1,\ldots,\Xi_l$ be the irreducible components of $(q')^{-1}(Z)$.
Then we have
$$\bigcup_{i}q'(\iota_U'(PTU)\cap \Xi_i)=\iota'(PTU)\cap Z.$$
There exists $\Xi_i$ such that $q'(\iota_U'(PTU)\cap \Xi_i)\subset Z$ is Zariski dense in $Z$.
We set $\Xi=\Xi_i$.
Since $q'$ is {\'e}tale, we have $\dim\Xi=\dim Z$.
Hence $\iota_U'(PTU)\cap \Xi\subset \Xi$ is Zariski dense in $\Xi$.

Since the composite of $U\overset{\iota_U}{\to} A\times U\overset{\Phi}{\to} A\times U$ is the graph of a constant map, the image of the composite of $PTU\overset{\iota_U'}{\to} U_{1,A}\overset{D\Phi}{\to}U_{1,A}$ is $U_{1,A}^*$.
Hence $\iota_U'(PTU)\subset (D\Phi)^{-1}U_{1,A}^*$, so $\Xi\cap \iota_U'(PTU)\subset (D\Phi)^{-1}U_{1,A}^*$.
Since $\Xi\cap \iota_U'(PTU)$ is Zariski dense in $\Xi$ and $(D\Phi)^{-1}U_{1,A}^*$ is Zariski closed, we have $\Xi\subset (D\Phi)^{-1}U_{1,A}^*$.
This concludes the proof.
\hspace{\fill} $\square$

\subsection{Statement of Proposition \ref{pro:320}}

We say that a smooth positive $(1,1)$-form $\omega_A$ on a semi-abelian variety $A$ is invariant if $\omega_A$ is invariant under the translation of $A$.

\begin{defn}\label{defn:20220601}
Given an infinite indexed family $\mathcal{F}=(f_i)_{i\in I}$ in $\mathrm{Hol}(\mathbb D,A)$, we denote by $\Pi(\mathcal{F})$ the set of all semi-abelian subvarieties $B\subset A$ satisfying the following property:
there is no infinite subset $J\subset I$ such that $\{ |(\varpi_B\circ f_i)'|_{\omega_{A/B}}\}_{i\in J}$ converges uniformly on compact subsets of $\mathbb D$ to $0$.
Here $\varpi_B: A\to A/B$ is the quotient map and $|\cdot|_{\omega_{A/B}}$ is a norm on $T(A/B)$ defined by an invariant (1,1)-form $\omega_{A/B}$ on $A/B$. 
\end{defn}

Note that every semi-abelian variety contains only countably many semi-abelain varieties (cf. \cite[Cor. 5.1.9]{NW}).
Hence $\Pi(\mathcal{F})$ is a countable set.
We note $A\not\in\Pi(\mathcal{F})$.
Hence if $B\in \Pi(\mathcal{F})$, then $E_{k,A,A/B}\subsetneqq P_{k,A}$ for all $k\geq 1$.

For a non-constant holomorphic map $f:\mathbb D\to A\times S$, we recall the notation $f_{S_{k,A}}:\mathbb D\to S_{k,A}$ from \eqref{eqn:20220206}.
We use the notation $f_{P_{k,A}}:\mathbb D\to P_{k,A}$ if $S=\{ \mathrm{pt}\}$, where $P_{k,A}=\{\mathrm{pt}\}_{k,A}$.
Given an infinite indexed family $\mathcal{F}=(f_i)_{i\in I}$ of non-constant holomorphic maps in $\mathrm{Hol}(\mathbb D,A)$, we define an infinite indexed family $\mathcal{F}_{P_{k,A}}$ in $\mathrm{Hol}(\mathbb D,P_{k,A})$ by
\begin{equation}\label{eqn:20230326}
\mathcal{F}_{P_{k,A}}=( (f_i)_{P_{k,A}})_{i\in I}.
\end{equation}
We consider the following assumption for an infinite indexed family $\mathcal{F}$ of non-constant holomorphic maps in $\mathrm{Hol}(\mathbb D,A)$.

\begin{ass}\label{ass:202103061}
$\mathrm{LIM}(\mathcal{F}_{P_{k,A}};\{ E_{k,A,A/B}\}_{B\in\Pi(\mathcal{F})})$ exists for all $k\geq 1$.
\end{ass}

If this assumption is satisfied, we write $T_k=\mathrm{LIM}(\mathcal{F}_{P_{k,A}};\{ E_{k,A,A/B}\}_{B\in\Pi(\mathcal{F})})\subset P_{k,A}$.
As we shall see later (cf. Lemma \ref{lem:20220523}), every infinite indexed family $\mathcal{F}$ in $\mathrm{Hol}(\mathbb D,A)$ contains an infinite  subfamily which satisfies this assumption.

Now we take an infinite indexed family $\mathcal{F}$ in $\mathrm{Hol}(\mathbb D,A)$ which satisfies Assumption \ref{ass:202103061}.
Let $k\geq 1$.
We claim 
\begin{equation}\label{eqn:20220219}
T_k\subset p(T_{k+1}),
\end{equation}
where $p:P_{k+1,A}\to P_{k,A}$ is the natural map.
Indeed, by Lemma \ref{lem:20220219}, we have $\mathcal{F}_{P_{k,A}}\Rightarrow p(T_{k+1})$.
Hence by Lemma \ref{lem:202105271}, we have either $T_k\subset p(T_{k+1})$ or $p(T_{k+1})\subsetneqq T_k$.
So assume contrary that $p(T_{k+1})\subsetneqq T_k$.
Then by the definition of $T_k$, we have $p(T_{k+1})\subset \cup_{B\in\Pi(\mathcal{F})} E_{k,A,A/B}$.
We have $p^{-1}(E_{k,A,A/B})\subset E_{k+1,A,A/B}$ (cf. Lemma \ref{lem:202206211}).
Hence
$$
T_{k+1}\subset p^{-1}(\cup_{B\in\Pi(\mathcal{F})} E_{k,A,A/B})\subset \cup_{B\in\Pi(\mathcal{F})}E_{k+1,A,A/B}.
$$
This is a contradiction.
Hence we have proved \eqref{eqn:20220219}.

We fix $k\geq 0$.
For each $l\geq 1$, let $p_l:P_{k+l,A}\to P_{k+1,A}$ be the natural map.
Then by \eqref{eqn:20220219}, we have a sequence 
$$p_1(T_{k+1})\subset p_{2}(T_{k+2})\subset p_{3}(T_{k+3})\subset\cdots.$$
By Remark \ref{rem:20210805} (2), each $p_l(T_{k+l})$ is irreducible.
Hence the sequence above stabilizes, i.e., there exists $l_0\geq 1$ such that $p_l(T_{k+l})=p_{l_0}(T_{k+l_0})$ for all $l\geq l_0$.
We set $Z=p_{l_0}(T_{k+l_0})$, namely
\begin{equation}\label{eqn:202202191}
Z=\cup_{l\geq 1}p_l(T_{k+l}).
\end{equation}

By Remark \ref{rem:20230312}, we have $Z\subset P_{k+1,A}\subset (P_{k,A})_{1,A}$.

\begin{pro}\label{pro:320}
Let $\mathcal{F}\subset \mathrm{Hol}(\mathbb D, A)$ be an infinite set of non-constant holomorphic maps which satisfies Assumption \ref{ass:202103061}.
Suppose that there exists $k\geq 0$ such that $Z\subset P_{k+1,A}\subset (P_{k,A})_{1,A}$ is horizontally integrable, where $Z$ is defined by \eqref{eqn:202202191}.
Then $\mathcal{F}$ is normal for every smooth equivariant compactification $\overline{A}$.
Namely, for every smooth equivariant compactification $\overline{A}$ and every sequence $(f_n)_{n\in\mathbb N}$ in $\mathcal{F}$, there exists a subsequence of $(f_n)_{n\in\mathbb N}$ which converges uniformly on compact subsets of $\mathbb D$ to $g:\mathbb D\to \overline{A}$.
\end{pro}

This proposition will be proved in Section \ref{sec:8} after the preparation in Sections \ref{sec:4}-\ref{sec:7}.

\section{Nevanlinna theory}\label{sec:4}

Let $V$ be an algebraic variety.
We denote by $\mathrm{Hol_m}(\mathbb D, V)$ all multi-valued holomorphic maps to $V$, i.e.,
\begin{equation*}
\begin{CD}
Y_f@>f>> V\\
@V\pi_fVV    \\
\mathbb D
\end{CD}
\end{equation*}
where $\pi_f:Y_f\to \mathbb D$ is a proper, surjective map.

We introduce the notion of Weil functions (cf., e.g., \cite[Def. 2.2.1]{Y2}).
Let $V$ be a projective variety and let $Z\subset V$ be a closed subscheme.
A Weil function $\lambda_Z$ for $Z$ is a continuous function $\lambda_Z:V-\supp Z\to \mathbb R$ which satisfies the following condition.
For each $x\in V$, there are a Zariski open neighborhood $U\subset V$ of $x$, holomorphic functions $g_1,\ldots ,g_l\in \Gamma (U,\mathcal{O}_V)$ which defines $Z\cap U$, and a continuous function $\alpha :U\to \mathbb R$ on $U$ such that
$$
\left\vert \lambda _Z(y)+\log \max_{1\leq i\leq l} \{  |g_i(y)|\} \right\vert \leq \alpha (y)
$$
for all $y\in U-\supp (Z\cap U)$.
We summarize the needed properties of Weil functions (cf. \cite[Section 2.2]{Y2}):

\begin{itemize}
\item
A Weil function $\lambda_Z$ exists for every closed subscheme $Z\subset V$.  
\item
If $\lambda_Z$ and $\lambda_Z'$ are Weil functions for $Z$, then there exists a positive constant $\gamma$ such that $|\lambda_Z(x)-\lambda_{Z}'(x)|\leq \gamma$ for all $x\in V-\supp Z$.
\item If $\lambda_Z$ is a Weil function for $Z$, then $\lambda_Z$ is bounded from below, namely there exists a constant $\gamma$ such that $\lambda_Z(x)>\gamma$ for all $x\in V-\supp Z$.
In particular, we may choose a Weil function $\lambda_Z$ such that $\lambda_Z\geq 0$.
\item 
Suppose that $p:\tilde{V}\to V$ is a morphism from another projective variety $\tilde{V}$.
Then $\lambda_Z\circ p$ is a Weil function for the pull-back $p^*Z\subset \tilde{V}$.
\item
Let $Z,Z'$ be closed subschemes of $V$.
Let $\lambda_Z$ and $\lambda_{Z'}$ be Weil functions for $Z$ and $Z'$, respectively.
Then $\lambda_Z+\lambda_{Z'}$ is a Weil function for $Z+Z'$, where $Z+Z'$ is defined by $\mathcal{I}_{Z+Z'}=\mathcal{I}_Z\cdot \mathcal{I}_{Z'}$.
\item
Let $Z,Z'$ be closed subschemes of $V$.
Assume that $\supp Z\subset \supp Z'$.
Then there exist positive constants $\gamma_1>0$ and $\gamma_2>0$ such that $\lambda_{Z}\leq \gamma_1\lambda_{Z'}+\gamma_2$.
\item 
Assume that $V$ is smooth.
Let $D$ be an effective Cartier divisor on $V$.
Let $L$ be a line bundle on $V$ associated to $D$, and let $h$ be a smooth Hermitian metric on $L$.
Let $\sigma$ be a section of $L$ associated to $D$ such that $h(\sigma (x),\sigma (x))\leq 1$ for all $x\in V$.
Then $\lambda_D(x)=-\log \sqrt{h(\sigma (x),\sigma (x))}$, where $x\in V-\supp D$, is a Weil function for $D$ with $\lambda_D\geq 0$.
\end{itemize}

\begin{rem}
In this paper, we always assume $\lambda_Z\geq 0$ for Weil functions unless otherwise specified.
\end{rem}

We introduce Nevanlinna theory.
Let $V$ be a projective variety.
Let $Z\subset V$ be a closed subscheme.
Let $f\in \mathrm{Hol_m}(\mathbb D,V)$ such that $f(Y_f)\not\subset \supp Z$.
For $0<s<r<1$, we set
$$
N_s(r,f,Z)=\frac{1}{\deg\pi_f}\int_s^r\left( \sum_{y\in Y_f(t)}\ord_yf^*Z\right)\frac{dt}{t},
$$
where $Y_f(t)=\pi_f^{-1}(\mathbb D(t))$.
Let $\lambda_Z$ be a Weil function for $Z$.
We set
$$
m(r,f,\lambda_Z)=\frac{1}{\deg\pi_f}\int_{y\in\partial Y_f(r)}\lambda_Z(f(y))\frac{d\arg \pi_f(y)}{2\pi}.
$$
We set
$$
T_s(r,f,Z,\lambda_Z)=N_s(r,f,Z)+m(r,f,\lambda_Z)-m(s,f,\lambda_Z).
$$

Assume $V$ is smooth.
Let $\omega$ be a smooth (1,1)-form on $V$.
Let $f\in \mathrm{Hol_m}(\mathbb D,V)$.
For $0\leq s<r<1$, we set
$$
T_s(r,f,\omega)=\frac{1}{\deg \pi_{f}}\int_s^r\frac{dt}{t}\int_{Y_f(t)}f^*\omega.
$$
Let $D\subset V$ be an effective Cartier divisor.
Let $L$ be the associated line bundle and let $h$ be a smooth Hermitian metric on $L$.
Let $\omega_{(L,h)}$ be the curvature form for the metrized line bundle $(L,h)$.
Let $\sigma$ be a section of $L$ such that $D=(\sigma=0)$.
Set $\lambda_D(x)=-\log \sqrt{h(\sigma (x),\sigma(x))}$.
The first main theorem states that
\begin{equation}\label{eqn:20211026}
T_s(r,f,\omega_{(L,h)})=N_s(r,f,D)+m(r,f,\lambda_D)-m(s,f,\lambda_D).
\end{equation}
We give a sketch of the proof.
By the Poinca{\'e}-Lelong formula, we have 
$$2dd^c\log (1/||\sigma\circ f||)=-\sum_{y\in Y_f}(\mathrm{ord}_yf^{*}D)\delta _y+f^{*}\omega_{(L,h)},$$
where $\delta_y$ is Dirac current suported on $y$.
Integrating over $Y_f(t)$, we get 
$$
2\int_{ Y_f(t)}dd^c\log (1/||\sigma\circ f||)=-\sum_{Y_f(t)}\mathrm{ord}_yf^{*}D+\int_{ Y_f(t)}f^{*}\omega_{(L,h)}.
$$
Hence, we get
\begin{equation*}
\begin{split}
-N_s(r,f,D)+T_s(r,f,\omega_{(L,h)})
&=\frac{2}{\deg \pi_f}\int_{s}^r\frac{dt}{t}\int_{Y_f(t)}dd^c\log \left( \frac{1}{||\sigma\circ f||}\right)\\
&=\frac{2}{\deg \pi_f}\int_{s}^r\frac{dt}{t}\int_{\partial Y_f(t)}d^c\log \left( \frac{1}{||\sigma\circ f||}\right)\\
&=m(r,f,\lambda_D)-m(s,f,\lambda_D).
\end{split}
\end{equation*}
This proves \eqref{eqn:20211026}.

As a corollary, we have the following:
Let $V$ be a projective variety, which is not necessarily smooth.
Let $D$ and $D'$ be linearly equivalent effective Cartier divisors on $V$.
Then we have
\begin{equation}\label{eqn:20211217}
|T_s(r,f,D,\lambda_D)-T_s(r,f,D',\lambda_{D'})|\leq c
\end{equation}
for all $0<s<r<1$, where $c$ is a positive constant which only depends on the choice of the Weil functions $\lambda_D$ and $\lambda_{D'}$.
Indeed, if $D$ and $D'$ are very ample, we reduce to the case of projective spaces by taking an embedding $\iota:V\hookrightarrow \mathbb P^k$ so that $D=\iota^*H$ and $D'=\iota^*H'$ for hyperplane sections $H$ and $H'$.
In this case, \eqref{eqn:20211217} follows from \eqref{eqn:20211026}.
In general, we take an effective Cartier divisor $E\subset V$ such that $D+E$ and $D'+E$ are very ample to reduce to the previous case.

\begin{lem}\label{lem:20211021}
Let $V$ be a projective variety.
Let $D\subset V$ be an effective Cartier divisor which is ample.
Let $Z\subset V$ be a closed subscheme.
Let $\lambda_D\geq 0$ and $\lambda_Z\geq 0$ be Weil functions for $D$ and $Z$, respectively.
Then there exist positive constants $c>0$, $c'>0$ such that for all $f\in \mathrm{Hol_m}(\mathbb D,V)$ with $f(Y_f)\not\subset \supp D\cup \supp Z$, we have
$$
T_s(r,f,Z,\lambda_Z)\leq cT_s(r,f,D,\lambda_D)+c'
$$
for all $0<s<r<1$.
\end{lem}

{\it Proof.}\
Set $V'=\mathrm{Bl}_{Z}V$.
Let $\varphi:V'\to V$ be the projection and let $\varphi^*Z\subset V'$ be the induced closed subscheme.
Then $\varphi^*Z$ is a Cartier divisor on $V'$.
We denote by $L$ the associated line bundle.
Let $M$ be the ample line bundle on $V$ associated to $D$.
Then there exists a positive integer $l$ such that $\varphi^*M^{\otimes l}\otimes L^{-1}$ is very ample on $V'$ (cf. \cite[II, Prop. 7.10 (b)]{H}).
There exists a closed immersion $\iota:V'\hookrightarrow \mathbb P^k$ such that $\iota^*\mathcal{O}_{\mathbb P^k}(1)=\varphi^*M^{\otimes l}\otimes L^{-1}$.
Let $H\subset \mathbb P^k$ be an effective divisor from $\mathcal{O}_{\mathbb P^k}(1)$.
Then $\varphi^*(lD)$ and $\varphi^*Z+\iota^*H$ are linearly equivalent Cartier divisors.
Let $h_{\mathcal{O}_{\mathbb P^k}(1)}$ be the Fubini-Study metric on $\mathcal{O}_{\mathbb P^k}(1)$.
Let $\sigma$ be a section of $\mathcal{O}_{\mathbb P^k}(1)$ such that $H=(\sigma=0)$.
Set $\lambda_H(x)=-\log \sqrt{h_{\mathcal{O}_{\mathbb P^k}(1)}(\sigma (x),\sigma(x))}\geq 0$.
By \eqref{eqn:20211217}, we have
\begin{equation*}
T_s(r,f,Z,\lambda_Z)+T_s(r,\iota\circ f,H,\lambda_H)\leq lT_s(r,f,D,\lambda_D)+\alpha,
\end{equation*}
where $\alpha$ is a positive constant which only depends on the choices of $\lambda_D$, $\lambda_Z$ and $\lambda_H$.
Set $\omega_{\mathbb P^k}=c_1(\mathcal{O}_{\mathbb P^k}(1),h_{\mathcal{O}_{\mathbb P^k}(1)})$.
By \eqref{eqn:20211026}, we have 
$$
T_s(r,\iota\circ f,H,\lambda_H)=
T_s(r,\iota\circ f,\omega_{\mathbb P^k})\geq 0.
$$
Hence we get
\begin{equation*}
T_s(r,f,Z,\lambda_Z)\leq lT_s(r,f,D,\lambda_D)+\alpha.
\end{equation*}
This conclude the proof.
\hspace{\fill} $\square$

\begin{lem}\label{lem:20211008}
Let $V$ be a smooth projective variety and let $\omega_V$ be a smooth positive (1,1)-form on $V$.
Let $Z\subset V$ be a closed subscheme with a Weil fuction $\lambda_Z\geq 0$.
Then there exists a positive constant $c>0$ such that for all $f\in\mathrm{Hol_m}(\mathbb D,V)$ with $f(Y_f)\not\subset \supp Z$, we have
$$
m(r,f,\lambda_Z)\leq cT_s(r,f,\omega_V)+m(s,f,\lambda_Z)+c
$$
for all $s\in (0,1)$ and $r\in (s,1)$.
\end{lem}

{\it Proof.}\
Since $N_s(r,f,Z)\geq 0$, we have
$$
m(r,f,\lambda_Z)-m(s,f,\lambda_Z)\leq T_s(r,f,Z,\lambda_Z).
$$
Let $D_1,\ldots,D_l\subset V$ be effective ample divisors such that $D_1\cap\cdots\cap D_l=\emptyset$.
Let $\lambda_{D_i}\geq 0$ be a Weil function for $D_i$.
We may take $D_i$ such that $f(Y_f)\not\subset \supp D_i$.
By Lemma \ref{lem:20211021}, we have
$$
m(r,f,\lambda_Z)-m(s,f,\lambda_Z)\leq c_iT_s(r,f,D_i,\lambda_{D_i})+c_i'.
$$
By \eqref{eqn:20211026}, we have
\begin{equation}\label{eqn:20211021}
T(r,f,D_i,\lambda_{D_i})\leq c_i''T_s(r,f,\omega_V).
\end{equation}
We set $\displaystyle{c=\max_{1\leq i\leq l}\max\{ c_ic_i'',c_i'\}}$ to complete the proof.
\hspace{\fill} $\square$

\medskip

We apply Lemma \ref{lem:20211008} for $f:\mathbb D\to V$ with $f(0)\not\in \supp Z$.
Then $m(s,f,\lambda_Z)\to \lambda_Z(f(0))$ when $s\to 0+$.
Hence we get
\begin{equation}\label{eqn:20211216}
m(r,f,\lambda_Z)\leq cT_0(r,f,\omega_V)+\lambda_Z(f(0))+c.
\end{equation}

\begin{lem}\label{lem:20230207}
Let $p:\Sigma\to V$ be a generically finite surjective morphism of projective varieties, where $V$ is smooth.
Let $D$ be an effective Cartier divisor on $\Sigma$ with a Weil function $\lambda_D\geq 0$.
Let $\omega_V$ be a smooth positive (1,1)-form on $V$.
Let $E\subset \Sigma$ be the exceptional locus of $\varphi:\Sigma\to \Sigma^+$, where $\Sigma\to\Sigma^+\to V$ is the Stein factorization.
Let $\lambda_E\geq 0$ be a Weil function for $E$.
Then there exists a positive constant $c>0$ such that for all $f\in\mathrm{Hol_m}(\mathbb D,\Sigma)$ with $f(Y_f)\not\subset \supp D\cup\supp E$, we have
$$
T_s(r,f,D,\lambda_D)\leq cT_s(r,p\circ f,\omega_V)+cm(s,f,\lambda_E)+c
$$
for all $0<s<r<1$.
\end{lem}

{\it Proof.}\
Let $\widetilde{D}\subset \Sigma^+$ be the scheme theoretic image of the composite $D\hookrightarrow \Sigma\to\Sigma^+$.
Then we have $\mathcal{I}_{\varphi^{*}\widetilde{D}}\subset \mathcal{I}_{D}\subset \mathcal{O}_{\Sigma}$.
Since $\mathcal{I}_{D}$ is an invertible sheaf, we have a closed subscheme $Z\subset \Sigma$ such that $\varphi^{*}\widetilde{D}=D+Z$.
Since $\Sigma-E\to \Sigma^+$ is an open immersion, we note that $\supp Z\subset E$.
Hence we have a positive constant $c'>0$ such that $\lambda_{Z}\leq c'\lambda_{E}+c'$.
Thus
\begin{equation*}
\begin{split}
T_s(r,f,D,\lambda_{D})&\leq
T_s(r,f,\varphi^{*}\widetilde{D},\varphi^{*}\lambda_{\widetilde{D}})+m(s,f,\lambda_{Z})
\\
&\leq
T_s(r,\varphi\circ f,\widetilde{D},\lambda_{\widetilde{D}})+c'm(s,f,\lambda_E)+c'.
\end{split}
\end{equation*}
Let $G_1,\ldots,G_{\nu}\subset V$ be effective ample divisors such that $G_1\cap\cdots\cap G_{\nu}=\emptyset$.
Then since $q:\Sigma^+\to V$ is finite, $q^*G_j\subset \Sigma^+$ are ample.
We may take $G_j$ such that $\varphi\circ f(Y_f)\not\subset \supp q^*G_j$.
By Lemma \ref{lem:20211021}, we get
\begin{equation*}
\begin{split}
T_s(r,\varphi\circ f,\widetilde{D},\lambda_{\widetilde{D}})
&\leq \gamma_{j}T_s(r,\varphi\circ f,q^*G_j,\lambda_{G_j}\circ q)+\gamma_{j}'
\\
&=\gamma_{j}T_s(r,p\circ f,G_j,\lambda_{G_j})+\gamma_{j}'
\end{split}
\end{equation*}
By \eqref{eqn:20211026}, we have 
$$T_s(r,p\circ f,G_j,\lambda_{G_j})\leq \mu_j T_s(r,p\circ f,\omega_V).$$
We set $\displaystyle{c=\max_{1\leq j\leq \nu}\max\{\gamma_{j}\mu_j,c'+\gamma_j'\}}$ to conclude the proof.
\hspace{\fill} $\square$

\begin{lem}\label{lem:202110212}
Let $p:\Sigma\to V$ be a generically finite surjective morphism of smooth projective varieties. 
Let $\omega_{\Sigma}$ and $\omega_V$ be smooth positive (1,1)-forms on $\Sigma$ and $V$, respectively.
Let $E\subset \Sigma$ be the exceptional locus of $\varphi:\Sigma\to \Sigma^+$, where $\Sigma\to\Sigma^+\to V$ is the Stein factorization.
Let $\lambda_E\geq 0$ be a Weil function for $E$.
Then there exists a positive constant $c>0$ such that for all $f\in\mathrm{Hol_m}(\mathbb D,\Sigma)$ with $f(Y_f)\not\subset \supp E$, we have
$$
T_s(r,f,\omega_{\Sigma})\leq cT_s(r,p\circ f,\omega_V)+cm(s,f,\lambda_E)+c
$$
for all $0<s<r<1$.
\end{lem}

{\it Proof.}\
Let $D_1,\ldots,D_l\subset \Sigma$ be linearly equivalent, effective ample divisors such that $D_1\cap\cdots\cap D_l=\emptyset$.
Let $L$ be the associated line bundle on $\Sigma$.
Then since $L$ is ample, there exists a smooth Hermitian metric $h$ on $L$ such that the associated curvature form $\omega_{(L,h)}$ is positive.
Then there exists a positive constant $\alpha>0$ such that $\omega_{\Sigma}<\alpha\omega_{(L,h)}$.
Let $\sigma_i$ be a section of $L$ such that $D_i=(\sigma_i=0)$.
Set $\lambda_{D_i}(x)=-\log \sqrt{h(\sigma_i(x),\sigma_i(x))}\geq 0$.
We may take $D_i$ such that $f(Y_f)\not\subset \supp D_i$.
Then by \eqref{eqn:20211026}, we have
$$
T_s(r,f,\omega_{\Sigma})\leq \alpha T_s(r,f,D_i,\lambda_{D_i}).
$$
By Lemma \ref{lem:20230207}, we have
$$
T_s(r,f,D_i,\lambda_{D_i})\leq c_iT_s(r,p\circ f,\omega_V)+c_im(s,f,\lambda_E)+c_i
$$
for all $0<s<r<1$.
We set $\displaystyle{c=\max_{1\leq i\leq l}\{\alpha c_i\}}$ to conclude the proof.
\hspace{\fill} $\square$

\begin{rem}
Let $V$ be a smooth projective variety and let $Z\subset V$ be a closed subscheme.
Suppose $\mathrm{Bl}_ZV$ is smooth.
Then there exists a positive constant $c>0$ such that for all $f:\mathbb D\to V$ with $f(\mathbb D)\not\subset \supp Z$, we have
\begin{equation}\label{eqn:202110212}
T_s(r,f,\omega_{\mathrm{Bl}_ZV})\leq cT_s(r,f,\omega_V)+cm(s,f,\lambda_{\supp Z})+c
\end{equation}
for all $s\in (0,1)$ and $r\in (s,1)$.
This follows from Lemma \ref{lem:202110212} applied to $\mathrm{Bl}_ZV\to V$.
Indeed we have $m(s,f,\lambda_E)\leq cm(s,f,\lambda_{\supp Z})$.
\end{rem}

\begin{rem}\label{rem:20210803}
Let $\mathcal{F}\subset \mathrm{Hol}(\mathbb D,V)$ be an infinite set, and let $Z\subset V$ be a Zariski closed set.
For $s\in (0,1)$, $\gamma>0$, and an open neighbourhood $U\subset V$ of $Z$, we set
$$
\mathcal{F}_{s,\gamma,U}=\{ f\in \mathcal{F};\ \text{$f(z)\in U$ for $\gamma$-almost all $z\in \mathbb D(s)$}
\}.
$$
Then if $s'\geq s$, $\gamma'\leq \gamma$ and $U'\subset U$, we have $\mathcal{F}_{s',\gamma',U'}\subset \mathcal{F}_{s,\gamma,U}$.
We remark the following two facts which directly follow from the definition.
\begin{enumerate}
\item
$\mathcal{F}\not\to Z$ if and only if there exist $s\in (0,1)$, $\gamma>0$, $Z\subset U\subset V$ such that $\mathcal{F}-\mathcal{F}_{s,\gamma,U}$ is infinite.
\item
$\mathcal{G}\not\to Z$ for all infinite subset $\mathcal{G}\subset \mathcal{F}$ if and only if there exist $s\in (0,1)$, $\gamma>0$, $Z\subset U\subset V$ such that $\mathcal{F}_{s,\gamma,U}$ is finite.
Moreover if $f(\mathbb D)\not\subset Z$ for all $f\in \mathcal{F}$, then we may take $s$, $\gamma$, $U$ such that $\mathcal{F}_{s,\gamma,U}=\emptyset$.
\end{enumerate}
\end{rem}

\begin{lem}\label{lem:20210803}
Let $\mathcal{F}\subset \mathrm{Hol}(\mathbb D,V)$ be an infinite set, and let $Z\subset V$ be a Zariski closed set.
Assume that $\mathcal{F}\not\to Z$.
Then there exists an infinite subset $\mathcal{G}\subset \mathcal{F}$ such that $f(\mathbb D)\not\subset Z$ for all $f\in \mathcal{G}$, and that $\mathcal{G}'\not\to Z$ for all infinite subset $\mathcal{G}'\subset \mathcal{G}$.
\end{lem}

{\it Proof.}\
We may take $s\in (0,1)$, $\gamma>0$, $Z\subset U\subset V$ such that $\mathcal{F}-\mathcal{F}_{s,\gamma,U}$ is infinite. 
Set $\mathcal{G}=\mathcal{F}-\mathcal{F}_{s,\gamma,U}$.
Then $\mathcal{G}_{s,\gamma,U}=\mathcal{G}\cap \mathcal{F}_{s,\gamma,U}=\emptyset$.
Hence $\mathcal{G}'\not\to Z$ for all infinite subset $\mathcal{G}'\subset \mathcal{G}$.
Note that $f(\mathbb D)\not\subset Z$ for all $f\in \mathcal{G}$.
\hspace{\fill} $\square$

\begin{lem}\label{lem:20201119}
Let $V$ be a smooth projective variety.
Let $\omega_V$ be a smooth positive (1,1)-form.  
Let $\mathcal{F}\subset \mathrm{Hol}(\mathbb D,V)$ be an infinite subset.
Let $Z\subset V$ be a closed subscheme with a Weil function $\lambda_Z$.
Assume that $\mathcal{G}\not\to \supp Z$ for all infinite subset $\mathcal{G}\subset \mathcal{F}$ and that $f(\mathbb D)\not\subset \supp Z$ for all $f\in \mathcal{F}$.
Then there exist $\sigma\in (0,1)$ and $\alpha>0$ such that for all $s\in (\sigma,1)$, $r\in(s,1)$ and $f\in \mathcal{F}$, we have
$$
m(r,f,\lambda_Z)
\leq \frac{\alpha}{(r-s)}
T_{s}(r,f,\omega_{V})
+\alpha.
$$
\end{lem}

To prove this lemma, we start from the following consideration.
Let $w\in\mathbb D(r)$.
We set
\begin{equation}\label{eqn:20210716}
\varphi_{w,r}(z)=r^2\frac{w-z}{r^2-\overline{w}z},
\end{equation}
which is an isomorphism $\varphi_{w,r}:\mathbb D(r)\to \mathbb D(r)$ such that $\varphi_{w,r}(0)=w$.
We have $\varphi_{w,r}'(z)=r^2\frac{|w|^2-r^2}{(r^2-\overline{w}z)^2}$.
Applying the maximum and minimum principles on $\mathbb D(r)$, we get
\begin{equation}\label{eqn:202107161}
\frac{r-|w|}{r+|w|}\leq |\varphi_{w,r}'(z)|\leq \frac{r+|w|}{r-|w|}
\end{equation}
for all $z\in \overline{\mathbb D(r)}$.

\begin{lem}\label{lem:20210711}
Let $0<\sigma<r<1$.
Then for all non-negative function $\Lambda$ on $\partial \mathbb D(r)$ and $w\in \mathbb D(\sigma)$, we have
$$
\frac{r-\sigma}{r+\sigma}\int_0^{2\pi}\Lambda(\varphi_{w,r}(re^{i\theta}))\frac{d\theta}{2\pi}
\leq 
\int_0^{2\pi}\Lambda(re^{i\theta})\frac{d\theta}{2\pi}
\leq \frac{r+\sigma}{r-\sigma} \int_0^{2\pi}\Lambda(\varphi_{w,r}(re^{i\theta}))\frac{d\theta}{2\pi}.
$$
\end{lem}

{\it Proof.}\
We have
$$
\int_0^{2\pi}\Lambda(re^{i\theta})\frac{d\theta}{2\pi}=\int_0^{2\pi}\Lambda(\varphi_{w,r}(re^{i\theta}))|\varphi_{w,r}'(re^{i\theta})|\frac{d\theta}{2\pi}.
$$
Hence by \eqref{eqn:202107161}, we have
$$
\frac{r-|w|}{r+|w|}\int_0^{2\pi}\Lambda(\varphi_{w,r}(re^{i\theta}))\frac{d\theta}{2\pi}
\leq 
\int_0^{2\pi}\Lambda(re^{i\theta})\frac{d\theta}{2\pi}
\leq \frac{r+|w|}{r-|w|} \int_0^{2\pi}\Lambda(\varphi_{w,r}(re^{i\theta}))\frac{d\theta}{2\pi}.
$$
Since $r\in (\sigma,1)$ and $w\in \overline{\mathbb D(\sigma)}$, we have $\frac{r+|w|}{r-|w|}<\frac{r+\sigma}{r-\sigma}$.
The proof is completed.
\hspace{\fill} $\square$

\medskip

In the proof of Lemma \ref{lem:20201119}, we need the following estimate from \cite[Lemma 6.17]{Hay}:
If $\mu$ is a mass distribution on $\mathbb C$ with finite total mass $M$ and $\gamma$ is a constant with $0<\gamma <1$, then we have
\begin{equation}\label{eqn:bce}
\int_{\mathbb C}\log \frac{1}{|z -w|}d\mu_{z}\leq \tau_{\gamma} M
\end{equation}
for $\gamma$-almost all $w\in \mathbb C$, where $\tau_{\gamma} >0$ is a positive constant which depends on $\gamma$.
For instance, we may take as $\tau_{\gamma} =\log (6/\gamma )$.
See also \cite[VIII, \S 3]{Lan}.

\medskip

{\it Proof of Lemma \ref{lem:20201119}.}\
By Remark \ref{rem:20210803}, there exist $s_0\in (0,1)$, $\gamma_0>0$ and an open neighbourhood $U_0\subset V$ of $\supp Z$ such that $\mathcal{F}_{s_0,\gamma_0,U_0}=\emptyset$.
We fix $\sigma\in (s_0,1)$.
Let $f\in \mathcal{F}$.
By Lemma \ref{lem:20211008}, we get
$$
m(r,f,\lambda_Z)\leq cT_{\sigma}(r,f,\omega_V)+m(\sigma,f,\lambda_Z)
$$
for all $r\in (\sigma,1)$, where $c>0$ is a positive constant which only depends on $\omega_V$ and $\lambda_Z$.
Given an isomorphism $\varphi_{w,\sigma}:\mathbb D(\sigma)\to \mathbb D(\sigma)$ such that $w\in \mathbb D(s_0)$ and $f(w)\not\in \supp Z$, we get (cf. \eqref{eqn:20211216})
$$
m(\sigma,f\circ \varphi_{w,\sigma},\lambda_Z)\leq cT_{0}(\sigma,f\circ \varphi_{w,\sigma},\omega_V)+ \lambda_Z(f(w)).
$$
By Lemma \ref{lem:20210711}, we have
$$
m(\sigma,f,\lambda_Z)\leq \frac{\sigma+s_0}{\sigma-s_0}m(\sigma,f\circ \varphi_{w,\sigma},\lambda_Z)\leq \frac{2}{\sigma-s_0}m(\sigma,f\circ \varphi_{w,\sigma},\lambda_Z).
$$
Hence we get
$$
m(\sigma,f,\lambda_Z)\leq \frac{2c}{\sigma-s_0}T_{0}(\sigma,f\circ \varphi_{w,\sigma},\omega_V)+\frac{2}{\sigma-s_0} \lambda_Z(f(w)),
$$
so that
\begin{equation}\label{eqn:202205251}
m(r,f,\lambda_Z)\leq cT_{\sigma}(r,f,\omega_V)+\frac{2c}{\sigma-s_0}T_{0}(\sigma,f\circ \varphi_{w,\sigma},\omega_V)+\frac{2}{\sigma-s_0} \lambda_Z(f(w)).
\end{equation}

Now we estimate the right hand side of this estimate.
We apply \eqref{eqn:bce}.
We choose $w\in\mathbb D(s_0)$ such that $f(w)\not\in U_0$ and
$$
\int_{\mathbb C}\log \frac{1}{|z-w|}d\mu_z\leq \tau_{\gamma_0}\int_{\mathbb D(\sigma)}f^*\omega_V,
$$
where $\mu=\mathbb{I}_{\mathbb D(\sigma)}f^*\omega_V$.
We set $\eta=\sup_{x\in V-U_0}\lambda_Z(x)$.
Then we have
\begin{equation}\label{eqn:202205252}
\lambda_Z(f(w))\leq \eta.
\end{equation}
We have
\begin{equation*}
\begin{split}
\int_0^{\sigma}\frac{dt}{t}\int_{\mathbb D(t)}(f\circ \varphi_{w,\sigma})^*\omega_V
&=\int_{\mathbb C}\log \frac{\sigma}{|\xi |}d(\varphi_{w,\sigma}^*\mu)_{\xi}=\int_{\mathbb C}\log\frac{\sigma}{|\varphi_{w,\sigma}^{-1}(z)|}d\mu_z
\\
&=\int_{\mathbb C}\log\frac{|\sigma-(\overline{w}/\sigma)z|}{|z -w|}d\mu_{z}\leq (\tau_{\gamma_0} +\log 2)\int_{\mathbb D(\sigma)}f^*\omega_V.
\end{split}
\end{equation*}
Hence 
\begin{equation}
T_{0}(\sigma,f\circ \varphi_{w,\sigma},\omega_V)\leq \frac{\tau_{\gamma_0} +\log 2}{r-s}T_s(r,f,\omega_V)
\end{equation}\label{eqn:202205253}
for $\sigma<s<r<1$.
We have
\begin{equation}\label{eqn:202205254}
T_{\sigma}(r,f,\omega_V)\leq \frac{1}{\sigma(r-s)}T_s(r,f,\omega_V)
\end{equation}
for $\sigma<s<r<1$.
We substitute \eqref{eqn:202205252}-\eqref{eqn:202205254} into \eqref{eqn:202205251} to conclude the proof.
Here we set
$$
\alpha=\max\left\{
\frac{c}{\sigma}+\frac{2c(\tau_{\gamma_0} +\log 2)}{\sigma-s_0},\frac{2\eta}{\sigma-s_0}
\right\},
$$
which is independent of the choice of $s\in (\sigma,1)$, $r\in(s,1)$ and $f\in \mathcal{F}$.
\hspace{\fill} $\square$

\section{Main proposition for the proof of Proposition \ref{pro:320}}

We recall the notation $f_{S_{k,A}}$ from \eqref{eqn:20220206}.
The following proposition plays an important role in the proof of Proposition \ref{pro:320}.

\begin{pro}\label{pro:mpro}
Let $\overline{A}$ be a smooth projective equivariant compactification of a semi-abelian variety $A$.
Let $S$ be a smooth projective variety.
Let $\tau:S_{1,A}\to S$ be the natural projection.
Let $Z\subset S_{1,A}$ be an irreducible Zariski closed set.
Assume that $Z$ is horizontally integrable.
Then there exists a proper Zariski closed set $W\subsetneqq \tau(Z)$ 
with the following property:
Let $\mathcal{F}\subset \mathrm{Hol}(\mathbb D,A\times S)$ be an infinite set of non-constant holomorphic maps such that $( f_{S_{1,A}})_{f\in \mathcal{F}}\Rightarrow Z$.
Let $\omega_{\overline{A}}$ and $\omega_S$ be smooth positive $(1,1)$ forms on $\overline{A}$ and $S$, respectively. 
Let $\lambda_W\geq 0$ be a Weil function for $W$.
Let $s\in (1/2,1)$ and $\delta>0$.
Then there exist positive constants $c_1>0$, $c_2>0$, $c_3>0$ such that for all $f\in \mathcal{F}$ with $f_S(\mathbb D)\not\subset W$, we have
\begin{equation}\label{eqn:20230208}
T_s(r,f_{A},\omega_{\overline{A}}) \leq 
c_1T_s(r,f_{S},\omega_{S})
+
c_2
m((s+r)/2,f_S,\lambda_{W})
+c_3
\end{equation}
for all $r\in (s,1 )$ outside some exceptional set $E\subset (s,1)$ with the linear measure $|E|<\delta$.
\end{pro}

The rest of this section is devoted to the proof of Proposition \ref{pro:mpro}.

\subsection{Outline of the proof of Proposition \ref{pro:mpro}}
We briefly describe an outline of the proof of Proposition \ref{pro:mpro}.
We recall that a semi-abelian variety $A$ is an algebraic group with an expression 
\begin{equation}\label{eqn:20230227}
0\to T\to A\overset{\rho}{\to} A_0\to 0,
\end{equation}
where $A_0$ is an abelian variety and $T\simeq \mathbb G_m^l$ is an algebraic torus (cf. Appendix \ref{sec:a}).

The starting point of the proof is the following estimate for the left hand side of \eqref{eqn:20230208}.

\begin{lem}\label{lem:202302022}
Let $\overline{A}$ be a smooth projective equivariant compactification of a semi-abelian variety $A$.
Let $\rho:A\to A_0$ be the canonical quotient as in \eqref{eqn:20230227}.
Let $\omega_{\overline{A}}$ be a smooth positive $(1,1)$-form on $\overline{A}$ and $\omega_{A_0}$ be an invariant positive $(1,1)$-form on $A_0$.
For an irreducible component $D$ of $\partial A$, let $\lambda_D\geq 0$ be a Weil function. 
Then there exists a positive constant $c>0$ such that for all $g\in\mathrm{Hol}(\mathbb D,A)$, we have
\begin{equation}\label{eqn:202302081}
T_s(r,g,\omega_{\overline{A}})\leq cT_s(r,\rho\circ g,\omega_{A_0})+c\sum_{D\subset \partial A}|m(r,g,\lambda_D)-m((s+r)/2,g,\lambda_D)|+c
\end{equation}
for all $s\in (1/2,1)$ and $r\in (s,1)$, where $D\subset \partial A$ runs over all irreducible components of $\partial A$.
\end{lem}

{\it Proof.}\
We first note that
\begin{equation}\label{eqn:20230226}
T_s(r,g,\omega_{\overline{A}})\leq 4T_{\sigma}(r,g,\omega_{\overline{A}})
\end{equation}
for all $1/2<s<r<1$, where $\sigma=(s+r)/2$.

Let $\overline{A}=(\overline{T}\times A)/T$ (cf. Lemma \ref{lem:20200906}).
Then $\overline{T}\subset \overline{A}$ is projective.
Since the Picard group of $T$ is trivial, we may take a (not necessarily effective) very ample divisor $D_0$ on $\overline{T}$ such that $\supp D_0\subset \partial T$.
Note that $D_0$ is $T$-invariant.
Set $D=(D_0\times A)/T$.
Then $D$ is a Cartier divisor on $\overline{A}$.
We take a Zariski open covering $\{ U_i\}$ of $A_0$ such that $\bar{\rho}^{-1}(U_i)=U_i\times \overline{T}$, where $\bar{\rho}:\overline{A}\to A_0$ is the extension of $\rho$.
Then we have $D\cap \bar{\rho}^{-1}(U_i)=D_0\times U_i$.
Hence $\mathcal{O}_{\overline{A}}(D)$ is very ample for $\bar{\rho}$ in the sense of \cite[Def. 13.52]{Gortz}.
Since $A_0$ is projective, we may take an ample line bundle $L$ on $A_0$.
Then by \cite[Prop. 13.65]{Gortz}, there exists a positive integer $n\geq 1$ such that $\mathcal{O}_{\overline{A}}(D)\otimes \bar{\rho}^*L^{\otimes k}$ is ample on $\overline{A}$.
We write as $D=D^+-D^-$ by the positive and negative parts of $D$.
By \eqref{eqn:20211026}, there exist positive constants $\alpha>0$ and $\alpha'>0$ such that for all $g\in \mathrm{Hol}(\mathbb D,A)$, we have
\begin{equation}\label{eqn:202110261}
T_s(r,g,\omega_{\overline{A}})\leq \alpha T_s(r,\rho\circ g,\omega_{A_0})+\alpha' T_{\sigma}(r,g,D^+,\lambda_{D^+})-\alpha' T_{\sigma}(r,g,D^-,\lambda_{D^-})+\alpha
\end{equation}
for all $1/2<s<r<1$ (cf. \eqref{eqn:20230226}).
For $g\in \mathrm{Hol}(\mathbb D,A)$, we have $N_{\sigma}(r,g,D^+)=0$,
hence 
$$T_{\sigma}(r,g,D^+,\lambda_{D^+})=m(r,g,\lambda_{D^+})-m(\sigma,g,\lambda_{D^+}).$$
Similarly, we have $T_{\sigma}(r,g,D^-,\lambda_{D^-})=m(r,g,\lambda_{D^-})-m(\sigma,g,\lambda_{D^-})$.
Hence
$$
T_{\sigma}(r,g,D^+,\lambda_{D^+})-T_{\sigma}(r,g,D^-,\lambda_{D^-})
\leq \sum_{D\subset \partial A}|m(r,g,\lambda_D)-m(\sigma,g,\lambda_D)|,
$$
where $D\subset \partial A$ runs over all irreducible components of $\partial A$.
Combining this estimate with \eqref {eqn:202110261}, we conclude the proof of the lemma.
Here we set $c=\max\{ \alpha,\alpha'\}$.
\hspace{\fill} $\square$

\medskip

The proof of Proposition \ref{pro:mpro} roughly goes as follows.
Since $Z\subset S_{1,A}$ is horizontally integrable, we may take $U\hookrightarrow A\times S$ as in Definition \ref{defn:20201122}.
Then as in \eqref{eqn:20220708}, we get an isomorphism $\Phi:A\times U\to A\times U$.
For each $f\in\mathcal{F}$, we may take a lift $\hat{f}
\in\mathrm{Hol_m}(\mathbb D,A\times \Sigma)$ of $f$ so that $\{(\Phi\circ \hat{f})_{U_{1,A}}\}_{f\in \mathcal{F}}$ becomes very close to $U_{1,A}^*$, thanks to Lemma \ref{lem:20220707}.
See Subsection \ref{subsec:202302221} for the detail of the construction of $\hat{f}$.
Then the derivative of $(\Phi\circ \hat{f})_A$ is bounded by the derivative of $(\Phi\circ \hat{f})_U$.
Thus the estimate \eqref{eqn:20230208} for $\Phi\circ \hat{f}$, instead of $f$, is easier to prove.

However, there are several problems to work out this argument rigorously.
\begin{itemize}
\item
$\{(\Phi\circ \hat{f})_{U_{1,A}}\}_{f\in \mathcal{F}}$ becomes very close to $U_{1,A}^*$ only on some brunch $\Omega_f'\to \Omega_f$ over a domain $\Omega_f\subset \mathbb D$. 
The situation will be explained in Definition \ref{defn:20230124} below.
\item We need to compactify $U$.
Moreover, to complete the argument rigorously, we need to take the compactification $\Sigma$ carefully.
This is the theme of Subsection \ref{subsec:goodcom}.
\item
After the compactification of $U$, $\Phi$ only extends to a rational map $\overline{A}\times \Sigma\dashrightarrow \overline{A}\times \Sigma$.
Hence $(\Phi\circ \hat{f})_A$ may hit the boundary of $\overline{A}$.
\item
We estimate the first and second terms of the right hand side of \eqref{eqn:202302081} separably.
See Lemmas \ref{lem:202302193} and \ref{lem:202302203}.
The first term is easier to estimate.
The reason is that since $(\rho\circ f_A)^*\omega_{A_0}$ is subharmonic, we could estimate $T_s(r,\rho\circ f_A,\omega_{A_0})$ directly from the information of $\rho\circ f_A$ over $\Omega_f$, based on \cite[Lemma 7.1]{Y}.
The main technique to estimate the second term will be discussed in Subsections \ref{subsec:20230222}--\ref{subsec:area}.
\end{itemize}

Now we introduce the following definition.

\begin{defn}\label{defn:20230124}
Let $V$ be a variety.
Let $\mathcal{G}=((g_i,\Omega_i,\Omega'_i))_{i\in I}$ be an infinite indexed family
of triples $(g_i,\Omega_i,\Omega'_i)$, where $g_i\in \mathrm{Hol_m}(\mathbb D,V)$, $\Omega_i\subset \mathbb D$ is a connected open subset, and $\Omega'_i\subset Y_{g_i}$ is a connected component of $\pi_{g_i}^{-1}(\Omega_i)$.
Let $T\subset V$ be a Zariski closed set.
We write $\mathcal{G} \rightsquigarrow T$ if the following holds:
For every $0<s<1$, $\gamma >0$, and an open subset $U\subset V$ such that $T\subset U$, there exists a finite subset $E\subset I$ such that for every $i\in I-E$, we have
\begin{itemize}
\item $z\in \Omega_i$ for $\gamma$-almost all $z\in\mathbb D(s)$,
\item $g_i(\Omega_i')\subset U$.
\end{itemize}
\end{defn}

\medskip

We prepare some notations in Nevanlinna theory in the context of the definition above.
Let $\Sigma$ be a projective variety and let $Z\subset \Sigma$ be a closed subscheme with a Weil function $\lambda_Z$.
Let $(g,\Omega,\Omega')$ be a triple as in Definition \ref{defn:20230124}, i.e., $g\in \mathrm{Hol_m}(\mathbb D,\Sigma)$, $\Omega\subset \mathbb D$ is a connected open subset, and $\Omega'\subset Y_{g}$ is a connected component of $\pi_{g}^{-1}(\Omega)$.
Assume $g(Y_g)\not\subset \mathrm{supp}Z$.
We set
\begin{equation}\label{eqn:202302247}
m(r,(g,\Omega,\Omega'),\lambda_Z)=\frac{1}{\deg (\pi_g|_{\Omega'})}\int_{y\in \Omega'\cap \partial Y_g(r)}\lambda_Z(g(y))\frac{d\arg \pi_g(y)}{2\pi},
\end{equation}
where $\deg (\pi_g|_{\Omega'})$ is the degree of the restriction $\pi_g|_{\Omega'}:\Omega'\to\Omega$.

Suppose $\Sigma$ is smooth.
Let $\omega$ be a smooth $(1,1)$-form on $\Sigma$.
We set
\begin{equation}\label{eqn:20230215}
T_s(r,(g,\Omega,\Omega'),\omega)=\frac{1}{\deg (\pi_g|_{\Omega'})}\int_s^r\frac{dt}{t}\int_{\Omega'\cap Y_g(t)}g^*\omega.
\end{equation}

In the rest of this section, $A$ is a semi-abelian variety, $S$ is a smooth projective variety, and $Z\subset S_{1,A}$ is an irreducible Zariski closed set which is horizontally integrable.

\subsection{Good liftings of $f\in\mathcal{F}$}\label{subsec:202302221}
Since $Z\subset S_{1,A}$ is horizontally integrable, we take $\iota:U\hookrightarrow A\times S$ as in Definition \ref{defn:20201122}.
We get the isomorphism 
$$
\Phi:A\times U\to A\times U
$$
from \eqref{eqn:20220708}.
We apply Lemma \ref{lem:20220707} to get an irreducible component $\Xi\subset U_{1,A}$ of $(q')^{-1}(Z)$ such that
\begin{equation}\label{eqn:202207081?}
D\Phi(\Xi)\subset U_{1,A}^*,
\end{equation}
where $q':U_{1,A}\to S_{1,A}$ is induced from $q:U\to S$.
Let $\Theta\subset U$ be the image of $\Xi$ under $U_{1,A}\to U$.
Since $U_{1,A}\to U$ is proper, $\Theta$ is a Zariski closed subset of $U$.
Let $\Sigma$ be a smooth compactification of $U$ such that $q:U\to S$ extends to a morphism $\bar{q}:\Sigma\to S$.
We denote by $\overline{\Theta}\subset \Sigma$ the Zariski closure of $\Theta$ in $\Sigma$.

\begin{rem}\label{rem:202205294}
(1)
Let $g\in\mathrm{Hol_m}(\mathbb D,A\times\Sigma)$ such that $g_{\Sigma}(Y_g)\not\subset \partial U$.
Then we may define 
$$\Phi\circ g:Y_{g}\to \overline{A}\times  \Sigma$$ 
as follows.
Set $Q_{g}=(g_{\Sigma})^{-1}(\partial U)$.
Then $Q_{g}\subset Y_{g}$ is a discrete subset.
Thus we get 
\begin{equation*}
\Phi\circ (g|_{Y_{g}\backslash Q_g}):Y_{g}\backslash Q_g\to A\times U\subset \overline{A}\times \Sigma.
\end{equation*}
Since $\overline{A}\times \Sigma$ is compact and $\Phi:\overline{A}\times\Sigma\dashrightarrow \overline{A}\times\Sigma$ is rational, $\Phi\circ (g|_{Y_{\hat{f}}\backslash Q_f})$ extends to a holomorphic map $\Phi\circ g:Y_{g}\to \overline{A}\times  \Sigma$.
By the construction, we have $(\Phi\circ g)_{\overline{A}}(Y_g)\not\subset \partial A$.

(2)
Let $h\in\mathrm{Hol_m}(\mathbb D,\overline{A}\times\Sigma)$ be non-constant such that $h_{\overline{A}}(Y_h)\not\subset \partial A$.
Then we may define $h_{\Sigma_{1,A}}:Y_{h}\to \Sigma_{1,A}$ as follows.
Set $R_{h}=(h_{\overline{A}})^{-1}(\partial A)$.
Then $R_{h}\subset Y_{h}$ is a discrete subset.
We get $h|_{Y_{h}\backslash R_h}:Y_{h}\backslash R_h\to A\times \Sigma$.
Since $h$ is non-constant, this induces $(h|_{Y_{h}\backslash Q_{h}})_{\Sigma_{1,A}}:Y_{h}\backslash Q_{h}\to \Sigma_{1,A}$.
Since $\Sigma_{1,A}$ is compact, we get $h_{\Sigma_{1,A}}:Y_{h}\to \Sigma_{1,A}$.
\end{rem}

With these notations, we state the following lemma.

\begin{lem}\label{lem:20230218}
Let $\mathcal{F}\subset \mathrm{Hol}(\mathbb D,A\times S)$ be an infinite set of non-constant holomorphic maps such that $( f_{S_{1,A}})_{f\in \mathcal{F}}\Rightarrow Z$.
Then there exists a finite subset $\mathcal{E}\subset \mathcal{F}$ with the following property:
For each $f\in\mathcal{F}\backslash\mathcal{E}$, there exist 
\begin{itemize}
\item
a lifting $\hat{f}\in \mathrm{Hol_m}(\mathbb D,A\times \Sigma)$ of $f\in \mathrm{Hol}(\mathbb D,A\times S)$,
\item
a connected open subset $\Omega_f\subset \mathbb D$,
\item
a connected component $\Omega_f'$ of $\pi_{\hat{f}}^{-1}(\Omega_f)\subset Y_{\hat{f}}$
\end{itemize}
such that 
\begin{equation}\label{eqn:20210504?}
\hat{f}_{\Sigma}(Y_{\hat{f}})\not\subset \partial U,
\end{equation}
\begin{equation}\label{eqn:20211024?}
( (\hat{f}_{\Sigma},\Omega_f,\Omega_f'))_{f\in \mathcal{F}\backslash\mathcal{E}} \rightsquigarrow\overline{\Theta},
\end{equation}
\begin{equation}\label{eqn:202102261?}
( ((\Phi\circ \hat{f})_{\Sigma_{1,A}},\Omega_f,\Omega_f'))_{f\in \mathcal{F}\backslash\mathcal{E}} \rightsquigarrow \Sigma_{1,A}^*
\end{equation}
and
\begin{equation}\label{eqn:202302034?}
\deg \pi_{\hat{f}}\leq [\mathbb C(S):\mathbb C(\Sigma)].
\end{equation}
\end{lem}

Before going to prove this lemma, we need some preparations.

\subsubsection{Preliminary lemmas for the proof of Lemma \ref{lem:20230218}}

\begin{lem}\label{lem:20220529}
Let $s\in (0,1)$ and $\gamma>0$.
Let $V\subset \mathbb D(s)$ be an open set such that $z\in V$ for $\gamma$-almost all $z\in\mathbb D(s)$.
Let $s'\in (0,s)$.
Then there exists an open subset $\Omega\subset V$ such that $\Omega$ is connected and $z\in \Omega$ for $2\gamma$-almost all $z\in\mathbb D(s')$.
\end{lem}

{\it Proof.}\
Set $K=\overline{\mathbb D(s')}\backslash V$.
Then $K$ is compact.
We first show that there exists a finite correction of closed discs $D_1,\ldots,D_n$ such that 
\begin{enumerate}
\item
$K\subset \cup_{i=1}^nD_i$, and 
\item
$\sum_{i=1}^nr_i<2\gamma$, where $r_i$ is the radius of $D_i$.
\end{enumerate}
Indeed, let $\{ E_i\}$ be a countable set of closed discs such that $\mathbb D(s)\backslash V\subset \cup E_i$ and the sum of the radii of $E_i$ is less than $\gamma$.
Let $O_i$ be the open disc whose center is equal to that of $E_i$ and radius is equal to the double of that of $E_i$.
We have $K\subset \cup O_i$.
Hence we may take $O_1,O_2,\ldots,O_n$ such that $K\subset O_1\cup \cdots\cup O_n$ and the sum of radii of $O_i$ is less than $2\gamma$.
Thus the closed discs $\overline{O_1},\ldots,\overline{O_n}$ satisfy our requirements.

Now let $n$ be the minimum such that there exist closed discs $D_1,\ldots,D_n$ with the properties (1) and (2) above.
We claim that $D_i\cap D_j=\emptyset$ for $i\not=j$.
To prove this, we suppose contrary that there exists $i\not=j$ such that $D_i\cap D_j\not=\emptyset$.
We may assume without loss of generality that $i=1$ and $j=2$.
Let $p_1,p_2$ be the centers of $D_1,D_2$, respectively.
Let $p$ be the point which divides the line segment $\overline{p_1p_2}$ internally in the ratio $\overline{p_1p}:\overline{pp_2}=r_2:r_1$.
Let $D$ be the closed disc whose center is $p$ and radius is equal to $r_1+r_2$.
Then $D_1\cup D_2\subset D$.
Hence the closed discs $D,D_3,\ldots,D_n$ satisfy the property (1) and (2) above.
This contradicts to the choice of $n$.
Hence we have proved $D_i\cap D_j=\emptyset$ for $i\not=j$.
Now we set $\Omega=\mathbb D(s')\backslash \cup_{i=1}^n D_i$.
Then $\Omega$ is open and connected.
We have $z\in \Omega$ for $2\gamma$-almost all $z\in\mathbb D(s')$.
\hspace{\fill} $\square$

\begin{lem}\label{lem:20210502}
Let $\tilde{V}$ and $V$ be smooth projective varieties.
Let $p:\tilde{V} \to V$ be a proper, surjective, generically finite morphism.
Let $\tilde{V}_0\subset \tilde{V}$ be a nonempty Zariski open set such that $p:\tilde{V} \to V$ is quasi-finite on $\tilde{V}_0$.
Let $Z\subset V$ be an irreducibe Zariski closed set.
Let $\tilde{Z}\subset \tilde{V}$ be an irreducible component of $p^{-1}(Z)$ such that $\tilde{Z}\cap \tilde{V}_0\not= \emptyset$.
Let $\mathcal{G}=(g_i)_{i\in I}$ be an infinite indexed family in $\mathrm{Hol}(\mathbb D,V)$ such that $\mathcal{G}\Rightarrow Z$.
Then for all but finitely many $i\in I$, there exist
\begin{itemize}
\item
a lift $\tilde{g}_i\in \mathrm{Hol_m}(\mathbb D,\tilde{V})$ of $g_i$, 
\item
a connected open subset $\Omega_i\subset \mathbb D$, and
\item
a connected component $\Omega_i'$ of $\pi_{\tilde{g}_i}^{-1}(\Omega_i)\subset Y_{\tilde{g}_i}$
\end{itemize}
such that 
\begin{itemize}
\item
$((\tilde{g}_i,\Omega_i,\Omega_i'))_{i\in I}\rightsquigarrow \tilde{Z}$,\item
$\deg\pi_{\tilde{g}_i}\leq [\mathbb C(V):\mathbb C(\tilde{V})]$, and
\item
$\tilde{g}_i(Y_{\tilde{g}_i})\not\subset \tilde{V}\backslash \tilde{V}_0$.
\end{itemize}
\end{lem}

{\it Proof.}\
We construct a $Z$-admissible modification $V'\to V$ and $\tilde{Z}$-admissible modification $\tilde{V}'\to \tilde{V}$ with the following properties (cf. Definition \ref{defn:20211007}):
\begin{itemize}
\item
$p$ induces a morphism $p':\tilde{V}'\to V'$ with the following commutative diagram
\begin{equation*}
\begin{CD}
\tilde{V}@<<<\tilde{V}'\\
@VpVV @VVp'V   \\
V@<<< V'
\end{CD}
\end{equation*}
\item
$p'$ is quasi-finite on some neighbourhood of $\tilde{Z}'$, where $\tilde{Z}'\subset \tilde{V}'$ is the minimal transform.
\item
$(p')^{-1}(Z')=\tilde{Z}'\amalg E$, where $Z'\subset V'$ is the minimal transform.
\end{itemize}
We construct these objects as follows.
We decompose $p^{-1}(Z)$ into irreducible components as $p^{-1}(Z)=\tilde{Z}\cup F_1\cup \cdots \cup F_k$ and set $F=F_1\cup \cdots \cup F_k$.
Replacing $\tilde{V}$ by $\mathrm{Bl}_{\tilde{Z}\cap F}\tilde{V}$, we may assume $\tilde{Z}\cap F_j=\emptyset$ if $p(F_j)=Z$ (cf. the proof of \eqref{eqn:20220522}).
By Lemma \ref{lem:202110273}, we may take a $Z$-admissible modification $V'\to V$ such that $\tilde{V}'|_{Z'}\to Z'$ is flat.
We may assume moreover that $V'$ is smooth.
Then $\tilde{V}'|_{Z'}\to Z'$ is an open map.
Hence every irreducible component of $(p')^{-1}(Z')$ dominates $Z'$.
Thus $\tilde{Z}'$ is a connected component of $(p')^{-1}(Z')$.
Hence $\tilde{Z}'\to Z'$ is flat, so finite.
Hence $p'$ is quasi-finite on some neighbourhood of $\tilde{Z}'$, and $(p')^{-1}(Z')=\tilde{Z}'\amalg E$.

Let $T\subset \tilde{V}'$ be a Zariski open neighbourhood of $\tilde{Z}'$ such that $p'|_T:T\to V'$ is quasi-finite and $T\cap E=\emptyset$.
Then, by Zariski's main theorem, there exist a compactification $T\subset \overline{T}$ and a finite map 
$$p'':\overline{T}\to V'$$
such that $p''|_T=p'|_T$.
We have $(p'')^{-1}(Z')=\tilde{Z}'\amalg E'$.
Let $\tilde{Z}'\subset O_{\tilde{Z}'}$ and $E'\subset O_{E'}$ be open neighbourhoods in $\overline{T}$ such that $O_{\tilde{Z}'}\cap O_{E'}=\emptyset$.
We may assume $O_{\tilde{Z}'}\subset T$.
Note that $p''\left(\overline{T} -(O_{\tilde{Z}'}\cup O_{E'})\right)\subset V'$ is compact and disjoint from $Z'$.
Hence we may take an open neighbourhood $Z'\subset O_{Z'}$ which is disjoint from $p''\left(\overline{T}-(O_{\tilde{Z}'}\cup O_{E})\right)$.
Then $(p'')^{-1}(O_{Z'})\subset O_{\tilde{Z}'}\cup O_{E'}$.
We replace $O_{\tilde{Z}'}$ by $(p'')^{-1}(O_{Z'})\cap O_{\tilde{Z}'}$.
Then $p''|_{O_{\tilde{Z}'}}:O_{\tilde{Z}'}\to O_{Z'}$ is a proper map.
Since $V'$ is smooth, by Remmert open mapping theorem, the finite map $p'':\overline{T}\to V'$ is an open map.
Hence $p''|_{O_{\tilde{Z}'}}:O_{\tilde{Z}'}\to O_{Z'}$ is proper and open.

By $\mathcal{G}\Rightarrow Z$, we have $\mathcal{G}\to Z'$.
Here by removing finite elements from $I$, we assume that every $g_i:\mathbb D\to V$ has a unique lift $\mathbb D\to V'$ (cf. Remark \ref{rem:202205291}), which we continue to use the same notation $g_i$.
Let $O_{Z'}\Supset U_1\Supset U_2\Supset\cdots$ be a sequence of open neighbourhoods of $Z'$ such that every open neighbourhood of $Z'$ in $O_{Z'}$ contains some $U_n$ for sufficiently large $n$.
Set $\gamma_n=\frac{1}{2^{n+1}}$ and 
$$
I_n=\{i\in I;\ \text{$g_i(z)\in U_n$ for $\gamma_n$-almost all $z\in \mathbb D(1-\gamma_n)$}\}.
$$
Then we get a descending sequence $I\supset I_1\supset I_2\supset \cdots$.
By $\mathcal{G}\to Z'$, $I\backslash I_n$ is a finite set for every $n$.
We set $I_{\infty}=\cap_nI_n$.

For $i\in I_1$, we are going to construct a lift $\tilde{g}_i\in \mathrm{Hol_m}(\mathbb D,\overline{T})$, a connected open set $\Omega_i\subset \mathbb D$ and a connected component $\Omega_i'\subset Y_{\tilde{g}_i}$ of $\pi_{\tilde{g}_i}^{-1}(\Omega_i)$. 
We first consider the case $i\not\in I_{\infty}$.
We fix $i\in I_1\backslash I_{\infty}$.
We take $n_i$ such that $i\in I_{n_i}$, but $i\not\in I_{n_i+1}$.
By Lemma \ref{lem:20220529}, we may take an open set $\Omega_i\subset g_i^{-1}(U_{n_i})\cap\mathbb D(1-2\gamma_{n_i})$ such that $\Omega_i$ is connected and $z\in \Omega_i$ for $2\gamma_{n_i}$-almost all $z\in \mathbb D(1-2\gamma_{n_i})$.
We take $h:\mathcal{Y}\to \overline{T}$ by pull-back, where $\mathcal{Y}$ is an one dimensional analytic space with finite map $q:\mathcal{Y}\to \mathbb D$.
\begin{equation*}
\begin{CD}
\mathcal{Y}@>h>>\overline{T}\\
@VqVV @VVp''V   \\
\mathbb D@>>g_i> V'
\end{CD}
\end{equation*}
Let $\mathcal{Y}_1,\ldots ,\mathcal{Y}_l$ be irreducible components of $\mathcal{Y}$.
We may take $\mathcal{Y}_j$, where $1\leq j\leq l$, such that $h(\mathcal{Y}_j)\cap O_{\tilde{Z}'}\not=\emptyset$.
Then $h|_{\mathcal{Y}_j}:\mathcal{Y}_j\to \overline{T}$ defines a lift $\tilde{g}_i\in \mathrm{Hol_m}(\mathbb D,\overline{T})$ of $g_i$ such that $\tilde{g}_i(Y_{\tilde{g}_i})\cap O_{\tilde{Z}'}\not=\emptyset$.
We take a connected component $\Omega_i'$ of $\pi_{\tilde{g}_i}^{-1}(\Omega_i)\subset Y_{\tilde{g}_i}$ such that $\tilde{g}_i(\Omega_i')\cap O_{\tilde{Z}'}\not=\emptyset$.
By $g_i(\Omega_i)\subset U_{n_i}\subset O_{Z'}$, we have $\tilde{g}_i(\Omega_i')\subset O_{\tilde{Z}'}\cup O_{E}$.
Hence, since $\Omega_i'$ is connected, we have 
\begin{equation}\label{eqn:202301221}
\tilde{g}_i(\Omega_i')\subset O_{\tilde{Z}'}
\end{equation}
When $i\in I_{\infty}$, we have $g_i(\mathbb D)\subset Z'$.
Since $\tilde{Z}'\to Z'$ is a finite map, we have a lift $\tilde{g}_i:Y_{\tilde{g}_i}\to \tilde{Z}'$ of $g_i$.
We set $\Omega_i=\mathbb D$ and $\Omega_i'=Y_{\tilde{g}_i}$.
Thus we have constructed $\tilde{g}_i\in \mathrm{Hol_m}(\mathbb D,\overline{T})$, $\Omega_i$ and $\Omega_i'$ for all $i\in I_1$.

We claim that, for all $i\in I_n$, we have $\tilde{g}_i(\Omega_i')\subset (p'')^{-1}(U_n)\cap O_{\tilde{Z}'}$.
This is obvious for $i\in I_{\infty}$ by $\tilde{g}_i(Y_{\tilde{g}_i})\subset \tilde{Z}'$.
When $i\not\in I_{\infty}$, we have $n_i\geq n$, hence $g_i(\Omega_i)\subset U_{n_i}\subset U_n$.
Thus $\tilde{g}_i(\Omega_i')\subset (p'')^{-1}(U_n)\cap O_{\tilde{Z}'}$ (cf. \eqref{eqn:202301221}).
This proves our claim.

Now we prove $((\tilde{g}_i,\Omega_i,\Omega_i'))_{i\in I_1}\rightsquigarrow \tilde{Z}'$.
Let $M\subset \overline{T}$ be an open neighbourhood of $\tilde{Z}'$.
Let $s\in (0,1)$ and $\delta>0$.
We take a sufficiently large $n$ such that $s<1-2\gamma_n$ and $2\gamma_n<\delta$.
Note that $O_{\tilde{Z}'}\backslash M$ is a closed subset of $O_{\tilde{Z}'}$.
By $(p'')^{-1}(Z')=\Xi'\amalg E'$, we have $(p'')^{-1}(Z')\cap O_{\tilde{Z}'}=\tilde{Z}'$.
Hence $Z'\cap p''(O_{\tilde{Z}'}\backslash M)=\emptyset$.
Since $p''|_{O_{\tilde{Z}'}}:O_{\tilde{Z}'}\to O_{Z'}$ is proper, $p''(O_{\tilde{Z}'}\backslash M)$ is a closed subset of $O_{Z'}$.
Hence there exists $n'$ such that $U_{n'}\cap p''(O_{\tilde{Z}'}\backslash M)=\emptyset$.
Hence we have $(p'')^{-1}(U_{n'})\cap O_{\tilde{Z}'}\subset M$.
We set $n''=\max\{ n,n'\}$.
Then for all $i\in I_{n''}$, we have $\tilde{g}_i(\Omega_i')\subset M$ and $z\in \Omega_i$ for $\delta$-almost all $z\in \mathbb D(s)$.
Hence $((\tilde{g}_i,\Omega_i,\Omega_i'))_{i\in I_1}\rightsquigarrow \tilde{Z}'$.

So far, we have considered that $\tilde{g}_i$ are maps into $\overline{T}$.
Since the birational map $\overline{T}\dashrightarrow \tilde{V}'$ is an isomorphism on $T$ and $O_{\tilde{Z}'}\subset T$, we have $((\tilde{g}_i,\Omega_i,\Omega_i'))_{i\in I_1}\rightsquigarrow \tilde{Z}'$ under the consideration $\tilde{g}_i\in \mathrm{Hol_m}(\mathbb D,\tilde{V}')$.
Thus we have $((\tilde{g}_i,\Omega_i,\Omega_i'))_{i\in I_1}\rightsquigarrow \tilde{Z}$ under the consideration $\tilde{g}_i\in \mathrm{Hol_m}(\mathbb D,\tilde{V})$.
By the construction, we have $\deg \pi_{\tilde{g}_i}\leq [\mathbb C(V),\mathbb C(\tilde{V} )]$.

Now let $\tilde{V}'_0\subset \tilde{V}'$ be the inverse image of $\tilde{V}_0$ under the map $\tilde{V}'\to\tilde{V}$.
We define $J\subset I_1$ to be the set of $i\in I_1$ such that $\tilde{g}_i(Y_{\tilde{g}_i})\subset \tilde{V}'\backslash \tilde{V}_0'$ under the consideration $\tilde{g}_i\in \mathrm{Hol_m}(\mathbb D,\tilde{V}')$.
To prove that $J$ is finite, we assume contrary that $J$ is infinite.
Then there exists an irreducible component $D$ of $\tilde{V}'\backslash \tilde{V}_0'$ such that $\tilde{g}_i(Y_{\tilde{g}_i})\subset D$ for infinitely many $i\in J$.
By $((\tilde{g}_i,\Omega_i,\Omega_i'))_{i\in I_1}\rightsquigarrow \tilde{Z}'$, we have $D\cap \tilde{Z}'\not=\emptyset$.
Hence $T\cap D\not=\emptyset$.
Hence $Z'\not\subset p'(D)$ and $g_i(\mathbb D)\subset p'(D)$ for infinitely many $i\in J$.
This contradicts to the assumption $\mathcal{G}\Rightarrow Z$.\hspace{\fill} $\square$

\subsubsection{Proof of Lemma \ref{lem:20230218}}
We first fix a modification of $\Sigma_{1,A}$ as follows.
The isomorphism $D\Phi:U_{1,A}\to U_{1,A}$ induces a rational map $\Sigma_{1,A}\dashrightarrow \Sigma_{1,A}$. 
We also have a rational map $\Sigma_{1,A}\dashrightarrow S_{1,A}$.
By taking a birational modification $\widetilde{\Sigma_{1,A}}\to\Sigma_{1,A}$, we may assume that
\begin{itemize}
\item
$p:\widetilde{\Sigma_{1,A}}\to S_{1,A}$ is generically finite and surjective, and 
\item
$D\Phi :\widetilde{\Sigma_{1,A}}\to \Sigma_{1,A}$ is holomorphic.
\end{itemize}
Note that the natural map $U_{1,A}\to S_{1,A}$ is {\'e}tale and $D\Phi$ is holomorphic on $U_{1,A}$.
Hence we may assume $U_{1,A}\subset \widetilde{\Sigma_{1,A}}$.
Let $\overline{\Xi}\subset \widetilde{\Sigma_{1,A}}$ be the closure of $\Xi\subset U_{1,A}$.
Then $\overline{\Xi}$ is an irreducible component of $p^{-1}(Z)$ such that 
\begin{equation}\label{eqn:20230203200}
\overline{\Xi}\cap U_{1,A}\not= \emptyset.
\end{equation}

Now we are given an infinite subset $\mathcal{F}\subset \mathrm{Hol}(\mathbb D,A\times S)$ of non-constant holomorphic maps so that $(f_{S_{1,A}})_{f\in\mathcal{F}}\Rightarrow Z$.
By \eqref{eqn:20230203200} and $(f_{S_{1,A}})_{f\in\mathcal{F}}\Rightarrow Z$, we may apply Lemma \ref{lem:20210502} for $p:\widetilde{\Sigma_{1,A}}\to S_{1,A}$ and $(f_{S_{1,A}})_{f\in\mathcal{F}}$.
Then there exists a finite subset $\mathcal{E}\subset \mathcal{F}$ such that for each $f\in \mathcal{F}\backslash\mathcal{E}$, there exist a lift $g\in \mathrm{Hol_m}(\mathbb D,\widetilde{\Sigma_{1,A}})$ of $f_{S_{1,A}}$, a connected open subset $\Omega_f\subset \mathbb D$ and a connected component $\Omega_f'$ of $\pi_g^{-1}(\Omega_f)\subset Y_g$ such that 
\begin{equation}\label{eqn:202205301}
( (g,\Omega_f,\Omega_f'))_{f\in \mathcal{F}\backslash\mathcal{E}}\rightsquigarrow \overline{\Xi},
\end{equation}
\begin{equation}\label{eqn:202302033}
\deg \pi_{g}\leq [\mathbb C(S_{1,A}):\mathbb C(\Sigma_{1,A})],
\end{equation}
and
\begin{equation}\label{eqn:202301181}
g(Y_g)\not\subset \widetilde{\Sigma_{1,A}}\backslash U_{1,A}.\end{equation}
Let $\tau':\widetilde{\Sigma_{1,A}} \to \Sigma$ be the natural projection.
Then we have the following commutative diagram:
\begin{equation}\label{eqn:202302036}
\begin{CD}
Y_g@>g>>\widetilde{\Sigma_{1,A}} @>\tau'>>\Sigma\\
@V\pi_{g}VV @VVpV   @VV\bar{q}V\\
\mathbb D@>>f_{S_{1,A}}> S_{1,A} @>>\tau> S
\end{CD}
\end{equation}
For $f\in \mathcal{F}\backslash\mathcal{E}$, we set $\hat{f}=(f_A,\tau'\circ g)\in \mathrm{Hol_m}(\mathbb D,A\times \Sigma)$.
Then we get the following commutative diagram:
\begin{equation}\label{eqn:202205303}
\begin{CD}
Y_{\hat{f}}@>\hat{f}>>A\times \Sigma @>>>\Sigma\\
@V\pi_{\hat{f}}VV @VVV   @VV\bar{q}V \\
\mathbb D@>>f> A\times S@>>> S
\end{CD}
\end{equation}
Here $Y_{\hat{f}}=Y_g$.
Note that $\tau'(\overline{\Xi})=\overline{\Theta}$.
Hence by \eqref{eqn:202205301}, we get \eqref{eqn:20211024?}.

Next we show \eqref{eqn:20210504?} and \eqref{eqn:202302034?}.
Let $f\in\mathcal{F}\backslash\mathcal{E}$.
By \eqref{eqn:202301181}, we get \eqref{eqn:20210504?}.
We note that $[\mathbb C(S):\mathbb C(\Sigma)]=[\mathbb C(S_{1,A}):\mathbb C(\Sigma_{1,A})]$.
Hence by \eqref{eqn:202302033}, we get \eqref{eqn:202302034?}.

Finally we prove \eqref{eqn:202102261?}.
We have $p\circ g=f_{S_{1,A}}\circ\pi_g$ (cf. \eqref{eqn:202302036}) and $p\circ (\hat{f})_{\widetilde{\Sigma_{1,A}}}=f_{S_{1,A}}\circ\pi_{\hat{f}}$ (cf. \eqref{eqn:202205303}).
Hence we get
\begin{equation}\label{eqn:202302037}
p\circ g=p\circ (\hat{f})_{\widetilde{\Sigma_{1,A}}}.
\end{equation}
By the definition of $\hat{f}$, we get $\tau'\circ g=\hat{f}_{\Sigma}$.
By \eqref{eqn:202205303}, we have $\tau'\circ (\hat{f})_{\widetilde{\Sigma_{1,A}}}=\hat{f}_{\Sigma}$.
Hence we get 
\begin{equation}\label{eqn:202302038}
\tau'\circ g=\tau'\circ (\hat{f})_{\widetilde{\Sigma_{1,A}}}.
\end{equation}
Since $U_{1,A}=S_{1,A}\times_SU$, the two relations \eqref{eqn:202302037} and \eqref{eqn:202302038} yield $(\hat{f})_{\widetilde{\Sigma_{1,A}}}=g$.
Hence by \eqref{eqn:202205301}, we have
\begin{equation}\label{eqn:20201207}
(((\hat{f})_{\widetilde{\Sigma_{1,A}}},\Omega_f,\Omega_f'))_{f\in \mathcal{F}\backslash\mathcal{E}}\to \bar{\Xi}.
\end{equation}
By \eqref{eqn:202207081?}, we get
$$
D\Phi (\overline{\Xi})\subset \Sigma_{1,A}^*.
$$
Thus by \eqref{eqn:20201207}, we get $( (D\Phi \circ \hat{f}_{\widetilde{\Sigma_{1,A}}},\Omega_f,\Omega_f'))_{f\in \mathcal{F}\backslash\mathcal{E}} \to \Sigma_{1,A}^*$.
By $D\Phi \circ \hat{f}_{\widetilde{\Sigma_{1,A}}}=(\Phi\circ \hat{f})_{\Sigma_{1,A}}$, we get \eqref{eqn:202102261?}.
\hspace{\fill} $\square$

\subsection{Good compactification of $U$}
\label{subsec:goodcom}

\begin{lem}\label{lem:202302021}
Let $S$ be a smooth projective variety.
Let $q:U \to S$ be an {\'e}tale morphism from a smooth variety $U$.
Let $\Theta\subset U$ be an irreducible Zariski closed set.
Then there exist
\begin{itemize}
\item a smooth compactification $\Sigma$ of $U$ such that $q:U\to S$ extends to a morphism $\bar{q}:\Sigma\to S$,
\item a Zariski open set $\Sigma^o\subset \Sigma$ such that $\overline{\Theta}\subset \Sigma^o$, where $\overline{\Theta}\subset \Sigma$ is the Zariski closure of $\Theta$ in $\Sigma$,
\item
a smooth semi-positive $(1,1)$-form $\eta\geq 0$ on $\Sigma$ such that $\eta>0$ on $\Sigma^o$,
\item
a proper Zariski closed subset $W\subset S$ with $W\subsetneqq\bar{q}(\overline{\Theta})$
\end{itemize}
such that the following properties hold:
\begin{enumerate}
\item
Let $\Delta\subset \Sigma$ be an irreducible component of $\partial U$ such that $\Delta\cap\overline{\Theta}\not=\emptyset$.
Then $\bar{q}(\Delta)\subset W$.
\item
Let $\omega_S$ be a smooth positive $(1,1)$-form on $S$ and let $\lambda_W\geq 0$ be a Weil function for $W$.
Then there exists a positive constant $c>0$ such that for all $g\in  \mathrm{Hol_m}(\mathbb D,\Sigma)$ with $\bar{q}\circ g(Y_g)\not\subset W$, we have
$$
T_s(r,g,\eta)\leq cT_s(r,\bar{q}\circ g,\omega_{S})+cm((s+r)/2,\bar{q}\circ g,\lambda_W)+c
$$
for all $s\in (1/2,1)$ and $r\in (s,1)$.
\end{enumerate}
\end{lem}

To prove this lemma, we start from the following lemma.

\begin{lem}\label{lem:202302231}
Let $X$ be a projective variety.
Let $V$ and $V'$ be non-empty Zariski open subsets of $X$ such that $V\cap V'$ is smooth.
Then there exists a proper birational morphism $p:X'\to X$ from a projective variety $X'$ such that $p^{-1}(V)$ is smooth and $p^{-1}(V')\to V'$ is isomorphic.
\end{lem}

{\it Proof.}\
By a strong desingularization of $V$, we have a sequence
$$
V=V_0\leftarrow V_1\leftarrow V_2\leftarrow \cdots\leftarrow V_k
$$
such that 
\begin{itemize}
\item $V_k$ is smooth, and
\item
each $V_{i+1}\to V_i$ is a blowing-up $V_{i+1}=\mathrm{Bl}_{C_i}V_i$ such that the center $C_i\subset V_i$ satisfies $p_i(C_i)\subset V\backslash (V\cap V')$, where $p_i:V_i\to V$ is the natural morphism.
\end{itemize}

We inductively construct a projective variety $X_i$ with an open immersion $V_i\subset X_i$ as follows.
Set $X_0=X$.
Then $V_0\subset X_0$.
Suppose we have constructed $X_i$ with $V_i\subset X_i$.
Let $\overline{C_i}\subset X_i$ be the (schematic) closure of $C_i\subset V_i$.
We set $X_{i+1}=\mathrm{Bl}_{\overline{C_i}}X_i$.
Then the inverse image of $V_i\subset X_i$ under the projection $X_{i+1}\to X_i$ is equal to $V_{i+1}\subset X_{i+1}$.
Since $X_i$ is projective, the blowing-up $X_{i+1}$ is also projective.
Thus we have constructed projective varieties $X_i$ with open immersions $V_i\subset X_i$ for all $i=0,1,\ldots,k$.

By the construction, we have extensions $\bar{p}_i:X_i\to X$ of $p_i$.
By $p_i(C_i)\subset V\backslash(V\cap V')$, we have $\bar{p}_i(\overline{C_i})\subset X\backslash V'$.
Hence $\bar{p}_i:X_i\to X$ are isomorphisms over $V'$ for all $i$.
Moreover we have $\bar{p}_i^{-1}(V)=V_i$ for all $i$.
We set $X'=X_k$ and $p=\bar{p}_k$ to conclude the proof.
\hspace{\fill} $\square$

\medskip

{\it Proof of Lemma \ref{lem:202302021}.}\
We first construct $\Sigma$ and $\Sigma^o$.
We apply Zariski's main theorem for the quasi-finite map $q:U\to S$.
Then we get an open immersion $U\hookrightarrow \Sigma_1$ and a finite map 
$$\bar{q}_1:\Sigma_1\to S.$$
Since $S$ is projective, $\Sigma_1$ is projective.
Since $U$ is smooth, we may assume that $\Sigma_1$ is normal.
Set 
\begin{equation}\label{eqn:202302235}
\Sigma_2=\mathrm{Bl}_{\overline{\Theta}_1\cap\partial U}\Sigma_1,
\end{equation}
where $\overline{\Theta}_1\subset \Sigma_1$ is the Zariski closure of $\Theta\subset U$ in $\Sigma_1$.
Then $\Sigma_2$ is projective.
Set 
$$p_1:\Sigma_2\to\Sigma_1.
$$
Then $p_1$ is isomorphic over $U\subset \Sigma_1$.
Hence $U\subset \Sigma_2$.
Let $\overline{\Theta}_2\subset \Sigma_2$ be the Zariski closure of $\Theta\subset U$ in $\Sigma_2$.
If $\Delta\subset \Sigma_2$ is an irreducible component of $\Sigma_2\backslash U$ such that $\Delta\cap \overline{\Theta}_2\not=\emptyset$, then 
\begin{equation}\label{eqn:202302223}
p_1(\Delta)\subsetneqq \overline{\Theta}_1.
\end{equation}
Let $V\subset \Sigma_2$ be a Zariski open set defined by $V=\Sigma_2\backslash \bigcup_{\Delta'}\Delta'$, where $\Delta'$ runs over all irreducible components of $\Sigma_2\backslash U$ such that $\Delta'\cap \overline{\Theta}_2=\emptyset$.
Then we have 
\begin{equation}\label{eqn:202302233}
\overline{\Theta}_2\subset V.
\end{equation}
 We define a Zariski closed set $D\subset\Sigma_2$ by $D=\bigcup_{\Delta}\Delta$, where $\Delta$ runs over all irreducible components of $\Sigma_2\backslash U$ such that $\Delta\cap \overline{\Theta}_2\not=\emptyset$.
We have $U= V\cap (\Sigma_2\backslash D)$, where $U$ is smooth.
By Lemma \ref{lem:202302231}, there exists a proper birational modification 
$$p_2:\Sigma_3\to \Sigma_2$$ 
such that $p_2^{-1}(V)$ is smooth and $p_2$ is isomorphic over $\Sigma_2\backslash D$.
We define $\Sigma$ by a smooth modification 
$$p_3:\Sigma\to \Sigma_3$$ 
which is an isomorphism over $p_2^{-1}(V)\subset \Sigma_3$.
We set $\Sigma^o=(p_3\circ p_2)^{-1}(V)$.
By \eqref{eqn:202302233}, we have $\overline{\Theta}\subset\Sigma^o$.

We construct $W\subset S$ with $W\subsetneqq\bar{q}(\overline{\Theta})$ and prove (1).
We denote the exceptional locus of $p_1:\Sigma_2\to \Sigma_1$ by $E_1\subset \Sigma_2$.
We set 
$$
W=\bar{q}_1\circ p_1(E_1\cup D)\subset S. 
$$
By \eqref{eqn:202302235}, we have $p_1(E_1)\subsetneqq \overline{\Theta}_1$.
By \eqref{eqn:202302223}, we have $p_1(D)\subsetneqq \overline{\Theta}_1$.
Hence we have $W\subsetneqq\bar{q}(\overline{\Theta})$.
By the definition of $D\subset \Sigma_2$, we obtain the assertion (1).

Let $E\subset \Sigma_3$ be the exceptional locus of $\Sigma_3\to \Sigma_1$.
Then we claim
\begin{equation}\label{eqn:202302237}
\bar{q}_1\circ p_1\circ p_2(E)\subset W.
\end{equation}
To prove this, we denote the exceptional locus of $p_2:\Sigma_3\to \Sigma_2$ by $E_2\subset \Sigma_3$.
Then we have $E\subset p_2^{-1}(E_1)\cup E_2$.
Hence $p_2(E)\subset E_1\cup p_2(E_2)$.
Since $p_2:\Sigma_3\to\Sigma_2$ is isomorphic over $\Sigma_2\backslash D$, we have $p_2(E_2)\subset D$.
Hence $p_2(E)\subset E_1\cup D$.
Thus we have proved \eqref{eqn:202302237}.

Next we construct $\eta$ and prove the property (2).
Since $\Sigma_3$ is projective, we have an immersion $\iota:\Sigma_3\hookrightarrow \mathbb P^n$.
Let $\omega_{\mathbb P^n}$ be the Fubini-Study metric on $\mathbb P^n$.
Let $\varphi:\Sigma\to \mathbb P^n$ be the composite of $\Sigma\to\Sigma_3\to\mathbb P^n$.
We set
$$
\eta=\varphi^*\omega_{\mathbb P^n}.
$$
Then $\eta$ is semi-positive on $\Sigma$ and positive on $\Sigma^o$, for the composite of $\Sigma^o\hookrightarrow\Sigma\to \Sigma_3$ is an open immersion.

Now we prove (2).
Let $H_1,\ldots,H_l\subset \mathbb P^n$ be hyperplanes such that $H_1\cap\cdots\cap H_l=\emptyset$.
Let $h$ be the Fubini-Study metric on $\mathcal{O}_{\mathbb P^n}(1)$.
Let $\tau_i$ be a section of $\mathcal{O}_{\mathbb P^n}(1)$ such that $H_i=(\tau_i=0)$.
Set $\lambda_{H_i}(x)=-\log \sqrt{h(\tau_i(x),\tau_i(x))}\geq 0$.
By Lemma \ref{lem:20230207}, there exists a positive constant $c_i>0$ such that for all $h\in  \mathrm{Hol_m}(\mathbb D,\Sigma_3)$ with $h(Y_h)\not\subset E\cup \iota^*H_i$, we have
\begin{equation}\label{eqn:20230224}
T_{\sigma}(r,\iota\circ h,H_i,\lambda_{H_i})\leq
c_iT_{\sigma}(r,\bar{q}_1\circ p_1\circ p_2\circ h,\omega_{S})+c_im(\sigma,h,\lambda_E)+c_i
\end{equation}
for all $0<\sigma<r<1$.

Let $g\in  \mathrm{Hol_m}(\mathbb D,\Sigma)$ with $\bar{q}\circ g(Y_g)\not\subset W$.
We may take $H_i$ such that $\varphi\circ g(Y_g)\not\subset \supp H_i$, where $\varphi:\Sigma\to\mathbb P^n$ is the morphism above.
Then by \eqref{eqn:20211026}, we have
$$
T_{\sigma}(r,g,\eta)=T_{\sigma}(r,\varphi\circ g,\omega_{\mathbb P^n})= T_{\sigma}(r,\varphi\circ g,H_i,\lambda_{H_i})
$$
for all $0<\sigma<r<1$.
By \eqref{eqn:202302237} and \eqref{eqn:20230224}, we have
\begin{equation*}
T_{\sigma}(r,\varphi\circ g,H_i,\lambda_{H_i})\leq
c_iT_{\sigma}(r,\bar{q}\circ g,\omega_{S})+c_im(\sigma,\bar{q}\circ g,\lambda_W)+c_i.
\end{equation*}
We set $\displaystyle{c'=\max_{1\leq i\leq l}\{c_i\}}$.
Then we have
$$
T_{\sigma}(r,g,\eta)\leq c'T_{\sigma}(r,\bar{q}\circ g,\omega_{S})+c'm(\sigma,\bar{q}\circ g,\lambda_W)+c'.
$$
Now for $1/2<s<r<1$ and $\sigma=(s+r)/2$, we have
$$
T_{s}(r,g,\eta)\leq 4T_{\sigma}(r,g,\eta)\leq 4c'T_{\sigma}(r,\bar{q}\circ g,\omega_{S})+4c'm(\sigma,\bar{q}\circ g,\lambda_W)+4c'.
$$
We set $c=4c'$.
By $T_{\sigma}(r,\bar{q}\circ g,\omega_{S})\leq T_{s}(r,\bar{q}\circ g,\omega_{S})$, we conclude the proof.
\hspace{\fill} $\square$

\begin{rem}\label{rem:202108051}
Let $A$ be a semi-abelian variety and let $S$ be a smooth projective variety.
Let $Z\subset S_{1,A}$ be an irreducible Zariski closed set.
Assume that $Z$ is horizontally integrable.
Then we take and fix the following objects:
\begin{itemize}
\item
An immersion $U\hookrightarrow A\times S$ as in Definition \ref{defn:20201122}.
\item
An isomorphism $\Phi:A\times U\to  A\times U$ from \eqref{eqn:20220708}.
\item
An irreducible component $\Xi\subset U_{1,A}$ of $(q')^{-1}(Z)$ from Lemma \ref{lem:20220707}, where $q':U_{1,A}\to S_{1,A}$ is induced from $q:U\to S$.
\item
The image $\Theta\subset U$ of $\Xi$ under $U_{1,A}\to U$.
Since $U_{1,A}\to U$ is proper, $\Theta$ is an irreducible Zariski closed subset of $U$.
\end{itemize}
We may apply Lemma \ref{lem:202302021} for $q:U\to S$ and $\Theta\subset U$ to get and fix the following objects:
\begin{itemize}
\item A smooth compactification $\Sigma$ of $U$ such that $q:U\to S$ extends to a morphism $\bar{q}:\Sigma\to S$.
\item A Zariski open set $\Sigma^o\subset \Sigma$ such that $\overline{\Theta}\subset \Sigma^o$, where $\overline{\Theta}\subset \Sigma$ is the Zariski closure of $\Theta$ in $\Sigma$.
\item
A smooth semi-positive $(1,1)$-form $\eta\geq 0$ on $\Sigma$ such that $\eta>0$ on $\Sigma^o$.
\item
A proper Zariski closed subset $W\subset S$ with $W\subsetneqq\bar{q}(\overline{\Theta})$.
\end{itemize}
By $\bar{q}(\overline{\Theta})=\tau(Z)$, we have $W\subsetneqq\tau(Z)$, where $\tau:S_{1,A}\to S$ is the natural map.
\end{rem}

\subsection{Estimate for the first term of RHS of \eqref{eqn:202302081}}
Let $A$ be a semi-abelian variety with the quotient $\rho:A\to A_0$ as in \eqref{eqn:20230227}.
Let $S$ be a smooth projective variety.
Let $Z\subset S_{1,A}$ be an irreducible Zariski closed set.
Assume that $Z$ is horizontally integrable.
We recall $W\subset S$ from Remark \ref{rem:202108051}.

\begin{lem}\label{lem:202302193}
Let $\mathcal{F}\subset \mathrm{Hol}(\mathbb D,A\times S)$ be an infinite set of non-constant holomorphic maps such that $( f_{S_{1,A}})_{f\in \mathcal{F}}\Rightarrow Z$.
Let $\omega_{A_0}$ be an invariant positive $(1,1)$ form on $A_0$ and let $\omega_S$ be a smooth positive $(1,1)$ form on $S$. 
Let $\lambda_W\geq 0$ be a Weil function for $W$.
Let $s\in (1/2,1)$ and $\delta>0$.
Then there exists a positive constant $\alpha>0$ such that for all $f\in \mathcal{F}$ with $f_S(\mathbb D)\not\subset W$, we have
\begin{equation*}
T_s(r,\rho\circ f_{A},\omega_{A_0}) \leq 
\alpha T_s(r,f_{S},\omega_{S})
+
\alpha
m((s+r)/2,f_S,\lambda_{W})
+\alpha
\end{equation*}
for all $r\in (s,1 )$ outside some exceptional set $E\subset (s,1)$ with the linear measure $|E|<\delta$.
\end{lem}

{\it Proof.}\
We use the objects fixed in Remark \ref{rem:202108051}.
Note that the isomorphism $\Phi:A\times U\to A\times U$ induces an isomorphism $\Phi_0:A_0\times \Sigma\to A_0\times \Sigma$ with the following commutative diagram:
\begin{equation}\label{eqn:202302243}
\begin{CD}
A\times U@>\Phi>>A\times U\\
@V(\rho,i)VV @VV(\rho,i)V   \\
A_0\times \Sigma@>>\Phi_0> A_0\times \Sigma
\end{CD}
\end{equation}
Here $i:U\to \Sigma$ is the open immersion.

We prove this claim.
Let $\varphi:U\to A$ be the composite of the immersion $\iota:U\hookrightarrow A\times U$ and the first projection $A\times S\to A$.
Then $\Phi:A\times U\to A\times U$ is defined by $\Phi (a,u)=(a-\varphi (u),u)$ (cf. \eqref{eqn:202302248}).
We have a rational map $\rho\circ\varphi:\Sigma\dashrightarrow A_0$.
Since $A_0$ is an abelian variety and $\Sigma$ is smooth, this rational map extends to a morphism $\rho\circ\varphi:\Sigma\to A_0$.
Then $\Phi_0$ is defined by $\Phi_0(b,u)=(b-\rho\circ\varphi(u),u)$ for $(b,u)\in A_0\times \Sigma$.
Then $\Phi_0$ fits in the commutative diagram \eqref{eqn:202302243}.

We set
\begin{equation}\label{eqn:202302041}
\delta'=\delta/801.
\end{equation}
We take an open neighbourhood $O_{\overline{\Theta}}\subset \Sigma$ of $\overline{\Theta}\subset \Sigma$ such that $\overline{O_{\overline{\Theta}}}\subset \Sigma^o$, where we recall $\Sigma^o\subset \Sigma$ from Remark \ref{rem:202108051}.
We take a smooth positive $(1,1)$-form $\omega_{\Sigma}$ on $\Sigma$ such that
\begin{equation}\label{eqn:202302156}
\omega_{\Sigma}\leq \eta
\end{equation}
on $O_{\overline{\Theta}}$, where we recall $\eta$ from Remark \ref{rem:202108051}.

We apply Lemma \ref{lem:20230218} to get a finite subset $\mathcal{E}\subset \mathcal{F}$ such that for each $f\in\mathcal{F}\backslash\mathcal{E}$, we get a lifting $\hat{f}\in \mathrm{Hol_m}(\mathbb D,A\times \Sigma)$ of $f\in \mathrm{Hol}(\mathbb D,A\times S)$, a connected open set $\Omega_f\subset \mathbb D$ and a connected component $\Omega_f'\subset Y_{\hat{f}}$ of $\pi_{\hat{f}}^{-1}(\Omega_f)$ with the properties \eqref{eqn:20210504?}--\eqref{eqn:202302034?} described in Lemma \ref{lem:20230218}.
Then we have (cf. \eqref{eqn:20211024?})
$$( (\hat{f}_{\Sigma},\Omega_f,\Omega_f'))_{f\in \mathcal{F}\backslash\mathcal{E}} \rightsquigarrow\overline{\Theta}.$$
Hence there exists a finite subset $\mathcal{E}_1\subset \mathcal{F}$, where $\mathcal{E}\subset \mathcal{E}_1$, such that for all $f\in \mathcal{F}\backslash \mathcal{E}_1$, we have $z\in \Omega_f$ for $\delta'$-almost all $z\in \mathbb D(1-\delta')$ and 
\begin{equation}\label{eqn:202302157}
\hat{f}_{\Sigma}(\Omega_f')\subset O_{\overline{\Theta}}.
\end{equation}

Let $\omega_A$ be an invariant positive $(1,1)$-form on $A$ such that 
\begin{equation}\label{eqn:20230216}
\rho^*\omega_{A_0}\leq \omega_A.
\end{equation}
We take an open neighbourhood $O_{\Sigma_{1,A}^*}\subset \Sigma_{1,A}$ of $\Sigma_{1,A}^*$ such that
\begin{equation}\label{eqn:202302029}
|v_{\mathrm{Lie}A}|_{\omega_A}\leq | v_{\Sigma}|_{\omega_{\Sigma}}
\end{equation}
for all $(v_{\Sigma},v_{\mathrm{Lie}A})\in T\Sigma\times \mathrm{Lie}A$ with $[(v_{\Sigma},v_{\mathrm{Lie}A})]\in O_{\Sigma_{1,A}^*}
\subset \Sigma_{1,A}$.
We have (cf. \eqref{eqn:202102261?})
$$( ((\Phi\circ \hat{f})_{\Sigma_{1,A}},\Omega_f,\Omega_f'))_{f\in \mathcal{F}\backslash\mathcal{E}} \rightsquigarrow \Sigma_{1,A}^*.$$
Hence, there exists a finite subset $\mathcal{E}_2\subset \mathcal{F}$, where $\mathcal{E}_1\subset \mathcal{E}_2$, such that for all $f\in \mathcal{F}\backslash \mathcal{E}_2$, we have 
\begin{equation}\label{eqn:202302159}
(\Phi\circ\hat{f})_{\Sigma_{1,A}}(\Omega_f')\subset O_{\Sigma_{1,A}^*}.
\end{equation}

Now we take $f\in \mathcal{F}\backslash \mathcal{E}_2$ such that $f_S(\mathbb D)\not\subset W$.
Since $z\in \Omega_f$ for $\delta'$-almost all $z\in \mathbb D(1-\delta')$, we may apply \cite[Lemma 7.1]{Y} to get 
\begin{equation}\label{eqn:202302152}
T_s(r,\rho\circ f_A,\omega_{A_0})\leq
8\int_s^r\frac{dt}{t}\int_{\mathbb D(t)\cap \Omega_f}(\rho\circ f_A)^*\omega_{A_0}
\end{equation}
for all $r\in (s,1-\delta')$ outside some exceptional set $E'\subset (s,1-\delta')$  whose linear measure satisfies
\begin{equation}\label{eqn:2023021510}
|E'|\leq 800\delta'.
\end{equation}
By $\hat{f}_A=f_A\circ\pi_{\hat{f}}$, we have
\begin{equation}\label{eqn:202302153}
\int_s^r\frac{dt}{t}\int_{\mathbb D(t)\cap \Omega_f}(\rho\circ f_A)^*\omega_{A_0}=T_s(r,(\rho\circ \hat{f}_A,\Omega_f,\Omega_f'),\omega_{A_0}).
\end{equation}
Here we recall \eqref{eqn:20230215} for the definition of the right hand side.

Let $\Phi_0:A_0\times\Sigma\to A_0\times\Sigma$ be the isomorphism above (cf. \eqref{eqn:202302243}).
Let $\nu_1:A_0\times \Sigma\to A_0$ be the first projection, and $\nu_2:A_0\times \Sigma\to \Sigma$ be the second projection.
Then there exists a positive constant $\alpha_1>0$ such that 
$$(\nu_1\circ\Phi_0^{-1})^*\omega_{A_0}\leq \alpha_1(\nu_1^*\omega_{A_0}+\nu_2^*\omega_{\Sigma}).$$
We define $\tilde{\rho}:A\times \Sigma \to A_0\times \Sigma$ by $\tilde{\rho}(a,s)=(\rho(a),s)$ for $(a,s)\in A\times \Sigma$.
Then, we get
\begin{multline}\label{eqn:202302154}
T_s(r,(\rho\circ \hat{f}_A,\Omega_f,\Omega_f'),\omega_{A_0})
\\
\leq \alpha_1T_s(r,((\Phi_0\circ\tilde{\rho}\circ\hat{f})_{A_0},\Omega_f,\Omega_f'),\omega_{A_0})
+\alpha_1T_s(r,(\hat{f}_{\Sigma},\Omega_f,\Omega_f'),\omega_{\Sigma}).
\end{multline}
By \eqref{eqn:202302029} and \eqref{eqn:202302159}, we have
$$
T_s(r,((\Phi\circ\hat{f})_A,\Omega_f,\Omega_f'),\omega_A)
\leq T_s(r,(\hat{f}_{\Sigma},\Omega_f,\Omega_f'),\omega_{\Sigma}).
$$
By the commutativity of \eqref{eqn:202302243} and \eqref{eqn:20230216}, we have
$$
T_s(r,((\Phi_0\circ\tilde{\rho}\circ\hat{f})_{A_0},\Omega_f,\Omega_f'),\omega_{A_0})
\leq T_s(r,((\Phi\circ\hat{f})_A,\Omega_f,\Omega_f'),\omega_A).
$$
Hence we get
\begin{equation}\label{eqn:202302155}
T_s(r,((\Phi_0\circ\tilde{\rho}\circ\hat{f})_{A_0},\Omega_f,\Omega_f'),\omega_{A_0})
\leq T_s(r,(\hat{f}_{\Sigma},\Omega_f,\Omega_f'),\omega_{\Sigma}).
\end{equation}
Hence by \eqref{eqn:202302152}, \eqref{eqn:202302153}, \eqref{eqn:202302154} and \eqref{eqn:202302155}, we get
$$
T_s(r,\rho\circ f_A,\omega_{A_0})\leq 16\alpha_1T_s(r,(\hat{f}_{\Sigma},\Omega_f,\Omega_f'),\omega_{\Sigma})
$$
for all $r\in (s,1-\delta')\backslash E'$.
Using \eqref{eqn:202302156} and \eqref{eqn:202302157}, we have
$$
T_s(r,(\hat{f}_{\Sigma},\Omega_f,\Omega_f'),\omega_{\Sigma})\leq T_s(r,(\hat{f}_{\Sigma},\Omega_f,\Omega_f'),\eta)\leq (\deg\pi_{\hat{f}})T_s(r,\hat{f}_{\Sigma},\eta).
$$
Hence we get
\begin{equation}\label{eqn:202302245}
T_s(r,\rho\circ f_A,\omega_{A_0})\leq \alpha_2T_s(r,\hat{f}_{\Sigma},\eta)
\end{equation}
for all $r\in (s,1-\delta')\backslash E'$, where we set $\alpha_2=16\alpha_1[\mathbb C(S):\mathbb C(\Sigma)]$ (cf. \eqref{eqn:202302034?}).
We set $\sigma=(r+s)/2$.
By Lemma \ref{lem:202302021}, we get 
\begin{equation}\label{eqn:202302246}
T_s(r,\hat{f}_{\Sigma},\eta)\leq \alpha_3 T_{s}(r,f_S,\omega_S)+\alpha_3 m(\sigma,f_S,\lambda_W)+\alpha_3,
\end{equation}
where $\alpha_3>0$ is a positive constant which is independent of the choice of $f\in\mathcal{F}\backslash\mathcal{E}_2$.
We set $\alpha=\alpha_2\alpha_3$ and $E=E'\cup (1-\delta',1)$.
By \eqref{eqn:202302041} and \eqref{eqn:2023021510}, we have $|E|\leq \delta$.
Then by \eqref{eqn:202302245} and \eqref{eqn:202302246}, we get the desired estimate for $f\in \mathcal{F}\backslash \mathcal{E}_2$ with $f_S(\mathbb D)\not\subset W$.

Finally we enlarge $\alpha>0$ so that
$$
\max_{f\in\mathcal{E}_2}T_s\left(1-\delta,\rho\circ f_A,\omega_{A_0}\right)
\leq \alpha
$$
to complete the proof of the lemma.
\hspace{\fill} $\square$

\subsection{Estimate of Weil functions}\label{subsec:20230222}

Let $\omega_A$ be an invariant positive $(1,1)$-form on $A$.
Given $x,y\in A$, we denote by $d(x,y)$ the distance with respect to $\omega_A$.
Let $\pi:\mathbb C^n\to A$ be a universal covering, where $n=\dim A$.
We denote by $d_{\mathbb C^n}$ the Euclidean distance on $\mathbb C^n$.
Note that $\pi^*\omega_A$ is an invariant positive $(1,1)$-form on $\mathbb C^n$ with respect to the additive structure of $\mathbb C^n$.
Hence there exists a positive constant $\alpha>1$ such that
\begin{equation}\label{eqn:202205261}
\frac{1}{\alpha}d(x,y)\leq d_{\mathbb C^n}(\pi^{-1}(x),\pi^{-1}(y))\leq \alpha d(x,y)
\end{equation}
for all $x,y\in A$.

\begin{lem}\label{lem:20220526}
Let $B\subset A$ be a semi-abelian subvariety.
Let $d_B$ be the distance function on $B$ with respect to some invariant positive $(1,1)$-form on $B$.
Then there exists a positive constant $\beta>1$ such that
$$
\frac{1}{\beta}d_B(x,y)\leq d(x,y)\leq \beta d_B(x,y)
$$
for all $x,y\in B$.
\end{lem}

{\it Proof.}\
Let $\mathbb C^k\to B$ be a universal covering, where $k=\dim B$.
Then $B\subset A$ induces an immersion $\mathbb C^k\subset \mathbb C^n$ into the universal covering of $A$.
This is a linear subspace.
Hence for $p,q\in \mathbb C^k$, we have $d_{\mathbb C^k}(p,q)=d_{\mathbb C^n}(p,q)$.
Hence by \eqref{eqn:202205261}, we obtain our lemma.
\hspace{\fill} $\square$

\begin{lem}\label{lem:torusdistance}
Let $\overline{A}$ be a smooth equivariant compactification.
Let $D\subset \partial A$ be an irreducible component with a Weil function $\lambda_{D}$.
There exists a positive constant $c>0$ such that
$$
\lambda_{D} (y)\leq \lambda_{D}(x)+cd(x,y)+c
$$
for all $x,y\in A$.
\end{lem}

{\it Proof.}
We first consider the case that $A$ is an algebraic torus, and then prove the general case.

{\it The case of algebraic tori.}
The proof is by induction on the dimension of $A$.
When $\dim A=1$, we have $A=\mathbb G_m$ and $\overline{A}=\mathbb P^1$.
We may assume $D=(\infty)$.
Then
$$
|\lambda_{D}(x)-\lambda_{D}(y)|\leq |\log^+|x|-\log^+|y||\leq |\log |x/y||.
$$
On the other hand, by \eqref{eqn:202205261}, we have
\begin{equation}\label{eqn:20230117}
|\log |x/y||\leq \alpha d(x,y).
\end{equation}
This shows our estimate in the one dimensional case.

Now we consider the case $\dim A\geq 2$ and assume that our lemma is true for algebraic tori whose dimension is less than $\dim A$.
Let $D\subset \partial A$ be an irreducible component.
Let $I\subset A$ be the isotoropy group of $D$.
Then $I=\mathbb{G}_m$.
Set $B=A/I$.
We take a compactification $\overline{B}$ of $B$ such that $\overline{B}=D$.
By Lemma \ref{lem:vect}, there exist a Zariski open neighborhood $U\subset \overline{A}$ of $D\subset \overline{A}$ and a canonical projection $p:U\to D$ which extends $A\to B$.
Note that $p:U\to D$ is a total space of a line bundle over $\overline{B}$ whose zero section is $D$. 
Let $||\cdot ||$ be a smooth Hermitian metric on the line bundle $p:U\to \overline{B}$.
We define $\mu: U\backslash D\to \mathbb R$ by $\mu (x)=\log (1/||x||)$ where $x\in U$. 
We prove the estimate of Lemma \ref{lem:torusdistance} in several steps.

\medskip

{\it Step 1.}\
We first show that there exists a positive constant $c_1>0$ such that
\begin{equation}\label{eqn:202205263}
|\mu (x)-\mu(y)|\leq c_1 d(x,y)+c_1
\end{equation}
for all $x,y\in A$ which satisfy $p(x)=p(y)$.
To prove this, we take a finite affine covering $\overline{B}=\cup_{i\in I}V_i$ such that $U|_{V_i}=V_i\times \mathbb C$.
Let $\tau_i:U|_{V_i}\to \mathbb C$ be the second projection.
We may take an open set $W_i\Subset V_i$ such that $\overline{B}=\cup_{i\in I}W_i$ and $|\mu(x)-\log |1/\tau_i(x)||<\gamma_i$ on $x\in p^{-1}(W_i)\backslash D$.
We set $\gamma=\max_{i\in I}\gamma_i$.
Now given $x,y\in A$ such that $p(x)=p(y)$, we take $i\in I$ such that $p(x)\in W_i$.
Then we have
$$
|\mu (x)-\mu(y)|\leq |\log |1/\tau_i(x)|-\log |1/\tau_i(y)||+2\gamma. 
$$
By the one dimensional case (cf. \eqref{eqn:20230117}), we have
$$
|\log |1/\tau_i(x)|-\log |1/\tau_i(y)||\leq \alpha d_{\mathbb G_m}(\tau_i(x),\tau_i(y)).
$$
Hence we get
$$
|\mu (x)-\mu(y)|\leq \alpha d_{\mathbb G_m}(\tau_i(x),\tau_i(y))+2\gamma.
$$
We take $g\in I$ such that $y=g\cdot x$.
Since $\tau_i$ is $\mathbb{G}_m$ equivariant, we have $\tau_i(y)=g\cdot \tau_i(x)$.
Hence we have $d_{\mathbb G_m}(\tau_i(x),\tau_i(y))=d_{\mathbb G_m}(e_{\mathbb {G}_m},g)$.
Hence we get
$$
|\mu (x)-\mu(y)|\leq \alpha d_{\mathbb G_m}(e_{\mathbb {G}_m},g)+2\gamma.
$$
Similarly, we have $d(x,y)=d(e_A,g)$.
By Lemma \ref{lem:20220526}, we have $d_{I}(e_A,g)\leq \beta d(e_A,g)$.
Hence we get \eqref{eqn:202205263} with $c_1=\max\{ \alpha\beta, 2\gamma\}$.

\medskip

{\it Step 2.}\
By \cite[Corollary, p. 115]{Borellinear}, the quotient $A\to B$ has a section $s:B\to A$ of group varieties.
We next show 
\begin{equation}\label{eqn:202205264}
|\mu (x)-\mu (y)|\leq c_2 d(x,y)+c_2
\end{equation}
for all $x,y\in s(B)$.
By the section $s:B\to A$, we get a rational section $\bar{s}:\overline{B}\dashrightarrow U$, which is holomorphic and zero-free on $B$.
Set $(s)=E-F$.
Then $\mu (s(b))=\lambda_E(b)-\lambda_F(b)$ for $b\in B$.
Then by the induction hypothesis, we have, for $b,b'\in B$ with $x=s(b)$ and $y=s(b')$,
\begin{equation*}
\begin{split}
|\mu (x)-\mu (y)|
&= | (\lambda_E(b)-\lambda_F(b))-(\lambda_E(b')-\lambda_F(b'))|
\\
&\leq |\lambda_E(b)-\lambda_E(b')|+|\lambda_F(b)-\lambda_F(b')|
\\
&\leq \gamma d_B(b,b')+\gamma.
\end{split}
\end{equation*}
By Lemma \ref{lem:20220526}, we have $d_B(b,b')\leq \alpha d(x,y)$.
Hence we get \eqref{eqn:202205264} with $c_2=\max\{ \alpha\gamma,\gamma\}$.

\medskip

{\it Step 3.}\
Now we take $x,y\in A$ such that $x\in s(B)$.
Then by \eqref{eqn:202205263} and \eqref{eqn:202205264}, we have
\begin{equation*}
\begin{split}
|\mu (x)-\mu(y)|
&\leq |\mu (x)-\mu(s(p(y)))|+|\mu (s(p(y)))-\mu (y)|
\\
&\leq c_2 d(x,s(p(y)))+c_1 d(s(p(y)),y)+c_1+c_2.
\end{split}
\end{equation*}
By
$$
d(s(p(y)),y)\leq d(x,s(p(y)))+d(x,y),
$$
we get
$$
|\mu (x)-\mu(y)|\leq (c_1+c_2) d(x,s(p(y)))+c_1 d(x,y)+c_1+c_2.
$$
Since $s((p(x))=x$, Lemma \ref{lem:20220526} yields $d(x,s(p(y)))\leq \alpha d_B(p(x),p(y))$.
Hence
$$
|\mu (x)-\mu(y)|\leq (c_1+c_2)\alpha d_B(p(x),p(y))+c_1 d(x,y)+c_1+c_2.
$$
Since $p|_A:A\to B$ is a group homomorphism, we have
\begin{equation}\label{eqn:202205266}
d_B(p(x),p(y))\leq \alpha' d(x,y).
\end{equation}
Indeed, $(p|_A)^*\omega_B$ is an invariant $(1,1)$-form on $A$, where $\omega_B$ is the invariant positive $(1,1)$-form on $B$ used to define $d_B$.
Hence there exists a positive constant $\alpha'>0$ such that $(p|_A)^*\omega_B\leq \alpha'\omega_A$.
To conclude, we get
$$
|\mu (x)-\mu(y)|\leq  c d(x,y)+c,
$$
where $c=\max\{ (c_1+c_2)\alpha\alpha'+c_1,c_1+c_2\}$.

\medskip

{\it Step 4.}\
Now we take $x,y\in A$ in general.
We take $g\in  I$ such that $s(p(x))=g\cdot x$.
Then we get
$$
|\mu(x)-\mu(y)|=|\mu(g\cdot x)-\mu(g\cdot y)|\leq c d(g\cdot x,g\cdot y)+c=c d(x,y)+c.
$$
Hence 
$$|\lambda_D(x)-\lambda_D(y)|\leq |\mu (x)-\mu(y)|\leq cd(x,y)+c,
$$
which concludes the induction step.
Hence we get our lemma in the case of algebraic tori.

\medskip

{\it The case of general semi-abelian varieties.}\
We treat the general case as in \eqref{eqn:20230227}:
\begin{equation*}
0\to T\to A\overset{\rho}{\to} A_0\to 0.
\end{equation*}
We pull-back this sequence by the universal covering $\pi: \mathbb C^n\to A_0$ to get
$$
0\to T\to A\times_{A_0}\mathbb C^n\overset{r}{\to} \mathbb C^n\to 0.
$$
We have a section $s:\mathbb C^n\to A\times_{A_0}\mathbb C^n$ of complex Lie groups.
Let $p:A\times_{A_0}\mathbb C^n\to A$ be the natural projection.
We set 
$$\psi=p\circ s:\mathbb C^n\to A.$$
Then $p$ and $\psi$ are morphisms of complex Lie groups.
We take a closed ball $\mathbb B\subset \mathbb C^n$ centered at the origin such that $\pi (\mathbb B)=A_0$.

We first show that there exists a positive constant $\alpha_1>0$ such that for all $g\in T$ and $z,w\in \mathbb B$, we have
\begin{equation}\label{eqn:202102271}
|\lambda_D(g\cdot \psi(z))-\lambda_D(g\cdot \psi(w))|\leq \alpha_1.
\end{equation}

We prove this.
Let $\overline{A}$ corresponds to a torus embedding $T\subset \overline{T}$ (cf. Lemma \ref{lem:20200906}).
Note that $\rho:A\to A_0$ extends to $\bar{\rho}:\overline{A}\to A_0$.
Then $r:A\times_{A_0}\mathbb C^n\to\mathbb C^n$ extends to $\bar{r}:\overline{A}\times_{A_0}\mathbb C^n\to\mathbb C^n$.
The section $s:\mathbb C^n\to A\times_{A_0}\mathbb C^n$ induces a (non-canonical) splitting $\overline{A}\times_{A_0}\mathbb C^n=\overline{T}\times \mathbb C^n$ such that the composite $q\circ s:\mathbb C^n\to T$ with the first projection $q:\overline{A}\times_{A_0}\mathbb C^n\to \overline{T}$ is the constant map identically equal to $e_T$.
\begin{equation*}
\begin{CD}
\overline{A}@<p<<\overline{A}\times_{A_0}\mathbb C^n  @>q>> \overline{T}\\
@V\bar{\rho} VV @VV\bar{r}V    \\
A_0@<<\pi< \mathbb C^n
\end{CD}
\end{equation*}
Note that $p^*\lambda_{D}(x)-(\lambda_D|_{\rho^{-1}(0)})(q(x))$ is continuous on $\overline{A}\times_{A_0}\mathbb C^n$.
Since $\bar{r}^{-1}(\mathbb B)$ is compact, there exists a positive constant $\alpha_1'>0$ such that $|p^*\lambda_{D}(x)-(\lambda_D|_{\rho^{-1}(0)})(q(x))|<\alpha_1'$ for all $x\in \bar{r}^{-1}(\mathbb B)$.
Note that $q(g\cdot s(z))=g$ for all $g\in T$ and $z\in\mathbb C^n$.
Hence for all $z\in \mathbb B$ and $g\in T$, we have
$$
|p^*\lambda_{D}(g\cdot s(z))-(\lambda_D|_{\rho^{-1}(0)})(g)|\leq \alpha_1'.
$$
Hence by
$$
|\lambda_D(g\cdot \psi(z))-\lambda_D(g\cdot \psi(w))|
\leq |p^*\lambda_{D}(g\cdot s(z))-(\lambda_D|_{\rho^{-1}(0)})(g)|+|p^*\lambda_{D}(g\cdot s(w))-(\lambda_D|_{\rho^{-1}(0)})(g)|,$$
we get \eqref{eqn:202102271} with $\alpha_1=2\alpha_1'$.

Next we prove that there exists a positive constant $\alpha_2>0$ such that for all $g\in T$ and $z,w\in \mathbb B$, we have
\begin{equation}\label{eqn:202102272}
d(g\cdot \psi(z),g\cdot \psi(w))\leq \alpha_2.
\end{equation}

We prove this.
Since $\psi:\mathbb C^n\to A$ is a group homomorphism, $\psi^*\omega_A$ is an invariant $(1,1)$-form on $\mathbb C^n$.
Hence we have $d(\psi(z),\psi(w))\leq \alpha_2'd_{\mathbb C^n}(z,w)$ for all $z,w\in\mathbb C^n$.
See the argument for \eqref{eqn:202205266}.
Hence
$$
d(g\cdot \psi(z),g\cdot \psi(w))=d(\psi(z),\psi(w))\leq \alpha_2'd_{\mathbb C^n}(z,w)
$$
for all $g\in T$ and $z,w\in\mathbb C^n$.
Since $\mathbb B\subset \mathbb C^n$ is compact, there exists a positive constant $\gamma>0$ such that $d_{\mathbb C^n}(z,w)<\gamma$ for all $z,w\in \mathbb B$.
Hence we get \eqref{eqn:202102272} with $\alpha_2=\alpha_2'\gamma$.

Now, for $x,y\in A$, we take $x',y'\in r^{-1}(\mathbb B)$ such that $p(x')=x$ and $p(y')=y$.
We take $g_x,g_y\in T$ such that $g_x\cdot s(r(x'))= x'$ and $g_y\cdot s(r(y'))= y'$.
Then by the torus case above, we have
$$
|\lambda_D(g_x\cdot \psi(0))-\lambda_D(g_y\cdot \psi(0))|\leq c_1 d_T(g_x\cdot \psi(0),g_y\cdot \psi(0))+c_1.
$$
By Lemma \ref{lem:20220526}, we have $d_T(g_x\cdot\psi(0),g_y\cdot \psi(0))\leq \beta d(g_x\cdot \psi(0),g_y\cdot \psi(0))$.
Hence
\begin{equation}\label{eqn:202102273}
|\lambda_D(g_x\cdot \psi(0))-\lambda_D(g_y\cdot \psi(0))|\leq c_1\beta d(g_x\cdot\psi(0),g_y\cdot \psi(0))+c_1.
\end{equation}
We have
\begin{multline*}
|\lambda_D(x)-\lambda_D(y)|
\leq |\lambda_D(g_x\cdot \psi(r(x')))-\lambda_D(g_x\cdot \psi(0))|
\\
+|\lambda_D(g_x\cdot \psi(0))-\lambda_D(g_y\cdot \psi(0))|+|\lambda_D(g_y\cdot \psi(0))-\lambda_D(g_y\cdot \psi(r(y')))|.
\end{multline*}
Hence by \eqref{eqn:202102271} and \eqref{eqn:202102273}, we get
$$
|\lambda_D(x)-\lambda_D(y)|\leq c_1\beta d(g_x\cdot \psi(0),g_y\cdot \psi(0))+c_1+2\alpha_1.
$$
By \eqref{eqn:202102272}, we have
\begin{equation*}
\begin{split}
d(g_x\cdot \psi(0),g_y\cdot \psi(0))&\leq d(g_x\cdot \psi(0),g_x\cdot \psi(r(x')))\\
&\qquad +d(g_x\cdot \psi(r(x')),g_y\cdot \psi(r(y')))+d(g_y\cdot \psi(r(y')),g_y\cdot \psi(0))\\
&\leq  d(x,y)+2\alpha_2.
\end{split}
\end{equation*}
Hence we get
$$
|\lambda_D(x)-\lambda_D(y)|\leq c d(x,y)+c,
$$
where $c=\max\{ c_1\beta,2c_1\alpha_2\beta+c_1+2\alpha_1\}$.
This conclude the proof.
\hspace{\fill} $\square$

\subsection{Application of Lemma \ref{lem:torusdistance}}

We recall the notation from \eqref{eqn:202302247}. 

\begin{lem}\label{lem:20230219}
Let $\overline{A}$ be a smooth equivariant compactification of a semi-abelian variety $A$, and let $\Sigma$ be a smooth projective variety.
Let $D\subset \partial A$ be an irreducible component with a Weil function $\lambda_D\geq 0$.
Let $\omega_A$ be an invariant positive $(1,1)$-form on $A$.
Then there exists a positive constant $c>0$ with the following property:
Let $(g,\Omega,\Omega')$ be a triple as in Definition \ref{defn:20230124}, where $g\in \mathrm{Hol_m}(\mathbb D,\overline{A})$ with $g(Y_g)\not\subset \partial A$.
We take $r,r'\in (0,1)$ such that $\partial \mathbb D(r)\subset \Omega$ and $\partial\mathbb D(r')\subset \Omega$.
We take $\theta\in\mathbb R$ such that the line segment $\gamma$ connecting $re^{i\theta}$ and $r'e^{i\theta}$ satisfies $\gamma\subset \Omega$.
Assume that $\partial\mathbb D(r')\cup \partial\mathbb D(r)\cup \gamma\subset \Omega$ does not contain the critical values of $\pi_g:Y_g\to\mathbb D$.
Then we have
\begin{multline}\label{eqn:202302125}
|m(r,(g,\Omega,\Omega'),\lambda_D)-m(r',(g,\Omega,\Omega'),\lambda_D)|\\
\leq \frac{c}{\deg(\pi_g|_{\Omega'})}\ell_{\omega_{A}}(g(\Omega'\cap\pi_g^{-1}(\partial\mathbb D(r)\cup\partial\mathbb D(r')\cup\gamma)))+c.
\end{multline}
\end{lem}

{\it Proof.}\
By Lemma \ref{lem:torusdistance}, there exists a positive constant $\alpha>0$ such that
\begin{equation}\label{eqn:20230212}
|\lambda_D(x)-\lambda_D(y)|\leq \alpha d(x,y)+\alpha
\end{equation}
for all $x,y\in A$, where $d(x,y)$ is the distance with respect to $\omega_A$.

Let $(g,\Omega,\Omega')$ be a triple as in Definition \ref{defn:20230124}, where $g\in \mathrm{Hol_m}(\mathbb D,\overline{A})$ with $g(Y_g)\not\subset \partial A$.
For each $z\in \Omega$, we set
$$
\mu_g(z)=\frac{1}{\deg(\pi_g|_{\Omega'})}\sum_{y\in \Omega'\cap\pi_g^{-1}(z)}\lambda_D(g(y)).
$$
We have
\begin{equation}\label{eqn:202302123}
m(r,(g,\Omega,\Omega'),\lambda_D)=\int_{z\in \Omega\cap\partial\mathbb D(r)}\mu_g(z)\frac{d\arg \pi_g(z)}{2\pi}
\end{equation}

Let $z,z'\in \Omega$ be two points connected by a smooth arc $\sigma$ in $\Omega$.
Suppose that $\sigma\subset \Omega$ does not contain the critical values of $\pi_g:Y_g\to\mathbb D$.
Then we claim
\begin{equation}\label{eqn:202302121}
|\mu_g(z)-\mu_g(z')|\leq 
\frac{\alpha}{\deg(\pi_g|_{\Omega'})}\ell_{\omega_{A}}(g(\Omega'\cap\pi_g^{-1}(\sigma)))+\alpha.
\end{equation}

To show this, we set $\Omega'\cap\pi_g^{-1}(z)=\{y_1,\ldots,y_k\}$ and $\Omega'\cap\pi_g^{-1}(z)=\{y_1',\ldots,y_k'\}$, where $k=\deg(\pi_g|_{\Omega'})$.
We may assume that $y_i$ and $y_i'$ are connected by a smooth arc $\sigma_i$ in $\Omega'$, where $\sigma_i$ is a lift of $\sigma$.
Then by \eqref{eqn:20230212}, we get 
$$
|\mu_g(z)-\mu_g(z')|\leq \frac{1}{k}\sum_{i=1}^k|\lambda_D(g(y_i))-\lambda_D(g(y_i'))|\leq \frac{\alpha}{k}\sum_{i=1}^kd(g(y_i),g(y_i'))+\alpha.
$$
By $d(g(y_i),g(y_i'))\leq \ell_{\omega_A}(g(\sigma_i))$, we get \eqref{eqn:202302121}.

Now we take $r,r'\in (0,1)$, $\theta\in\mathbb R$ and $\gamma$ as in Lemma \ref{lem:20230219}.
We have $\partial \mathbb D(r)\subset \Omega$, $\partial\mathbb D(r')\subset \Omega$ and $\gamma\subset \Omega$.
Moreover $\partial\mathbb D(r')\cup \partial\mathbb D(r)\cup \gamma\subset \Omega$ does not contain the critical values of $\pi_g:Y_g\to\mathbb D$.
By \eqref{eqn:202302121}, we have
$$
|\mu_g(re^{i\theta'})-\mu_g(re^{i\theta})|\leq 
\frac{\alpha}{\deg(\pi_g|_{\Omega'})}\ell_{\omega_{A}}(g(\Omega'\cap\pi_g^{-1}(\partial\mathbb D(r))))+\alpha$$
for all $\theta'\in\mathbb R$.
Hence by \eqref{eqn:202302123}, we have
$$
|m(r,(g,\Omega,\Omega'),\lambda_D)-\mu_g(re^{i\theta})|
\leq 
\frac{\alpha}{\deg(\pi_g|_{\Omega'})}\ell_{\omega_{A}}(g(\Omega'\cap\pi_g^{-1}(\partial\mathbb D(r))))+\alpha.$$
Similarly we have
$$
|m(r',(g,\Omega,\Omega'),\lambda_D)-\mu_g(r'e^{i\theta})|
\leq 
\frac{\alpha}{\deg(\pi_g|_{\Omega'})}\ell_{\omega_{A}}(g(\Omega'\cap\pi_g^{-1}(\partial\mathbb D(r'))))+\alpha.$$
By \eqref{eqn:202302121}, we have
$$
|\mu_g(re^{i\theta})-\mu_g(r'e^{i\theta})|
\leq 
\frac{\alpha}{\deg(\pi_g|_{\Omega'})}\ell_{\omega_{A}}(g(\Omega'\cap\pi_g^{-1}(\gamma)))+\alpha.
$$
The last three estimates yield \eqref{eqn:202302125}, where $c=3\alpha$.
\hspace{\fill} $\square$

\subsection{Application of the area-length method}\label{subsec:area}

For $r\in (\frac{1}{2},1)$ and $\theta\in [0,2\pi]$, we denote by $\gamma_{r,\theta}\subset \mathbb D$ the line segment connecting $\frac{1}{2}e^{i\theta}$ and $re^{i\theta}$.
Let $\eta$ be a smooth semi-positive $(1,1)$-form on a smooth projective variety $\Sigma$. 
We denote by $\ell_{\eta}$ the length of curves in $\Sigma$ with respect to $\eta$. 
The next lemma is an application of the area-length method.

\begin{lem}\label{lem:20210423}
Let $f\in \mathrm{Hol_m}(\mathbb D,\Sigma)$ and $s\in (\frac{1}{2},1)$.
Let $\delta >0$.

(1)
There exists a subset $E_1\subset (s,1)$ with $|E_1|\leq \delta$ such that, for all $r\in (s,1)\backslash E_1$, we have
$$
\ell_{\eta}(f(\pi_f^{-1}(\partial\mathbb D(r))))\leq c_1 (\deg \pi_f)^2T_s(r,f,\eta)+c_1(\deg \pi_f)^2,
$$
where $c_1>0$ is a positive constant which only depends on $\delta$.

(2)
There exists a subset $E_2\subset (s,1)$ with $|E_2|\leq \delta$ such that, for all $r\in (s,1)\backslash E_2$, we have
$$
\ell_{\eta}(f(\pi_f^{-1}(\gamma_{r,\theta})))\leq c_2 (\deg \pi_f)T_s(r,f,\eta)+c_2 (\deg \pi_f)
$$
for all $\theta\in (0,2\pi )$ outside some exceptional set $E_3\subset  (0,2\pi )$ with linear measure $|E_3|<\delta$.
Here $c_2>0$ is a positive constant which only depends on $\delta$.
\end{lem}

{\it Proof.}
Let $f^*\eta=\varphi^2\pi_f^*(dx\wedge dy)$, where $z=x+iy$.
Set
$$
A(r)=\int_{Y_f(r)}f^*\eta
$$
and
$$
T(r)=T_s(r,f,\eta).
$$
Then we have 
$$
A(r)=r(\deg \pi_f)T'(r)\leq (\deg \pi_f)T'(r).
$$
By
$$
A(r)=\int_0^rdt \int_{\pi_f^{-1}(\partial\mathbb D(t))}\varphi^2td\arg \pi_f,
$$
we have
$$
A'(r)=\int_{\pi_f^{-1}(\partial\mathbb D(r))}\varphi^2rd\arg \pi_f .
$$
Set $L(r)=\ell_{\eta}(f(\pi_f^{-1}(\partial\mathbb D(r))))$.
Then
$$
L(r)=\int_{\pi_f^{-1}(\partial\mathbb D(r))}\varphi rd\arg \pi_f .
$$
Hence by the Cauchy-Schwarz inequality, we get
\begin{equation*}
L(r)^2
\leq 2\pi r (\deg \pi_f)\int_{\pi_f^{-1}(\partial\mathbb D(r))}\varphi^2rd\arg \pi_f \leq  2\pi (\deg \pi_f)A'(r).
\end{equation*}
We estimate the right hand side.
Note that since $\eta$ is semi-positive, we have $A'(r)\geq 0$ for all $r\in (0,1)$.
Hence we may apply Borel growth lemma (Lemma \ref{lem:growth} below) for $A(r)$.
Letting $\varepsilon=\sqrt{2}-1$ and $\delta'=\delta/2$, we have
$$
A'(r)\leq \frac{2}{\varepsilon\delta'}\max\{ 1, A(r)^{1+\varepsilon}\}\leq \frac{2}{\varepsilon\delta'}(\deg \pi_f)^{1+\varepsilon}\max\{ 1,(T'(r))^{1+\varepsilon}\}
$$
for all $r\in (s,1)$ outside some exceptional set $E_1'\subset (s,1)$ of linear measure less than $\delta'$.
Again by Lemma \ref{lem:growth}, we have
$$
(T'(r))^{1+\varepsilon}\leq \left(\frac{2}{\varepsilon\delta'}\right)^{1+\varepsilon}\max\left\{ 1, T(r)^{(1+\varepsilon)^2}\right\}
$$
for all $r\in (s,1)$ outside some exceptional set $E_1''\subset (s,1)$ of linear measure less than $\delta'$.
Hence
$$
A'(r)\leq \left(\frac{2}{\varepsilon\delta'}\right)^{2+\varepsilon}(\deg \pi_f)^{1+\varepsilon}\max\{ 1,T(r)^{(1+\varepsilon)^2}\}
$$
for all $r\in (s,1)$ outside $E_1=E_1'\cup E_1''$, where $|E_1|<\delta$.
Then we get
$$
L(r)^2\leq 2\pi\left(\frac{2}{\varepsilon\delta'}\right)^{2+\varepsilon}(\deg \pi_f)^{2+\varepsilon}\max\{ 1,T(r)^{2}\},
$$
thus
\begin{equation*}
\begin{split}
L(r)
&\leq \sqrt{2\pi}\left(\frac{2}{\varepsilon\delta'}\right)^{1+\varepsilon/2}(\deg \pi_f)^{1+\varepsilon/2}\max\{ 1,T(r)\}\\
&\leq c_1(\deg \pi_f)^2T(r)+c_1(\deg \pi_f)^2
\end{split}
\end{equation*}
for all $r\in (s,1)$ outside $E_1$.
Here we set $c_1=\sqrt{2\pi}\left(\frac{2}{\varepsilon\delta'}\right)^{1+\varepsilon/2}$, which only depends on $\delta$.
This is the first estimate.

Next we prove the second estimate.
Set $L_{\gamma}(r,\theta )=\ell_{\eta}(f(\pi_f^{-1}(\gamma_{r,\theta})))$.
Then
$$
L_{\gamma}(r,\theta )=\int_{\pi_f^{-1}(\gamma_{r,\theta})}\varphi d|\pi_f|.$$
Since $|\pi_f|>\frac{1}{2}$ on $\pi_f^{-1}(\gamma_{r,\theta})$, we have
$$
L_{\gamma}(r,\theta )\leq 2\int_{\pi_f^{-1}(\gamma_{r,\theta})}\varphi |\pi_f|d|\pi_f|.
$$
By the Cauchy-Schwarz inequality, we get
$$
L_{\gamma}(r,\theta )^2\leq 4(\deg \pi_f)\int_{\pi_f^{-1}(\gamma_{r,\theta})}\varphi^2 |\pi_f|d|\pi_f|,
$$
hence
$$
\int_0^{2\pi}L_{\gamma}(r,\theta )^2d\theta \leq 4(\deg \pi_f)A(r)\leq 4(\deg \pi_f)^2T'(r).
$$
By Lemma \ref{lem:growth}, letting $\varepsilon=1$, we have
$$
\int_0^{2\pi}L_{\gamma}(r,\theta )^2d\theta \leq \frac{8(\deg \pi_f)^2}{\delta}\max\{ 1,T(r)^2\}
$$
for all $r\in (s,1)\backslash E_2$ with $|E_2|< \delta$.
Then, for each $r\in (s,1)\backslash E_2$, we get
$$
L_{\gamma}(r,\theta )^2 \leq \frac{8(\deg \pi_f)^2}{\delta^2}\max\{ 1,T(r)^2\}
$$
for all $\theta \in (0,2\pi )$ outside some exceptional set $E_3$ with $|E_3|<\delta$.
Hence
$$
L_{\gamma}(r,\theta )\leq \frac{2\sqrt{2}\deg \pi_f}{\delta}\max\{ 1,T(r)\}
$$
for all $\theta \in (0,2\pi )\backslash E_3$.
Thus we obtain the second estimate by letting $c_2=2\sqrt{2}/\delta$, which only depends on $\delta$.
\hspace{\fill} $\square$

\begin{lem}\label{lem:growth}
Let $g$ be a continuously differentiable, increasing function on $[s,1 )$ with $g(s)\geq 0$.
Let $\delta >0$ and $0<\varepsilon\leq 1$.
Then we have 
$$g'(r)\leq \frac{2}{\varepsilon\delta}\max \{ 1,g(r)^{1+\varepsilon}\}$$ 
for all $r\in (s,1)$ outside a set $E$ with $|E|<\delta$.
\end{lem}

{\it Proof.}\
Set
$$
E=\left\{ r\in (s,1); g'(r)> \frac{2}{\varepsilon\delta}\max \{ 1,g(r)^{1+\varepsilon}\}  \right\} .
$$
If $E=\emptyset$, then our assertion is trivial.
Suppose $E\not=\emptyset$.
We have
$$
|E|< \frac{\varepsilon\delta}{2}\int_{E}\frac{g'(r)}{\max\{ 1,g(r)^{1+\varepsilon}\}}dr\leq \frac{\varepsilon\delta}{2}\int_{s}^{1}\frac{g'(r)}{\max\{ 1,g(r)^{1+\varepsilon}\}}dr.
$$
We have the following three cases.

{\it Case 1:}\ $g(r)\geq 1$ for all $r\in [s,1)$.
Then we have
$$
\int_{s}^{1}\frac{g'(r)}{\max\{ 1,g(r)^{1+\varepsilon}\}}dr=\int_{s}^{1}\frac{g'(r)}{g(r)^{1+\varepsilon}}dr= \lim_{t\to 1-0}\left[ \frac{-1}{\varepsilon g(r)^{\varepsilon}}\right]_{s}^{t}\leq \frac{1}{\varepsilon}.
$$

{\it Case 2:}\ $g(r)\leq 1$ for all $r\in [s,1)$.
Then we have
$$
\int_{s}^{1}\frac{g'(r)}{\max\{ 1,g(r)^{1+\varepsilon}\}}dr=\int_s^{1}g'(r)dr\leq 1.
$$

{\it Case 3:}\ Otherwise, we have $g(s)<1$ and $\lim_{r\to 1-0}g(r)>1$.
We set $\kappa =\sup\{ r\in [s,1); \ g(r)\leq 1\}$.
Then we have $s<\kappa <1$ and $g(\kappa )=1$.
Hence we have
$$
\int_{s}^{1}\frac{g'(r)}{\max\{ 1,g(r)^{1+\varepsilon}\}}dr=\int_s^{\kappa}g'(r)dr+\int_{\kappa}^{1}\frac{g'(r)}{g(r)^{1+\varepsilon}}dr\leq 1+\lim_{t\to 1-0}\left[ \frac{-1}{\varepsilon g(r)^{\varepsilon}}\right]_{\kappa}^{t}\leq \frac{2}{\varepsilon}.
$$

Thus in all cases, we have proved $|E|<\delta$.
\hspace{\fill} $\square$

\subsection{Estimate for the second term of RHS of \eqref{eqn:202302081}}
Let $A$ be a semi-abelian variety with a smooth projective equivariant compactification $\overline{A}$.
Let $S$ be a smooth projective variety.
Let $Z\subset S_{1,A}$ be an irreducible Zariski closed set.
Assume that $Z$ is horizontally integrable.
We recall $W\subset S$ from Remark \ref{rem:202108051}.

\begin{lem}\label{lem:202302203}
Let $D\subset \overline{A}$ be an irreducible component of $\partial A$ with a Weil function $\lambda_D\geq 0$.
Let $\mathcal{F}\subset \mathrm{Hol}(\mathbb D,A\times S)$ be an infinite set of non-constant holomorphic maps such that $( f_{S_{1,A}})_{f\in \mathcal{F}}\Rightarrow Z$.
Let $\omega_S$ be a smooth positive $(1,1)$ form on $S$ and let $\lambda_W\geq 0$ be a Weil function for $W$.
Let $s\in (1/2,1)$ and $\delta>0$.
Then there exists a positive constant $\alpha>0$ such that for all $f\in \mathcal{F}$ with $f_S(\mathbb D)\not\subset W$, we have
\begin{equation*}
|m(r,f_A,\lambda_D)-m((s+r)/2,f_A,\lambda_D)|\leq \alpha T_{s}(r,f_S,\omega_S)+\alpha m((s+r)/2,f_S,\lambda_W)+\alpha
\end{equation*}
for all $r\in (s,1)\backslash E$, where $E\subset (s,1)$ is an exceptional set with $|E|<\delta$.
\end{lem}

{\it Proof.}\
We recall the objects fixed in Remark \ref{rem:202108051}.
Set $\partial U=\Sigma-U$ with a Weil function $\lambda_{\partial U}\geq 0$.
Then by Lemma \ref{lem:202302021} (1), we get
\begin{equation}\label{eqn:202105291}
\bigcup_{\Delta}\bar{q}(\Delta)\subset W,
\end{equation}
where $\Delta$ runs over all irreducible components $\Delta\subset \partial U$ such that $\Delta\cap \overline{\Theta}\not=\emptyset$.
Here $\bar{q}:\Sigma\to S$ is the extension of $q:U\to S$.
Let $K\subset \Sigma$ be a compact neighbourhood of $\overline{\Theta}\subset \Sigma$.
We assume that $K\subset \Sigma^o$ and
\begin{equation}\label{eqn:202302164}
\Delta'\cap K=\emptyset
\end{equation}
for all irreducible components $\Delta'\subset \partial U$ such that $\Delta'\cap\overline{\Theta}=\emptyset$.
By \eqref{eqn:202105291} and \eqref{eqn:202302164}, there exists a positive constant $\alpha_1>0$ suct that
\begin{equation}\label{eqn:202302165}
\lambda_{\partial U}(x)\leq \alpha_1\lambda_W(q(x))+\alpha_1
\end{equation}
for all $x\in K$.

Since $\overline{A}$ is an equivariant compactification, the isomorphism $\Phi:A\times U\to A\times U$ extends to a morphism
$$
\overline{\Phi}:\overline{A}\times U\to\overline{A}\times U
$$
by the definition \eqref{eqn:202302248}.
Then the inverse $\overline{\Phi}^{-1}:\overline{A}\times U\to\overline{A}\times U$ induces a rational map
$$
\overline{\Phi}^{-1}:\overline{A}\times \Sigma\dashrightarrow\overline{A}\times \Sigma,
$$
which is holomorphic over $\overline{A}\times U\subset \overline{A}\times \Sigma$.
Let $\nu_1:\overline{A}\times \Sigma\to\overline{A}$ be the first projection and let $\nu_2:\overline{A}\times \Sigma\to\Sigma$ be the second projection.
We claim that there exists a positive constant $\alpha_2>0$ such that 
\begin{equation}\label{eqn:202302166}
|\lambda_{D}(\nu_1(x))-\lambda_{D}(\nu_1\circ\Phi^{-1}(x))|\leq \alpha_2\lambda_{\partial U}(\nu_2(x))+\alpha_2
\end{equation}
for all $x\in A\times U\subset \overline{A}\times \Sigma$.
Indeed, since $D\subset \overline{A}$ is $A$-invariant, we have $\overline{\Phi}(D\times U)=D\times U$.
Hence we have 
$$\left(\nu_1\circ\overline{\Phi}^{-1}\right)^{-1}(D)\vert_{\overline{A}\times U}=\nu_1^{-1}(D)\vert_{\overline{A}\times U}$$ 
over $\overline{A}\times U$.
Thus by \cite[Prop. 2.2.9 (7)]{Y2}, we get \eqref{eqn:202302166}.

Now we are given an infinite subset $\mathcal{F}\subset \mathrm{Hol}(\mathbb D,A\times S)$ of non-constant holomorphic maps so that $(f_{S_{1,A}})_{f\in\mathcal{F}}\Rightarrow Z$.
We apply Lemma \ref{lem:20230218} to get a finite subset $\mathcal{E}\subset \mathcal{F}$ such that for each $f\in\mathcal{F}\backslash\mathcal{E}$, we get a lifting $\hat{f}\in \mathrm{Hol_m}(\mathbb D,A\times \Sigma)$ of $f\in \mathrm{Hol}(\mathbb D,A\times S)$, a connected open set $\Omega_f\subset \mathbb D$ and a connected component $\Omega_f'\subset Y_{\hat{f}}$ of $\pi_{\hat{f}}^{-1}(\Omega_f)$ with the properties \eqref{eqn:20210504?}--\eqref{eqn:202302034?} described in Lemma \ref{lem:20230218}.

Let $s\in (1/2,1)$, $\delta>0$ be given.
Set
$$
\delta'=\delta/11.
$$
We have (cf. \eqref{eqn:20211024?})
\begin{equation*}
( (\hat{f}_{\Sigma},\Omega_f,\Omega_f'))_{f\in \mathcal{F}\backslash\mathcal{E}} \rightsquigarrow\overline{\Theta}.
\end{equation*}
Hence, there exists a finite subset $\mathcal{E}_1\subset \mathcal{F}$ with $\mathcal{E}\subset \mathcal{E}_1$ such that for all $f\in \mathcal{F}\backslash \mathcal{E}_1$, we have $z\in \Omega_f$ for $\delta'$-almost all $z\in \mathbb D(1-\delta')$ and 
\begin{equation}\label{eqn:202302163}
\hat{f}_{\Sigma}(\Omega_f')\subset K.
\end{equation}

Let $\omega_{\Sigma}$ be a smooth positive $(1,1)$-form on $\Sigma$ such that 
\begin{equation}\label{eqn:202302249}
\omega_{\Sigma}\leq \eta
\end{equation}
on $K$.
Here we recall $\eta$ from Remark \ref{rem:202108051}.
We take an open neighbourhood $O\subset \Sigma_{1,A}$ of $\Sigma_{1,A}^*$ such that
\begin{equation}\label{eqn:20230202}
|v_{\mathrm{Lie}A}|_{\omega_A}\leq | v_{\Sigma}|_{\omega_{\Sigma}}
\end{equation}
for all $(v_{\Sigma},v_{\mathrm{Lie}A})\in T\Sigma\times \mathrm{Lie}A$ with $[(v_{\Sigma},v_{\mathrm{Lie}A})]\in O\subset \Sigma_{1,A}$.
We have (cf. \eqref{eqn:202102261?})

\begin{equation*}
( ((\Phi\circ \hat{f})_{\Sigma_{1,A}},\Omega_f,\Omega_f'))_{f\in \mathcal{F}\backslash\mathcal{E}} \rightsquigarrow \Sigma_{1,A}^*.
\end{equation*}
Hence, we may take a finite subset $\mathcal{E}_2\subset \mathcal{F}$ with $\mathcal{E}_1\subset \mathcal{E}_2$ such that for all $f\in\mathcal{F}\backslash\mathcal{E}_2$, we have 
\begin{equation}\label{eqn:202302201}
(\Phi\circ \hat{f})_{\Sigma_{1,A}}(\Omega_f')\subset O.
\end{equation}

Let $f\in\mathcal{G}\backslash\mathcal{E}_2$ with $f_S(\mathbb D)\not\subset W$.
Then by \eqref{eqn:202302165} and \eqref{eqn:202302163}, we have
\begin{equation}\label{eqn:202302044}
m(r,(\hat{f}_{\Sigma},\Omega_f,\Omega_f'),\lambda_{\partial U})\leq \alpha_1m(r,f_S,\lambda_W)+\alpha_1
\end{equation}
for all $r\in (0,1)$.
We set
$$
E_1=\{r\in (s,1-\delta'); (\mathbb D\backslash \Omega_f)\cap\partial\mathbb D(r)\not=\emptyset\}\cup \{r\in (s,1-\delta'); (\mathbb D\backslash \Omega_f)\cap\partial\mathbb D(\sigma)\not=\emptyset\},
$$
where $\sigma=(s+r)/2$.
By adding finite points to $E_1$, we may assume that if $r\in (s,1-\delta')\backslash E_1$, then $\partial \mathbb D(r)\cup \partial \mathbb D(\sigma)$ is contained in $\Omega_f$ and does not contain the critical values of $\pi_{\hat{f}}:Y_{\hat{f}}\to\mathbb D$.
Then we have
$$
|E_1|\leq 6\delta'.
$$
For all $r\in(s,1-\delta')\backslash E_1$, we have
$$
m(r,f_{A},\lambda_D)=m(r,(\hat{f}_A,\Omega_f,\Omega_f'),\lambda_D)
$$
and
$$
m(\sigma,f_{A},\lambda_D)=m(\sigma,(\hat{f}_A,\Omega_f,\Omega_f'),\lambda_D).
$$
By \eqref{eqn:202302166} and \eqref{eqn:202302044}, we get
$$
|m(r,((\Phi\circ \hat{f})_A,\Omega_f,\Omega_f'),\lambda_D)-m(r,(\hat{f}_A,\Omega_f,\Omega_f'),\lambda_D)|\leq \alpha_3m(r,f_S,\lambda_{W})+\alpha_3,
$$
where $\alpha_3=\alpha_2\alpha_1+\alpha_2$.
Similarly, we have
$$
|m(\sigma,((\Phi\circ \hat{f})_A,\Omega_f,\Omega_f'),\lambda_D)-m(\sigma,(\hat{f}_A,\Omega_f,\Omega_f'),\lambda_D)|\leq \alpha_3m(\sigma,f_S,\lambda_{W})+\alpha_3.
$$
Thus for $r\in (s,1-\delta')\backslash E_1$, we have
\begin{align*}
|m(r,f_A,\lambda_D)-m(\sigma,f_A,\lambda_D)|&=
|m(r,(\hat{f}_A,\Omega_f,\Omega_f'),\lambda_D)-m(\sigma,(\hat{f}_A,\Omega_f,\Omega_f'),\lambda_D)|
\\
&\leq |m(r,((\Phi\circ \hat{f})_A,\Omega_f,\Omega_f'),\lambda_D)-m(\sigma,((\Phi\circ \hat{f})_A,\Omega_f,\Omega_f'),\lambda_D)|
\\
&\qquad
+\alpha_3m(r,f_S,\lambda_{W})+\alpha_3m(\sigma,f_S,\lambda_{W})+2\alpha_3.
\end{align*}
Hence
\begin{multline}\label{eqn:202302195}
|m(r,f_A,\lambda_D)-m(\sigma,f_A,\lambda_D)|
\leq |m(r,((\Phi\circ \hat{f})_A,\Omega_f,\Omega_f'),\lambda_D)-m(\sigma,((\Phi\circ \hat{f})_A,\Omega_f,\Omega_f'),\lambda_D)|
\\
+\alpha_3m(r,f_S,\lambda_{W})+\alpha_3m(\sigma,f_S,\lambda_{W})+2\alpha_3.
\end{multline}

Next we claim
\begin{equation}\label{eqn:202302192}
|m(r,((\Phi\circ \hat{f})_A,\Omega_f,\Omega_f'),\lambda_D)-m(\sigma,((\Phi\circ \hat{f})_A,\Omega_f,\Omega_f'),\lambda_D)|\leq \alpha_4T_{s}(r,\hat{f}_{\Sigma},\eta)+\alpha_4
\end{equation}
for all $r\in (s,1)\backslash E_2$, where $E_2\subset (0,1)$ is an exceptional set with $|E_2|<11\delta'$.
Here $\alpha_4>0$ is a positive constant which does not depend on the choice of $f\in\mathcal{F}\backslash\mathcal{E}_2$.

We prove this.
Set $\gamma_{r,\theta}: te^{i\theta}, 1/2\leq t\leq r$.
If $\partial\mathbb D(\sigma)+\gamma_{r,\theta}+\partial\mathbb D(r)\subset \Omega_f$, then by \eqref{eqn:20230202} and \eqref{eqn:202302201}, we have
$$
\ell_{\omega_{A}}((\Phi\circ\hat{f})_{\overline{A}}(\Omega_f'\cap\pi_{\hat{f}}^{-1}(\partial\mathbb D(r)\cup\partial\mathbb D(\sigma)\cup\gamma_{r,\theta})))\leq 
\ell_{\omega_{\Sigma}}(\hat{f}_{\Sigma}(\Omega_f'\cap\pi_{\hat{f}}^{-1}(\partial\mathbb D(r)\cup\partial\mathbb D(\sigma)\cup\gamma_{r,\theta}))).
$$
By \eqref{eqn:202302163} and \eqref{eqn:202302249}, we have
$$
\ell_{\omega_{\Sigma}}(\hat{f}_{\Sigma}(\Omega_f'\cap\pi_{\hat{f}}^{-1}(\partial\mathbb D(r)\cup\partial\mathbb D(\sigma)\cup\gamma_{r,\theta})))
\leq \ell_{\eta}(\hat{f}_{\Sigma}(\Omega_f'\cap\pi_{\hat{f}}^{-1}(\partial\mathbb D(r)\cup\partial\mathbb D(\sigma)\cup\gamma_{r,\theta}))).
$$
Thus by Lemma \ref{lem:20230219}, we get
\begin{multline}\label{eqn:202104232}
|m(r,((\Phi\circ\hat{f})_{\overline{A}},\Omega_f,\Omega_f'),\lambda_D)-m(\sigma,((\Phi\circ\hat{f})_{\overline{A}},\Omega_f,\Omega_f'),\lambda_D)|
\\
\leq \frac{c}{\deg(\pi_{\hat{f}}|_{\Omega_{\hat{f}}'})}\ell_{\eta}(\hat{f}_{\Sigma}(\pi_{\hat{f}}^{-1}(\partial\mathbb D(r)\cup\partial\mathbb D(\sigma)\cup\gamma_{r,\theta})))+c,
\end{multline}
provided that $\partial\mathbb D(r)\cup \partial\mathbb D(\sigma)\cup \gamma_{r,\theta}\subset \Omega_f$ does not contain the critical values of $\pi_{\hat{f}}:Y_{\hat{f}}\to\mathbb D$.
Here $c>0$ is a positive constant which appears in Lemma \ref{lem:20230219}, hence independent of the choice of $f\in\mathcal{F}\backslash\mathcal{E}_2$, $r$ and $\theta$.

Now we apply Lemma \ref{lem:20210423} to get
$$
\ell_{\eta}(\hat{f}_{\Sigma}(\pi_{\hat{f}}^{-1}(\partial\mathbb D(r))))\leq \alpha_5 (\deg \pi_{\hat{f}})^2T_s(r,\hat{f}_{\Sigma},\eta)+\alpha_5(\deg \pi_{\hat{f}})^2,
$$
for all $r\in (s,1)$ outside some exceptional set $E_3\subset (s,1)$ with 
$$|E_3|<\delta'.$$
Here $\alpha_5>0$ only depends on $\delta'$.
We define $E_4\subset (s,1)$ by $r\in E_4$ iff $(s+r)/2\in E_3$.
We have 
$$|E_4|\leq 2| E_3|\leq 2\delta'.$$
Then for $r\in (s,1-\delta')\backslash (E_1\cup E_4)$, we have
\begin{equation*}
\ell_{\eta}(\hat{f}_{\Sigma}(\pi_{\hat{f}}^{-1}(\partial\mathbb D(\sigma))))\leq \alpha_5(\deg\pi_{\hat{f}})^2T_s(\sigma,\hat{f}_{\Sigma},\eta)+\alpha_5(\deg\pi_{\hat{f}})^2,
\end{equation*}
and $\partial\mathbb D(\sigma)\subset \Omega_f$.
Thus for $r\in (s,1-\delta')\backslash (E_1\cup E_3\cup E_4)$, the estimate \eqref{eqn:202104232} yields
\begin{multline*}
|m(r,((\Phi\circ\hat{f})_{\overline{A}},\Omega_f,\Omega_f'),\lambda_D)-m(\sigma,((\Phi\circ\hat{f})_{\overline{A}},\Omega_f,\Omega_f'),\lambda_D)|
\\
\leq
2c\alpha_5(\deg\pi_{\hat{f}})T_s(r,\hat{f}_{\Sigma},\eta)+2c\alpha_5(\deg\pi_{\hat{f}})+c+
\frac{c}{\deg(\pi_{\hat{f}}|_{\Omega_{\hat{f}}'})}\ell_{\eta}(\hat{f}_{\Sigma}(\pi_{\hat{f}}^{-1}(\gamma_{r,\theta})))
\end{multline*}
provided $\gamma_{r,\theta}$ is contained in $\Omega_f$ and does not contain the critical values of $\pi_{\hat{f}}:Y_{\hat{f}}\to\mathbb D$.
By Lemma \ref{lem:20210423} (2), there exists $E_5\subset (s,1)$ with 
$$|E_5|<\delta'$$ 
such that for each $r\in (s,1-\delta')\backslash E_5$, we may choose $\theta\in (0,2\pi)$ such that $\gamma_{r,\theta}\subset \Omega_f$, $\gamma_{r,\theta}$ does not contain the critical values of $\pi_{\hat{f}}:Y_{\hat{f}}\to\mathbb D$, and 
$$
\ell_{\eta}(\hat{f}_{\Sigma}(\pi_{\hat{f}}^{-1}(\gamma_{r,\theta})))\leq \alpha_6 (\deg \pi_{\hat{f}})T_s(r,\hat{f}_{\Sigma
},\eta)+\alpha_6 (\deg \pi_{\hat{f}}).
$$
Here $\alpha_6>0$ only depends on $\delta'$.
Hence, we get
\begin{multline*}
|m(r,((\Phi\circ\hat{f})_{\overline{A}},\Omega_f,\Omega_f'),\lambda_D)-m(\sigma,((\Phi\circ\hat{f})_{\overline{A}},\Omega_f,\Omega_f'),\lambda_D)|
\\
\leq c(2\alpha_5+\alpha_6)(\deg\pi_{\hat{f}})T_s(r,\hat{f}_{\Sigma},\eta)+c(2\alpha_5+\alpha_6)(\deg\pi_{\hat{f}})+c
\end{multline*}
for all $r\in (s,1)\backslash E_2$, where 
\begin{equation}\label{eqn:2023022410}
E_2=E_1\cup E_3\cup E_4\cup E_5\cup (1-\delta',1).
\end{equation}
We have
$$
|E_2|<11\delta'=\delta.
$$
This conclude the proof for \eqref{eqn:202302192}, where we set $\alpha_4=c(2\alpha_5+\alpha_6)[\mathbb C(S):\mathbb C(\Sigma)]+c$.

Now by \eqref{eqn:202302195} and \eqref{eqn:202302192}, we get
$$
|m(r,f_A,\lambda_D)-m(\sigma,f_A,\lambda_D)|\leq \alpha_7T_{s}(r,\hat{f}_{\Sigma},\eta)+\alpha_7m(r,f_S,\lambda_{W})+\alpha_7m(\sigma,f_S,\lambda_{W})+\alpha_7
$$
for all $r\in (s,1)\backslash E_2$, where $\alpha_7=2\alpha_3+\alpha_4$.
Here we note that $E_1\cup(1-\delta',1)\subset E_2$ (cf. \eqref{eqn:2023022410}).
By Lemma \ref{lem:202302021}, we have
$$
T_s(r,\hat{f}_{\Sigma},\eta)\leq \alpha_8T_s(r,\bar{q}\circ \hat{f}_{\Sigma},\omega_{S})+\alpha_8m(\sigma,f_S,\lambda_W)+\alpha_8
$$
for all $r\in (1/2,1)$.
Here $\alpha_8>0$ is a positive constant which does not depend on the choice of $f\in\mathcal{F}\backslash\mathcal{E}_2$.
Hence we get
$$
|m(r,f_{A},\lambda_D)-m(\sigma,f_A,\lambda_D)|
\leq \alpha_9T_{s}(r,f_S,\omega_S)+\alpha_9m(r,f_S,\lambda_{W})+\alpha_9m(\sigma,f_S,\lambda_{W})+\alpha_9
$$
for all $r\in (s,1)\backslash E_2$, where $\alpha_9=\alpha_7\alpha_8+\alpha_7$.
By Lemma \ref{lem:20211008}, we get 
$$
m(r,f_S,\lambda_W)\leq \alpha_{10}T_{\sigma}(r,f_S,\omega_S)+m(\sigma,f_S,\lambda_W)+\alpha_{10}
$$
for all $r\in (1/2,1)$.
Here $\alpha_{10}>0$ is a positive constant which does not depend on the choice of $f\in\mathcal{F}\backslash\mathcal{E}_2$.
Hence we get
\begin{equation}\label{eqn:202302196}
|m(r,f_A,\lambda_D)-m(\sigma,f_A,\lambda_D)|\leq \alpha T_{s}(r,f_S,\omega_S)+\alpha m(\sigma,f_S,\lambda_{W})+\alpha
\end{equation}
for all $r\in (s,1)\backslash E_2$, where $\alpha=\max\{2\alpha_9,\alpha_9\alpha_{10}+\alpha_9\}$.
This proves the desired estimate for $f\in \mathcal{F}\backslash \mathcal{E}_2$ with $f_S(\mathbb D)\not\subset W$.

Finally we enlarge $\alpha>0$ so that
$$
2\sup_{r\in (s,1-\delta)}m(r,f,\lambda_D)\leq \alpha
$$
for all $f\in \mathcal{E}_2$.
Then \eqref{eqn:202302196} is valid for all $f\in\mathcal{F}$ with $f_S(\mathbb D)\not\subset W$.
\hspace{\fill} $\square$

\subsection{Proof of Proposition \ref{pro:mpro}}\label{subsec:proofmpro}
We take $W\subset S$ as in Remark \ref{rem:202108051}.
Then we have $W\subsetneqq\tau(Z)$ (cf. Remark \ref{rem:202108051}).
Let $\mathcal{F}\subset\mathrm{Hol}(\mathbb D,A\times S)$, $\omega_{\overline{A}}$, $\omega_S$, $\lambda_W$, $s\in (1/2,1)$, $\delta>0$ be given as in Proposition \ref{pro:mpro}.
Let $\rho:A\to A_0$ be the canonical quotient (cf. \eqref{eqn:20230227}) and let $\omega_{A_0}$ be an invariant positive $(1,1)$ form on $A_0$.
Note that we are given a smooth projective equivariant compactification $\overline{A}$.
For an irreducible component $D\subset\overline{A}$ of $\partial A$, let $\lambda_D\geq 0$ be a Weil function.

Now let $f\in\mathcal{F}$ such that $f_S(\mathbb D)\not\subset W$.
By Lemma \ref{lem:202302022}, we get
\begin{equation*}
T_{s}(r,f_A,\omega_{\overline{A}})
\leq c'T_{s}(r,\rho\circ f_A,\omega_{A_0})+c'\sum_{D\subset \partial A}|m(r,f_A,\lambda_D)-m(\sigma,f_A,\lambda_D)|+c'
\end{equation*}
for all $r\in (s,1)$,
where $c'>0$ is independent of the choice of $f\in\mathcal{F}$.
By Lemma \ref{lem:202302193}, we get 
\begin{equation*}
T_s(r,\rho\circ f_{A},\omega_{A_0}) \leq 
\alpha T_s(r,f_{S},\omega_{S})
+
\alpha
m((s+r)/2,f_S,\lambda_{W})
+\alpha
\end{equation*}
for all $r\in (s,1)\backslash E'$ with $|E'|<\delta/2$.
Here $\alpha>0$ is independent of the choice of $f\in\mathcal{F}$.
Let $D_1,\ldots,D_l$ be the irreducible components of $\partial A$.
By Lemma \ref{lem:202302203}, we get
\begin{equation*}
|m(r,f_A,\lambda_{D_i})-m((s+r)/2,f_A,\lambda_{D_i})|\leq \beta_i T_{s}(r,f_S,\omega_S)+\beta_i m((s+r)/2,f_S,\lambda_W)+\beta_i
\end{equation*}
for all $r\in (s,1)\backslash E_i$, where $E_i\subset (s,1)$ is an exceptional set with $|E_i|<\delta/2l$.
Here $\beta_i>0$ is independent of the choice of $f\in\mathcal{F}$.
We set $c_1=c_2=c'(\alpha+\sum_{i=1}^l\beta_i)$, $c_3=c'+c'(\alpha+\sum_{i=1}^l\beta_i)$, and $E=E'\cup \cup_{i=1}^l E_i$.
Then $|E|<\delta$.
This shows \eqref{eqn:20230208}, where $c_1>0$, $c_2>0$, $c_3>0$ are independent of the choice of $f\in\mathcal{F}$.
This completes the proof.
\hspace{\fill} $\square$

\section{Application of logarithmic tautological inequality}

\subsection{Logarithmic tautological inequality}
Let $X$ be a smooth projective variety with a smooth, positive (1,1)-form $\omega_X$.
Let $D\subset X$ be a simple normal crossing divisor.
We set $\overline{T}X(-\log D)=P(TX(-\log D)\oplus \mathcal{O}_X)$, which is a smooth compactification of $TX(-\log D)$.
Let $\partial TX(-\log D)$ be the Cartier divisor on the boundary which  corresponds to a section of $\mathcal{O}_{\overline{T}X(-\log D)}(1)$.
If $f\in \mathrm{Hol}(\mathbb D,X)$ is holomorphic with $f(\mathbb D)\not\subset D$, then the derivative of $f$ induces a holomorphic map $f':\mathbb D\to \overline{T}X(-\log D)$.

\begin{lem}\label{lem:llld}
Let $\varepsilon >0$, $\delta>0$ and $s\in (0,1)$.
Let $\lambda_{\partial TX(-\log D)}$ and $\lambda_D$ be Weil functions for $\partial TX(-\log D)\subset \overline{T}X(-\log D)$ and $D\subset X$, respectively.
Then there exists a positive constant $\mu>0$ such that for all $f\in \mathrm{Hol}(\mathbb D,X)$ with $f(\mathbb D)\not\subset D$, we have
$$
m(r,f',\lambda_{\partial TX(-\log D)})\leq \varepsilon T_s(r,f,\omega_X)+
\varepsilon m(r,f,\lambda_D)+
\mu
$$
for all $r\in (s,1)$ outside some exceptional set of linear measure less than $\delta$.
\end{lem}

This lemma is a variant of the estimate for entire curves $f:\mathbb C\to X$ due to R.~Kobayashi \cite{Kjapan} and McQuillan \cite{M}.
We refer Vojta \cite[Thm A.2]{Vojta0} for the precise statement and simplified proof.
We remark that \cite[Thm A.2]{Vojta0} implies classical Nevanlinna's lemma on logarithmic derivatives when applied to entire curves $f:\mathbb C\to\mathbb P^1$ and $D=(0)+(\infty)$, while Vojta's proof of \cite[Thm A.2]{Vojta0} is based on Nevanlinna's lemma.
In the following, we follow another method described in \cite[Thm 4.8]{Ykob}, which does not use Nevanlinna's lemma.
See also \cite[Sec. 2]{Brunella1}.

\medskip

{\it Proof of Lemma \ref{lem:llld}.}\
The proof divides into three steps.
The following proof is similar to the proof of \cite[Thm 4.8]{Ykob}.
\par

{\it Step 1.}\
Let $\overline{T}X=P(TX\oplus \mathcal{O}_X)$ be a smooth compactification of $TX$ and let $\partial TX$ be the Cartier divisor on the boundary which  corresponds to a section of $\mathcal{O}_{\overline{T}X}(1)$.
Let $\lambda_{\partial TX}$ be the Weil function of $\partial TX$ defined by $\lambda_{\partial TX}(v)=\log \sqrt{1+|v|^2_{\omega_X}}$ for $v\in TX$, where $|\cdot|_{\omega_{X}}$ is a norm on $TX$ defined by $\omega_{X}$.
We prove the following estimate for all $g\in \mathrm{Hol_m}(\mathbb D,X)$ with $g(Y_g)\not\subset D$:
\begin{equation}\label{eqn:loglld11}
m(r,g',\lambda_{\partial TX(-\log D)})\leq m(r,g,\lambda_D)+m(r,g',\lambda_{\partial TX})+\mu_1
\end{equation}
for all $r\in(0,1)$.
Here $\mu_1>0$ is a positive constant which only depends on the choices of Weil functions.

We prove this.
The natural morphism $\iota_1 :TX(-\log D)\to TX$ induces a birational map 
$$\psi :\overline{T}X\dashrightarrow \overline{T}X(-\log D).$$
Let $Z\subset \overline{T}X$ be the indeterminacy locus of $\psi$.
Let $p:\overline{T}X\to X$ be the projection.
Then we have
\begin{equation}\label{eqn:madamada}
(\psi |_{\overline{T}X\backslash Z})^*\partial TX(-\log D)=(p^*D+\partial TX)|_{\overline{T}X\backslash Z}.
\end{equation}
Let $\alpha :\widetilde{\overline{T}X}\to \overline{T}X$ be a modification such that $\psi$ induces a morphism $\widetilde{\psi}:\widetilde{\overline{T}X}\to \overline{T}X(-\log D)$.
Then there exists an effective Cartier divisor $E\subset \widetilde{\overline{T}X}$ such that 
\begin{equation}\label{eqn:mouii}
\widetilde{\psi}^*\partial TX(-\log D)=\alpha^*(p^*D+\partial TX)-E.
\end{equation}
Indeed, by \eqref{eqn:madamada}, we have
\begin{equation}\label{eqn:20220601}
(\psi |_{\overline{T}X\backslash Z})^*\mathcal{O}_{\overline{T}X(-\log D)}(1)=( p^*\mathcal{O}_X(D)\otimes \mathcal{O}_{\overline{T}X}(1))|_{\overline{T}X\backslash Z}.
\end{equation}
Let $p^{\mathrm{log}}:\overline{T}X(-\log D)\to X$ be the projection.
There exists a natural surjection
\begin{equation}\label{eqn:202206011}
(p^{\mathrm{log}})^*(p^{\mathrm{log}})_*\mathcal{O}_{\overline{T}X(-\log D)}(1)\to \mathcal{O}_{\overline{T}X(-\log D)}(1).
\end{equation}
We pull-back this by $\psi|_{\overline{T}X\backslash Z}:\overline{T}X\backslash Z\to \overline{T}X(-\log D)$ and combine with \eqref{eqn:20220601}.
Then we have the following surjection on $\overline{T}X\backslash Z$:
$$
p^*(p^{\mathrm{log}})_*\mathcal{O}_{\overline{T}X(-\log D)}(1)|_{\overline{T}X\backslash Z}\to p^*\mathcal{O}_X(D)\otimes  \mathcal{O}_{\overline{T}X}(1)|_{\overline{T}X\backslash Z}.
$$
Since $Z$ has codimension greater than one, this morphism extends over whole $\overline{T}X$.
We denote by $K$ the kernel of the following surjection obtained from \eqref{eqn:202206011}:
$$
\widetilde{\psi}^*(p^{\mathrm{log}})^*(p^{\mathrm{log}})_*\mathcal{O}_{\overline{T}X(-\log D)}(1)\to \widetilde{\psi}^*\mathcal{O}_{\overline{T}X(-\log D)}(1).
$$
Then the composition of
$$
K\to \alpha^*p^*(p^{\mathrm{log}})_*\mathcal{O}_{\overline{T}X(-\log D)}(1)\to \alpha^*\left( p^*\mathcal{O}_X(D)\otimes \mathcal{O}_{\overline{T}X}(1)  \right)
$$
is a zero map on $\alpha^{-1}(\overline{T}X\backslash Z)$, hence on $\widetilde{\overline{T}X}$.
Hence we get a morphism
$$
\widetilde{\psi}^*\mathcal{O}_{\overline{T}X(-\log D)}(1)\to \alpha^*\left( p^*\mathcal{O}_X(D)\otimes \mathcal{O}_{\overline{T}X}(1)  \right).
$$
Hence there exists an effective Cartier divisor $E$ on $\widetilde{\overline{T}X}$ such that \eqref{eqn:mouii} is valid.
Hence, we get \eqref{eqn:loglld11}.

\medskip

{\it Step 2.}\
We estimate $m(r,g',\lambda_{\partial TX})$, where
$$
m(r,g',\lambda_{\partial TX})=
\frac{1}{\deg \pi_g}\int_{y\in \partial Y_g(r)}\log\sqrt{1+ | g'(y)|^2_{\omega_{X}}}\frac{d\arg \pi_g(y)}{2\pi}
$$
Using concavity of $\log$, we have
\begin{equation*}
m(r,g',\lambda_{\partial TX})
\leq \frac{1}{2}\log \left( 1+\frac{1}{\deg \pi _g}\int_{y\in \partial Y_g} | g'(y)|^2_{\omega_{X}}\frac{d\arg \pi _g(y)}{2\pi}\right).
\end{equation*}
We set
$$
\tau(r)=\frac{1}{\deg\pi_g}\int_s^rdt\int_{Y_g(t)}g^*\omega_{X}.
$$
Then we have
$$
\frac{1}{2\pi r}\frac{d^2}{dr^2}\tau(r) = \frac{1}{\deg\pi_g}\int_{y\in \partial Y_g(r)} | g'(y)|^2_{\omega_{X}}\frac{d\arg\pi_g(y)}{2\pi}.
$$
Hence
$$
m(r,g',\lambda_{\partial TX})\leq  \frac{1}{2}\log \left( 1+\frac{1}{2\pi r}\frac{d^2}{dr^2}\tau(r)
\right) . 
$$
Hence for $r>s$, we have
$$
m(r,g',\lambda_{\partial TX})\leq  \frac{1}{2}\log \left( 1+\frac{1}{2\pi s}\frac{d^2}{dr^2}\tau(r)
\right) . 
$$
Now we apply Lemma \ref{lem:growth} twice.
We have
\begin{equation*}
\begin{split}
\frac{d^2}{dr^2}\tau(r)
 &\leq  \frac{4}{\delta}\max \left\{ 1,\left( \tau'(r)\right)^{2}
\right\}  \\
&\leq   \frac{4}{\delta}\max \left\{ 1,\left( \frac{4}{\delta}\max\{ 1, \tau(r)^{2}\}\right)^{2}
\right\}    \\
&=  \frac{4^{3}}{ \delta^{3}}\max \left\{ 1, \tau(r)^{4}
\right\} 
\\
&\leq \frac{4^{3}}{ \delta^{3}}+\frac{4^{3}}{ \delta^{3}}\tau(r)^{4}
\end{split}
\end{equation*}
for $r\in (s ,1)$ outside some exceptional set $E$ with $|E|<\delta$.
Hence by $\tau(r)\leq T_s(r,g,\omega_X)$, we get
\begin{equation}\label{eqn:20201217}
m(r,g',\lambda_{\partial TX})
\leq \frac{1}{2}\log \left(
1+c
+cT_s(r,g,\omega_X)^4
\right)
\end{equation}
for $r\in (s ,1)$ outside $E$, where $c=\frac{4^3}{2\pi s \delta^3}$.

\par

\medskip

{\it Step 3.}\
We take a positive integer $l$ such that $\frac{1}{l}<\varepsilon$.
There exists a ramification covering $\varphi:X'\to X$ such that (1) $X'$ is smooth, (2) $D':=(\varphi^*D)_{\mathrm{red}}$ is normal crossing, (3) $lD'\subset \varphi^*D$ (cf. \cite[Thm 17]{Ka2}).
We take a holomorphic map $g:Y\to X'$, where $Y$ is a Riemann surface with a proper surjective holomorphic map $\pi_g:Y\to \mathbb D$, with the following commutative diagram.
\begin{equation*}
\begin{CD}
Y@>g>>  X'\\
@V\pi_gVV    @VV\varphi V\\
\mathbb D @>>f> X
\end{CD}
\end{equation*}
Then we have
$$
lm(r,g,\lambda_{D'})\leq m(r,g,\lambda_{\varphi^*D})=m(r,f,\lambda_{D}).
$$
The morphism $TX'(-\log D')\to TX(-\log D)$ induces a rational map $\Phi:\overline{T}X'(-\log D')\dashrightarrow \overline{T}X(-\log D)$.
Let $Z'\subset \overline{T}X'(-\log D')$ be the indeterminacy locus of $\Phi$.
Then we have 
\begin{equation}\label{eqn:202206012}
(\Phi|_{\overline{T}X'(-\log D')\backslash Z'})^*(\partial TX(-\log D))=\partial TX'(-\log D')|_{\overline{T}X'(-\log D')\backslash Z'}.
\end{equation}
Let $\beta :\widetilde{\overline{T}X'(-\log D')}\to \overline{T}X'(-\log D')$ be a modification such that $\Phi$ induces a morphism $\widetilde{\Phi}:\widetilde{\overline{T}X'(-\log D')}\to \overline{T}X(-\log D)$.
Then by a similar argument used to justify \eqref{eqn:mouii}, we have 
\begin{equation}\label{eqn:202206015}
\widetilde{\Phi}^{*}(\partial TX(-\log D))= \beta^*(\partial TX'(-\log D'))-E',
\end{equation}
where $E'\subset \widetilde{\overline{T}X'(-\log D')}$ is an effective Cartier divisor.

We prove this.
By \eqref{eqn:202206012}, we have the following on $\overline{T}X'(-\log D')\backslash Z'$:
\begin{equation}\label{eqn:20220604}
(\Phi |_{\overline{T}X'(-\log D')\backslash Z'})^*\mathcal{O}_{\overline{T}X(-\log D)}(1)=\mathcal{O}_{\overline{T}X'(-\log D')}(1)|_{\overline{T}X'(-\log D')\backslash Z'}.
\end{equation}
We pull-back \eqref{eqn:202206011} by $\Phi|_{\overline{T}X'(-\log D')\backslash Z'}:\overline{T}X'(-\log D')\backslash Z'\to \overline{T}X(-\log D)$ and combine with \eqref{eqn:20220604}.
Then we have the following surjection on $\overline{T}X'(-\log D')\backslash Z'$:
$$
((p')^{\log})^*(p^{\mathrm{log}})_*\mathcal{O}_{\overline{T}X(-\log D)}(1)|_{\overline{T}X'(-\log D')\backslash Z'}\to \mathcal{O}_{\overline{T}X'(-\log D')}(1)|_{\overline{T}X'(-\log D')\backslash Z'},
$$
where $(p')^{\log}:\overline{T}X'(-\log D')\to X'$ is the projection.
Since $Z'$ has codimension greater than one, this morphism extends over whole $\overline{T}X'(-\log D')$.
We denote by $K$ the kernel of the following surjection obtained from \eqref{eqn:202206011}:
$$
\widetilde{\Phi}^*(p^{\mathrm{log}})^*(p^{\mathrm{log}})_*\mathcal{O}_{\overline{T}X(-\log D)}(1)\to \widetilde{\Phi}^*\mathcal{O}_{\overline{T}X(-\log D)}(1).
$$
Then the composition of
$$
K\to \beta^*((p')^{\log})^*(p^{\mathrm{log}})_*\mathcal{O}_{\overline{T}X(-\log D)}(1)\to \beta^*\mathcal{O}_{\overline{T}X'(-\log D')}(1)
$$
is a zero map on $\beta^{-1}(\overline{T}X'(-\log D')\backslash Z')$, hence on $\widetilde{\overline{T}X'(-\log D')}$.
Hence we get a morphism
$$
\widetilde{\Phi}^*\mathcal{O}_{\overline{T}X(-\log D)}(1)\to  \beta^*\mathcal{O}_{\overline{T}X'(-\log D')}(1).
$$
Hence there exists an effective Cartier divisor $E'$ on $\widetilde{\overline{T}X'(-\log D')}$ such that \eqref{eqn:202206015} is valid.

Now by \eqref{eqn:202206015}, we have
$$
m(r,f',\lambda_{\partial TX(-\log D)})\leq m(r,g',\lambda_{\partial TX'(-\log D')})+\mu_2.
$$
Here $\mu_2>0$ is a positive constant which only depends on the choices of Weil functions.
Using \eqref{eqn:loglld11} for the pair $(X',D')$ instead of $(X,D)$, we get
\begin{equation*}
\begin{split}
m(r,f',\lambda_{\partial TX(-\log D)})&\leq m(r,g',\lambda_{\partial TX'(-\log D')})+\mu_2\\
&\leq m(r,g,\lambda_{D'})+m(r,g',\lambda_{\partial TX'}) +\mu_1+\mu_2\\
&\leq \frac{1}{l}m(r,f,\lambda_{D})+m(r,g',\lambda_{\partial TX'})+\mu_1+\mu_2\\
&\leq \varepsilon m(r,f,\lambda_{D})+m(r,g',\lambda_{\partial TX'})+\mu_1+\mu_2.
\end{split}
\end{equation*}
By \eqref{eqn:20201217} for the pair $(X',\omega_{X'})$ instead of $(X,\omega_X)$, we get
$$
m(r,f',\lambda_{\partial TX(-\log D)})\leq \varepsilon m(r,f,\lambda_{D})
+
\frac{1}{2}\log \left(
1+c
+cT_s(r,g,\omega_{X'})^4
\right)+\mu_1+\mu_2
$$
for all $r\in (s,1)$ outside some exceptional set of linear measure less than $\delta$.
Since $\varphi:X'\to X$ is finite, there exist positive constants $c'>0$ and $c''>0$ such that
$$
T_s(r,g,\omega_{X'})\leq c'T_s(r,g,\varphi^*\omega_X)+c''=c'T_s(r,f,\omega_X)+c''
$$ 
for all $r\in (s,1)$, where $c'$ and $c''$ are independent of the choice of $g$ (cf. Lemma \ref{lem:202110212}).
Hence we get
$$
m(r,f',\lambda_{\partial TX(-\log D)})\leq \varepsilon m(r,f,\lambda_{D})
+
\frac{1}{2}\log \left(
1+c
+c(c'T_s(r,f,\omega_{X})+c'')^4
\right)+\mu_1+\mu_2
$$
for all $r\in (s,1)$ outside some exceptional set of linear measure less than $\delta$.
We take a positive constant $\mu >0$ such that 
$$
\frac{1}{2}\log \left(
1+c
+c(c'x+c'')^4
\right)+\mu_1+\mu_2
\leq \varepsilon x+\mu
$$
for $x\geq 0$.
Then we obtain our estimate.
\hspace{\fill} $\square$

\subsection{The case of semi-abelian varieties}

Let $\overline{A}$ be an equivariant compactification of a semi-abelian variety $A$.
Let $\mathcal{G}$ be an infinite indexed family in $\mathrm{Hol}(\mathbb D,A)$.
We consider the following assumption.

\begin{ass}\label{ass:20201109}
Let $D\subset \partial{A}$ be an irreducible component.
Then $\mathcal{G}'\not\to D$ for every infinite subfamily $\mathcal{G}'$ of $\mathcal{G}$.
\end{ass}

We recall $\Pi(\mathcal{G})$ from Definition \ref{defn:20220601}.
In this subsection, we prove the following lemma.

\begin{lem}\label{lem:tcor}
Let $A$ be a semi-abelian variety and let $S$ be a smooth projective variety.
Let $\overline{A}$ be a smooth projective equivariant compactification.
Let $\omega_{\overline{A}}$ and $\omega_S$ be smooth, positive (1,1)-forms on $\overline{A}$ and $S$, respectively.
Let $\mathcal{F}\subset \mathrm{Hol}(\mathbb D,A\times S)$ be an infinite set such that $\mathcal{F}_A=( f_A)_{f\in\mathcal{F}}$ satisfies Assumption \ref{ass:20201109}.
Assume that $\{0\}\in \Pi(\mathcal{F}_A)$ and that $f_A$ is non-constant for all $f\in \mathcal{F}$.

Let $k\in\mathbb Z_{\geq 1}$.
Let $\omega_{S_{k,A}}$ be a smooth, positive (1,1)-form on $S_{k,A}$.
Then there exists $\sigma\in (0,1)$ with the following property:
Let $s\in (\sigma,1)$, $\varepsilon>0$ and $\delta>0$.
Then there exist positive constants $\mu_1>0$, $\mu_2>0$ such that for all $f\in\mathcal{F}$, the estimate
\begin{equation*}
T_{s}(r,f_{S_{k,A}},\omega_{S_{k,A}})
\leq \varepsilon 
T_{s}(r,f_A,\omega_{\overline{A}})
+\mu_1
T_{s}(r,f_S,\omega_S)
+\mu_2
\end{equation*}
holds for all $r\in (s,1)$ outside some exceptional set of linear measure less than $\delta$.
\end{lem}

To prove this, we prepare several lemmas.

\begin{lem}\label{lem:2020110211}
Let $\overline{A}$ be a smooth projective equivariant compactification.
Let $\mathcal{F}=(f_i)_{i\in I}$ be an infinite indexed family in $\mathrm{Hol}(\mathbb D,A)$ such that the assumption \ref{ass:20201109} is satisfied.
Then there exists $\sigma\in (0,1)$ with the following property:
Let $s\in (\sigma,1)$ and $\delta>0$.
Then there exists a positive constant $c>0$ such that for all $i\in I$, we have
$$
m(r,f_i,\lambda_{\partial A})\leq cT_{s}(r,f_i,\omega_{\overline{A}})+c
$$
for all $r\in (s+\delta,1)$.
\end{lem}

{\it Proof.}\
We first consider the case that the subset
$$
\mathcal{F}_o=\{ f_i;\ i\in I\}\subset \mathrm{Hol}(\mathbb D,A)
$$
is finite.
In this case, our estimate follows from Lemma \ref{lem:20211008}.
Thus we assume that $\mathcal{F}_o$ is infinite.
Then $\mathcal{F}_o$ satisfies the assumption \ref{ass:20201109}.
By replacing $\mathcal{F}$ by $\mathcal{F}_o$, we may assume that $\mathcal{F}$ is a subset of $\mathrm{Hol}(\mathbb D,A)$.

Let $D_1,\ldots,D_n$ be the irreducible components of $\partial A$.
Then for each $j=1,\ldots,n$, by the assumption \ref{ass:20201109}, we apply Lemma \ref{lem:20201119} to get $\sigma_j\in (0,1)$ and $\alpha_j>0$ such that, for $s\in (\sigma_j,1)$ and $\delta>0$, we have
$$
m(r,f,\lambda_{D_j})\leq \frac{\alpha_{j}}{\delta}T_s(r,f,\omega_{\overline{A}})+\alpha_{j}
$$
for all $r\in (s+\delta,1)$ and all $f\in \mathcal{F}$.
We set $\sigma=\max\{ \sigma_1,\ldots,\sigma_n\}$ and $\alpha=\alpha_1+\cdots+\alpha_n$.
Let $s\in (\sigma,1)$ and $\delta>0$.
Then for all $f\in \mathcal{F}$, we have
$$
m(r,f,\lambda_{\partial A})=\sum_{j=1}^{n}m(r,f,\lambda_{D_j})\leq \frac{\alpha}{\delta}T_s(r,f,\omega_{\overline{A}})+\alpha
$$
for all $r\in (s+\delta,1)$.
We set $c=\max\{ \alpha,\alpha/\delta\}$ to conclude the proof.
\hspace{\fill} $\square$

\medskip

Let $V_1$ and $V_2$ be smooth algebraic varieties.
Let $p_1:V_1\times V_2\to V_1$ be the first projection and let $p_2:V_1\times V_2\to V_2$ be the second projection.
Let $\omega_{V_1}$ and $\omega_{V_2}$ be smooth $(1,1)$-forms on $V_1$ and $V_2$, respectively.
We set 
\begin{equation}\label{eqn:20230306}
\omega_{V_1\times V_2}=p_1^*\omega_{V_1}+p_2^*\omega_{V_2}.
\end{equation}

Let $w\in\mathbb D(r)$.
We recall $\varphi_{w,r}:\mathbb D(r)\to \mathbb D(r)$ from \eqref{eqn:20210716}.

\begin{lem}\label{lem:20201102}
Let $S$ be a smooth projective variety.
Let $\overline{A}$ be a smooth projective equivariant compactification.
Let $\omega_{\overline{A}}$ and $\omega_S$ be smooth, positive (1,1)-forms on $\overline{A}$ and $S$, respectively.
Let $\omega_A$ be an invariant positive $(1,1)$-form on $A$.
Let $\mathcal{F}\subset \mathrm{Hol}(\mathbb D,A\times S)$ be an infinite set of holomorphic maps such that $\mathcal{F}_A=( f_A)_{f\in\mathcal{F}}$ satisfies Assumption \ref{ass:20201109}.
Then there exists $\sigma\in (0,1)$ with the following property:
Let $\varepsilon>0$, $s\in(\sigma,1)$ and $\delta>0$.
Then there exists a positive constant $\mu>0$ such that for all $f\in \mathcal{F}$ and $w\in\mathbb D(\sigma)$,  we have
\begin{equation*}
\int_0^{2\pi}\log | (f\circ\varphi_{w,r})'(re^{i\theta})|_{\omega_{A\times S}}\frac{d\theta}{2\pi}
\leq \varepsilon T_s(r,f,\omega_{\overline{A}\times S})+
\mu
\end{equation*}
for all $r\in (s,1)$ outside some exceptional set of linear measure less than $\delta$.
\end{lem}

{\it Proof.}\
We first apply Lemma \ref{lem:2020110211} for $\mathcal{F}_A=(f_A)_{f\in\mathcal{F}}$ to get $\sigma\in (0,1)$.
Let $\varepsilon>0$, $s\in(\sigma,1)$, $\delta>0$ be given.

Set $D=\partial A\times S\subset \overline{A}\times S$.
We first note that by Lemma \ref{lem:202111111}, we have 
$$T(\overline{A}\times S)(-\log D)=T\overline{A}(-\log\partial A)\times TS=\overline{A}\times \mathrm{Lie}A\times TS.$$
Hence $\omega_{A\times S}$ defines a Hermitian metric on $T(\overline{A}\times S)(-\log D)$.
Hence we may define the Weil function by
\begin{equation}\label{eqn:202111112}
\lambda_{\partial T(\overline{A}\times S)(-\log D)}(v)=\log \sqrt{1+|v|^2_{\omega_{A\times S}}},
\end{equation}
where $v\in T(\overline{A}\times S)(-\log D)$.

By Lemma \ref{lem:2020110211}, there exists $c>0$ such that
\begin{equation*}
m(r,f_A,\lambda_{\partial A})\leq cT_s(r,f_A,\omega_{\overline{A}})+c
\end{equation*}
for all $f\in\mathcal{F}$ and $r\in (s+\delta',1)$, where $\delta'=\delta/2$.
Set $\varepsilon'=\frac{(s-\sigma)\varepsilon}{(s+\sigma)(1+c)}$.
By Lemma \ref{lem:llld}, there exists $\mu_0>0$ such that
\begin{equation*}
m(r,f',\lambda_{\partial T(\overline{A}\times S)(-\log D)})\leq \varepsilon'T_s(r,f,\omega_{\overline{A}\times S})+\varepsilon'm(r,f,\lambda_D)+\mu_0
\end{equation*}
for all $f\in\mathcal{F}$ and $r\in (s,1)$ outside some exceptional set of linear measure less than $\delta'$.
Hence we have 
\begin{equation}\label{eqn:20210717}
m(r,f',\lambda_{\partial T(\overline{A}\times S)(-\log D)})\leq \varepsilon'(1+c)T_s(r,f,\omega_{\overline{A}\times S})+\varepsilon' c+\mu_0
\end{equation}
for all $f\in\mathcal{F}$ and $r\in (s,1)$ outside some exceptional set of linear measure less than $2\delta'$.

Now we have
$$
| (f\circ \varphi_{r,w})'(z)|_{\omega_{A\times S}}=| f'\circ \varphi_{r,w}(z)|_{\omega_{A\times S}}\times|\varphi_{r,w}'(z)|.
$$
Hence by \eqref{eqn:202107161}, we have
$$
\int_0^{2\pi}\log | (f\circ\varphi_{w,r})'(re^{i\theta})|_{\omega_{A\times S}}\frac{d\theta}{2\pi}\leq \int_0^{2\pi}\log | f'\circ\varphi_{w,r}(re^{i\theta})|_{\omega_{A\times S}}\frac{d\theta}{2\pi}+\frac{s+\sigma}{s-\sigma}
$$
for all $r\in (s,1)$.
Using \eqref{eqn:202111112}, we have
\begin{equation*}
\begin{split}
\int_0^{2\pi}\log | f'\circ\varphi_{r,w}(re^{i\theta})|_{\omega_{A\times S}}\frac{d\theta}{2\pi}&\leq
\int_0^{2\pi}\log\sqrt{1+ | f'\circ\varphi_{r,w}(re^{i\theta})|^2_{\omega_{A\times S}}}\frac{d\theta}{2\pi}
\\
&=m(r,f'\circ\varphi_{r,w},\lambda_{\partial T(\overline{A}\times S)(-\log D)}).
\end{split}
\end{equation*}
By Lemma \ref{lem:20210711}, we have
$$
m(r,f'\circ\varphi_{r,w},\lambda_{\partial T(\overline{A}\times S)(-\log D)})\leq \frac{s+\sigma}{s-\sigma} m(r,f',\lambda_{\partial T(\overline{A}\times S)(-\log D)})
$$
for all $r\in (s,1)$.
Hence
$$
\int_0^{2\pi}\log | (f\circ\varphi_{w,r})'(re^{i\theta})|_{\omega_{A\times S}}\frac{d\theta}{2\pi}
\leq \frac{s+\sigma}{s-\sigma} m(r,f',\lambda_{\partial T(\overline{A}\times S)(-\log D)})+\frac{s+\sigma}{s-\sigma}
$$
for all $r\in (s,1)$.
Combining this with \eqref{eqn:20210717}, we get ou lemma.
Here we set $\mu=\frac{s+\sigma}{s-\sigma}(\varepsilon'c+\mu_0+1)$.
\hspace{\fill} $\square$

Let $\mu$ be a non-negative mass on $\mathbb D(r)$.
For $|w|<s<r$, we define
$$
T_s^w(r,\mu)=\int_s^r\frac{dt}{t}\int_{\mathbb D(t)}\varphi_{w,r}^*\mu.
$$

\begin{lem}\label{lem:20210711a}
Let $0<\sigma<s<r<1$.
Then for all non-negative mass $\mu$ on $\mathbb D(r)$ and $w\in \mathbb D(\sigma)$, we have
$$
\frac{s(s-\sigma)}{s+\sigma}T_s^w(r,\mu)\leq T_s(r,\mu)\leq \frac{s+\sigma}{s(s-\sigma)} T_s^w(r,\mu).
$$
\end{lem}

To prove this lemma, we prepare the following two estimates:
Let $|w|<\sigma<s<r<1$, then we have 
\begin{equation}\label{eqn:202206021}
|\varphi_{w,r}(z)| \geq \frac{s+\sigma}{s-\sigma}(|z|-r)+r,
\end{equation}
\begin{equation}\label{eqn:202206022}
|\varphi_{w,r}(z)|\leq \frac{s-\sigma}{s+\sigma}(|z|-r)+r,
\end{equation}
for all $z\in\mathbb D(r)$.
Indeed, we have
\begin{equation}\label{eqn:20220602}
r^2\frac{|z|-|w|}{r^2-|w||z|}\leq |\varphi_{w,r}(z)|\leq r^2\frac{|w|+|z|}{r^2+|w||z|}
\end{equation}
for all $z\in \mathbb D(r)$.
We have
$$
r^2\frac{|z|-|w|}{r^2-|w||z|}=r\frac{r+|w|}{r^2-|w||z|}(|z|-r)+r.
$$
For $z\in \mathbb D(r)$, we have
$$
\frac{r+|w|}{r^2-|w||z|}\leq \frac{r+|w|}{r(r-|w|)}\leq \frac{s+\sigma}{r(s-\sigma)}.
$$
Hence combining these with \eqref{eqn:20220602}, we get \eqref{eqn:202206021}.
To prove \eqref{eqn:202206022}, we note
$$
r^2\frac{|w|+|z|}{r^2+|w||z|}=r\frac{r-|w|}{r^2+|w||z|}(|z|-r)+r.
$$
For $z\in \mathbb D(r)$, we have
$$
\frac{r-|w|}{r^2+|w||z|}\geq \frac{r-|w|}{r(r+|w|)}\geq \frac{s-\sigma}{r(s+\sigma)}.
$$
Hence combining these with \eqref{eqn:20220602}, we get \eqref{eqn:202206022}.

{\it Proof of Lemma \ref{lem:20210711a}.}\
We set
$$
U_s^w(r,\mu)=\int_s^rdt\int_{\mathbb D(t)}\varphi_{w,r}^*\mu.
$$
Then we have
$$
U_s^w(r,\mu)\leq T_s^w(r,\mu)\leq \frac{1}{s}U_s^w(r,\mu)
$$
for all $r\in (s,1)$.

We fix $r\in (s,1)$ and $w\in \mathbb D(\sigma)$.
For $x\in (0,r)$, we set
$$
\Lambda(x)=\min\left\{ r-x,r-s\right\}.
$$
Then we have
$$
U_s^w(r,\mu)=\int_{\mathbb D(r)}\Lambda(|z|)\varphi_{w,r}^*\mu_z.
$$
In particular,
$$
U_s^0(r,\mu)=\int_{\mathbb D(r)}\Lambda(|z|)\mu_z.
$$

First by \eqref{eqn:202206021}, we have
$$
 r-|\varphi_{w,r}(z)| \leq \frac{s+\sigma}{s-\sigma}(r-|z|)
$$
for $|z|\leq r$.
Hence we have
$$
\Lambda(|\varphi_{w,r}(z)|)\leq r-|\varphi_{w,r}(z)|
\leq
\frac{s+\sigma}{s-\sigma}\Lambda(|z|)
$$
for $s\leq |z|\leq r$.
By $\Lambda(|\varphi_{w,r}(z)|)\leq r-s$, this holds for all $|z|\leq r$.
Hence we get
\begin{equation*}
U_s^0(r,\mu)=\int_{\mathbb D(r)}\Lambda(|z|)\mu_z
=\int_{\mathbb D(r)}\Lambda(|\varphi_{w,r}(z)|)\varphi_{w,r}^*\mu_z
\leq \frac{s+\sigma}{s-\sigma}U_s^w(r,\mu).
\end{equation*}
Thus
$$
T_s(r,\mu)\leq \frac{1}{s}U_s^0(r,\mu)\leq \frac{s+\sigma}{s(s-\sigma)}T_s^w(r,\mu).
$$

Next by \eqref{eqn:202206022}, we have
$$
\frac{s-\sigma}{s+\sigma}(r-|z|)
\leq r-|\varphi_{w,r}(z)|
$$
for $|z|\leq r$.
Hence we have
$$
\Lambda(|\varphi_{w,r}(z)|)\geq  \frac{s-\sigma}{s+\sigma}\Lambda(|z|)
$$
for $|z|\leq r$, provided $s\leq |\varphi_{w,r}(z)|\leq r$.
By $\Lambda(|z|)\leq r-s$, this holds for all $|z|\leq r$.
Hence we get
\begin{equation*}
U_s^0(r,\mu)=\int_{\mathbb D(r)}\Lambda(|z|)\mu_z
=\int_{\mathbb D(r)}\Lambda(|\varphi_{w,r}(z)|)\varphi_{w,r}^*\mu_z
\geq \frac{s-\sigma}{s+\sigma}U_s^w(r,\mu).
\end{equation*}
Thus
$$
T_s(r,\mu)\geq U_s^0(r,\mu)\geq \frac{s(s-\sigma)}{s+\sigma}T_s^w(r,\mu).
$$
This concludes the proof of our lemma.
\hspace{\fill} $\square$

\medskip

{\it Proof of Lemma \ref{lem:tcor}.}
We prove our lemma by the induction on $k$.
Thus we first consider the case $k=1$.
We set $D=\partial (A\times S)$.
Then $D\subset \overline{A}\times S$ is a simple normal crossing divisor.
By Lemma \ref{lem:202111111}, we note that $\omega_{A\times S}$ induces a Hermitian metric on the tautological line bundle $\mathcal{O}_{PT(\overline{A}\times S)(-\log D)}(1)$ on $PT(\overline{A}\times S)(-\log D)$.
Let $\omega_{\mathcal{O}_{PT(\overline{A}\times S)(-\log D)}(1)}$ be the associated curvature form for this tautological line bundle.
Similarly, we note that $\omega_{A\times S}$ induces a Hermitian metric on the tautological line bundle $\mathcal{O}_{P(TS\times\mathrm{Lie}A)}(1)$ on $S_{1,A}$.
Let $\omega_{\mathcal{O}_{P(TS\times\mathrm{Lie}A)}(1)}$ be the associated curvature form for this tautological line bundle.
By $PT(\overline{A}\times S)(-\log D)=\overline{A}\times S_{1,A}$, we have
$q^*\omega_{\mathcal{O}_{P(TS\times\mathrm{Lie}A)}(1)}=\omega_{\mathcal{O}_{PT(\overline{A}\times S)(-\log D)}(1)}$ where $q:\overline{A}\times S_{1,A}\to S_{1,A}$ is the projection.
There exist positive constants $\alpha_1$ and $\alpha_2$ such that
\begin{equation}\label{eqn:20230307}
\omega_{S_{1,A}}\leq \alpha_1\tau^*\omega_{S}+\alpha_{2}\omega_{\mathcal{O}_{P(TS\times\mathrm{Lie}A)}(1)},
\end{equation}
where $\tau:S_{1,A}\to S$.

We first  choose $\sigma\in (0,1)$ which appears in Lemma \ref{lem:20201102}.
We modify $\sigma$ as follows.
Let $\omega_{A}$ be an invariant positive $(1,1)$ form on $A$.
For each $n\in \mathbb Z_{\geq 2}$, let $\mathcal{E}_n\subset \mathcal{F}$ be the set of $f\in \mathcal{F}$ such that 
$$\sup_{z\in \mathbb D(1-\frac{1}{n})}|f_A'(z)|_{\omega_A}\leq \frac{1}{n}.$$
Then we have $\mathcal{E}_2\supset \mathcal{E}_3\supset\cdots$.
By the assumption $\{ 0\}\in\Pi (\mathcal{F}_A)$, we may take $n_0\in \mathbb Z_{\geq 2}$ such that $\mathcal{E}_{n_0}$ is finite.
Since $f_A$ is non-constant for all $f\in\mathcal{F}$, we may take $n_1\geq n_0$ such that $\mathcal{E}_{n_1}=\emptyset$.
By enlarging $\sigma$ if necessary, we may assume that $1-\frac{1}{n_1}\leq \sigma<1$.
Then for all $f\in\mathcal{F}$, we have
$$
\sup_{z\in \mathbb D(\sigma)}|f_A'(z)|_{\omega_A}\geq 1-\sigma.
$$
For each $f\in\mathcal{F}$, we take $w_f\in \overline{\mathbb D(\sigma)}$ such that $$|f_A'(w_f)|_{\omega_A}=\sup_{z\in \mathbb D(\sigma)}|f_A'(z)|_{\omega_A}.$$

Let $s\in (\sigma,1)$ and let $f\in \mathcal{F}$.
For $r\in (s,1)$, we take $\varphi_{w_f,r}:\mathbb D(r)\to \mathbb D(r)$ as in \eqref{eqn:20210716} and set $g=f\circ\varphi_{w_f,r}$.
Then by \eqref{eqn:20230307}, we have
\begin{equation*}
T_s(r,g_{S_{1,A}},\omega_{S_{1,A}})
\leq 
\alpha_1
T_s(r,g_S,\omega_S) +\alpha_2T_s(r,g_{[1]},\omega_{\mathcal{O}_{PT(\overline{A}\times S)(-\log D)}(1)}).
\end{equation*}
By Lemma \ref{lem:20210711a}, we have
\begin{equation*}
T_s(r,f_{S_{1,A}},\omega_{S_{1,A}})
\leq \frac{s+\sigma}{s(s-\sigma)}T_s(r,g_{S_{1,A}},\omega_{S_{1,A}})
\end{equation*}
and
\begin{equation*}
T_s(r,f_{S},\omega_{S})
\geq \frac{s(s-\sigma)}{s+\sigma}T_s(r,g_{S},\omega_{S}).
\end{equation*}
Thus we get
\begin{equation}\label{eqn:202206027}
T_s(r,f_{S_{1,A}},\omega_{S_{1,A}})
\leq 
\frac{\alpha_1(s+\sigma)^2}{s^2(s-\sigma)^2}
T_s(r,f_S,\omega_S) +\frac{\alpha_2(s+\sigma)}{s(s-\sigma)}T_s(r,g_{[1]},\omega_{\mathcal{O}_{PT(\overline{A}\times S)(-\log D)}(1)}).
\end{equation}

Next we estimate the second term of the right hand side.
We claim
\begin{equation}\label{eqn:202206024}
T_s(r,g_{[1]},\omega_{\mathcal{O}_{PT(\overline{A}\times S)(-\log D)}(1)})
\leq
\int_0^{2\pi}\log | g'(re^{i\theta})|_{\omega_{A\times S}}\frac{d\theta}{2\pi}-\int_0^{2\pi}\log | g'(se^{i\theta})|_{\omega_{A\times S}}\frac{d\theta}{2\pi},
\end{equation}
where $|\cdot |_{\omega_{A\times S}}$ is the length with respect to $\omega_{A\times S}$.
This is obtained as follows.
The metric $\omega_{A\times S}$ defines a Hermitian metric $|\cdot |_{\omega_{A\times S}}$ on $\mathcal{O}_{PT(\overline{A}\times S)(-\log D)}(-1)$, whose curvature form is $-\omega_{\mathcal{O}_{PT(\overline{A}\times S)(-\log D)}(1)}$.
By the Poincar{\'e}-Lelong formula, we have
$$
-(g_{[1]})^*
\omega_{\mathcal{O}_{PT(\overline{A}\times S)(-\log D)}(1)}
=[(g')^*Z]-2dd^c\log |g'|_{\omega_{A\times S}}
$$
as currents on $\mathbb D$, where $Z$ is the zero section of $\mathcal{O}_{PT(\overline{A}\times S)(-\log D)}(-1)$.
By the Jensen formula, we get \eqref{eqn:202206024}.
Compare with the proof of \eqref{eqn:20211026}.

Now let $\varepsilon>0$ and $\delta>0$.
We set $\varepsilon'=\frac{s(s-\sigma)\varepsilon}{\alpha_2(s+\sigma)}$.
By Lemma \ref{lem:20201102}, we get
\begin{equation*}
\int_0^{2\pi}\log | g'(re^{i\theta})|_{\omega_{A\times S}}\frac{d\theta}{2\pi}
\leq \varepsilon' T_s(r,f,\omega_{\overline{A}\times S})+
\mu
\end{equation*}
for $r\in (s,1)$ outside some exceptional set of linear measure less than $\delta$.
Here $\mu>0$ is a positive constant which appears in Lemma \ref{lem:20201102}.
In particular $\mu$ is independent of the choices of $f\in\mathcal{F}$ and $w_f$.
Since $\log |g_A'|_{\omega_A}$ is subharmonic, we have
$$
\int_0^{2\pi}\log | g'(se^{i\theta})|_{\omega_{A\times S}}\frac{d\theta}{2\pi}\geq \int_0^{2\pi}\log | g_A'(se^{i\theta})|_{\omega_{A}}\frac{d\theta}{2\pi}\geq
\log | g_A'(0)|_{\omega_{A}}.
$$
Note that by \eqref{eqn:202107161}, we have
$$
 | g_A'(0)|_{\omega_{A}}=| f_A'(w_f)|_{\omega_{A}} |\varphi_{w_f,r}'(0)|\geq (1-\sigma)\frac{s-\sigma}{s+\sigma}.
$$
Hence, taking into account these estimate in \eqref{eqn:202206024}
, we get
\begin{equation}\label{eqn:202107081}
T_s(r,g_{[1]},\omega_{\mathcal{O}_{PT(\overline{A}\times S)(-\log D)}(1)})
\leq \varepsilon'
T_s(r,f,\omega_{\overline{A}\times S})
+\mu'
\end{equation}
for $r\in (s,1)$ outside some exceptional set of linear measure less than $\delta$.
Here we set $\mu'=\mu+\log \frac{s+\sigma}{(1-\sigma)(s-\sigma)}$.

Combining \eqref{eqn:202206027} with \eqref{eqn:202107081}, we get
\begin{equation*}
T_s(r,f_{S_{1,A}},\omega_{S_{1,A}})
\leq \varepsilon
T_s(r,f_A,\omega_{\overline{A}})
+\mu_1
T_s(r,f_S,\omega_S) +\mu_2
\end{equation*}
for $r\in (s,1)$ outside some exceptional set of linear measure less than $\delta$.
Here we set $\mu_1=\frac{\alpha_1(s+\sigma)^2}{s^2(s-\sigma)^2}+\frac{\alpha_2(s+\sigma)\varepsilon'}{s(s-\sigma)}$, $\mu_2=\frac{\alpha_2(s+\sigma)\mu'}{s(s-\sigma)}$ to conclude the proof of the case $k=1$.

\par

We assume that $k\geq 2$ and the lemma is true for $k-1$.
By the construction of $S_{k,A}$, we have $S_{k,A}\subset (S_{k-1,A})_{1,A}$ (cf. \eqref{eqn:202112033}).
Hence there exists a positive constant $\alpha>0$ such that $\omega_{S_{k,A}}\leq \alpha \omega_{(S_{k-1,A})_{1,A}}$ on $S_{k,A}$.
Hence we have
\begin{equation}\label{eqn:202303071}
T_s(r,f_{S_{k,A}},\omega_{S_{k,A}})\leq \alpha T_s(r,(f_{[k-1]})_{(S_{k-1,A})_{1,A}},\omega_{(S_{k-1,A})_{1,A}}).
\end{equation}
By the first step applied to $A\times S_{k-1,A}$ and $\{f_{[k-1]};f\in\mathcal{F}\}\subset \mathrm{Hol}(\mathbb D,A\times S_{k-1,A})$, we get $\sigma'\in (0,1)$ with the following property:
For $s\in(\sigma',1)$, $\varepsilon>0$ and $\delta>0$, there exist $\mu_1'=\mu_1'(s,\varepsilon,\delta)>0$ and $\mu_2'=\mu_2'(s,\varepsilon,\delta)>0$ such that for all $f\in\mathcal{F}$, we have
\begin{equation}
T_s(r,(f_{[k-1]})_{(S_{k-1,A})_{1,A}},\omega_{(S_{k-1,A})_{1,A}})
\leq \frac{\varepsilon}{2\alpha} 
T_s(r,f_A,\omega_{\overline{A}})
+\mu_1'
T_s(r,(f_{[k-1]})_{S_{k-1,A}},\omega_{S_{k-1,A}})
 +\mu_2'
\end{equation}
for all $r\in (s,1)$ outside some exceptional set of linear measure less than $\delta/2$.
Now by the induction hypothesis, we get $\sigma\in (\sigma',1)$ with the following property:
For $s\in(\sigma,1)$, $\varepsilon>0$ and $\delta>0$, there exist $\mu_1''>0$ and $\mu_2''>0$ such that for all $f\in\mathcal{F}$, we have
\begin{equation}\label{eqn:202303072}
T_s(r,(f_{[k-1]})_{S_{k-1,A}},\omega_{S_{k-1,A}})
\leq \frac{\varepsilon}{2\alpha\mu_1'(s,\varepsilon,\delta)} 
T_s(r,f_A,\omega_{\overline{A}})
+\mu_1''
T_s(r,f_S,\omega_S)
 +\mu_2''
\end{equation}
holds for $r\in (s,1)$ outside some exceptional set of linear measure less than $\delta/2$.
Hence by \eqref{eqn:202303071}--\eqref{eqn:202303072},
our $\sigma$ satisfies the required property to conclude the induction step.
Here we set $\mu_1=\alpha\mu_1'\mu_1''$ and $\mu_2=\alpha(\mu_1'\mu_2''+\mu_2')$.
\hspace{\fill} $\square$

\section{Nevanlinna theory and blowing-ups}\label{sec:7}

In this section, we shall establish Lemma \ref{lem:20210302}.
For this purpose, we start from a general estimate (cf. Lemma \ref{lem:1.1}).

\subsection{A general estimate}

\begin{lem}\label{lem:202105294}
Let $S$ be a smooth projective variety.
Let $W\subset S$ be an irreducible Zariski closed set.
Let $\mathcal{F}\subset \mathrm{Hol}(\mathbb D,S)$ be an infinite subset such that $\mathcal{F}\not\Rightarrow W$ and $f(\mathbb D)\not\subset W$ for all $f\in \mathcal{F}$ with finite exception.
Then there exist an infinite subset $\mathcal{F}'\subset \mathcal{F}$ and a $W$-admissible modification $\varphi:S'\to S$ such that the following two properties hold:
\begin{enumerate}
\item
$S'=\mathrm{Bl}_YS$ for some closed subscheme $Y\subset S$ such that $\supp Y\subsetneqq W$ and $f(\mathbb D)\not\subset \supp Y$ for all $f\in \mathcal{F}'$.
\item
$\mathcal{F}'\not\to W'$, where $W'\subset S'$ is the strict transform.
\end{enumerate}
Moreover we may take $S'$ to be smooth. 
\end{lem}

{\it Proof.}\
We consider the following two cases.

{\it Case 1.}\
There exists a Zariski closed set $Z\subset S$ such that $W\not\subset Z$ and $f(\mathbb D)\subset Z$ for infinitely many $f\in \mathcal{F}$.
In this case, we take an infinite subset $\mathcal{F}'\subset \mathcal{F}$ such that $f(\mathbb D)\subset Z$ for all $f\in \mathcal{F}'$.
Set $Y=W\cap Z$ as a closed subscheme.
Then $\supp Y\subsetneqq W$.
Set $S'=\mathrm{Bl}_YS$.
Then $Z'\cap W'=\emptyset$ in $S'$, where $Z'$ and $W'$ are the strict transforms of $Z$ and $W$, respectively (cf. \eqref{eqn:20220522}).
By $\supp Y\subset W$, we have $f(\mathbb D)\not\subset \supp Y$ for all $f\in \mathcal{F}'$ with finite exception.
Hence by removing these finite elements from $\mathcal{F}'$, we may assume $f(\mathbb D)\not\subset \supp Y$ for all $f\in \mathcal{F}'$.
Hence we may consider as $f(\mathbb D)\subset Z'$ for all $f\in \mathcal{F}'$.
Hence $\mathcal{F}'\not\to W'$.

{\it Case 2.}\
Otherwise, there exists a $W$-admissible modification $\hat{S}\to S$ such that $\mathcal{F}\not\to \hat{W}$, where $\hat{W}\subset \hat{S}$ is the minimal transform.
By \cite[p. 171, Exercise 7.11 (c)]{H}, there exists a closed subscheme $Z\subset S$ such that $\hat{S}=\mathrm{Bl}_ZS$ and $W\not\subset \supp Z$.
Set $Y=W\cap Z$ as closed subschemes of $S$.
Then $\supp Y\subsetneqq W$.
Since we are in case 2, $W\not\subset \supp Z$ implies that only finitely many $f\in \mathcal{F}$ satisfies $f(\mathbb D)\subset \supp Z$.
By removing these finite elements from $\mathcal{F}$, we get $\mathcal{F}'\subset \mathcal{F}$ so that $f(\mathbb D)\not\subset \supp Z$ for all $f\in \mathcal{F}'$.

Now we set $S'=\mathrm{Bl}_YS$ with $\varphi:S'\to S$.
There exists a closed subscheme $Z^{\dagger}\subset S'$ such that $\mathcal{I}_{\varphi^*Z}=\mathcal{I}_{\varphi^*Y}\cdot\mathcal{I}_{Z^{\dagger}}$.
Indeed, since $\varphi^*Y\subset S'$ is a Cartier divisor, $\mathcal{I}_{\varphi^*Y}$ is an invertible sheaf.
By $\mathcal{I}_{\varphi^*Z}\subset \mathcal{I}_{\varphi^*Y}\subset \mathcal{O}_{S'}$, we take $Z^{\dagger}\subset S'$ such that $\mathcal{I}_{Z^{\dagger}}=\mathcal{I}_{\varphi^*Z}\otimes (\mathcal{I}_{\varphi^*Y})^{-1}\subset \mathcal{O}_{S'}$.
Similarly there exists $W^{\dagger}\subset S'$ such that $\mathcal{I}_{\varphi^*W}=\mathcal{I}_{\varphi^*Y}\cdot\mathcal{I}_{W^{\dagger}}$.
Then by $\mathcal{I}_{\varphi^*Z}+\mathcal{I}_{\varphi^*W}=\mathcal{I}_{\varphi^*Y}$, we have $\mathcal{I}_{\varphi^*Y}\cdot (\mathcal{I}_{Z^{\dagger}}+\mathcal{I}_{W^{\dagger}})=\mathcal{I}_{\varphi^*Y}$ so that $\mathcal{I}_{Z^{\dagger}}+\mathcal{I}_{W^{\dagger}}=\mathcal{O}_{S'}$.
Hence $Z^{\dagger}\cap W^{\dagger}=\emptyset$.
Note that $W'\subset W^{\dagger}$, where $W'\subset S'$ is the strict transform.
Hence we get
$$
Z^{\dagger}\cap W'=\emptyset.
$$

We set $\tilde{S}=\mathrm{Bl}_{Z^{\dagger}}S'$.
Then since $(\varphi\circ p)^*Z$ is a Cartier divisor, $\tilde{S}\to \hat{S}$ exists.
\begin{equation*}
\begin{CD}
S'@=\mathrm{Bl}_{Y}S@<p<< \mathrm{Bl}_{Z^{\dagger}}S'@= \tilde{S}\\
 @.  @V{\varphi}VV    @VVV @.\\
 @. S @<<< \mathrm{Bl}_{Z}S @= \hat{S}
\end{CD}
\end{equation*}
Let $\tilde{W}\subset \tilde{S}$ be the strict transform.
Then by $\mathcal{F}\not\to \hat{W}$, we have $\mathcal{F}'\not\to \tilde{W}$.
On the other hand, by $Z^{\dagger}\cap W'=\emptyset$, the morphism $p:\tilde{S}\to S'$ is an isomorphism on a neighbourhood of $\tilde{W}$.
Hence $\mathcal{F}'\not\to W'$.

Now suppose $S'$ is not smooth.
Then we take a birational modification $S''\to S'$ such that $S''$ is smooth.
Since $S$ is smooth, we may assume that $S''\to S$ is an isomorphism over $S\backslash \supp Y$.
 By \cite[p. 171, Exercise 7.11 (c)]{H}, we may take a closed subscheme $Y'\subset S$ such that $S''=\mathrm{Bl}_{Y'}S$ and $\supp Y'\subset \supp Y$.
We replace $S'$ by $S''$.
The proof of the lemma is completed.
\hspace{\fill} $\square$

\begin{lem}\label{lem:1.1}
Let $S$ be a smooth projective variety with a smooth positive $(1,1)$-form $\omega_S$.
Let $W\subset S$ be a Zariski closed set with a Weil function $\lambda_W\geq 0$.
Let $\mathcal{F}\subset \mathrm{Hol}(\mathbb D,S)$ be an infinite set of holomorphic maps and let $\{ E_i\}_{i\in I}$ be a countable set of Zariski closed subsets of $S$.
Assume that $\mathrm{LIM}(\mathcal{F},\{ E_i\}_{i\in I})$ exists and $\mathrm{LIM}(\mathcal{F},\{ E_i\}_{i\in I})\not\subset W$.
Then there exist $\sigma\in (0,1)$, $\beta>0$, $k\in \mathbb Z_{\geq 0}$, a Zariski closed set $E\subset \cup E_i$ with a Weil function $\lambda_E\geq 0$ and an infinite subset $\mathcal{G}\subset \mathcal{F}$ with the following two properties.
\begin{enumerate}
\item
For all $f\in \mathcal{G}$, we have $f(\mathbb D)\not\subset E\cup W$.
\item
Let $s\in (\sigma,1)$ and $f\in \mathcal{G}$.
Then we have
\begin{equation*}
m(r,f,\lambda_W)
\leq \frac{\beta}{(r-s)^k}
T_{s}(r,f,\omega_S)
+\frac{\beta}{(r-s)^k}m\left(s,f,\lambda_E\right)
+\frac{\beta}{(r-s)^k},
\end{equation*}
for all  $r\in (s,1)$.
\end{enumerate}
\end{lem}

{\it Proof.}\
The proof is by Noetherian induction on $W$.
Let $\mathcal{P}$ be the set of all Zariski closed set $W\subset S$ such that our lemma is false for $W$.
To show $\mathcal{P}=\emptyset$, we assume contrary that $\mathcal{P}\not=\emptyset$.
We take a minimal element $W\in \mathcal{P}$.
In the following, we shall show that our lemma is true for $W$.

Let $\mathcal{F}\subset \mathrm{Hol}(\mathbb D,S)$ and $\{ E_i\}_{i\in I}$ be the objects which appear in our lemma.
Thus $\mathrm{LIM}(\mathcal{F},\{ E_i\}_{i\in I})$ exists and $\mathrm{LIM}(\mathcal{F},\{ E_i\}_{i\in I})\not\subset W$.

We first observe that $W$ is irreducible.
Otherwise, we have $W=W_1\cup W_2$.
Let $\lambda_{W_1}\geq 0$ and $\lambda_{W_2}\geq 0$ be Weil functions such that $\lambda_W\leq \lambda_{W_1}+\lambda_{W_2}$.
By $W_1\subsetneqq W$, we have $W_1\not\in\mathcal{P}$.
Note that $\mathrm{LIM}(\mathcal{F},\{ E_i\}_{i\in I})\not\subset W_1$.
Hence we may take
$\sigma_1\in (0,1)$, $\beta_1>0$, $k_1\in \mathbb N$, a Zariski closed set $E_1\subset \cup E_i$ with a Weil function $\lambda_{E_1}\geq 0$ and an infinite subset $\mathcal{G}_1\subset \mathcal{F}$ such that for all $f\in\mathcal{G}_1$, we have $f(\mathbb D)\not\subset E_1\cup W_1$ and
$$
m(r,f,\lambda_{W_1})
\leq \frac{\beta_1}{(r-s)^{k_1}}
T_{s}(r,f,\omega_S)
+\frac{\beta_1}{(r-s)^{k_1}}m\left(s,f,\lambda_{E_1}\right)
+\frac{\beta_1}{(r-s)^{k_1}},
$$
where $\sigma_1<s<r<1$.
Now $\mathrm{LIM}(\mathcal{G}_1,\{ E_i\}_{i\in I})$ exists and 
$$\mathrm{LIM}(\mathcal{G}_1,\{ E_i\}_{i\in I})=\mathrm{LIM}(\mathcal{F},\{ E_i\}_{i\in I})\not\subset W_2.$$
Hence by $W_2\not\in\mathcal{P}$, we may take $\sigma_2\in (0,1)$, $\beta_2>0$, $k_2\in \mathbb N$, a Zariski closed set $E_2\subset \cup E_i$ with a Weil function $\lambda_{E_2}\geq 0$ and an infinite subset $\mathcal{G}_2\subset \mathcal{G}_1$ such that for all $f\in\mathcal{G}_2$, we have $f(\mathbb D)\not\subset E_2\cup W_2$ and
$$
m(r,f,\lambda_{W_2})
\leq \frac{\beta_2}{(r-s)^{k_2}}
T_{s}(r,f,\omega_S)
+\frac{\beta_2}{(r-s)^{k_2}}m\left(s,f,\lambda_{E_2}\right)
+\frac{\beta_1}{(r-s)^{k_2}},
$$
where $\sigma_2<s<r<1$.
We set $\sigma=\max\{\sigma_1,\sigma_2\}$, $\beta=\beta_1+\beta_2$, $k=\max\{ k_1,k_2\}$, $E=E_1\cup E_2$, $\mathcal{G}=\mathcal{G}_2$.
We take a Weil function $\lambda_E$ such that $\lambda_E\geq \max\{\lambda_{E_1},\lambda_{E_2}\}$.
Then by $m(r,f,\lambda_W)\leq m(r,f,\lambda_{W_1})+m(r,f,\lambda_{W_2})$, the two properties of our lemma is satisfied.
Hence $W\not\in\mathcal{P}$, a contradiction.
Hence we may assume that $W$ is irreducible.

We consider two cases.

\medskip

{\it Case 1:}\
$W\subset \cup E_i$.
In this case, we set $\sigma=1/2$, $E=W$ and $\lambda_E=\lambda_W$.
By $\mathrm{LIM}(\mathcal{F};\{ E_i\})\not\subset W$, only finitely many $f\in \mathcal{F}$ satisfies $f(\mathbb D)\subset W$.
We remove these $f$ from $\mathcal{F}$ to get an infinite subset $\mathcal{G}\subset \mathcal{F}$.
Then $f(\mathbb D)\not\subset E\cup W$ for all $f\in \mathcal{G}$.
By Lemma \ref{lem:20211008}, there exists a positive constant $c>0$ such that
\begin{equation*}
m(r,f,\lambda_W)
\leq c
T_{s}(r,f,\omega_S)
+c
m(s,f,\lambda_W)
\end{equation*}
for all $s\in (\sigma,1)$, $r\in (s,1)$ and $f\in \mathcal{G}$.
By $\lambda_W=\lambda_E$,
our Lemma is true for $W$.
Here we set $\beta=c$ and $k=0$.

\medskip

{\it Case 2:}\
$W\not\subset \cup E_i$.
By $\mathrm{LIM}(\mathcal{F};\{ E_i\})\not\subset W$ and $W\not\subset \cup E_i$, we have $\mathcal{F}\not\Rightarrow W$ (cf. Lemma \ref{lem:202105271}).
By Lemma \ref{lem:202105294}, there exist an infinite subset $\mathcal{F}'\subset \mathcal{F}$ and a closed subscheme $\mathcal{Z}\subset S$ such that 
\begin{enumerate}
\item
$\supp \mathcal{Z}\subsetneqq W$,
\item
$f(\mathbb D)\not\subset \supp \mathcal{Z}$ for all $f\in \mathcal{F}'$,
\item
$\mathcal{F}'\not\to W'$, where $W'\subset \mathrm{Bl}_{\mathcal{Z}}S$ is the strict transform.
\end{enumerate}
Set $S'=\mathrm{Bl}_{\mathcal{Z}}S$.
We may assume that $S'$ is smooth.
By Lemma \ref{lem:20210803}, replacing $\mathcal{F}'$ by its infinite subset, we may assume that $\mathcal{F}''\not\to W'$ for all infinite subsets $\mathcal{F}''\subset \mathcal{F}'$ and
\begin{equation}\label{eqn:202110084}
f(\mathbb D)\not\subset W'
\end{equation}
for all $f\in \mathcal{F}'$.
Let $\lambda_{W'}\geq 0$ be a Weil function for $W'$ and let $\omega_{S'}$ be a smooth positive $(1,1)$-form on $S'$.
By Lemma \ref{lem:20201119}, there exist $\sigma_0\in (0,1)$ and $\alpha>0$ such that for all $\rho\in (\sigma_0,1)$, $r\in (\rho,1)$ and $f\in\mathcal{F}'$, we have
\begin{equation}\label{eqn:202110081}
m(r,f,\lambda_{W'})\leq \frac{\alpha}{r-\rho}T_{\rho
}(r,f,\omega_{S'})+\alpha.
\end{equation}
Set $Z=\supp \mathcal{Z}$ with a Weil function $\lambda_Z\geq 0$.
By $Z\subsetneqq W$, we have $Z\not\in\mathcal{P}$.
Note that $\mathrm{LIM}(\mathcal{F}';\{ E_i\})$ exists and $\mathrm{LIM}(\mathcal{F}';\{ E_i\})=\mathrm{LIM}(\mathcal{F};\{ E_i\})$.
Hence $\mathrm{LIM}(\mathcal{F}';\{ E_i\})\not\subset Z$.
Thus by the induction hypothesis, we get $\sigma_1\in (0,1)$, $\beta_1>0$, $k_1\in \mathbb Z_{\geq 0}$, $E\subset \cup E_i$ with $\lambda_E\geq 0$ and $\mathcal{G}\subset \mathcal{F}'$ such that for all $s\in(\sigma_1,1)$, $r\in (s,1)$ and $f\in\mathcal{G}$, we have
\begin{equation}\label{eqn:202110086}
m(r,f,\lambda_Z)
\leq \frac{\beta_1}{(r-s)^{k_1}}
T_{s}(r,f,\omega_S)
+\frac{\beta_1}{(r-s)^{k_1}}m\left(s,f,\lambda_E\right)
+\frac{\beta_1}{(r-s)^{k_1}}
\end{equation}
and $f(\mathbb D)\not\subset E\cup Z$.
By \eqref{eqn:202110084}, we have
\begin{equation}\label{eqn:202110087}
f(\mathbb D)\not\subset  E\cup W
\end{equation}
for all $f\in\mathcal{G}$.

Now we set $\sigma=\max\{\sigma_0,\sigma_1\}$.
Let $s\in(\sigma,1)$ and $r\in (s,1)$.
Set $\rho=(s+r)/2$.
Then by \eqref{eqn:202110212}, there exists a positive constant $c>0$ such that 
$$
T_{\rho}(r,f,\omega_{S'})\leq cT_{\rho}(r,f,\omega_{S})+cm(\rho,f,\lambda_{Z})+c
$$
for all $f\in \mathcal{G}$.
Then by \eqref{eqn:202110086} applied to $r=\rho$, we get
$$
T_{\rho}(r,f,\omega_{S'})\leq \frac{\beta_2}{(r-s)^{k_1}}
T_{s}(r,f,\omega_S)
+\frac{\beta_2}{(r-s)^{k_1}}m\left(s,f,\lambda_E\right)
+\frac{\beta_2}{(r-s)^{k_1}},
$$
where $\beta_2=\max\{ c+2^{k_1}c\beta_1, 2^{k_1}c\beta_1\}$.
Combining this with \eqref{eqn:202110081}, we get
\begin{equation}\label{eqn:20220603}
m(r,f,\lambda_{W'})\leq \frac{2\alpha\beta_2}{(r-s)^{k_1+1}}
T_{s}(r,f,\omega_S)
+\frac{2\alpha\beta_2}{(r-s)^{k_1+1}}m\left(s,f,\lambda_E\right)
+\frac{\alpha+2\alpha\beta_2}{(r-s)^{k_1+1}}.
\end{equation}
There exists a positive constant $c'>0$ such that $m(r,f,\lambda_W)\leq c'm(r,f,\lambda_{Z})+m(r,f,\lambda_{W'})$ for all $f\in \mathcal{G}$.
Combining this with \eqref{eqn:202110086} and \eqref{eqn:20220603}, we get for all $f\in\mathcal{G}$,
\begin{equation*}
m(r,f,\lambda_W)
\leq \frac{\beta}{(r-s)^k}
T_{s}(r,f,\omega_S)
+\frac{\beta}{(r-s)^k}m\left(s,f,\lambda_E\right)
+\frac{\beta}{(r-s)^k}.
\end{equation*}
Here we set $\beta=\alpha+c'\beta_1+2\alpha\beta_2$ and $k=k_1+1$.
This and \eqref{eqn:202110087} imply $W\not\in\mathcal{P}$.

Now in both cases above, we have $W\not\in \mathcal{P}$.
This is a contradiction.
Thus $\mathcal{P}=\emptyset$.
We have proved our lemma.
\hspace{\fill} $\square$

\subsection{The case of semi-abelian varieties}

We recall $E_{k,A,A/B}\subset P_{k,A}$ from Definition \ref{defn:20201225} and $\Pi(\mathcal{G})$ from Definition \ref{defn:20220601}.
We use the convention from \eqref{eqn:20230306} in the proof.

\begin{lem}\label{lem:1.2cor}
Let $\overline{A}$ be a smooth projective equivariant compactification. 
Let $\omega_{\overline{A}}$ be a smooth positive $(1,1)$-form on $\overline{A}$.
Let $\mathcal{F}\subset \mathrm{Hol}(\mathbb D,A)$ be an infinite set such that $\mathcal{F}$ satisfies Assumption \ref{ass:20201109}.
Let $B\subset A$ be a semi-abelian subvariety such that $B\in \Pi(\mathcal{F})$ and that $\varpi_{B}\circ f$ is non-constant for all $f\in\mathcal{F}$, where $\varpi_B:A\to A/B$ is the quotient.
Let $k\in\mathbb Z_{\geq 1}$.
Let $\lambda_{E_{k,A,A/B}}\geq 0$ be a Weil function for $E_{k,A,A/B}\subset P_{k,A}$.
Then there exists $\sigma\in (0,1)$ with the following property:
Let $s\in (\sigma,1)$, $\varepsilon>0$ and $\delta>0$.
Then there exists a positive constant $\beta>0$ such that for all $f\in \mathcal{F}$, the estimate
\begin{equation*}
m(r,f_{P_{k,A}},\lambda_{E_{k,A,A/B}})
\leq 
\varepsilon T_{s}(r,f,\omega_{\overline{A}})
+\beta
\end{equation*}
holds for all $r\in (s,1)$ outside some exceptional set of linear measure less than $\delta$.
\end{lem}

{\it Proof.}\
We first take $\sigma_0\in (0,1)$ as follows.
Let $\omega_{A/B}$ be an invariant positive $(1,1)$ form on $A/B$.
For each $n\in \mathbb Z_{\geq 2}$, let $\mathcal{E}_n\subset \mathcal{F}$ be the set of $f\in \mathcal{F}$ such that 
$$\sup_{z\in \mathbb D(1-\frac{1}{n})}|(\varpi_B\circ f)'(z)|_{\omega_{A/B}}\leq \frac{1}{n}.$$
Then we have $\mathcal{E}_2\supset \mathcal{E}_3\supset\cdots$.
By the assumption $B\in\Pi (\mathcal{F})$, we may take $n_0\in \mathbb Z_{\geq 2}$ such that $\mathcal{E}_{n_0}$ is finite.
Since $\varpi_{B}\circ f$ is non-constant for all $f\in\mathcal{F}$, we may take $n_1\geq n_0$ such that $\mathcal{E}_{n_1}=\emptyset$.
We take $\sigma_0\in (0,1)$ such that $1-\frac{1}{n_1}\leq \sigma_0<1$.
Then we have
$$
\sup_{z\in \mathbb D(\sigma_0)}|(\varpi_B\circ f)'(z)|_{\omega_{A/B}}\geq 1-\sigma_0$$
for all $f\in \mathcal{F}$.
For each $f\in\mathcal{F}$, we take $w_f\in \overline{\mathbb D(\sigma_0)}$ such that 
\begin{equation}\label{eqn:202303075}
|(\varpi_B\circ f)'(w_f)|_{\omega_{A/B}}\geq 1-\sigma_0.
\end{equation}

Next we take $\sigma\in (\sigma_0,1)$ as follows.
Let $\omega_{P_{k-1,A}}$ be a smooth positive $(1,1)$-form on $P_{k-1,A}$.
We note that by $B\in \Pi(\mathcal{F})$, we have $\{ 0\}\in \Pi(\mathcal{F})$.
We first apply Lemma \ref{lem:tcor} to get $\sigma_1\in (\sigma_0,1)$ with the following property:
For $s\in (\sigma_1,1)$, $\varepsilon>0$ and $\delta>0$, there exists a positive constant $\mu_1=\mu_1(s,\varepsilon,\delta)>0$ such that for all $f\in\mathcal{F}$, we have
\begin{equation}\label{eqn:20230314}
T_{s}(r,f_{P_{k-1,A}},\omega_{P_{k-1,A}})
\leq \frac{\varepsilon}{2} 
T_{s}(r,f,\omega_{\overline{A}})
+\mu_1
\end{equation}
for all $r\in (s,1)$ outside some exceptional set of linear measure less than $\delta/2$.
Let $\omega_{A}$ be an invariant positive $(1,1)$ form on $A$.
By Lemma \ref{lem:20201102}, we get $\sigma\in (\sigma_1,1)$ with the following property:
For $s\in (\sigma,1)$, $\varepsilon>0$ and $\delta>0$, there exists a positive constant $\mu_2=\mu_2(s,\varepsilon,\delta)>0$ such that for all $f\in\mathcal{F}$ and $w\in\mathbb D(\sigma)$, we have
\begin{equation}\label{eqn:202303141}
\int_0^{2\pi}\log |(f_{[k-1]}\circ\varphi_{w,r})'
(re^{i\theta})|_{\omega_{A\times P_{k-1,A}}}\frac{d\theta}{2\pi}
\leq 
\frac{\varepsilon}{2}T_s(r,f_{[k-1]},\omega_{\overline{A}\times P_{k-1,A}})+\mu_2
\end{equation}
for all $r\in (s,1)$ outside some exceptional set of linear measure less than $\delta/2$.
Here we recall $\varphi_{w,r}:\mathbb D(r)\to \mathbb D(r)$ from \eqref{eqn:20210716}.

Let $s\in (\sigma,1)$, $\varepsilon>0$ and $\delta>0$.
We set $\varepsilon'=\varepsilon\frac{s-\sigma}{s+\sigma}$.
Let $f\in\mathcal{F}$.
Given $r\in (s,1)$, we set $g=f\circ\varphi_{w_f,r}$.
Then we have $g_{[k-1]}=f_{[k-1]}\circ\varphi_{w_f,r}$.
Hence by \eqref{eqn:20230314} and \eqref{eqn:202303141} both applied to $\varepsilon'$, we get
\begin{equation}\label{eqn:20211016}
\int_0^{2\pi}\log |(g_{[k-1]})'
(re^{i\theta})|_{\omega_{A\times P_{k-1,A}}}\frac{d\theta}{2\pi}
\leq 
\varepsilon'
T_s(r,f,\omega_{\overline{A}}) 
+\mu_1'+\mu_2'
\end{equation}
for all $r\in (s,1)$ outside some exceptional set of linear measure less than $\delta$.
Here we set $\mu_1'=\mu_1(s,\varepsilon',\delta)$ and $\mu_2'=\mu_2(s,\varepsilon',\delta)$.

Let $p:A\times P_{k-1,A}\to A/B$ be the composite of the first projection and $\varpi_B$.
Then there exists a positive constant $c>0$ such that for all $v\in T(A\times P_{k-1,A})$ with $[v]\in A\times P_{k,A}\subset PT(A\times P_{k-1,A})$, we have 
$$
\lambda_{E_{k,A,A/B}}([v]_{P_{k,A}})\leq \log \left(\frac{|v|_{\omega_{A\times P_{k-1,A}}}
}
{
| p'(v)|_{\omega_{A/B}}}\right)+c,
$$
where $[v]_{P_{k,A}}\in P_{k,A}$ is the image of $[v]\in A\times P_{k,A}$ under the second projection $A\times P_{k,A}\to P_{k,A}$.
Hence for all $f\in\mathcal{F}$ and $r\in (s,1)$, we have
\begin{equation*}
m(r,g_{P_{k,A}},\lambda_{E_{k,A,A/B}})
\leq
\int_0^{2\pi}
\log \frac{|(g_{[k-1]})'
(re^{i\theta})|_{\omega_{A\times P_{k-1,A}}}
}
{
| (p\circ g_{[k-1]})'(re^{i\theta})|_{\omega_{A/B}}}
\frac{d\theta}{2\pi}+c.
\end{equation*}
Since $\lambda_{E_{k,A,A/B}}\geq 0$ and $g_{[k]}=f_{[k]}\circ\varphi_{w_f,r}$, Lemma \ref{lem:20210711} yields that
\begin{equation*}
\begin{split}
m(r,f_{P_{k,A}},\lambda_{E_{k,A,A/B}})
&\leq \frac{r+\sigma}{r-\sigma} m(r,g_{P_{k,A}},\lambda_{E_{k,A,A/B}})\\
&\leq \frac{s+\sigma}{s-\sigma} m(r,g_{P_{k,A}},\lambda_{E_{k,A,A/B}})
\end{split}
\end{equation*}
for all $r\in (s,1)$.
Hence we have
\begin{equation*}
m(r,f_{P_{k,A}},\lambda_{E_{k,A,A/B}})
\leq
\frac{s+\sigma}{s-\sigma}\int_0^{2\pi}
\log \frac{|(g_{[k-1]})'
(re^{i\theta})|_{\omega_{A\times P_{k-1,A}}}
}
{
| (p\circ g_{[k-1]})'(re^{i\theta})|_{\omega_{A/B}}}
\frac{d\theta}{2\pi}+c\frac{s+\sigma}{s-\sigma}
\end{equation*}
for all $r\in (s,1)$.
Combining this estimate with \eqref{eqn:20211016}, we get
\begin{multline}\label{eqn:202303142}
m(r,f_{P_{k,A}},\lambda_{E_{k,A,A/B}})
\leq \varepsilon
T_s(r,f,\omega_{\overline{A}}) +\frac{s+\sigma}{s-\sigma}\int_0^{2\pi}
\log \frac{1
}
{
| (p\circ g_{[k-1]})'(re^{i\theta})|_{\omega_{A/B}}}
\frac{d\theta}{2\pi}\\
+\frac{s+\sigma}{s-\sigma}
(\mu_1'+\mu_2'+c)
\end{multline}
for all $r\in (s,1)$ outside some exceptional set of linear measure less than $\delta$.

Next, by the subharmonicity of $\log |(\varpi_B\circ g)'|_{\omega_{A/B}}$, we get
\begin{equation*}
\begin{split}
\int_0^{2\pi}
\log \frac{1
}
{
| (p\circ g_{[k-1]})'(re^{i\theta})|_{\omega_{A/B}}}
\frac{d\theta}{2\pi}
=& 
\int_0^{2\pi}
\log \frac{1
}
{
| (\varpi_B\circ g)'(re^{i\theta})|_{\omega_{A/B}}}
\frac{d\theta}{2\pi}
\\
\leq& 
\log \frac{1
}
{
| (\varpi_B\circ g)'(0)|_{\omega_{A/B}}}
\end{split}
\end{equation*}
for all $r\in (s,1)$.
By \eqref{eqn:202107161} and \eqref{eqn:202303075}, we have
$$
 | (\varpi_B\circ g)'(0)|_{\omega_{A/B}}=| (\varpi_B\circ f)'(w_f)|_{\omega_{A/B}}|\varphi_{w_f,r}'(0)|\geq \frac{(1-\sigma)(s-\sigma)}{s+\sigma}.
$$
Hence
\begin{equation}\label{eqn:202206034}
\int_0^{2\pi}
\log \frac{1
}
{
| (p\circ g_{[k-1]})'(re^{i\theta})|_{\omega_{A/B}}}
\frac{d\theta}{2\pi}\leq \log\left(\frac{s+\sigma}{(1-\sigma)(s-\sigma)}\right)
\end{equation}
for all $r\in (s,1)$.
Using \eqref{eqn:202303142} and \eqref{eqn:202206034} and letting
$$
\beta=\frac{s+\sigma}{s-\sigma}\left(\mu_1'+\mu_2'+\log\left(\frac{s+\sigma}{(1-\sigma)(s-\sigma)}\right)+c\right),
$$
we complete the proof.
\hspace{\fill} $\square$

\begin{lem}\label{lem:20210302}
Let $k\in\mathbb Z_{\geq 1}$.
Let $\overline{A}$ be a smooth projective equivariant compactification. 
Let $\omega_{\overline{A}}$ be a smooth positive $(1,1)$-form on $\overline{A}$.
Let $W\subset P_{k,A}$ be a closed subscheme with a Weil function $\lambda_W\geq 0$.
Let $\mathcal{F}\subset \mathrm{Hol}(\mathbb D,A)$ be an infinite set  such that $\mathcal{F}$ satisfies Assumption \ref{ass:20201109} and $\{0\}\in \Pi(\mathcal{F})$.
Let $\Pi'\subset \Pi(\mathcal{F})$ be a subset such that $\mathrm{LIM}(\mathcal{F}_{P_{k,A}};\{ E_{k,A,A/B}\}_{B\in \Pi'})$ exists, where $\mathcal{F}_{P_{k,A}}=(f_{P_{k,A}})_{f\in\mathcal{F}}$.
Assume that $\mathrm{LIM}(\mathcal{F}_{P_{k,A}};\{ E_{k,A,A/B}\}_{B\in \Pi'})\not\subset \supp W$.
Then there exist $\sigma\in (0,1)$ and an infinite subset $\mathcal{G}\subset \mathcal{F}$ with the following property:
Let $s\in (\sigma,1)$, $\varepsilon >0$, $\delta>0$.
Then there exists a positive constant $\beta>0$ such that, for all $f\in \mathcal{G}$, we have
$$
m(r,f_{P_{k,A}},\lambda_W)
\leq \varepsilon T_{s}(r,f_A,\omega_{\overline{A}})+
\beta
$$
for all $r\in (s ,1)$ outside some exceptional set of linear measure less than $\delta$.
\end{lem}

{\it Proof.}\
We first consider the case that the subset
\begin{equation}\label{eqn:20220423}
\{ f_{P_{k,A}}; \ f\in \mathcal{F}\}\subset \mathrm{Hol}(\mathbb D, P_{k,A})
\end{equation}
is finite.
We may choose an infinite subset $\mathcal{G}\subset \mathcal{F}$ by removing finite elements from $\mathcal{F}$ such that $f_{P_{k,A}}(\mathbb D)\not\subset \supp W$ for all $f\in\mathcal{G}$.
Then our estimate is valid for this $\mathcal{G}$ and any $\sigma\in (0,1)$.
Indeed we just need to set $\beta=\max_{f\in \mathcal{G}}\sup_{r\in (s,1-\delta)}m(r,f_{P_{k,A}},\lambda_W)$.
Hence in the following, we assume that the set \eqref{eqn:20220423} is infinite.
Replacing $\mathcal{F}$ by its infinite subset, we may assume that the set \eqref{eqn:20220423} is infinite and the map $\mathcal{F}\ni f\mapsto f_{P_{k,A}}\in \mathrm{Hol}(\mathbb D,P_{k,A})$ is injective.
Hence we may assume that the infinite indexed family $\mathcal{F}_{P_{k,A}}$ is an infinite subset of $\mathrm{Hol}(\mathbb D,P_{k,A})$.

There exists a positive constant $c>0$ such that $\lambda_W\leq c\lambda_{\supp W}+c$.
Hence we may assume that $W$ is reduced.
Let $\omega_{P_{k,A}}$ be a smooth positive $(1,1)$-form on $P_{k,A}$.
We apply Lemma \ref{lem:1.1} to obtain $\sigma\in (0, 1)$, $\beta_0>0$, $l\in\mathbb Z_{\geq 0}$, $\mathcal{G}\subset \mathcal{F}$ and a Zariski closed set $E\subset \cup_{B\in \Pi'}E_{k,A,A/B}$ with a Weil function $\lambda_E\geq 0$ such that
\begin{equation}\label{eqn:202206035}
m(r,f_{P_{k,A}},\lambda_W)\leq \frac{\beta_0}{(r-\rho)^l} T_{\rho}(r,f_{P_{k,A}},\omega_{P_{k,A}})+\frac{\beta_0}{(r-\rho)^l} m(\rho,f_{P_{k,A}},\lambda_E)+\frac{\beta_0}{(r-\rho)^l}
\end{equation}
for all $\sigma<\rho<r<1$ and $f\in\mathcal{G}$.
We may take a finite subset $\Lambda\subset \Pi'$ such that letting $E_{\Lambda}=\cup_{B\in \Lambda}E_{k,A,A/B}$, we have $E\subset E_{\Lambda}$.
Note that $\mathcal{G}$ satisfies Assumption \ref{ass:20201109}.
By removing finite elements from $\mathcal{G}$, we may assume that $\varpi_B\circ f$ is non-constant for all $B\in \Lambda$ and $f\in\mathcal{G}$.
By enlarging $\sigma\in (0,1)$, if necessary, we may assume that $\sigma$ fits in Lemma \ref{lem:1.2cor} for each $B\in \Lambda$.
Note that $\{0\}\in\Pi(\mathcal{G})$.
Hence, again by enlarging $\sigma\in (0,1)$, we may assume that $\sigma$ fits in Lemma \ref{lem:tcor}.

Now we fix $s\in (\sigma,1)$, $\varepsilon>0$ and $\delta>0$.
We set $\varepsilon'=\frac{\varepsilon\delta^l}{2^{2l+1}\beta_0}$.
We take $f\in \mathcal{G}$.
We apply Lemma \ref{lem:1.2cor} to get
\begin{equation*}
m(\rho,f_{P_{k,A}},\lambda_{E_{\Lambda}})
\leq 
\varepsilon' T_{s}(\rho,f,\omega_{\overline{A}})
+\beta_1
\end{equation*}
for all $\rho\in (s,1)$ outside some exceptional set of linear measure less than $\delta/4$.
Hence we may take $\rho\in (s,s+\delta/4)$ which satisfies this property.
Hence by \eqref{eqn:202206035}, we get
$$
m(r,f_{P_{k,A}},\lambda_{W})
\leq
\frac{4^l\beta_0}{\delta^l}
T_{s}(r,f_{P_{k,A}},\omega_{P_{k,A}})+\varepsilon'\frac{4^l\beta_0}{\delta^l} T_{s}(r,f,\omega_{\overline{A}})+\frac{4^l\beta_0(1+\beta_1)}{\delta^l}
$$
for all $r\in (s+\delta/2,1)$.
We apply Lemma \ref{lem:tcor} to get
$$
T_{s}(r,f_{P_{k,A}},\omega_{P_{k,A}})\leq \varepsilon'T_s(r,f_A,\omega_{\overline{A}})+\mu
$$
for all $r\in (s,1)$ outside some exceptional set of linear measure less than $\delta/2$.
Combining these two estimates, we get our claim.
Here we set $\beta=\frac{4^l\beta_0(1+\beta_1+\mu)}{\delta^l}$.
\hspace{\fill} $\square$

\section{Proof of Proposition \ref{pro:320}}\label{sec:8}

Let $\overline{A}$ be a smooth equivariant compactification.
Given an infinite indexed family $\mathcal{F}$ in $\mathrm{Hol}(\mathbb D,A)$, we set
$$
\mathcal{I}(\mathcal{F})=\{ D; \text{irreducible component of $\partial A$ s.t. $\exists \mathcal
{F}'$ subfamily of $\mathcal{F}$ such that $\mathcal{F}'\to D$}
\}
$$
For an irreducible component $D\subset \partial A$, let $I_D\subset A$ be the isotropy group for $D$.
Then $I_D=\mathbb G_m$.
Set $I_{\mathcal{F}}=I_{D_1}\cdots I_{D_l}\subset A$ where $\{ D_1,\ldots,D_l\}=\mathcal{I}(\mathcal{F})$.
In the next lemma, for $\tau\in A$, we use the same notation for the map $\tau:\overline{A}\to\overline{A}$ defined by $a\mapsto a+\tau$.

\begin{lem}\label{lem:202011091}
Let $\mathcal{F}=(f_i)_{i\in I}$ be an infinite indexed family in $\mathrm{Hol}(\mathbb D,A)$.
For $i\in I$, there exists $\tau_i\in I_{\mathcal{F}}$ such that $(\tau_i\circ f_i)_{i\in I}$ satisfies Assumption \ref{ass:20201109}.
\end{lem}

In the following proof, we use the notation of the 1-dimensional unlimited Hausdorff content defined as follows.
Let $K\subset \mathbb C$, we set
\begin{equation}\label{20220604}
C_H^1(K)=\inf\left\{ \sum_jr_j; \text{$\exists$a countable cover of $K$ by closed discs with radii $r_j>0$}\right\}.
\end{equation}

{\it Proof of Lemma \ref{lem:202011091}.}\
If $\mathcal{I}(\mathcal{F})=\emptyset$, then our lemma is trivial.
Hence we assume $\mathcal{I}(\mathcal{F})\not=\emptyset$.
We take $D\in \mathcal{I}(\mathcal{F})$.
For each $i\in I$, we shall show that there exists $\tau_{i,D}\in I_D$ such that the indexed family $\mathcal{F}'=( \tau_{i,D}\circ f_i)_{i\in I}$ satisfies the following two properties:
\begin{itemize}
\item $D\not\in \mathcal{I}(\mathcal{F}')$.
\item If $E\not\in\mathcal{I}(\mathcal{F})$, then $E\not\in\mathcal{I}(\mathcal{F}')$.
\end{itemize}

We are going to construct $\tau_{i,D}$.
For each $E\not\in \mathcal{I}(\mathcal{F})$, there exist $s_E\in (0,1)$, $\gamma_E>0$, an open neighbourhood $U_E\subset \overline{A}$ of $E$, and a finite subset $J_E\subset I$ such that
\begin{equation}\label{eqn:202206047}
C_H^1(\mathbb D(s_E)\backslash f_i^{-1}(U_E))\geq \gamma_E
\end{equation}
for all $i\in I\backslash J_E$ (cf. Remark \ref{rem:20210803}).
We may assume moreover that if $E\cap D=\emptyset$, then $U_E\Subset \overline{A}-D$.
We set $s=\max_{E\not\in\mathcal{I}(\mathcal{F})}\{ s_E\}$, $\gamma=\min_{E\not\in\mathcal{I}(\mathcal{F})}\{ \gamma_E\}$ and $\gamma'=\gamma/2$.

We apply Lemma \ref{lem:vect} to get $p:W\to D$, where $W\subset \overline{A}$ is a Zariski open neighbourhood of $D$ and $W$ is a total space of a line bundle over $D$.
Let $||\cdot||$ be a smooth Hermitian metric on this line bundle $p:W\to D$.
We set $U_D=\{ x\in W;||x||<1\}$.
Then $U_D\subset \overline{A}$ is an open neighbourhood of $D$.
By replacing $||\cdot||$ if necessary, we may assume that $U_D\cap U_E=\emptyset$ if $D\cap E=\emptyset$.
Note that $I_D$ acts on $W$ by the scholar product (cf. Lemma \ref{lem:vect}).
Hence for $a\in I_D$, we may define its norm $|a|$ by $|a|=||a\cdot x||/||x||$ for $x\in W-D$.
We set $U_D'=\{ x\in W;||x||<1/2\}$.
Then $U_D'\Subset U_D$.

Let $i\in I$.
We choose $\tau_{i,D}\in I_D$ in the two cases below.

\medskip

{\it Case 1.}\
$f_i(z)\in U_D'$ for $\gamma'$-almost all $z\in \mathbb D(s)$.
We set
$$
\sigma_i=\sup\{ |\tau| ;\  \text{$\tau\in I_D$ and $\tau\circ f_i(z)\in U_D'$ for $\gamma'$-almost all $z\in \mathbb D(s)$}\}.
$$
Note that if $\tau\in I_D$ satisfies $|\tau|\min_{z\in\overline{\mathbb D(s)}}||f_i(z)||\geq 1/2$, then $\tau\circ f_i(z)\not\in U_D'$ for all $z\in \mathbb D(s)$.
Hence $\sigma_i<\infty$.
We take $\tau_{i,D}\in I_D$ such that $|\tau_{i,D}|=\sigma_i$.
Then we have
\begin{equation}\label{eqn:202206046}
|\tau_{i,D}|\geq 1.
\end{equation}
This follows from $\sigma_i\geq |e_{I_D}|$, where $e_{I_D}$ is the identity element of $I_D$.
\medskip

{\it Case 2.}\
Otherwise, we set $\tau_{i,D}=e_{I_D}$.
\medskip

We set $\Omega_i=\mathbb D(s)\cap (\tau_{i,D}\circ f_i)^{-1}(U_D)$.
We define $G\subset I$ by the set of $i\in I$ such that $f_i$ belongs to the case 1.
By the definition of $\tau_{i,D}$, we have 
\begin{equation}\label{eqn:20210727}
C_H^1(\mathbb D(s)\backslash \Omega_i)<\gamma'
\end{equation}
for all $i\in G$.
We show that $\mathcal{F}'=(\tau_{i,D}\circ f_i)_{i\in I}$ satisfies our requirement.
We first observe that $D\not\in\mathcal{I}(\mathcal{F}')$.
To show this, let $U_D''=\{ x\in W;||x||<1/4\}$ so that $U_D''\Subset U_D'$.
Then we have $C_H^1(\mathbb D(s)\backslash (\tau_{i,D}\circ f_i)^{-1}(U_D''))\geq \gamma'$ for all $i\in I$.
Hence $D\not\in\mathcal{I}(\mathcal{F}')$.

Next we take $E\not\in\mathcal{I}(\mathcal{F})$.
We first consider the case $E\cap D=\emptyset$.
By $U_E\cap U_D=\emptyset$ and \eqref{eqn:20210727}, we have $C_H^1(\mathbb D(s)\cap (\tau_{i,D}\circ f_i)^{-1}(U_E))< \gamma'$, provided $i\in G$.
By \eqref{eqn:202206047}, we have $C_H^1(\mathbb D(s_E)\backslash (\tau_{i,D}\circ f_i)^{-1}(U_E))\geq \gamma_E$ for all $i\in I\backslash (G\cup J_E)$.
Hence $E\not\in\mathcal{I}(\mathcal{F}')$ as desired.

In the following, we consider the case $E\cap D\not=\emptyset$.
We take an open neighbourhood $U_E'\subset U_E$ of $E$ so that $U_E'\cap U_D\subset U_D-p^{-1}(p(\overline{U_D}-U_E))$.
Then for $\tau\in I_D$ with $|\tau|\geq 1$, we have 
\begin{equation}\label{eqn:202206043}
U_E'\cap U_D\subset \tau(U_E).
\end{equation}
By $\tau_{i,D}\circ f_i(\Omega_i)\subset U_D$, we have
$$
(\tau_{i,D}\circ f_i)^{-1}(U_E')\cap \Omega_i
\subset (\tau_{i,D}\circ f_i)^{-1}(U_E'\cap U_D).
$$
By \eqref{eqn:202206046} and \eqref{eqn:202206043}, we have
\begin{equation*}
(\tau_{i,D}\circ f_i)^{-1}(U_E'\cap U_D)\subset (\tau_{i,D}\circ f_i)^{-1}(\tau_{i,D}(U_E))=
f_i^{-1}(U_E).
\end{equation*}
Hence we get
\begin{equation}\label{eqn:20210728}
(\tau_{i,D}\circ f_i)^{-1}(U_E')\cap \Omega_i\subset f_i^{-1}(U_E).
\end{equation}
By \eqref{eqn:202206047} and the definitions of $s$ and $\gamma$, we have $C_H^1(\mathbb D(s)\backslash f_i^{-1}(U_E))\geq \gamma$ for all $i\in I\backslash J_E$.
Then \eqref{eqn:20210727} yields that $C_H^1(\Omega_i\backslash f_i^{-1}(U_E))\geq \gamma'$ for all $i\in G\backslash J_E$.
Hence by \eqref{eqn:20210728}, we have $C_H^1(\Omega_i\backslash (\tau_{i,D}\circ f_i)^{-1}(U_E'))\geq \gamma'$, thus $C_H^1(\mathbb D(s)\backslash (\tau_{i,D}\circ f_i)^{-1}(U_E'))\geq \gamma'$ for all $i\in G\backslash J_E$.
By \eqref{eqn:202206047}, the same holds for all $i\in I\backslash J_E$.
Thus $E\not\in\mathcal{I}(\mathcal{F}')$ as desired.

Now starting from $\mathcal{F}$, we set $\mathcal{F}_1=\mathcal{F}'$.
Inductively, we get $\mathcal{F}_2,\mathcal{F}_3,\ldots$ by $\mathcal{F}_{k+1}=\mathcal{F}_k'$.
Then by the above consideration, we have $\mathcal{I}(\mathcal{F}_k)\supsetneqq \mathcal{I}(\mathcal{F}_{k+1})$.
Hence we get $\mathcal{I}_{\mathcal{F}_n}=\emptyset$ for some $n$.
Then $\mathcal{F}_n$ satisfies the assumption \ref{ass:20201109}.
\hspace{\fill} $\square$

\begin{lem}\label{lem:20210303}
Let $p:W\to V$ be a vector bundle, where $V$ is a smooth projective variety.
We consider $V\subset W$ by the zero-section.
Let $\mathcal{F}=(f_i)_{i\in I}$ be an infinite indexed family in $\mathrm{Hol}(\mathbb D,W)$.
Assume that $\mathcal{F}\to V$ and that $\{ p\circ f_i\}_{i\in I}$ converges uniformly on compact subsets of $\mathbb D$ to $g:\mathbb D\to V$.
Then $\mathcal{F}$ converges uniformly on compact subsets of $\mathbb D$ to $g$.
\end{lem}

{\it Proof.}\
We equip $W$ with a smooth Hermitian metric.
We denote by $||\cdot||_W$ the associated Euclidean norm on each fiber of $p:W\to V$.
Let $O_1,\ldots,O_n\subset V$ be an open set in the standard topology of $V$ such that $V\subset O_1\cup\cdots\cup O_n$.
Assume that $W|_{O_j}=\mathbb C^k\times O_j$ and that letting $\varphi_j:W|_{O_j}\to \mathbb C^k$ be the first projection, there exists a positive constant $c_j>1$ such that
\begin{equation}\label{eqn:20210729}
\frac{1}{c_j}||\varphi_j(x)||\leq ||x||_W\leq c_j||\varphi_j(x)||
\end{equation}
for all $x\in W|_{O_j}$.
Let $\delta$ be a distance function on $V$ which determines the standard topology on $V$.
Let $O_j'\Subset O_j$ be an open set such that $V\subset O_1'\cup\cdots\cup O_l'$.
There exists a positive constant $\alpha>0$ such that for all $j=1,\ldots,n$, we have
\begin{equation}\label{eqn:202107313}
\delta(\overline{O_j'},V-O_j)>\alpha.
\end{equation}

Let $\tau:V\times W\to \mathbb R_{\geq 0}$ be defined by 
$$
\tau(x,y)=\delta(x,p(y))+||y||_W.
$$
Let $r\in (0,\alpha/2)$ and $s\in (0,1)$.
We claim that there exists a finite set $E\subset I$ such that, for all $i\in I\backslash E$ and $z\in \mathbb D(s)$, we have
\begin{equation}\label{eqn:20220607}
\tau(g(z),f_i(z))<r.
\end{equation}
We prove this.
Let $s'\in (s,1)$.
Since $\{ p\circ f_i\}_{i\in I}$ converges uniformly on compact subsets of $\mathbb D$ to $g$, there exists a finite subset $E_1\subset I$ such that for all $i\in I\backslash E_1$ and $z\in \mathbb D(s')$, we have
\begin{equation}\label{eqn:202107312}
\delta(g(z),p\circ f_i(z))<r/2.
\end{equation}
Set $c=\max_{1\leq j\leq n}\{ c_j\}$ and
$$
U=\{ x\in W;\ ||x||_W<r/2c^2\}.
$$
Then since $||\cdot ||_W:W\to \mathbb R_{\geq 0}$ is continuous, $U\subset W$ is an open neighbourhood of $V\subset W$.
We take $\gamma\in (0,\frac{s'-s}{2})$ such that if $\Delta\subset \mathbb D(s')$ is a disc of radius $\gamma$, then for all $z,z'\in \Delta$, we have 
\begin{equation}\label{eqn:202107311}
\delta(g(z),g(z'))<\alpha/2.
\end{equation}
By $\mathcal{F}\to V$, there exists a finite subset $E_2\subset I$ such that, for all $i\in I\backslash E_2$, we have 
\begin{equation}\label{eqn:20220606}
C_H^1(\mathbb D(s')\backslash f_i^{-1}(U))<\gamma',
\end{equation}
where $\gamma'=\gamma/4$.
Set $E=E_1\cup E_2$.
For all $i\in I\backslash E$ and $z\in\mathbb D(s)$, we claim
\begin{equation}\label{eqn:202107315}
||f_i(z)||_W<r/2.
\end{equation}
To show this, we take $z\in\mathbb D(s)$ and $i\in I\backslash E$.
Given $\rho\in(0,\gamma)$, we denote by $\Delta_{\rho}$ the disc of radius $\rho$ centered at $z$.
By $\gamma<\frac{s'-s}{2}$, we have $\overline{\Delta_{\rho}}\subset \mathbb D(s')$.
We take $O_j'$ such that $g(z)\in O_j'$.
To show $p\circ f_i(\Delta_{\rho})\subset O_j$, we take $w\in \Delta_{\rho}$.
Then by \eqref{eqn:202107311}, we have $\delta(g(z),g(w))<\alpha/2$.
By $w\in \mathbb D(s')$, the estimate \eqref{eqn:202107312} yields $\delta(g(w),p\circ f_i(w))<\alpha/4$.
Hence we have $\delta(g(z),p\circ f_i(w))<\alpha$.
Hence by \eqref{eqn:202107313}, we have $p\circ f_i(w)\in O_j$.
This shows $p\circ f_i(\Delta_{\rho})\subset O_j$, thus $f_i(\Delta_{\rho})\subset W|_{O_j}$.
By \eqref{eqn:20220606}, there exists $\rho\in (0,\gamma)$ such that $f_i(\partial \Delta_{\rho})\subset U$.
By \eqref{eqn:20210729}, we have $||\varphi_j\circ f_i(w)||<r /2c_j$  for all $w\in \partial \Delta_{\rho}$.
Since $||\varphi_j\circ f_i||$ is subharmonic on $\Delta_{\rho}$, the maximum principle yields $||\varphi_j\circ f_i(z)||<r /2c_j$. 
By \eqref{eqn:20210729}, we get \eqref{eqn:202107315}.
Thus by \eqref{eqn:202107312} and \eqref{eqn:202107315}, we get \eqref{eqn:20220607}.

For $r>0$ and $x\in V$, we set $B_x(r)=\{ y\in W;\ \tau (x,y)<r\}$.
Note that $\tau(x,y)$ is continuous with respect to $y\in W$.
Hence $B_x(r)\subset W$ is an open subset.
Let $d$ be a distance function on $W$ which induces the topology on $W$.
Let $\varepsilon>0$.
Then we may take a positive constant $r_{\varepsilon}>0$ such that 
\begin{equation}\label{eqn:202206062}
\sup_{y\in B_x(r_{\varepsilon})}d(x,y)<\varepsilon
\end{equation}
for all $x\in V$.
To prove this, we note that for each $x\in V$, we may take $r_{\varepsilon,x}>0$ such that $\sup_{y\in B_x(r_{\varepsilon,x})}d(x,y)<\varepsilon/4$.
We consider the open covering $\{ B_x(r_{\varepsilon,x}/2)\}_{x\in V}$ of $V\subset W$.
Since $V$ is compact, there exist $x_1,\ldots,x_l\in V$ such that $V\subset B_{x_1}(r_{\varepsilon,x_1}/2)\cup \cdots\cup B_{x_l}(r_{\varepsilon,x_l}/2)$.
We set $r_{\varepsilon}=\min_{1\leq j\leq l}\{ r_{\varepsilon,x_j}/2\}$.
Let $x\in V$.
Then there exists $x_j$ such that $x\in B_{x_j}(r_{\varepsilon,x_j}/2)$.
Let $y\in B_{x}(r_{\varepsilon})$.
Then we have
$$\tau(x_j,y)\leq \delta(x_j,x)+\tau(x,y)<r_{\varepsilon,x_j}.$$
Hence $y\in B_{x_j}(r_{\varepsilon,x_j})$.
Thus $d(x,y)\leq d(x_j,x)+d(x_j,y)<\varepsilon/2$.
This proves \eqref{eqn:202206062}.

Now we prove that $\mathcal{F}$ converges uniformly on compact subsets of $\mathbb D$ to $g$.
Let $s\in (0,1)$ and let $\varepsilon>0$.
Then by \eqref{eqn:20220607} applied to $r=\min\{\alpha/3,r_{\varepsilon}\}$, there exists a finite subset $E\subset I$ such that for all $i\in I\backslash E$, we have
$$
\tau(g(z),f_i(z))<r_{\varepsilon}
$$
for all $z\in \mathbb D(s)$.
By \eqref{eqn:202206062}, we have $d(g(z),f_i(z))<\varepsilon$ for all $z\in\mathbb D(s)$ and $i\in I\backslash E$.
Thus $\mathcal{F}$ converges uniformly on compact subsets of $\mathbb D$ to $g$.
\hspace{\fill} $\square$

\begin{lem}\label{lem:compact}
Let $\overline{A}$ be a smooth equivariant compactification of $A$ and let $\omega_{\overline{A}}$ be a smooth positive $(1,1)$-form on $\overline{A}$.
Let $\mathcal{F}\subset \mathrm{Hol}(\mathbb D, A)$ be an infinite set of holomorphic maps which satisfies Assumption \ref{ass:20201109}.
Assume that for all $0<r<1$, we have
$$
\sup_{f\in\mathcal{F}}\left\{ \int_{\mathbb D(r)}f^*\omega_{\overline{A}}\right\}
 <\infty.
$$
Then there exists an infinite subset $\mathcal{G}\subset \mathcal{F}$ such that $\mathcal{G}$ converges uniformly on compact subsets of $\mathbb D$ to some $g:\mathbb D\to A$.
\end{lem}

{\it Proof.}\
The proof devides into three steps.

{\it Step 1.}\
We first consider the case that $A$ is an abelian variety.
In this case, we have $\overline{A}=A$.
We may assume that $\omega_{A}$ is a positive invariant $(1,1)$-form.
Let $\sigma\in (0,1)$.
We shall show that $\mathcal{F}$ is equi-continuous on $\mathbb D(\sigma)$.
Namely we prove
$$
\sup_{w\in \mathbb D(\sigma)}\sup_{f\in\mathcal{F}}\left\{ |f'(w)|_{\omega_A}\right\}<\infty.
$$
Let $s\in (\sigma,1)$.
We set 
$$
c=\sup_{f\in\mathcal{F}}\left\{ \int_{\mathbb D(s)}f^*\omega_{A}\right\}.
$$
Let $w\in \mathbb D(\sigma)$.
Recall $\varphi_{w,s}:\mathbb D(s)\to \mathbb D(s)$ from \eqref{eqn:20210716}.
Then
\begin{equation*}
\int_{\mathbb D(s)}(f\circ\varphi_{w,s})^*\omega_{A}
=\int_{\mathbb D(s)}f^*\omega_{A}\leq c.
\end{equation*}
For $r\in (0,s)$, we set
$$
u(r)=\int_{0}^{2\pi}|(f\circ\varphi_{w,s})'(re^{i\theta})|_{\omega_A}^2d\theta.
$$
Since $|(f\circ\varphi_{w,s})'(z)|_{\omega_A}^2$ is subharmonic on $\mathbb D(s)$, we have $2\pi |(f\circ\varphi_{w,s})'(0)|_{\omega_A}^2\leq u(r)$.
Hence 
$$
\int_{\mathbb D(s)}(f\circ\varphi_{w,s})^*\omega_{A}=\int_0^su(r)rdr\geq \pi s^2 |(f\circ\varphi_{w,s})'(0)|_{\omega_A}^2.
$$
This shows
$$
|(f\circ\varphi_{w,s})'(0)|_{\omega_A}\leq \sqrt{\frac{c}{\pi s^2}}.
$$
Hence by \eqref{eqn:202107161}, we have
$$
|f'(w)|_{\omega_A}\leq \frac{s+\sigma}{s-\sigma}\sqrt{\frac{c}{\pi s^2}}.
$$
Hence $\mathcal{F}$ is equi-continuous on $\mathbb D(\sigma)$.
By Arzel{\'a}-Ascoli theorem, there exists a subsequence $\mathcal{G}\subset \mathcal{F}$ such that $\mathcal{G}$ converges uniformly on compact subsets of $\mathbb D$ to some $g:\mathbb D\to A$.

{\it Step 2.}\
We consider the case $A=(\mathbb{G}_m)^k$.
We are given an equivariant compactification $\overline{A}$.
Then by Lemma \ref{lem:2020110211}, there exists $s_0\in (0,1)$
such that for all $r\in (s_0,1)$, we have
\begin{equation}\label{eqn:20201112}
\sup_{f\in \mathcal{F}}\left\{ m(r,f,\lambda_{\partial A})\right\}<\infty.
\end{equation}
Let $p_i:A\to \mathbb G_m$ be the $i$-th projection.
Let $\tau_i:W_i\to \bar{A}$ be a smooth modification such that $p_i:A\to \mathbb G_m$ extends to $\overline{p_i}:W_i\to \mathbb P^1$.
Then we have 
$$\supp \overline{p_i}^*\left( (0)+(\infty)\right)\subset \supp \tau_i^*\partial A.$$
Hence by \eqref{eqn:20201112}, we get
\begin{equation}\label{eqn:202108011}
\sup_{f\in \mathcal{F}}\left\{ m(r,p_i\circ f,\lambda_{(0)+(\infty)})\right\}<\infty
\end{equation}
for all $r\in (s_0,1)$.
Note that $m(r,p_i\circ f,\lambda_{(0)+(\infty)})$ is an increasing function on $r$.
Hence \eqref{eqn:202108011} holds for all $r\in (0,1)$.
Thus by Montel's theorem (cf. \cite[Thm 1.6, p. 230]{Lan} or \cite[p. 233]{Dra}), there exists an infinite subset $\mathcal{G}_1\subset \mathcal{F}$ such that $( p_1\circ f)_{f\in \mathcal{G}_1}$ converges uniformly on compact subsets of $\mathbb D$ to $g_1:\mathbb D\to \mathbb C$.
By \eqref{eqn:202108011}, we have $g_1\not\equiv 0$.
Hence we have $g_1:\mathbb D\to\mathbb G_m$.
By the same argument, there exists an infinite subset $\mathcal{G}_2\subset \mathcal{G}_1$ such that $(p_2\circ f)_{f\in \mathcal{G}_2}$ converges uniformly on compact subsets of $\mathbb D$ to $g_2:\mathbb D\to \mathbb C$.
Again we have $g_2:\mathbb D\to\mathbb G_m$.
In this way, we get $\mathcal{G}_1\supset\mathcal{G}_2\supset\cdots\supset\mathcal{G}_k$.
We set $\mathcal{G}=\mathcal{G}_k$.
Then $\mathcal{G}$ is an infinite set and converges uniformly on compact subsets of $\mathbb D$ to $(g_1,\ldots,g_k):\mathbb D\to (\mathbb G_m)^k$.

{\it Step 3.}\
We consider the general case with $0\to T\to A\overset{p}{\to} A_0\to 0$.
By Lemma \ref{lem:20200906}, there exists an equivariant compactification $\overline{T}$ such that $\overline{A}=(\overline{T}\times A)/T$.
We take the universal covering $\tilde{A_0}\to A_0$.
We consider $0\to T\to A\times_{A_0}\tilde{A_0}\overset{\tilde{p}}{\to} \tilde{A_0}\to 0$, which has a splitting.
Then we have $\overline{A}\times_{A_0}\tilde{A_0}=\overline{T}\times \tilde{A_0}$.
Let $q:\overline{A}\times_{A_0}\tilde{A_0}\to \overline{T}$ be the first projection and $\pi:\overline{A}\times_{A_0}\tilde{A_0}\to \overline{A}$ be the natural map.
We continue to write $p:\overline{A}\to A_0$ and $\tilde{p}:\overline{A}\times_{A_0}\tilde{A_0}\to \tilde{A_0}$.

\begin{equation*}
\begin{CD}
\overline{A}@<\pi<<\overline{A}\times_{A_0}\tilde{A_0}  @>q>> \overline{T}\\
@Vp VV @VV\tilde{p}V    \\
A_0@<<< \tilde{A_0}
\end{CD}
\end{equation*}

Now we claim that there exists an infinite subset $\mathcal{F}'\subset \mathcal{F}$ such that $(p\circ f)_{f\in\mathcal{F}'}$ converges uniformly on compact subsets of $\mathbb D$.
Indeed if $\{p\circ f;f\in\mathcal{F}\}\subset \mathrm{Hol}(\mathbb D,A_0)$ is a finite subset, we take infinite subset $\mathcal{F}'\subset \mathcal{F}$ such that $p\circ f$ is all the same for all $f\in\mathcal{F}'$; otherwise we apply the first step for the infinite subset $\{p\circ f;f\in\mathcal{F}\}\subset \mathrm{Hol}(\mathbb D,A_0)$ to get an infinite subset $\mathcal{F}'\subset \mathcal{F}$ such that $(p\circ f)_{f\in\mathcal{F}'}$ converges uniformly on compact subsets of $\mathbb D$.
For each $f\in \mathcal{F}'$, we take a lifting $\tilde{f}:\mathbb D\to A\times_{A_0}\tilde{A_0}$.
We may assume that $( \tilde{p}\circ \tilde{f})_{f\in \mathcal{F}'}$ converges uniformly on compact subsets of $\mathbb D$ to $h:\mathbb D\to \tilde{A_0}$.
For each fixed $r\in (0,1)$, there exists a compact set $K_r\subset \tilde{A_0}$ such that $h(\overline{\mathbb D(r)})\subset K_r^o$, where $K_r^o$ is the interior of $K_r$.
Then for all $f\in \mathcal{F}'$ with finite exception, we have $\tilde{p}\circ\tilde{f}(\overline{\mathbb D(r)})\subset K_r^o$.
Let $\omega_{\overline{T}}$ be a smooth positive $(1,1)$-form on $\overline{T}$.
Note that there exists a positive constant $c_r>0$ such that $q^*\omega_{\overline{T}}\leq c_r\pi^*\omega_{\overline{A}}$ on $\tilde{p}^{-1}(K_r)$.
Hence we have
$$
\sup_{f\in\mathcal{F}'}\left\{ \int_{\mathbb D(r)}(q\circ\tilde{f})^*\omega_{\overline{T}}\right\}
 <\infty.
$$

Next we claim that $(q\circ \tilde{f})_{f\in \mathcal{F}'}$ satisfies the assumption \ref{ass:20201109}.
Let $D\subset \partial T$ be an irreducible component.
Then $(D\times A)/T$ is an irreducible component of $\partial A$ with $\pi(q^{-1}(D))=(D\times A)/T$.
Since $\mathcal{F}'$ satisfies Assumption \ref{ass:20201109}, there exist an open neighbourhood $(D\times A)/T\subset U_0\subset \overline{A}$, $s_0\in (0,1)$, $\gamma_0>0$ and a finite subset $\mathcal{E}\subset \mathcal{F}'$ such that 
$$C_H^1(\mathbb D(s_0)\backslash f^{-1}(U_0))\geq \gamma_0$$
for all $f\in \mathcal{F}'\backslash \mathcal{E}$ (cf. Remark \ref{rem:20210803}).
Note that $q(\tilde{p}^{-1}(K_{s_0})\backslash \pi^{-1}(U_0))\subset \overline{T}$ is compact.
Set $W=\overline{T}\backslash q(\tilde{p}^{-1}(K_{s_0})\backslash \pi^{-1}(U_0))$.
Then $W$ is an open neighbourhood of $D$.
We have $q^{-1}(W)\cap \tilde{p}^{-1}(K_{s_0})\subset \pi^{-1}(U_0))$.
Hence for all $f\in \mathcal{F}'\backslash\mathcal{E}$ with $\tilde{p}\circ\tilde{f}(\overline{\mathbb D(s_0)})\subset K^o_{s_0}$, we have
$$
C_H^1(\mathbb D(s_0)\backslash (q\circ \tilde{f})^{-1}(W))\geq \gamma_0.
$$
Thus $( q\circ \tilde{f})_{f\in \mathcal{F}'}$ satisfies Assumption \ref{ass:20201109}.
Hence there exists an infinite subset $\mathcal{G}\subset \mathcal{F}'$ such that $(q\circ \tilde{f})_{f\in \mathcal{G}}$ converges uniformly on compact subsets of $\mathbb D$.
Indeed if $\{q\circ \tilde{f};f\in\mathcal{F}'\}\subset \mathrm{Hol}(\mathbb D,T)$ is a finite subset, we take infinite subset $\mathcal{G}\subset \mathcal{F}'$ such that $q\circ \tilde{f}$ is all the same for all $f\in\mathcal{G}$; otherwise we apply the second step for the infinite subset $\{q\circ \tilde{f};f\in\mathcal{F}'\}\subset \mathrm{Hol}(\mathbb D,T)$ to get an infinite subset $\mathcal{G}\subset \mathcal{F}'$ such that $(q\circ \tilde{f})_{f\in\mathcal{G}}$ converges uniformly on compact subsets of $\mathbb D$.
Then $\mathcal{G}$ converges uniformly on compact subsets of $\mathbb D$.
\hspace{\fill} $\square$

\medskip

{\it Proof of Proposition \ref{pro:320}.}\
We note that $Z$ is an irreducible Zariski closed set of $P_{k+1,A}$.
In the following, we fix $l\geq 1$ such that the natural map $T_{k+l}\to Z$ is surjective under the map $P_{k+l,A}\to P_{k+1,A}$.
We divide the proof into the following several steps.

{\it Step 1.}\
We are given a sequence $(f_n)_{n\in\mathbb N}$ in $\mathcal{F}$ and a smooth equivariant compactification $\overline{A}$.
For a while, we assume furthermore that  $\overline{A}$ is projective.
Replacing $(f_n)_{n\in\mathbb N}$ by its subsequence, we may assume that the sequence consists of distinct elements of $\mathcal{F}$, for otherwise the existence of a convergent subsequence is obvious. 
Replacing $(f_n)_{n\in\mathbb N}$ by its subsequence, we may assume that for every irreducible component $D\subset \partial A$, we have either (1) $(f_{n})_{n\in\mathbb N}\to D$, or (2) $(f_{n_k})_{k\in\mathbb N}\not\to D$ for every subsequence $(f_{n_k})_{k\in\mathbb N}$ of $(f_{n})_{n\in\mathbb N}$.
This is achieved as follows.
Let $D_1,\ldots ,D_k$ be all irreducible components of $\partial A$.
We define a subsequence $\mathcal{G}_1$ of $(f_n)$ as follows.
If $(f_n)$ contains a subsequence $(f_{n'})$ such that $(f_{n'})\to D_1$, then we set $\mathcal{G}_1=(f_{n'})$.
Otherwise, we set $\mathcal{G}_1=(f_n)$.
If $\mathcal{G}_{1}$ contains subsequence $(f_{n''})$ such that $(f_{n''})\to D_2$, then we set $\mathcal{G}_2=(f_{n''})$.
Otherwise we set $\mathcal{G}_2=\mathcal{G}_1$.
Continue this process to get $\mathcal{G}_k$.
Then $\mathcal{G}_k$ satisfies our requirement.
We replace $(f_n)$ by $\mathcal{G}_k$.

We denote by $\mathcal{I}$ the set of all irreducible components $D_i$ of $\partial A$ such that $(f_n)\to D_i$. 
We set $V=\cap_{D_i\in\mathcal{I}}D_i$.
When $\mathcal{I}=\emptyset$, we read $V=\overline{A}$.
By Lemma \ref{lem:20220626}, we have
\begin{equation}\label{eqn:20220608}
(f_n)\to V.
\end{equation}

{\it Step 2.}\
We may assume that $\{0\}\in\Pi((f_n)_{n\in\mathbb N})$, for otherwise the existence of a convergent subsequence is obvious.
We apply Lemma \ref{lem:202011091} to get an element $\tau_n$ of the isotropy group for $V$ such that $(\tau_n\circ f_n)_{n\in\mathbb N}$ satisfies Assumption \ref{ass:20201109}.
We take $p:W\to V$ from Lemma \ref{lem:vect}, where $W\subset \overline{A}$ is the total space of the vector bundle over $V$.
Then we have
\begin{equation}\label{eqn:202103031}
p\circ \tau_n\circ f_n=p\circ f_n
\end{equation}
for all $n\in \mathbb N$ in $\mathrm{Hol}(\mathbb D,V)$.

We set $S=P_{k,A}$.
Then we have $P_{k+1,A}\subset S_{1,A}$.
Since all elements of $\mathcal{F}$ are non-constant, $\tau_n\circ f_n:\mathbb D\to A$ is non-constant for all $n\in\mathbb N$.
Hence we get $(\tau_n\circ f_n)_{[k]}:\mathbb D\to A\times S$, which yields 
$$((\tau_n\circ f_n)_{[k]})_{S_{1,A}}:\mathbb D\to S_{1,A}.$$
Then under the inclusion $\iota:P_{k+1,A}\hookrightarrow  S_{1,A}$, we have $((\tau_n\circ f_n)_{[k]})_{S_{1,A}}=\iota\circ(\tau_n\circ f_n)_{P_{k+1,A}}$.
We note $Z\subset P_{k+1,A}\subset S_{1,A}$.

Now we claim that $(((\tau_n\circ f_n)_{[k]})_{S_{1,A}})_{n\in\mathbb N}\Rightarrow Z$.
Indeed, we have 
\begin{equation}\label{eqn:20230315}
(\tau_n\circ f_n)_{P_{k+l,A}}=(f_n)_{P_{k+l,A}}
\end{equation}
for all $n\in\mathbb N$ (cf. \eqref{eqn:202112035}).
By $\mathcal{F}_{P_{k+l,A}}\Rightarrow T_{k+l}$ and Remark \ref{rem:20220523}, we have $((f_n)_{P_{k+l,A}})_{n\in\mathbb N}\Rightarrow T_{k+l}$.
Hence $((\tau_n\circ f_n)_{P_{k+l,A}})_{n\in\mathbb N}\Rightarrow T_{k+l}$.
Thus by Lemma \ref{lem:20220219}, we get $((\tau_n\circ f_n)_{P_{k+1,A}})_{n\in\mathbb N}\Rightarrow Z$.
Hence again by Lemma \ref{lem:20220219}, we get $(((\tau_n\circ f_n)_{[k]})_{S_{1,A}})_{n\in\mathbb N}\Rightarrow Z$.

In the following, we are going to prove that there exists a subsequence $(\tau_{n_k}\circ f_{n_k})_{k\in\mathbb N}$ which converges uniformly on compact subsets on $\mathbb D$ to a holomorphic map $\mathbb D\to A$.
Set 
$$\mathcal{G}=\{ \tau_n\circ f_{n};\ n\in\mathbb N\}\subset \mathrm{Hol}(\mathbb D,A).
$$
If $\mathcal{G}$ is finite, the existence of such $(\tau_{n_k}\circ f_{n_k})_{k\in\mathbb N}$ is obvious.
Hence in the following, we assume that $\mathcal{G}$ is infinite.
We note the followings:
\begin{itemize}
\item
$\{0\}\in\Pi(\mathcal{G})$.
\item
$\mathcal{G}$ satisfies Assumption \ref{ass:20201109}.
\item
$((\varphi_{[k]})_{S_{1,A}})_{\varphi\in\mathcal{G}}
\Rightarrow Z$, where $\varphi_{[k]}:\mathbb D\to A\times S$.
\end{itemize}
These properties follow from the discussion of this step (step 2).

{\it Step 3.}
Let $\omega_{\overline{A}}$ and $\omega_S$ be smooth positive $(1,1)$-forms on $\overline{A}$ and $S$, respectively.
We apply Lemma \ref{lem:tcor} to get $\sigma_1\in (0,1)$ with the following property:
Let $s\in (\sigma_1,1)$, $\varepsilon>0$, $\delta>0$.
Then there exists $\mu>0$ such that, for all $\varphi\in\mathcal{G}$, we have
\begin{equation}\label{eqn:202108053}
T_{s}(r,\varphi_S,\omega_{S})
\leq \varepsilon 
T_{s}(r,\varphi,\omega_{\overline{A}})
+\mu
\end{equation}
for all $r\in (s,1 )$ outside some exceptional set with the linear measure less than $\delta$.

Next by \eqref{eqn:20230315} and Remark \ref{rem:20210805} (3), we have 
$$\mathrm{LIM}(\mathcal{G}_{P_{k+l,A}},\{ E_{k+l,A,A/B}\}_{B\in \Pi(\mathcal{F})})=T_{k+l}.$$
Since $Z\subset S_{1,A}$ is horizontally integrable, we may take a Zariski closed subset $W\subsetneqq t(Z)\subset S$ which appears in Proposition \ref{pro:mpro}, where $t:S_{1,A}\to S$ is the natural projection.
Let $\pi:P_{k+l,A}\to S$ be the natural projection.
By $\pi(T_{k+l})=t(Z)$, we have $T_{k+l}\not\subset \pi^{-1}(W)$.
Let $\lambda_{W}\geq 0$ be a Weil function for $W$.
Then $\lambda_W\circ \pi$ is a Weil function for $\pi^{*}W\subset P_{k+l,A}$.
By $\Pi(\mathcal{G})\supset \Pi(\mathcal{F})$, we may apply Lemma \ref{lem:20210302} to get an infinite subset $\mathcal{G}'\subset \mathcal{G}$ and $\sigma_2\in (0,1)$ with the following property:
Let $s\in (\sigma_2,1)$, $\varepsilon>0$, $\delta>0$.
Then there exists $\beta>0$ such that, for all $\varphi\in \mathcal{G}'$, we have
\begin{equation}\label{eqn:202108051}
m(r,\varphi_{P_{k+l,A}},\lambda_{W}\circ \pi)
\leq \varepsilon T_{s}(r,\varphi,\omega_{\overline{A}})+
\beta
\end{equation}
for all $r\in (s,1 )$ outside some exceptional set with the linear measure less than $\delta$.

We set $\sigma=\{ 1/2,\sigma_1,\sigma_2\}$.
Note that we have $\varphi_S(\mathbb D)\not\subset W$ for all $\varphi\in \mathcal{G}'$, which follows from \eqref{eqn:202108051}.

{\it Step 4.}\
Now we fix $s\in (\sigma,1)$ arbitrary.
We set 
$$\delta=(1-s)/5.$$
By Proposition \ref{pro:mpro} applied to $\{\varphi_{[k]};\varphi\in\mathcal{G}\}\subset\mathrm{Hol}(\mathbb D,A\times S)$, there exist $c_1>0$, $c_2>0$, $c_3>0$ such that for all $\varphi\in \mathcal{G}$ with $\varphi_S(\mathbb D)\not\subset W$, we have
\begin{equation*}
T_s(r,\varphi,\omega_{\overline{A}})
 \leq 
c_1
T_s(r,\varphi_{S},\omega_{S})
+
c_2
m((s+r)/2,\varphi_{S},\lambda_{W})
+c_3
\end{equation*}
for all $r\in (s,1 )$ outside some exceptional set with the linear measure less than $\delta$.
Letting $\varepsilon=1/3c_2$ in \eqref{eqn:202108051}, we obtain that, for all $\varphi\in \mathcal{G}'$, we have
$$
m((s+r)/2,\varphi_{S},\lambda_{W})
\leq \frac{1}{3c_2} T_{s}(r,\varphi,\omega_{\overline{A}})+
\beta
$$
for all $r\in (s,1)$ outside some exceptional set of the linear measure less than $2\delta$.
Hence for all $\varphi\in \mathcal{G}'$, we get 
$$
T_s(r,\varphi,\omega_{\overline{A}})
 \leq 
 \frac{1}{3}
 T_{s}(r,\varphi,\omega_{\overline{A}})
 +
c_1
T_s(r,\varphi_{S},\omega_{S})
+c_2\beta_2
+c_3
$$
for all $r\in (s,1)$ outside an exceptional set of the linear measure less than $3\delta$.
By \eqref{eqn:202108053}, letting $\varepsilon=1/3c_1$, we obtain that for all $\varphi\in \mathcal{G}'$, we have
$$
T_s(r,\varphi,\omega_{\overline{A}})
 \leq
 \frac{2}{3}
 T_{s}(r,\varphi,\omega_{\overline{A}})
 +
c_1\mu
+c_2\beta_2
+c_3
$$
for all $r\in (s,1)$ outside an exceptional set of the linear measure less than $4\delta$.
Thus for all $\varphi\in \mathcal{G}'$, we  get
$$
T_s(r,\varphi,\omega_{\overline{A}})
 \leq c
$$
for all $r\in (s,1 )$ outside some exceptional set of linear measure less than $4\delta$, where $c=3( c_1\mu
+c_2\beta_2
+c_3)$.
We may apply this estimate for some $r\in (s+\delta,1)$.
Hence we get
$$
T_s(s+\delta,\varphi,\omega_{\overline{A}})
 \leq c,
$$
thus
$$
\int_{\mathbb D(s)}\varphi^*\omega_{\overline{A}}
 \leq \frac{(s+\delta)c}{\delta}
$$
for all $\varphi\in \mathcal{G}'$.
Hence by Lemma \ref{lem:compact}, there exists a subsequence $(\tau_{n_k}\circ f_{n_k})_{k\in\mathbb N}$ which converges to a holomorphic map $g:\mathbb D\to A$.

{\it Step 5.}\
Now we show that $(f_{n_k})_{k\in\mathbb N}$ converges uniformly on compact subsets of $\mathbb D$ to $p\circ g$, where $p:W\to V$.
Note that $(p\circ \tau_{n_k}\circ f_{n_k})_{k\in\mathbb N}$ converges uniformly on compact subsets of $\mathbb D$ to $p\circ g$.
By \eqref{eqn:202103031}, $(p\circ f_{n_k})_{k\in\mathbb N}$ also converges uniformly on compact subsets of $\mathbb D$ to $p\circ g$.
We apply Lemma \ref{lem:20210303}.
Then by \eqref{eqn:20220608}, $(f_{n_k})_{k\in\mathbb N}$ converges uniformly on compact subsets of $\mathbb D$ to $p\circ g$.

{\it Step 6.}\
So far, we have assumed that $\overline{A}$ is projective.
Now let $\overline{A}$ be arbitrary smooth equivariant compactification.
We may take an equivariant modification $q:\hat{A}\to \overline{A}$ such that $\hat{A}$ is smooth and projective (cf. Lemma \ref{lem:202206031}).
Then we may take a subsequence $(f_{n_k})_{k\in\mathbb N}$ which converges uniformly on compact subsets of $\mathbb D$ to $g:\mathbb D\to \hat{A}$.
Then $(f_{n_k})_{k\in\mathbb N}$ converges uniformly on compact subsets of $\mathbb D$ to $q\circ g:\mathbb D\to \overline{A}$.
This conclude the proof.
\hspace{\fill} $\square$

\begin{cor}\label{cor:20220424}
Let $\mathcal{F}=(f_i)_{i\in I}$ be an infinite indexed family of non-constant holomorphic maps in $\mathrm{Hol}(\mathbb D,A)$ which satisfies Assumption \ref{ass:202103061}.
Suppose that there exists $k\geq 0$ such that $Z\subset P_{k+1,A}\subset (P_{k,A})_{1,A}$ is horizontally integrable, where $Z$ is defined by \eqref{eqn:202202191}.
Let $\bar{A}$ be a smooth equivariant compactification. 
Then there exists an infinite subfamily $\mathcal{G}$ of $\mathcal{F}$ such that $\mathcal{G}$ converges uniformly on compact subsets of $\mathbb D$ to $g:\mathbb D\to \bar{A}$.
\end{cor}

{\it Proof.}\
Set $\mathcal{F}_o=\{f_i;\ i\in I\}\subset \mathrm{Hol}(\mathbb D,A)$.
We may assume that the map $I\to \mathcal{F}_o$ is finite-to-one mapping.
Indeed otherwise, we may take an infinite subset $J\subset I$ such that $f_i$ are all the same for all $i\in J$, hence the assertion is obvious.
Thus we assume the map $I\to \mathcal{F}_o$ is finite-to-one mapping.
In particular, $\mathcal{F}_o$ is an infinite set.
We note that $\Pi(\mathcal{F}_o)=\Pi(\mathcal{F})$.
By taking a section of $I\to \mathcal{F}_o$, we take an infinite subset $I_o\subset I$ such that $I_o\to\mathcal{F}_o$ is bijective.
By Remark \ref{rem:20210805} (3), we have
$$
\mathrm{LIM}((\mathcal{F}_o)_{P_{k,A}};\{ E_{k,A,A/B}\}_{B\in\Pi(\mathcal{F}_o)})
=
\mathrm{LIM}(\mathcal{F}_{P_{k,A}};\{ E_{k,A,A/B}\}_{B\in\Pi(\mathcal{F})}).
$$  
Hence the assumption of Proposition \ref{pro:320} is satisfied for $\mathcal{F}_o$.
Hence by Proposition \ref{pro:320} for $\mathcal{F}_o$, we may take an infinite subset $\mathcal{G}_o\subset \mathcal{F}_o$ such that $\mathcal{G}_o$ converges uniformly on compact subsets of $\mathbb D$.
Then we may take an infinite subset $J\subset I_o$ such that $J\to\mathcal{G}_o$ is bijective.
We set $\mathcal{G}=(f_i)_{i\in J}$ to conclude the proof.
\hspace{\fill} $\square$

\section{Families of closed subschemes and regular jets}\label{sec:8.5}

We set $\Lambda_k =\mathrm{Spec}\ \mathbb C[\varepsilon ]/(\varepsilon^{k+1})$.
Then we have a sequence of closed immersions
\begin{equation}\label{eqn:202303151}
\Lambda_0\hookrightarrow \Lambda_1\hookrightarrow \Lambda_2\hookrightarrow \cdots.
\end{equation}
Let $S$ be a smooth variety.
Given a morphism $\eta:\Lambda_k\to S$ of schemes, we obtain the derivative $\eta':\Lambda_{k-1}\to TS$ of $\eta$.
This map satisfies the following:
Let $\varphi$ be a local holomorphic function on $S$.
Let $d:\mathcal{O}_S\to \Omega^1_S$ be the derivation.
Then we have 
\begin{equation}\label{eqn:202203221}
(\eta')^*d\varphi=\frac{d}{d\varepsilon}\eta^*\varphi,
\end{equation}
where $\frac{d}{d\varepsilon}:\mathbb C[\varepsilon]/(\varepsilon^{k+1})\to \mathbb C[\varepsilon]/(\varepsilon^{k})$ is the derivation.
A regular $k$-jet is a morphism $\eta:\Lambda_{k}\to S$ such that $\eta'(0)\not\in0_{TS}$, where $0_{TS}\subset TS$ is the zero section.
Hence by the composite of $\eta'$ and $TS-0_{TS}\to S_1$, we get $\eta_{[1]}:\Lambda_{k-1}\to S_1$.
Inductively, we obtain $\eta_{[l]}:\Lambda_{k-l}\to S_l$.
In particular, we get a point $\eta_{[k]}(0)\in S_k$.

Let $Z\subset S$ be a closed subscheme.
We define a closed subscheme $\mathcal{D}Z\subset S_1$ as follows.
Let $W\subset S$ be an affine open set where $Z\cap W$ is defined by $\varphi_1,\ldots ,\varphi_n \subset \Gamma (W,\mathcal{O}_W)$.
Then we define the closed subscheme $\widetilde{Z\cap W}\subset TW$ by $\varphi_1,\ldots ,\varphi_n,d\varphi_1,\ldots ,d\varphi_n$.
Then this definition of $\widetilde{Z\cap W}$ does not depend on the choice of generators $\varphi_1,\ldots ,\varphi_n$, so well defined over $W$.
In general, we cover $S$ by open affines $\{ W_i\}$
and make closed subschemes $\widetilde{Z\cap W_i}\subset TW_i$. 
Then we glue these subschemes and define the subscheme $\tilde{Z}\subset TS$.
By the construction, $\tilde{Z}$ is invariant under the $\mathbb C^*$-action on the fibers of $TS\to S$.
Thus we get a closed subscheme $\mathcal{D}Z=\tilde{Z}/\mathbb C^* \subset S_1$.

\par

\subsection{Family of closed subschemes}
Let $V$ and $S$ be smooth varieties and let $X\subset V\times S$ be a closed subscheme.
For $k\geq 0$, we define a closed subscheme $\mathcal{P}_kX\subset V\times S_k$ inductively as follows.
We have $V\times TS=T_{(V\times S)/V}\subset T(V\times S)$.
Hence we get $V\times S_1\subset (V\times S)_1$.
We define $\mathcal{P}_1X\subset V\times S_1$ by the restriction of $\mathcal{D}X\subset (V\times S)_1$ onto $V\times S_1\subset (V\times S)_1$.
Now suppose we get a closed subscheme $\mathcal{P}_kX\subset V\times S_k$.
We have $S_{k+1}\subset (S_k)_1$.
Thus we define $\mathcal{P}_{k+1}X\subset V\times S_{k+1}$ by the restriction of 
\begin{equation}\label{eqn:202206236}
\mathcal{P}_1(\mathcal{P}_kX)\subset V\times (S_k)_1
\end{equation}
onto $ V\times S_{k+1}\subset V\times (S_k)_1$.

\begin{lem}\label{lem:20210324}
Let $V$ and $S$ be smooth varieties and let $X\subset V\times S$ be a closed subscheme.
Let $\xi:\Lambda_k\to S$ be a regular $k$-jet.
Let $s_0=\xi(0)\in S$ and $s_k=\xi_{[k]}(0)\in S_k$.
Assume that the natural map $(\mathcal{P}_kX)_{s_k}\to X_{s_0}$ is an isomorphism as schemes.	
Then $X_{\Lambda_k}=X_{s_0}\times \Lambda_k$ as closed subschemes of $V\times \Lambda_k$.
Here $X_{\Lambda_k}$ is the pull-back of $X\to S$ by $\xi$.
\end{lem}

Before going to prove this lemma, we start from preliminary observation.
Let $W\subset V$ be an affine open subset and set $Z=X\cap (W\times S)$.
For $\nu=0,\ldots,k$, let
$I_{\nu}\subset \mathbb C[W]\otimes_{\mathbb C}\mathbb C[\Lambda_{k-\nu}]$
be the defining ideal of $(\mathcal{P}_{\nu}Z)_{\Lambda_{k-\nu}}$, where $(\mathcal{P}_{\nu}Z)_{\Lambda_{k-\nu}}$ is the pull-back of $\mathcal{P}_{\nu}Z\subset W\times S_{\nu}$ by $\xi_{[\nu]}:\Lambda_{k-\nu}\to S_{\nu}$.
We denote $\frac{\partial}{\partial\varepsilon}:\mathbb C[W]\otimes_{\mathbb C}\mathbb C[\Lambda_{k-\nu}]\to \mathbb C[W]\otimes_{\mathbb C}\mathbb C[\Lambda_{k-\nu-1}]$ by $\mathrm{id}_{\mathbb C[W]}\otimes \frac{d}{d\varepsilon}$.

\medskip

{\it Claim.}\
If $h\in I_{\nu}$, then $\frac{\partial}{\partial \varepsilon}h\in I_{\nu+1}$.

\medskip

We prove this claim.
Let $U\subset S_{\nu}$ be an affine open such that $\xi_{[\nu]}(0)\in U$.
Suppose that $\mathcal{P}_{\nu}Z$ is defined by $f_1,\ldots,f_n\in\Gamma(W\times U,\mathcal{O}_{W\times U})$ on $W\times U$.
Then $(\mathcal{P}_{\nu}Z)_{\Lambda_{k-\nu}}\subset W\times \Lambda_{k-\nu}$ is defined by $g_1,\ldots,g_n\in \mathbb C[W]\otimes_{\mathbb C}\mathbb C[\Lambda_{k-\nu}]$, where $g_1,\ldots,g_n$ are the images of $f_1,\ldots,f_n$ under the natural map $\mathbb C[W]\otimes_{\mathbb C}\mathbb C[U]\to \mathbb C[W]\otimes_{\mathbb C}\mathbb C[\Lambda_{k-\nu}]$.
We note that \eqref{eqn:202203221} yields the following commutative diagram of $\mathbb C[W]$-modules:
\begin{equation*}
	\begin{CD}
		\mathbb C[W]\otimes_{\mathbb C}\mathbb C[U] @>\mathrm{id}_{\mathbb C[W]}\otimes \xi_{[\nu]}^*>>   \mathbb C[W]\otimes_{\mathbb C}\mathbb C[\Lambda_{k-\nu}]\\
		@V\mathrm{id}_{\mathbb C[W]}\otimes d_UVV  @VV\frac{\partial}{\partial\varepsilon}V \\
		\mathbb C[W]\otimes_{\mathbb C}\mathbb C[TU] @>>\mathrm{id}_{\mathbb C[W]}\otimes (\xi_{[\nu]}')^*> \mathbb C[W]\otimes_{\mathbb C}\mathbb C[\Lambda_{k-\nu-1}]
	\end{CD}
\end{equation*}
Hence $(\mathcal{P}_{\nu+1}Z)_{\Lambda_{k-\nu-1}}$ is defined by 
\begin{equation}\label{eqn:20220320}
	\bar{g}_1,\ldots,\bar{g}_n,\frac{\partial}{\partial\varepsilon}g_1,\ldots,\frac{\partial}{\partial\varepsilon}g_n.
\end{equation}
Here $\bar{g}_1,\ldots,\bar{g}_k$ are the images of $g_1,\ldots,g_k$ under $\mathbb C[W]\otimes_{\mathbb C}\mathbb C[\Lambda_{k-\nu}]\to \mathbb C[W]\otimes_{\mathbb C}\mathbb C[\Lambda_{k-\nu-1}]$, where $\mathbb C[\Lambda_{k-\nu-1}]=\mathbb C[\Lambda_{k-\nu}]/(\varepsilon^{k-\nu})$ (cf. \eqref{eqn:202303151}).

Now let $h\in I_{\nu}$.
Then there exists $b_1,\ldots,b_n\in  \mathbb C[W]\otimes_{\mathbb C}\mathbb C[\Lambda_{k-\nu}]$ such that $h=b_1g_1+\cdots+b_ng_n$.
Then we have $\frac{\partial}{\partial \varepsilon}h=\bar{g}_1\frac{\partial}{\partial \varepsilon}b_1+\bar{b}_1\frac{\partial}{\partial \varepsilon}g_1+\cdots+\bar{g}_n\frac{\partial}{\partial \varepsilon}b_n+\bar{b}_n\frac{\partial}{\partial \varepsilon}g_n$.
Hence $\frac{\partial}{\partial \varepsilon}h\in I_{\nu+1}$.
This proves our claim.

\medskip

{\it Proof of Lemma \ref{lem:20210324}.}\
We may assume that $S$ is affine.
Let $X\subset V\times S$ be defined locally on an affine open $W\subset V$ by $f_1,\ldots ,f_n$.
Then $X_{\Lambda_k}\subset V\times \Lambda_k$ is defined locally by $g_1,\ldots,g_n\in \mathbb C[W]\otimes_{\mathbb C}\mathbb C[\Lambda_k]$, where $g_1,\ldots,g_n$ are the images of $f_1,\ldots,f_n$ under the natural map $\mathbb C[W]\otimes_{\mathbb C}\mathbb C[S]\to \mathbb C[W]\otimes_{\mathbb C}\mathbb C[\Lambda_k]$.
Let $\mathfrak{X}\subset V\times \Lambda_k$ be the constant family $\mathfrak{X}=(X_{s_0})\times \Lambda_k$ over $\Lambda_k$.
Then $\mathfrak{X}$ is defined locally by $g_1|_{\varepsilon=0},\ldots,g_n|_{\varepsilon=0}$.
Hence it is enough to show 
$$
(g_1,\ldots,g_n)=(g_1|_{\varepsilon=0},\ldots,g_n|_{\varepsilon=0})
$$
as ideals of $\mathbb C[W]\otimes_{\mathbb C}\mathbb C[\Lambda_k]$.

We first claim
\begin{equation}\label{eqn:20210405}
	\frac{\partial^{\nu} g_i}{\partial\varepsilon^{\nu}}\vert_{\varepsilon =0} \in (g_1|_{\varepsilon=0},\ldots,g_n|_{\varepsilon=0})
\end{equation}
for $\nu =0,\ldots ,k$, where $(g_1|_{\varepsilon=0},\ldots,g_n|_{\varepsilon=0})\subset \mathbb C[W]$.
To prove this, we set $s_{\nu}=\xi_{[\nu]}(0)\in S_{\nu}$.
Then we have $(\mathcal{P}_{k}X)_{s_k}\subset (\mathcal{P}_{k-1}X)_{s_{k-1}}\subset \cdots\subset (\mathcal{P}_{1}X)_{s_1}\subset X_{s_0}$.
Thus by the assumption $(\mathcal{P}_kX)_{s_k}= X_{s_0}$, we have $(\mathcal{P}_{\nu}X)_{s_{\nu}}= X_{s_0}$.
By the claim above, we have $\frac{\partial^{\nu} g_i}{\partial\varepsilon^{\nu}}\in I_{\nu}$.
Hence $\frac{\partial^{\nu} g_i}{\partial\varepsilon^{\nu}}|_{\varepsilon=0}\in\mathbb C[W]$ is contained in the defining ideal of $(\mathcal{P}_kX)_{s_k}\cap W$, hence that of $X_{s_0}\cap W$.
Note that $X_{s_0}\cap W$ is defined on $W$ by $(g_i\vert_{\varepsilon =0})_{i=1,\ldots ,n}$.
This proves \eqref{eqn:20210405}.

Now we have
\begin{equation}\label{eqn:202206081}
g_i=g_i|_{\varepsilon=0}+\varepsilon\frac{\partial g_i}{\partial\varepsilon}\vert_{\varepsilon =0}+\cdots+\frac{1}{k!}\varepsilon^k\frac{\partial^{k} g_i}{\partial\varepsilon^{k}}\vert_{\varepsilon =0}.
\end{equation}
Hence by \eqref{eqn:20210405}, we get $(g_1,\ldots,g_n)\subset (g_1|_{\varepsilon=0},\ldots,g_n|_{\varepsilon=0})$.

Next we show $(g_1|_{\varepsilon=0},\ldots,g_n|_{\varepsilon=0})\subset (g_1,\ldots,g_n)$.
For this we only show $g_1|_{\varepsilon=0}\in (g_1,\ldots,g_n)$, for the other  indices are treated in the same manner.
Let $J\subset (g_1,\ldots,g_n)$ be the subset of the elements of the form
$$
A=g_1|_{\varepsilon=0}+\varepsilon A_1+\cdots+\varepsilon^kA_k,
$$
where $A_1,\ldots,A_k\in (g_1|_{\varepsilon=0},\ldots,g_n|_{\varepsilon=0})\subset \mathbb C[W]$.
By \eqref{eqn:20210405} and \eqref{eqn:202206081}, we have $g_1\in J$.
Hence $J\not=\emptyset$.
For each element $A\in J$, we set 
$$\mu_A=\min\{ j; A_j\not=0\},$$
where we set $\mu_A=k+1$ if $A_j=0$ for all $j=1,\ldots,k$.
We put $\mu=\max_{A\in J} \mu_A$.
What we want to show is that $\mu=k+1$. 
So suppose $\mu\leq k$.
Take $A\in J$ such that $\mu_A=\mu$.
Then
$$
A=g_1|_{\varepsilon=0}+\varepsilon^{\mu}A_{\mu}+\cdots +\varepsilon^kA_k.
$$
By $A_{\mu}\in (g_1|_{\varepsilon=0},\ldots,g_n|_{\varepsilon=0})$, there exist $b_1,\ldots,b_n\in \mathbb C[W]$ such that 
$$
A_{\mu}=b_1g_1|_{\varepsilon=0}+\cdots +b_ng_n|_{\varepsilon=0}.
$$
Then by \eqref{eqn:20210405} and \eqref{eqn:202206081}, we have
$$
A'=A-\varepsilon^{\mu}b_1g_1-\cdots-\varepsilon^{\mu}b_ng_n\in J.
$$
Moreover we have $\mu_{A'}>\mu$.
This is a contradiction.
Hence we have $\mu=k+1$, thus $g_1|_{\varepsilon=0}\in (g_1,\ldots,g_n)$.
This completes the proof of our lemma.
\hspace{\fill} $\square$

\subsection{Regular jets and Demailly jet spaces}

Let $S$ be a smooth variety.
We are going to introduce the jet space $J_kS$ as in \cite[Sec. 4.2]{Ykob}.
This definition is equivalent to the usual definition described in \cite[Sec. 4.6.1]{NW}.
See Remark \ref{rem:20220329} below.

Set $J_0S=S$ and $J_1S=TS$.
For $k\geq 1$, the space $J_kS$ is a smooth variety with an embedding
\begin{equation*}\label{eqn:bigin}
J_{k}S\overset{\iota_k}{\hookrightarrow} TJ_{k-1}S,
\end{equation*}
where we set $\iota_1:J_1S\to TS$ to be the identity map.
We define $J_kS$ inductively as follows.
So we suppose that the smooth variety $J_kS$ and the embedding $\iota_k:J_{k}S\hookrightarrow TJ_{k-1}S$ are given.
Let $\varpi_k:J_kS\to J_{k-1}S$ be the composite of $\iota_k:J_{k}S\hookrightarrow TJ_{k-1}S$ and the projection $TJ_{k-1}S\to J_{k-1}S$.
Then $\varpi_{k}:J_kS\to J_{k-1}S$ induces the following commutative diagram
\begin{equation*}
\begin{CD}
 TJ_kS@>p_k>>   J_kS\\
@V(\varpi_{k})_*VV  @VV\varpi_{k}V \\
TJ_{k-1}S@>>p_{k-1}>    J_{k-1}S \end{CD}
\end{equation*}
This induces the morphism
\begin{equation*}
\mu_k:TJ_kS\to J_kS\times_{J_{k-1}S}TJ_{k-1}S.
\end{equation*}
The graph of $\iota_{k}:J_kS\to TJ_{k-1}S$ defines the closed immersion
\begin{equation*}
\hat{\iota}_k:J_kS\hookrightarrow J_kS\times_{J_{k-1}S}TJ_{k-1}S.
\end{equation*}
Then $\iota_{k+1}:J_{k+1}S\hookrightarrow TJ_{k}S$ is defined by the base change of $\hat{\iota}_k$ by $\mu_k$.
\begin{equation*}
\begin{CD}
J_{k+1}S @>\iota_{k+1}>> TJ_kS\\
 @V(\mu_k)'VV  @VV\mu_kV \\
 J_{k}S @>>\hat{\iota}_k>J_kS\times_{J_{k-1}S}TJ_{k-1}S
\end{CD}
\end{equation*}
We get a map $(\mu_k)':J_{k+1}S\to J_kS$ from this base change.
Let $\varpi_{k+1}:J_{k+1}S\to J_{k}S$ be the composite of $\iota_{k+1}:J_{k+1}S\hookrightarrow TJ_{k}S$ and the projection $p_k:TJ_{k}S\to J_{k}S$.
Then we note
\begin{equation}\label{eqn:202206082}
(\mu_k)'=\varpi_{k+1}.
\end{equation}
Indeed we have $(\mu_k)'=\mathrm{id}_{J_kS}\circ (\mu_k)'=r_1\circ\hat{\iota}_k\circ (\mu_k)'$, where $r_1:J_kS \times_{J_{k-1}S}TJ_{k-1}S\to J_kS$ is the first projection.
On the other hand, we have $r_1\circ\hat{\iota}_k\circ (\mu_k)'=r_1\circ \mu_k\circ\iota_{k+1}=p_k\circ \iota_{k+1}$.
Hence we get \eqref{eqn:202206082}.
By \cite[Cor. 4.4]{Ykob}, the map $\varpi_{k+1}:J_{k+1}S\to J_kS$ is smooth.
Hence $J_{k+1}S$ is smooth.

\begin{rem}\label{rem:20220329}
Assume that an affine open $W\subset S$ admits coordinate functions $\chi_1,\ldots,\chi_n$ such that $TW$ splits as $TW=W\times \mathbb C^n$ and $(d\chi_1,\ldots,d\chi_n)$ defines the second projection, where $n=\dim S$.
Then by \cite[Cor. 4.6]{Ykob}, we have $J_kW=W\times \mathbb C^{kn}$ and $(d\chi_1,\ldots,d\chi_n,\ldots,d^k\chi_1,\ldots,d^k\chi_n)$ defines the second projection.
This shows that our definition of $J_kS$ coincides with the usual definition described in \cite[Sec. 4.6.1]{NW}.
\end{rem}

The space $\iota_{k+1}:J_{k+1}S\hookrightarrow TJ_{k}S$ is characterized as follows:
\begin{lem}\label{lem:20220330}
A map $\eta:W\to TJ_{k}S$ from a scheme $W$ factors $\iota_{k+1}:J_{k+1}S\hookrightarrow TJ_{k}S$ if and only if $\iota_k\circ p_k\circ \eta=(\varpi_k)_*\circ \eta$.
In particular, if we consider the case $W=\mathrm{Spec}\ \mathbb C$, we have 
\begin{equation}\label{eqn:20211103}
J_{k+1}S=\{ x\in TJ_kS\ ;\ \iota_k\circ p_k(x)=(\varpi_{k})_*(x)\} .
\end{equation}
\end{lem}

{\it Proof.}\
To prove this, we suppose that $\eta:W\to TJ_{k}S$ satisfies $\iota_k\circ p_k\circ \eta=(\varpi_k)_*\circ \eta$.
We claim that 
\begin{equation}\label{eqn:202206091}
\mu_k\circ\eta=\hat{\iota}_k\circ p_k\circ\eta.
\end{equation}
To show this we note that 
\begin{equation}\label{eqn:20220609}
r_1\circ\mu_k\circ\eta=p_k\circ \eta=r_1\circ \hat{\iota}_k\circ p_k\circ\eta.
\end{equation}
We also have 
\begin{equation}\label{eq:202303152}
r_2\circ\mu_k\circ\eta=(\varpi_k)_*\circ \eta,
\end{equation}
where $r_2:J_kS \times_{J_{k-1}S}TJ_{k-1}S\to TJ_{k-1}S$ is the second projection of the fiber product:
\begin{equation*}
\begin{CD}
J_kS \times_{J_{k-1}S}TJ_{k-1}S@>r_2>>       TJ_{k-1}S \\
 @Vr_1VV  @VVp_{k-1}V \\
 J_{k}S  @>>\varpi_{k}>    J_{k-1}S\end{CD}
\end{equation*}
By our assumption $\iota_k\circ p_k\circ \eta=(\varpi_k)_*\circ \eta$, we have 
\begin{equation}\label{eq:202303153}
(\varpi_k)_*\circ \eta=\iota_k\circ p_k\circ\eta=r_2\circ \hat{\iota}_k\circ p_k\circ\eta.
\end{equation}
Hence by \eqref{eq:202303152} and \eqref{eq:202303153}, we have $r_2\circ\mu_k\circ\eta=r_2\circ \hat{\iota}_k\circ p_k\circ\eta$.
Combining this with \eqref{eqn:20220609}, we get \eqref{eqn:202206091} by the universal property of the fiber product.
Hence by \eqref{eqn:202206091}, the universal property of the fiber product yields the map $(\eta,p_k\circ\eta):W\to J_{k+1}S$ such that the composition with $\iota_{k+1}:J_{k+1}S\to TJ_kS$ is equal to $\eta$.
Hence $\eta$ factors $\iota_{k+1}:J_{k+1}S\hookrightarrow TJ_{k}S$.

To prove the converse, we claim
\begin{equation}\label{eqn:20220404}
\iota_k\circ p_k\circ \iota_{k+1}=(\varpi_k)_*\circ \iota_{k+1}.
\end{equation}
Indeed by \eqref{eqn:202206082}, we have $(\mu_k)'=p_k\circ\iota_{k+1}$.
Hence we have the following equalities
$$
\iota_k\circ p_k\circ \iota_{k+1}=r_2\circ\hat{\iota}_k\circ p_k\circ \iota_{k+1}
=r_2\circ\hat{\iota}_k\circ(\mu_k)'
=r_2\circ \mu_k\circ \iota_{k+1}=(\varpi_k)_*\circ \iota_{k+1}.
$$
Hence if $\eta:W\to TJ_{k}S$ factors $\iota_{k+1}:J_{k+1}S\hookrightarrow TJ_{k}S$, we have $\iota_k\circ p_k\circ \eta=(\varpi_k)_*\circ \eta$.
\hspace{\fill} $\square$

\medskip

Let $l\geq k$.
Given a morphism $\eta:\Lambda_l\to S$, we obtain $j_k\eta:\Lambda_{l-k}\to J_kS$ such that 
\begin{equation}\label{eqn:202203291}
\iota_k\circ j_{k}\eta =(j_{k-1}\eta)'.
\end{equation}
Indeed, to show the existence of  $j_k\eta$ by induction on $k$, we assume the existence for $k$.
Then we have
$$
\iota_k\circ p_k\circ (j_k\eta)'=\iota_k\circ j_k\eta=(j_{k-1}\eta)'=(\varpi_k)_*\circ (j_k\eta)'.
$$
Hence, by Lemma \ref{lem:20220330}, the map $(j_k\eta)':\Lambda_{l-k-1}\to TJ_kS$ factors $\iota_{k+1}:J_{k+1}S\hookrightarrow TJ_kS$.
Thus we get $j_{k+1}\eta:\Lambda_{l-k-1}\to J_{k+1}S$ such that $\iota_{k+1}\circ j_{k+1}\eta =(j_{k}\eta)'$.

Let $J_k^{\mathrm{reg}}S\subset J_kS$ be the Zariski open which is the inverse image of $TS-\{ 0_{TS}\}$ under the map $J_kS\to TS$, provided $k\geq 1$.
We set $J_0^{\mathrm{reg}}S=S$.
Note that the map $J_{k+1}S\to J_kS$ induces $J_{k+1}^{\mathrm{reg}}S\to J_k^{\mathrm{reg}}S$.

\begin{lem}\label{lem:202111041}
Let $S$ and $M$ be smooth algebraic varieties.
Let $V\subset TM$ be an algebraic vector subbundle and let $\tilde{M}=P(V)$.
Let $k\geq 0$.
Let $\varphi:J_k^{\mathrm{reg}}S\to M$ be a morphism such that the induced map $\varphi_*:TJ_{k}^{\mathrm{reg}}S\to TM$ satisfies $\varphi_*(\iota_{k+1}(J_{k+1}^{\mathrm{reg}}S))\subset V\backslash \{0_V\}$.
Let $\Phi:J_{k+1}^{\mathrm{reg}}S\to \tilde{M}$ be the composite of the following morphisms:
$$
J_{k+1}^{\mathrm{reg}}S\overset{\varphi_*\circ\iota_{k+1}}{\longrightarrow} V\backslash\{ 0_V\}\overset{\tau}{\longrightarrow} \tilde{M}.
$$
Then the induced map $\Phi_*:TJ_{k+1}^{\mathrm{reg}}S\to T\tilde{M}$ satisfies $\Phi_*(\iota_{k+2}(J_{k+2}^{\mathrm{reg}}S))\subset \tilde{V}\backslash \{0_V\}$, where $\tilde{V}\subset T\tilde{M}$ is defined by \eqref{eqn:20211104}.
\end{lem}

{\it Proof.}\
We first show that the following diagram is commutative:
\begin{equation}\label{eqn:202111041}
\begin{CD}
J_{k+1}^{\mathrm{reg}}S @>\Phi>>   \tilde{M}\\
@V\varpi_{k+1}VV  @VV\pi V \\
J_{k}^{\mathrm{reg}}S @>>\varphi> M
\end{CD}
\end{equation}
Indeed, $\pi\circ\Phi$ is the composite of 
$$
J_{k+1}^{\mathrm{reg}}S\overset{\iota_{k+1}}{\longrightarrow}TJ_{k}^{\mathrm{reg}}S\overset{\varphi_*}{\longrightarrow} TM\overset{q}{\longrightarrow} M.
$$
By the definition of the derivation, we have $q\circ\varphi_*= \varphi\circ p_k$, i.e., the following diagram is commutative.
\begin{equation*}
\begin{CD}
TJ_{k}^{\mathrm{reg}}S @>\varphi_*>>   TM\\
@Vp_kVV  @VVqV \\
 J_{k}^{\mathrm{reg}}S @>>\varphi> M
\end{CD}
\end{equation*}
Hence we get 
$$
\pi\circ\Phi=q\circ \varphi_*\circ\iota_{k+1}=\varphi\circ p_k\circ\iota_{k+1}=\varphi\circ\varpi_{k+1}.
$$
This shows that \eqref{eqn:202111041} is commutative.

Now by \eqref{eqn:202111041}, we get the following commutative diagram:
\begin{equation}\label{eqn:20220401}
\begin{CD}
J_{k+2}^{\mathrm{reg}}S @>\iota_{k+2}>>TJ_{k+1}^{\mathrm{reg}}S @>\Phi_*>>   T\tilde{M}@>\tilde{q}>> \tilde{M}\\
@V\varpi_{k+2}VV @V(\varpi_{k+1})_*VV  @VV\pi_*V @VV\pi V\\
J_{k+1}^{\mathrm{reg}}S@>>\iota_{k+1}> TJ_{k}^{\mathrm{reg}}S @>>\varphi_*> TM @>>q> M
\end{CD}
\end{equation}
Indeed the only non-trivial part is the relation $(\varpi_{k+1})_*\circ \iota_{k+2}=\iota_{k+1}\circ \varpi_{k+2}$. 
To show this, we note $(\varpi_{k+1})_*\circ \iota_{k+2}=\iota_{k+1}\circ p_{k+1}\circ \iota_{k+2}$, which follows from Lemma \ref{lem:20220330} (cf. \eqref{eqn:20220404}).
Thus by $p_{k+1}\circ \iota_{k+2}=\varpi_{k+2}$, we get $(\varpi_{k+1})_*\circ \iota_{k+2}=\iota_{k+1}\circ \varpi_{k+2}$. 
Thus the above diagram \eqref{eqn:20220401} is commutative.

We take $y\in J_{k+2}^{\mathrm{reg}}S$.
We want to show 
\begin{equation}\label{eqn:2021110315}
\Phi_*\circ \iota_{k+2}(y)\in \tilde{V}_{\tilde{q}\circ \Phi_*\circ \iota_{k+2}(y)}.
\end{equation}
Let $\tau:V\backslash\{0_V\}\to \tilde{M}$ be the projection.
By the definition of $\Phi$, we have
\begin{equation}\label{eqn:202204011}
\Phi\circ p_{k+1}\circ\iota_{k+2}=\Phi\circ\varpi_{k+2} =\tau\circ\varphi_{*}\circ\iota_{k+1}\circ\varpi_{k+2}=\tau\circ\pi_*\circ\Phi_*\circ\iota_{k+2},
\end{equation}
where the last equality follows from the commutativity of \eqref{eqn:20220401}.
Since $\Phi_*$ is the derivation of $\Phi:J_{k+1}^{\mathrm{reg}}S\to \tilde{M}$, we have $\tilde{q}\circ \Phi_*=\Phi\circ p_{k+1}$.
Combining this with \eqref{eqn:202204011}, we have
$$
\tilde{q}\circ \Phi_*\circ \iota_{k+2}=
\Phi\circ p_{k+1}\circ \iota_{k+2}=
\tau\circ\pi_*\circ\Phi_*\circ\iota_{k+2}.
$$
Thus we observe that $\tilde{q}(\Phi_*\circ \iota_{k+2}(y))\in \tilde{M}$ is the image of $\pi_*(\Phi_*\circ \iota_{k+2}(y))$ under the map $\tau:V\backslash\{0\}\to \tilde{M}$.
Hence by the definition of $\tilde{V}\subset T\tilde{M}$ (cf. \eqref{eqn:20211104}), we get \eqref{eqn:2021110315}.
By the assumption of $\varphi$, we have $\pi_*(\Phi_*\circ \iota_{k+2}(y))\not\in 0_{TM}$.
Hence $\Phi_*\circ \iota_{k+2}(y)\not\in 0_{T\tilde{M}}$.
Hence $\Phi_*\circ \iota_{k+2}(y)\in \tilde{V}\backslash \{0_V\}$.
The proof is completed.
\hspace{\fill} $\square$

\medskip

For each $k\geq 0$, we define a morphism $\varphi_k:J_k^{\mathrm{reg}}S\to S_k$ inductively as follows.
We set $\varphi_0:S\to S$ to be the identity map, where $S_0=S$ and $J_0^{\mathrm{reg}}S=S$.
Then we have $(\varphi_{0})_*(\iota_{1}(J_{1}^{\mathrm{reg}}S))\subset V_0\backslash \{0\}$, where $V_0=TS$.
Suppose we have constructed $\varphi_k:J_k^{\mathrm{reg}}S\to S_k$ which satisfies $(\varphi_{k})_*(\iota_{k+1}(J_{k+1}^{\mathrm{reg}}S))\subset V_k\backslash \{0\}$.
Then we define $\varphi_{k+1}:J_{k+1}^{\mathrm{reg}}S\to S_{k+1}$ by the composite of the following morphisms:
 $$
J_{k+1}^{\mathrm{reg}}S\overset{(\varphi_{k})_*\circ\iota_{k+1}}{\longrightarrow} V_k\backslash\{0\}\to S_{k+1}.
$$
Then by Lemma \ref{lem:202111041}, we have $(\varphi_{k+1})_*(\iota_{k+2}(J_{k+2}^{\mathrm{reg}}S))\subset V_{k+1}\backslash \{0\}$.
Hence we have constructed $\varphi_k:J_k^{\mathrm{reg}}S\to S_k$ inductively.
By the construction, we have
\begin{equation}\label{eqn:202111043}
(\varphi_{k})_*(\iota_{k+1}(J_{k+1}^{\mathrm{reg}}S))\subset V_k\backslash \{0\}
\end{equation}
for all $k\geq 0$.
By \eqref{eqn:202111041}, the following diagram commutes for all $k\geq 1$:
\begin{equation*}
\begin{CD}
J_k^{\mathrm{reg}}S @>\varphi_k>>   S_k\\
@V\varpi_{k}VV  @VV\pi_kV \\
J_{k-1}^{\mathrm{reg}}S @>>\varphi_{k-1}> S_{k-1}
\end{CD}
\end{equation*}
Hence for each $s\in J_{k-1}^{\mathrm{reg}}S$, we get the restriction map
\begin{equation}\label{eqn:20211105}
\varphi_k|_{\varpi_k^{-1}(s)}:\varpi_k^{-1}(s)\to \pi_k^{-1}(\varphi_{k-1}(s)).
\end{equation}
Here $\varpi_k^{-1}(s)\simeq \mathbb C^{\dim S}$ and $\pi_k^{-1}(\varphi_{k-1}(s))\simeq \mathbb P^{\dim S-1}$.
We set $S_k^{\mathrm{reg}}=S_k\backslash S_k^{\mathrm{sing}}$ (cf. \eqref{eqn:202206093}).

\begin{lem}\label{lem:202111031}
Let $k\geq 1$.
Then for each $s\in J_{k-1}^{\mathrm{reg}}S$, the map \eqref{eqn:20211105} is smooth and $\varphi_k(\varpi_k^{-1}(s))=\pi_k^{-1}(\varphi_{k-1}(s))\cap S_k^{\mathrm{reg}}$.
\end{lem}

{\it Proof.}\
We prove this by induction on $k$.
Thus we assume the lemma for $k$ and prove the lemma for $k+1$.
We take $s\in J_k^{\mathrm{reg}}S$.
We first note that
\begin{equation}\label{eqn:202111047}
(\varphi_{k})_{*,s}(T_{J_k^{\mathrm{reg}}S/J_{k-1}^{\mathrm{reg}}S,s})=T_{S_k/S_{k-1},\varphi_k(s)},
\end{equation}
where $(\varphi_{k})_{*,s}:T_sJ_k^{\mathrm{reg}}S\to T_{\varphi_k(s)}S_k$ is the induced map.
We prove this from the induction hypothesis as follows.
We have
$$
T_{J_k^{\mathrm{reg}}S/J_{k-1}^{\mathrm{reg}}S,s}=T_{s}(\varpi_k^{-1}(\varpi_k(s)))\overset{(\varphi_{k})_{*,s}}{\longrightarrow} T_{\varphi_k(s)}(\pi_k^{-1}(\varphi_{k-1}(\varpi_k(s))))=T_{S_k/S_{k-1},\varphi_k(s)},
$$
where we note the smoothness of $\varpi_k$ and $\pi_k$ on the first and last equality.
By the induction hypothesis, \eqref{eqn:20211105} is smooth for $k$.
Hence we get \eqref{eqn:202111047}.

Now we consider the following commutative diagram, where $(\varpi_k)_*\circ\iota_{k+1}=\iota_k\circ\varpi_{k+1}$ follows from \eqref{eqn:20220404} (cf. \eqref{eqn:20220401}).
\begin{equation}\label{eqn:202206094}
\begin{CD}
J_{k+1}^{\mathrm{reg}}S @>\iota_{k+1}>>  TJ^{\mathrm{reg}}_kS @>(\varphi_k)_*>> TS_k @<<< V_k\backslash \{0_{V_k}\}  @>>> S_{k+1}\\
@V\varpi_{k+1}VV  @VV(\varpi_k)_* V @VV(\pi_k)_*V @. @VV\pi_{k+1}V\\
J_{k}^{\mathrm{reg}}S @>>\iota_k> TJ^{\mathrm{reg}}_{k-1}S @>>(\varphi_{k-1})_*> TS_{k-1} @<<< V_{k-1}\backslash \{0_{V_{k-1}}\}  @>>> S_{k}
\end{CD}
\end{equation}
By \eqref{eqn:20211103}, we have
\begin{equation}\label{eqn:20220328}
\iota_{k+1}(\varpi_{k+1}^{-1}(s))=\{ x\in T_sJ^{\mathrm{reg}}_kS\ ;\ (\varpi_{k})_*(x)=\iota_k(s)\},
\end{equation}
where $\iota_k(s)\in T_{\varpi_k(s)}J^{\mathrm{reg}}_{k-1}S$.
We focus on the following two linear maps from \eqref{eqn:202206094}:
\begin{equation*}
\begin{CD}
T_{s}J^{\mathrm{reg}}_{k}S @>(\varphi_k)_{*,s}>>   T_{\varphi_k(s)}S_k\\
@V(\varpi_{k})_{*,s}VV  \\
T_{\varpi_k(s)}J^{\mathrm{reg}}_{k-1}S 
\end{CD}
\end{equation*}
By \eqref{eqn:20220328}, we have $(\varpi_k)_{*,s}^{-1}(\iota_k(s))=\iota_{k+1}(\varpi_{k+1}^{-1}(s))$.
Hence by $(\varpi_k)_{*,s}^{-1}(0)=T_{J^{\mathrm{reg}}_kS/J^{\mathrm{reg}}_{k-1}S,s}$, we note that $\iota_{k+1}(\varpi_{k+1}^{-1}(s))$ is a translate of the linear subspace $T_{J^{\mathrm{reg}}_kS/J^{\mathrm{reg}}_{k-1}S,s}$ in the linear space $T_{s}J^{\mathrm{reg}}_{k}S$.
Hence by \eqref{eqn:202111047}, we observe that $(\varphi_k)_{*,s}(\iota_{k+1}(\varpi_{k+1}^{-1}(s)))$ is a translate of $T_{S_k/S_{k-1},\varphi_k(s)}$.
By \eqref{eqn:202111043}, we have 
\begin{equation}\label{eqn:202203284}
(\varphi_k)_{*,s}(\iota_{k+1}(\varpi_{k+1}^{-1}(s)))\subset V_{k,\varphi_k(s)}\backslash \{0\}.
\end{equation}
On the other hand, we have $T_{S_k/S_{k-1},\varphi_k(s)}\subset V_{k,\varphi_k(s)}$ (cf. \eqref{eqn:202203281}) and $P V_{k,\varphi_k(s)}=\pi_{k+1}^{-1}(\varphi_k(s))$.
By $\varpi_k(s)\in J^{\mathrm{reg}}_{k-1}S$, the induction hypothesis yields that $\varphi_k(s)\in S_{k}^{\mathrm{reg}}$.
Hence the hyperplane $PT_{S_k/S_{k-1},\varphi_k(s)}\subset PV_{k,\varphi_k(s)}$ is equal to $\pi_{k+1}^{-1}(\varphi_k(s))\cap S_{k+1}^{\mathrm{sing}}$ (cf. \eqref{eqn:202203282}).
Since $(\varphi_k)_{*,s}(\iota_{k+1}(\varpi_{k+1}^{-1}(s)))$ is the translate of $T_{S_k/S_{k-1},\varphi_k(s)}$ in the linear space $V_{k,\varphi_k(s)}$,
\eqref{eqn:202203284} yields that 
$$(\varphi_k)_{*,s}(\iota_{k+1}(\varpi_{k+1}^{-1}(s)))\to \pi_{k+1}^{-1}(\varphi_k(s))\cap S_{k+1}^{\mathrm{reg}}$$ 
is an isomorphism under the restriction of the projectivization $V_{k,\varphi_k(s)}\backslash\{ 0\}\to PV_{k,\varphi_k(s)}$.
Now we look the two morphisms
$$
\iota_{k+1}(\varpi_{k+1}^{-1}(s))\overset{(\varphi_k)_{*,s}}{\longrightarrow} (\varphi_k)_{*,s}(\iota_{k+1}(\varpi_{k+1}^{-1}(s)))\to \pi_{k+1}^{-1}(\varphi_k(s))\cap S_{k+1}^{\mathrm{reg}},
$$
where the first map is a translate of the linear map $T_{J_kS/J_{k-1}S,s}\to T_{S_k/S_{k-1},\varphi_k(s)}$ which is surjective by \eqref{eqn:202111047}.
Hence $\varpi_{k+1}^{-1}(s)\to \pi_{k+1}^{-1}(\varphi_{k}(s))\cap S_{k+1}^{\mathrm{reg}}$ is smooth and $\varphi_{k+1}(\varpi_{k+1}^{-1}(s))=\pi_{k+1}^{-1}(\varphi_{k}(s))\cap S_{k+1}^{\mathrm{reg}}$.
This completes the induction step.
\hspace{\fill} $\square$

\begin{lem}\label{lem:202111045}
For each $k\geq 1$ we have $\varphi_k(J_k^{\mathrm{reg}}S)=S_k^{\mathrm{reg}}$.
\end{lem}

{\it Proof.}\
The proof is by induction on $k$.
The case $k=1$ is trivial.
So we assume the case $k-1$ and prove the case $k$.
Let $x\in S_k^{\mathrm{reg}}$.
Then $\pi_k(x)\in S_{k-1}^{\mathrm{reg}}$.
Hence by the induction hypothesis, there exists $s\in J_{k-1}^{\mathrm{reg}}S$ such that $\varphi_{k-1}(s)=\pi_k(x)$.
Then by Lemma \ref{lem:202111031}, there exists $s'\in \varpi_k^{-1}(s)$ such that $\varphi_k(s')=x$.
\hspace{\fill} $\square$

\begin{lem}\label{lem:202203296}
Let $\eta:\Lambda_l\to S$ be a regular $l$-jet.
Then we have $\varphi_k\circ j_k\eta=\eta_{[k]}$ as elements in $\mathrm{Hom}(\Lambda_{l-k},S_{k})$.
\end{lem}

{\it Proof.}\
The proof is by induction on $k$.
The case $k=0$ is trivial.
So we assume the case $k-1$ and prove the case $k$.
By the induction hypothesis, we have $(\varphi_{k-1})_*\circ (j_{k-1}\eta)'=(\eta_{[k-1]})'$.
By \eqref{eqn:202203291}, we have $(\varphi_{k-1})_*\circ (j_{k-1}\eta)'=(\varphi_{k-1})_*\circ\iota_k\circ j_k\eta$.
Thus we get $(\varphi_{k-1})_*\circ\iota_k\circ j_k\eta=(\eta_{[k-1]})'$.
We composite these maps with $V_{k-1}\backslash \{ 0\}\to S_{k}$.
Then by the definitions of $\varphi_k$ and $\eta_{[k]}$, we have $\varphi_k\circ j_k\eta=\eta_{[k]}$.
This completes the induction step.
\hspace{\fill} $\square$

\begin{lem}\label{lem:202203292}
Let $w\in J^{\mathrm{reg}}_kS$.
Then there exists a regular $k$-jet $\eta:\Lambda_k\to S$ such that $j_k\eta(0)=w$.
\end{lem}

{\it Proof.}\
We first consider the case $S=\mathbb A^1$.
Let $x$ be the coordinate of $\mathbb A^1$.
Then $J_k\mathbb A^1=\mathbb A^{k+1}$, where $x,dx,\ldots,d^kx$ are the coordinate functions of $J_k\mathbb A^1$.
Let $w=(w_0,w_1,\ldots,w_k)$.
We define $\eta:\Lambda_k\to \mathbb A^1$ by
$$
x=w_0+w_1\varepsilon+\frac{w_2}{2!}\varepsilon^2+\cdots+\frac{w_k}{k!}\varepsilon^k.
$$
Then we have $j_k\eta(0)=w$.
By $w\in J^{\mathrm{reg}}_kS$, we have $w_1\not=0$.
Hence $\eta$ is regular.
This proves our lemma when $S=\mathbb A^1$.

Next we consider the case $S=\mathbb A^n$.
In this case, we have the natural splitting $J_k\mathbb A^n=(J_k\mathbb A^1)^n$.
Let $p_i:J_k\mathbb A^n\to J_k\mathbb A^1$ be the $i$-th projection.
We take $\eta_i:\Lambda_k\to \mathbb A^1$ such that $j_k\eta_i(0)=p_i(w)$. 
We set $\eta=(\eta_1,\ldots,\eta_k)$.
Then we have $j_k\eta(0)=w$.
By $w\in J^{\mathrm{reg}}_k\mathbb A^n$, there exists $i$ such that $p_i(w)\in J^{\mathrm{reg}}_k\mathbb A^1$.
Then $\eta_i$ is regular.
Hence $\eta$ is regular.
This proves our lemma when $S=\mathbb A^n$.

In general, we may assume that $S$ is affine and has local coordinate functions $\chi_1,\ldots,\chi_n$ described in Remark \ref{rem:20220329}.
Then $\chi=(\chi_1,\ldots,\chi_n):S\to \mathbb A^n$ is {\'e}tale.
This induces $\chi_*:J_kS\to J_k\mathbb A^n$ so that $d^jx_i\circ \chi_*=d^j\chi_i$.
We have $\chi_*(w)\in J_k^{\mathrm{reg}}\mathbb A^n$.
By the previous step, there exists a regular $k$-jet $\xi:\Lambda_k\to \mathbb A^n$ such that $j_k\xi(0)=\chi_*(w)$.
We take $\eta:\Lambda_k\to S$ such that $\chi\circ\eta=\xi$ and $\eta(0)$ is the image of $w$ under $J_kS\to S$.
Then by $\chi_*\circ j_k\eta=j_k\xi$, we have $j_k\eta(0)=w$.
\hspace{\fill} $\square$

\begin{cor}\label{cor:20220320}
Let $w\in S_k^{\mathrm{reg}}$.
Then there exists a regular $k$-jet $\eta:\Lambda_k\to S$ such that $\eta_{[k]}(0)=w$.
\end{cor}

{\it Proof.}\
By Lemma \ref{lem:202111045}, there exists $w'\in J^{\mathrm{reg}}_kS$ such that $\varphi_k(w')=w$.
By Lemma \ref{lem:202203292}, there exists a regular $k$-jet $\eta:\Lambda_k\to S$ such that $j_k\eta(0)=w'$.
By Lemma \ref{lem:202203296}, we have $\eta_{[k]}(0)=w$.
This concludes the proof of the corollary.
See also \cite[Thm 6.8]{Dem}.
\hspace{\fill} $\square$

\subsection{One lemma for regular jets}

\begin{lem}\label{lem:20201224}
Let $S$ be a smooth variety.
Let $Z\subset S$ be a closed subscheme.
Then there exists a positive integer $k$ with the following property:
Let $\eta:\Lambda_k\to S$ be a regular $k$-jet with non-zero first derivative $\eta'(0)=v\in T_{\eta(0)}S$, where $v\not=0$.
Assume that $\eta$ factors $Z\subset S$.
Then there exists an irreducible component $Z'\subset Z_{\mathrm{red}}$ such that $[v]\in \mathcal{D}Z'\subset S_1$.
\end{lem}

{\it Proof.}\
We may assume that $S$ is affine.
We first prove a weaker statement $[v]\in \mathcal{D}(Z_{\mathrm{red}})$.
Let $Z_{\mathrm{red}}\subset S$ be defined by $\varphi_1,\ldots,\varphi_l$.  
We consider the map $\Phi:S\to \mathbb A^l$ defined by $\varphi_1,\ldots,\varphi_l$.
We have $\Phi^{-1}(0)=Z_{\mathrm{red}}$.
There exists $k$ such that $0_k=\mathrm{Spec} (\mathcal{O}_{0,\mathbb A^l}/\mathfrak{m}^k)\subset \mathbb A^l
$ satisfies $Z\subset \Phi^*0_k$.
Note that $\varphi_1,\ldots,\varphi_l,d\varphi_1,\ldots,d\varphi_l$ defines $\mathcal{D}(Z_{\mathrm{red}})\subset S_1$.
Let $\eta:\Lambda_k\to S$ be a regular $k$-jet with non-zero first derivative $\eta'(0)=v\in T_{\eta(0)}S$, where $v\not=0$, such that $\eta$ factors $Z\subset S$.
To show $[v]\in \mathcal{D}(Z_{\mathrm{red}})$, we assume contrary that $[v]\not\in \mathcal{D}(Z_{\mathrm{red}})$.
Then $\Phi\circ \eta:\Lambda_k\to \mathbb A^l$ is regular.
Since $\Phi\circ \eta$ factors $0_k\subset \mathbb A^l$, we have $((\Phi\circ\eta)^*x_1)^k=0,\ldots, ((\Phi\circ\eta)^*x_l)^k=0$, where $(\Phi\circ\eta)^*x_1,\ldots,(\Phi\circ\eta)^*x_l\in \mathbb C[\varepsilon]/(\varepsilon^{k+1})$.
Since $\Phi\circ \eta$ is regular, we may take $i=1,\ldots,l$ such that  $(\Phi\circ\eta)^*x_i=a_1\varepsilon+a_2\varepsilon^2+\cdots+a_k\varepsilon^k$ with $a_1\not=0$.
Then $((\Phi\circ\eta)^*x_i)^k=a_1^k\varepsilon^k\not=0$.
This is a contradiction.
Hence $[v]\in \mathcal{D}(Z_{\mathrm{red}})$.

Now let $Z_1,\ldots, Z_l$ be irreducible closed subschemes of $S$ such that $\supp Z_1,\ldots ,\supp Z_l$ are the irreducible components of $\supp Z$ and $Z\subset Z_1+\cdots+Z_l$.
Let $k_1,\ldots,k_l$ be positive integers which are obtained in the previous step for $Z_1,\ldots,Z_l$.
Set $k=k_1+\cdots+k_l$.
Let $\eta:\Lambda\to S$ be a regular $k$-jet with non-zero first derivative $\eta'(0)=v\in T_{\eta(0)}S$, where $v\not=0$, such that $\eta$ factors $Z\subset S$, hence $Z_1+\cdots+Z_l$.
Let $\eta_i:\Lambda_{k_i}\to S$, where $\Lambda_{k_i}=\mathrm{Spec}\mathbb C[\varepsilon]/(\varepsilon^{k_i+1})$, be induced from $\eta$.
Then there exists $i$ such that $\eta_i$ factors $Z_i$.
Indeed, to see this, we consider closed subschemes $\eta^*Z_i\subset \Lambda_k$.
We may write the defining ideal of $\eta^*Z_i$ as $(\varepsilon^{m_i})\subset \mathbb C[\varepsilon]/(\varepsilon^{k+1})$.
By $Z\subset Z_1+\cdots+Z_l$, we have $\eta^*Z_1+\cdots+\eta^*Z_l=\Lambda_k$.
The defining ideal of $\eta^*Z_1+\cdots+\eta^*Z_l$ is $(\varepsilon^{m_1+\cdots+m_l})\subset \mathbb C[\varepsilon]/(\varepsilon^{k+1})$.
Hence $m_1+\cdots+m_l\geq k+1$.
We may take $i$ such that $m_i\geq k_i+1$.
Hence $\eta_i:\Lambda_{k_i}\to S$ factors $Z_i$.
Then by the previous step, we have $[v]\in \mathcal{D}((Z_i)_{\mathrm{red}})$.
Note that $(Z_i)_{\mathrm{red}}$ is an irreducible component of $Z_{\mathrm{red}}$.
\quad $\square$

\section{Sufficient condition for horizontal integrability}

The main result of this section is Lemma \ref{lem:alglem}.
This gives a sufficient condition for $Z\subset S_{1,A/B}$ to be horizontally integrable.
This lemma is used in the proof of Proposition \ref{cor:20210306}.

Let $\overline{A}$ be an equivariant compactification of a semi-abelian variety $A$. 
Let $S$ be a smooth variety.
Let $X\subset \overline{A}\times S$ be a closed subscheme.
Let $\mathcal{X}\subset \overline{A}\times A\times S$ be the pull-back of $X$ by the action $m:\overline{A}\times A\times S\to \overline{A}\times S$ defined by $(x,a,s)\mapsto (x+a,s)$ so that $\mathcal{X}_{(a,s)}=X_{s}-a$.
Suppose that $X\subset \overline{A}\times S$ is $B$-invariant.
Then by Lemma \ref{lem:quotient} we have a closed subscheme $\mathcal{X}_B\subset \overline{A}\times (A/B)\times S$ such that $\mathcal{X}$ is the pull-back of $\mathcal{X}_B$ by the quotient $\overline{A}\times A\times S\to \overline{A}\times (A/B)\times S$ on the second factor.
We get the closed subscheme $\mathcal{P}_k\mathcal{X}_B\subset \overline{A}\times((A/B)\times S)_k$.
By the isomorphism \eqref{eqn:20211203}, we have 
$$S_{k,A/B}=\{ 0_{A/B}\}\times S_{k,A/B}\subset (A/B)\times S_{k,A/B}=((A/B)\times S)_k.$$
Using this immersion $S_{k,A/B}\subset ((A/B)\times S)_k$, we define a closed subscheme $X_{B,k}\subset \overline{A}\times S_{k,A/B}$ by
\begin{equation}\label{eqn:202109111}
X_{B,k}=(\mathcal{P}_k\mathcal{X}_B)\vert_{\overline{A}\times S_{k,A/B}}.
\end{equation}

Let $p_{k,B}:\overline{A}\times S_{k,A/B}\to S_{k,A/B}$ be the second projection.
Let $p_{k,B}|_{X_{B,k}}:X_{B,k}\to S_{k,A/B}$ be the composite of the closed immersion $X_{B,k}\hookrightarrow \overline{A}\times S_{k,A/B}$ and $p_{k,B}$.
Given $y\in S_{k,A/B}$, we denote by $(X_{B,k})_y\subset \overline{A}$ the scheme theoretic fiber of $p_{k,B}|_{X_{B,k}}:X_{B,k}\to S_{k,A/B}$ over $y\in S_{k,B}$.
For $k\geq l$, we have a natural morphism $\mathcal{P}_k\mathcal{X}_B\to \mathcal{P}_l\mathcal{X}_B$.
This induces the following commutative diagram:
\begin{equation*}
\begin{CD}
X_{B,k} @>>>   X_{B,l}\\
@Vp_{k,B}|_{X_{B,k}}VV  @VVp_{l,B}|_{X_{B,l}}V \\
S_{k,A/B} @>>> S_{l,A/B}
\end{CD}
\end{equation*}
For $y\in S_{k,A/B}$, let $y'\in S_{l,A/B}$ be the image of $y$ under the map $S_{k,A/B}\to S_{l,A/B}$.
Then we have the map $(X_{B,k})_{y}\to (X_{B,l})_{y'}$ of the closed subschemes of $\overline{A}$.

Let $V\subset \overline{A}$ be a closed subscheme.
We define $\mathrm{Stab}(V)\subset A$ by $a\in \mathrm{Stab}(V)$ if and only if $a+V=V$ as closed subschemes of $\overline{A}$.
Let $\mathrm{Stab^0}(V)$ be the connected component of $\mathrm{Stab}(V)$ containing the identity element of $A$.
Then $\mathrm{Stab^0}(V)\subset A$ is a semi-abelian subvariety.

Let $k\geq 2$.
We recall $S_{k,A}^{\mathrm{sing}}\subset S_{k,A}$ from \eqref{eqn:20220129}.
Set $S_{k,A}^{\mathrm{reg}}=S_{k,A}\backslash S_{k,A}^{\mathrm{sing}}$.
Then we have $(A\times S)_k^{\mathrm{reg}}=A\times S_{k,A}^{\mathrm{reg}}$ (cf. \eqref{eqn:20220611}).
The purpose of this section is to prove the following lemma.

\begin{lem}\label{lem:alglem}
Let $\overline{A}$ be an equivariant compactification of $A$, where $\overline{A}$ is projective.
Let $B\subsetneqq A$ be a proper semi-abelian subvariety.
Let $Z\subset S_{1,A/B}$ be an irreducible Zariski closed set.
Let $X\subset \overline{A}\times S$ be a $B$-invariant closed subscheme such that $Z\subset p_{1,B}(X_{B,1})$, where $p_{1,B}:\overline{A}\times S_{1,A/B}\to S_{1,A/B}$ is the second projection.
Assume that for every integer $k\geq 2$, there exists a non-empty Zariski open subset $O_k\subset Z$ such that for every $y\in O_k$ the followings hold:
\begin{enumerate}
\item
$\mathrm{Stab^0}((X_{B,1})_y)=B$.
\item 
The natural map $(X_{B,1})_y\to X_{\tau (y)}$ is an isomorphism as schemes, where  $\tau :S_{1,A/B}\to S$ is the induced map.
\item 
There exists $y_k\in S_{k,A/B}^{\mathrm{reg}}$ such that the image of $y_k$ under $S_{k,A/B}\to S_{1,A/B}$ is $y$, and that the map $(X_{B,k})_{y_k}\to (X_{B,1})_{y}$ is an isomorphism as schemes.
\end{enumerate}
Then $Z\subset S_{1,A/B}$ is horizontally integrable.
\end{lem}

According to Definition \ref{defn:20201122}, the conclusion of Lemma \ref{lem:alglem} reads the existence of an immersion $U\hookrightarrow (A/B)\times S$ with the following properties:
\begin{itemize}
\item 
$q:U\to S$ is {\'e}tale.
\item
$p'(PTU)\cap Z\subset Z$ is Zariski dense in $Z$, where $p':(A/B)\times S_{1,A/B}\to S_{1,A/B}$ is the second projection.
\end{itemize}

Before going to prove Lemma \ref{lem:alglem}, we start from algebro-geometric lemmas.

\begin{lem}\label{lem:20220405}
Let $\Sigma$ and $S$ be algebraic varieties such that $\dim \Sigma =\dim S $.
Let $p:\Sigma \to S$ be unramified.
Assume that $S$ is smooth.
Then $p:\Sigma\to S$ is smooth, hence {\'e}tale.
\end{lem}

{\it Proof.}\
Set $d=\dim \Sigma =\dim S$.
We first show that $\Sigma$ is smooth.
We have the exact sequence (cf. \cite[II, Prop 8.11]{H})
$$
p^*\Omega_S\to \Omega_{\Sigma}\to \Omega_{\Sigma /S}\to 0.
$$
Since $p:\Sigma\to S$ is unramified, we have $\Omega_{\Sigma /S}=0$ (cf. \cite[p. 221]{Liu}).
Hence the morphism $p^*\Omega_S\to \Omega_{\Sigma}$ is surjective.
Since $p^*\Omega_S$ is locally free of rank $d$, we have $\dim \Omega_{\Sigma}\otimes \mathbb C (x)\leq d$ for all $x\in \Sigma$.
For general $x\in \Sigma$, we have $\dim \Omega_{\Sigma}\otimes \mathbb C (x)= d$, so by the upper semicontinuity, this holds for all $x\in \Sigma$.
Since $\Sigma$ is reduced, this shows $\Omega_{\Sigma}$ is locally free of rank $d$ (cf. \cite[II, Ex. 5.8]{H}).
Hence $\Sigma$ is smooth.

Since $p^*\Omega_S\to \Omega_{\Sigma}$ is surjective, the induced map $\Omega_S\otimes \mathbb C(p(x))\to \Omega_{\Sigma}\otimes \mathbb C(x)$ is surjective for all $x\in \Sigma$.
Note that $\Omega_S\otimes \mathbb C(p(x))\to \Omega_{\Sigma}\otimes \mathbb C(x)$ is the dual of the induce map on the tangent spaces $p_*:T_{\Sigma,x}\to T_{S,p(x)}$.
Hence $p_*:T_{\Sigma,x}\to T_{S,p(x)}$ is injective for all $x\in\Sigma$.
Hence by $\dim \Sigma =\dim S$, we obtain that $p_*:T_{\Sigma,x}\to T_{S,p(x)}$ is surjective for all $x\in\Sigma$.
By \cite[III, Prop 10.4]{H}, $p:\Sigma\to S$ is smooth.
See also \cite[p. 141]{Liu}.
\hspace{\fill} $\square$

\medskip

Next we apply the previous lemma to prove the following lemma.

\begin{lem}\label{lem:20210325}
Let $S$ and $W$ be smooth algebraic varieties.
Let $p:W\to S$ be smooth morphism.
Let $\Sigma_0\subset W$ be an irreducible Zariski closed set and $x\in \Sigma_0$.
Assume that $p|_{\Sigma_0}:\Sigma_0 \to S$ is unramified at $x\in \Sigma_0$.
Then there exists $\Sigma\subset W$ such that $\Sigma_0\subset \Sigma$  and $p|_{\Sigma}:\Sigma \to S$ is {\'e}tale at $x$.
\end{lem}

{\it Proof.}\
Let $h_1,\ldots ,h_k$ be a local defining functions $\Sigma_0\subset W$ around $x$.
Let $L=p^{-1}(p(x))$ and set $s=\dim L$.
Since $p$ is smooth, $L\subset W$ is a smooth subvariety.
Since $p|_{\Sigma_0}$ is unramified at $x\in \Sigma_0$, we have 
$$\{ v\in T_xL;\ (dh_1|_{T_xL})(v)=\cdots =(dh_k|_{T_xL})(v)=0\}=\{ 0\}.$$
Since $T_xL$ is a $s$-dimensional vector space, we may assume that 
$$\{ v\in T_xL;\ (dh_1|_{T_xL})(v)=\cdots =(dh_s|_{T_xL})(v)=0\}=\{ 0\}.$$
We take $\Sigma\subset W$ such that $\Sigma$ is defined by $h_1=\cdots =h_s=0$ around $x$.
Then $\dim \Sigma \geq \dim S$ and $p|_{\Sigma}:\Sigma \to S$ is unramified at $x$.
Hence $\dim \Sigma = \dim S$.
Since unramified is an open condition, there exists a non-empty Zariski open set $\Sigma'\subset \Sigma$ such that $x\in\Sigma'$ and $\Sigma'\to S$ is unramified.
Since $S$ is smooth, $\Sigma'\to S$ is {\'e}tale (cf. Lemma \ref{lem:20220405}).
Hence $\Sigma\to S$ is {\'e}tale at $x$.
Since $\Sigma_0$ is irreducible, we have $\Sigma_0\subset \Sigma$.
\hspace{\fill} $\square$

\begin{lem}\label{lem:202206142}
Let $M$, $\Sigma$ and $H$ be algebraic varieties.
Let $c:M\times \Sigma \to H$ be a dominant morphism.
Let $V\subset \Sigma$ be a closed subvariety.
Assume that there exists a non-empty Zariski open set $W\subset M\times V$ such that the restriction $c|_W:W\to \overline{c(M\times V)}$ is {\'e}tale.
Let $a_0\in M$ satisfies $(\{a_0\}\times V)\cap W\not=\emptyset$.
Let $F\subset M\times \Sigma$ be an irreducible component of $\supp c^{-1}(\overline{c(\{a_0\}\times V)})$ such that $\{a_0\}\times V\subset F$.
Then we have the followings:
\begin{enumerate}
\item
For $t\in V$ with $(a_0,t)\in W$, the induced map $F\to \Sigma$ is unramified at $(a_0,t)\in F$. 
\item
Suppose moreover that the restriction $c|_{M\times V}:M\times V\to H$ is dominant and generically finite.
Then $F\to \Sigma$ is dominant and generically finite.
\end{enumerate}
\end{lem}

{\it Proof.}\
We take $t\in V$ such that $(a_0,t)\in W$.
Let $F_t$ be the scheme theoretic fiber of the map $F\to \Sigma$ over $t\in V$.
Then $F_t=F\cap (M\times\{ t\})$, where the intersection is taken scheme theoretically.
We are going to prove $\mathcal{O}_{F_t,(a_0,t)}=\mathbb C$.
Since $c|_W:W\to \overline{c(M\times V)}$ is {\'e}tale, the scheme theoretic intersection $F\cap (M\times V)$ coincides with $\{a_0\}\times V$ on some Zariski open neighbourhood $U\subset W$ of $(a_0,t)\in W$.
Namely $F\cap U=\{a_0\}\times V$.
Hence 
$$F_t\cap U= U\cap (\{a_0\}\times V)\cap (M\times\{ t\})=\{ (a_0,t)\},
$$
where $\{ (a_0,t)\}$ is a reduced scheme.
Hence $F_t\cap U=\mathrm{Spec}\ \mathbb C$.
Since $F_t\hookrightarrow M\times \Sigma$ factors $F_t\hookrightarrow M\times V\hookrightarrow M\times \Sigma$, we have $\mathcal{O}_{F_t,(a_0,t)}=\mathbb C$.
Hence $F\to \Sigma$ is unramified at $(a_0,t)$.

Next suppose $c|_{M\times V}:M\times V\to H$ is dominant and generically finite.
Note that all the irreducible components of fibers of $M\times \Sigma\to H$ have dimension greater than or equal to $\dim(M\times \Sigma)-\dim H$.
Hence we have
$$
\dim F\geq \dim( M\times \Sigma)-\dim H+\dim \overline{c(\{ a_0\}\times V)}.
$$
Since $c|_{M\times V}:M\times V\to H$ is dominant and generically finite, we have $\dim H=\dim (M\times V)$.
By the choice of $a_0$, we have $\dim \overline{c(\{ a_0\}\times V)}=\dim V$.
Hence
\begin{equation*}
\dim F\geq \dim( M\times \Sigma)-\dim (M\times V)+\dim V
=\dim \Sigma.
\end{equation*}
Thus $\dim F\geq \dim \Sigma$.
Since $F\to \Sigma$ is unramified at $(a_0,t)\in F$, where $(a_0,t)\in W$, the map $F\to \Sigma$ is dominant and generically finite.
\hspace{\fill} $\square$

\begin{lem}\label{lem:20220614}
Let $M$, $\Sigma$ and $H$ be algebraic varieties.
Let $c:M\times \Sigma \to H$ be a dominant morphism.
Assume that for generic $s\in\Sigma$, the restriction $c|_{M\times \{ s\}}:M\times \{ s\}\to H$ is quasi-finite.
Then there exists a closed subvariety $V\subset \Sigma$ such that $c|_{M\times V}:M\times V\to H$ is dominant and generically finite.
\end{lem}

{\it Proof.}\
Set $d=\dim (M\times \Sigma)-\dim H$.
Let $s\in \Sigma$ and $m\in M$ satisfy the followings:
\begin{itemize}
\item
$c|_{M\times \{ s\}}:M\times \{ s\}\to H$ is quasi-finite.
\item
Set $h=c((m,s))\in H$ and $c^{-1}(h)=X$.
Then all irreducible components of $X$ have dimension equal to $d$ (cf. \cite[II, Ex. 3.22]{H}).
\end{itemize}
Let $p:M\times \Sigma\to\Sigma$ be the second projection.
By $(m,s)\in X$, we have $s\in p(X)$.
Set $\overline{p(X)}=Y$.
Then all irreducible components of $Y$ have dimension equal to or less than $d$.
We take a closed subvariety $V\subset \Sigma$ of codimension equal to $d$ such that $\{ s\}=\supp (V\cap Y)$ on some Zariski open neighbourhood $U\subset \Sigma$ of $s\in \Sigma$. 
By $p^{-1}(s)\cap X=c|_{M\times\{s\}}^{-1}(h)$, the set $p^{-1}(s)\cap X$ is finite.
Hence
$$
c^{-1}(h)\cap (M\times V)\cap (M\times U)=p^{-1}(s)\cap X
$$
consists of finite points.
We note that this set is non-empty for it contains the point $(m,s)$.

Now we consider the map $c|_{M\times V}:M\times V\to H$.
Then $(c|_{M\times V})^{-1}(h)$ contains zero-dimensional irreducible components.
Moreover $\dim (M\times V)=\dim H$.
Hence $c|_{M\times V}:M\times V\to H$ is dominant and generically finite.
\hspace{\fill} $\square$

\medskip

{\it Proof of Lemma \ref{lem:alglem}.}\
The proof divides into several steps.
In the following argument, we fix an projective embedding $\overline{A}\subset \mathbb P^N$.

\medskip

{\it Step 1.}\
Let $Z_o=\tau (Z)\subset S$.
Let $X_{Z_{o}}\subset \overline{A}\times Z_{o}$ be the base change.
Let $P$ be the Hilbert polynomial of generic fibers of $X_{Z_{o}}\to Z_{o}$.
Let $Z_o^P\subset Z_o$ be a nonempty Zariski open set such that all the fibers over the points of $Z_o^P$ have Hilbert polynomial $P$.
Then since $Z_o^P$ is integral, $X_{Z_o^P}\to Z_o^P$ is flat (cf. \cite[III, Thm 9.9]{H}).
We denote by $\{\mathcal{S}_1,\ldots,\mathcal{S}_l\}$ the flattening stratification of the coherent sheaf $\mathcal{O}_X$ on $\mathbb P^N\times S$.
(See \cite[Thm 4.2.11]{Ser} for the existence of such stratification.)
Namely, each $\mathcal{S}_i$ is a locally closed subscheme of $S$ such that
\begin{itemize}
\item $S=\coprod \mathcal{S}_i$, and
\item a map $T\to S$ from a scheme $T$ factors $\coprod \mathcal{S}_i\to S$ if and only if $X_T\to T$ is flat.
\end{itemize}
We may choose $S^P\in\{\mathcal{S}_i\}$ such that $Z_{o}^P\subset S$ factors as 
$$Z_o^P\subset  S^P.$$
Then $S^P\subset S$ is a locally closed subscheme.

There exists a non-empty Zariski open set $S^o\subset S$ such that $S^P\hookrightarrow S$ factors $S^P\hookrightarrow S^o\subset S$, where $S^P\hookrightarrow S^o$ is a closed immersion.
Replacing $S$ by $S^o$, we may assume that $S^P\subset S$ is a closed subscheme.
Under this reduction, the relation $Z_o^P\subset  S^P$ implies $Z_o\subset S^P$.
Hence $Z_o^P=Z_o$.

\medskip

{\it Step 2.}
Let $\mathcal{X}\subset \overline{A}\times A\times S$ be the pull-back of $X$ by the action $m:\overline{A}\times A\times S\to \overline{A}\times S$, where $(x,a,s)\mapsto (x+a,s)$ so that 
\begin{equation}\label{eqn:20210818}
\mathcal{X}_{(a,s)}=X_{s}-a.
\end{equation}
Since $X$ is $B$-invariant, $\mathcal{X}$ is $B$-invariant under the $B$-action $\overline{A}\times A\times S\to \overline{A}\times A\times S$ defined by $(x,a,s)\mapsto (x,a+b,s)$ for $b\in B$.
Hence there exists a closed subscheme $\mathcal{X}_B\subset \overline{A}\times (A/B)\times S$ such that $\mathcal{X}$ is the pull-back of $\mathcal{X}_B$ by the quotient $\overline{A}\times A\times S\to \overline{A}\times (A/B)\times S$ on the second factor (cf. Lemma \ref{lem:quotient}).

Let $T\to (A/B)\times S$ be a map from a scheme $T$.
Let $(\mathcal{X}_B)_T\to T$ be the pull-back of $\mathcal{X}_B\to (A/B)\times S$ by this map $T\to (A/B)\times S$.
Let $X_T\to T$ be the pull-back of $X\to S$ by the composition of $T\to (A/B)\times S$ and the second projection $(A/B)\times S\to S$.

\medskip

{\it Claim 1.}\
$(\mathcal{X}_B)_T\to T$ is flat if and only if $X_T\to T$ is flat.

\medskip

We prove this.
Let $\mathcal{X}'\subset \overline{A}\times A\times S$ be the pull-back of $X\subset \overline{A}\times S$ by the map $\overline{A}\times A\times S\to \overline{A}\times S$ defined by $(x,a,s)\mapsto (x,s)$.
Let $\mathcal{X}_B'$ be the pull-back of $X$ by the map $\overline{A}\times (A/B)\times S\to \overline{A}\times S$ defined by $(x,a',s)\mapsto (x,s)$.
Then $\mathcal{X}'$ is the pull-back of $\mathcal{X}_B'$ by the quotient $\overline{A}\times A\times S\to \overline{A}\times (A/B)\times S$ on the second factor.
Note that the isomorphism $\kappa:\overline{A}\times A\times S\to \overline{A}\times A\times S$ defined by $(x,a,s)\mapsto (x+a,a,s)$ induces the isomorphism $\kappa|_{\mathcal{X}}:\mathcal{X}\to\mathcal{X'}$ over $A\times S$.

Let $T'\to A\times S$ be the pull-back of $T\to (A/B)\times S$ by the quotient map $A\times S\to (A/B)\times S$.
Then the induced map $T'\to T$ is faithfully flat, for $A\times S\to (A/B)\times S$ is faithfully flat (cf. Remark \ref{rem:20220613}).
Let $\mathcal{X}_{T'}\to T'$ be the pull-back of $\mathcal{X}\to A\times S$ by the map $T'\to A\times S$.
Then we have the following Cartesian diagram:
\begin{equation*}
\begin{CD}
\mathcal{X}_{T'} @>>>   (\mathcal{X}_B)_T \\
@VVV  @VVV \\
T'@>>> T
\end{CD}
\end{equation*}
Hence 
\begin{equation}\label{eqn:20230320}
\text{
$(\mathcal{X}_B)_T\to T$ is flat if and only if $\mathcal{X}_{T'}\to T'$ is flat.}
\end{equation}
Indeed if $\mathcal{X}_{T'}\to T'$ is flat, the composite map $\mathcal{X}_{T'}\to T$ is flat.
Note that the faithful flatness of $T'\to T$ yields that of $\mathcal{X}_{T'}\to (\mathcal{X}_B)_T$, for the faithful flatness is stable under the fiber product.
Hence by \cite[Cor. 14.12]{Gortz}, we obtain that $ (\mathcal{X}_B)_T\to T$ is flat.
Conversely, the flatness of $ (\mathcal{X}_B)_T\to T$ yields that of $\mathcal{X}_{T'}\to T'$, for the flatness is stable under the fiber product.
Thus we obtain \eqref{eqn:20230320}.

The isomorphism $\kappa$ induces an isomorphism $\mathcal{X}_{T'}\to (\mathcal{X}')_{T'}$ over $T'$.
Hence 
$$
\text{
$\mathcal{X}_{T'}\to T'$ is flat if and only if $(\mathcal{X}')_{T'}\to T'$ is flat.}
$$
By the same reason as \eqref{eqn:20230320}, we note that 
$$
\text{
$(\mathcal{X}')_{T'}\to T'$ is flat if and only if $(\mathcal{X}_B')_{T}\to T$ is flat.}
$$
Note that $(\mathcal{X}_B')_{T}=X_T$.
Hence $(\mathcal{X}_B)_T\to T$ is flat if and only if $X_T\to T$ is flat.
This conclude the proof of the claim.

Now, we note that $X_{S^P}\to S^P$ is flat.
Hence by the claim above, $(\mathcal{X}_B)_{(A/B)\times S^P}\to (A/B)\times S^P$ is flat.
Let $\mathrm{Hilb}$ be the Hilbert scheme of closed subschemes of $\overline{A}$ with Hilbert polynomials $P$ (cf. \cite[Thm 4.3.4]{Ser}).
By $(\mathcal{X}_B)_{(A/B)\times S^P}\to (A/B)\times S^P$, we have the classification map
$$
c:(A/B)\times S^P\to \mathrm{Hilb}.
$$

\medskip

{\it Step 3.}\
For each $(a,s)\in (A/B)\times S^P$, we set $E_{(a,s)}=c^{*}(c(a,s))$, which is the scheme theoretic fiber of $c$ over the point $c(a,s)\in \mathrm{Hilb}$.
Then $E_{(a,s)}\subset (A/B)\times S^P$ is a closed subscheme.
We claim that 
\begin{equation}\label{eqn:202206132}
E_{(a,s)}=a+E_{(0,s)}.
\end{equation}
We prove this.
For $t:T\to (A/B)\times S$, we get $t-a:T\to (A/B)\times S$ by the translation.
We denote by $(\mathcal{X}_B)_t\subset \overline{A}\times T$ the pull-back of $\mathcal{X}_B$ by $t:T\to (A/B)\times S$.
Similarly for $(\mathcal{X}_B)_{t-a}$.
Let $a'\in A$ be a point whose image under $A\to A/B$ is equal to $a$.
Then we note that 
\begin{equation}\label{eqn:202108181}
(\mathcal{X}_B)_t=(\mathcal{X}_B)_{t-a}-a'.
\end{equation}
To show this, it is enough to consider the case that $t=\mathrm{id}:(A/B)\times S\to (A/B)\times S$.
Let $t':A\times S\to A\times S$ be the identity map and $t'-a':A\times S\to A\times S$ be the map defined by $(x,s)\mapsto (x-a',s)$.
Then 
$$\mathcal{X}_{t'}=\mathcal{X}_{t'-a'}-a'.
$$
Note that $\mathcal{X}_{t'}$ (resp. $\mathcal{X}_{t'-a'}-a'$) is the pull-back of $(\mathcal{X}_B)_t$ (resp. $(\mathcal{X}_B)_{t-a}-a'$) under the quotient $\overline{A}\times A\times S\to \overline{A}\times (A/B)\times S$ on the second factor.
Moreover the induced maps $\mathcal{X}_{t'}\to (\mathcal{X}_B)_t$ and $\mathcal{X}_{t'-a'}-a'\to (\mathcal{X}_B)_{t-a}-a'$ are the categorical quotients (cf. Lemma \ref{lem:quotient}).
Hence we get \eqref{eqn:202108181}.

Now we note that  $t:T\to (A/B)\times S^P$ factors $E_{(a,s)}$ if and only if 
$$
(\mathcal{X}_B)_t=(\mathcal{X}_B)_{(a,s)}\times T
$$
as closed subschemes of $\overline{A}\times T$.
We note that
$$(\mathcal{X}_B)_{(a,s)}\times T=\mathcal{X}_{(a',s)}\times T=(X_s-a')\times T=(X_s\times T)-a'.$$ 
Hence $t:T\to (A/B)\times S^P$ factors $E_{(a,s)}$ if and only if
$$
(\mathcal{X}_B)_t=(X_s\times T)-a'
$$
as closed subschemes of $\overline{A}\times T$.
Using this, we note that $t-a:T\to (A/B)\times S^P$ factors $E_{(0,s)}$ if and only if
$$
(\mathcal{X}_B)_{t-a}=X_s\times T
$$
as closed subschemes of $\overline{A}\times T$.
Hence \eqref{eqn:202108181} yields that $t:T\to (A/B)\times S^P$ factors $E_{(a,s)}$ if and only if $t-a:T\to (A/B)\times S$ factors $E_{(0,s)}$.
Thus $E_{(a,s)}=a+E_{(0,s)}$.
This completes the proof of \eqref{eqn:202206132}.

\medskip

{\it Step 4.}\
For each $[v]\in Z$, we have $\tau([v])\in Z_o\subset S^P$.
Hence for each $a\in A/B$, we may consider $E_{(a,\tau([v]))}$.
We prove the following:

\medskip

{\it Claim 2.}\
For all $[v]\in (\cap_{k\geq 2} O_k)\subset Z$ and all $a\in A/B$, there exists an irreducible component $E$ of $\supp E_{(a,\tau([v]))}$ such that 
\begin{equation}\label{eqn:20210330}
(a,[v])\in \mathcal{D}E\subset (A/B)\times S_{1,A/B}
\end{equation}
under the identification $(A/B)\times S_{1,A/B}=((A/B)\times S)_1$.
Here $E\subset (A/B)\times S$ is a Zariski closed subset.

\medskip

We prove this.
Denoting $S'=(A/B)\times S$, we have $(S')_{k}=(A/B)\times S_{k,A/B}$.
For each $k\geq 2$, by our assumptions (2), (3) and Corollary \ref{cor:20220320},  we may take a regular $k$-jet $\eta:\Lambda_k\to S'$ with $\eta'(0)=(0,v)$ such that $(\mathcal{X}_B)_{(0,\tau{([v])})}=(\mathcal{P}_k\mathcal{X}_B)_{\eta_{[k]}(0)}$, where $\Lambda_k=\mathrm{Spec}\ \mathbb C[\varepsilon]/(\varepsilon^{k+1})$.
Hence by Lemma \ref{lem:20210324}, we have 
\begin{equation}\label{eqn:20220411}
(\mathcal{X}_B)_{\eta}=(\mathcal{X}_B)_{(0,\tau{([v])})}\times \Lambda_k.
\end{equation}
Hence $(\mathcal{X}_B)_{\eta}\to \Lambda_k$ is flat.
Let $\eta_2:\Lambda_k\to S$ be the compositions with $\eta$ and the second projection $(A/B)\times S\to S$.
Then by Claim 1 in the step 2, $X_{\eta_2}\to\Lambda_k$ is flat.
By $\eta_2(0)=\tau([v])\in Z_o$, $\eta_2:\Lambda_k\to S$ factors $S^P$.
Thus $\eta:\Lambda_k\to S'$ factors $(A/B)\times S^P$, hence $E_{(0,\tau([v]))}$ (cf. \eqref{eqn:20220411}).
This holds for all $k$.
Hence by Lemma \ref{lem:20201224}, we have
$$
(0,[v])\in \mathcal{D}E'
$$
for some irreducible component $E'$ of $\supp E_{(0,\tau([v]))}$.
By \eqref{eqn:202206132}, we have $E_{(a,\tau([v]))}=a+E_{(0,\tau([v]))}$.
Hence $E=a+E'$ is an irreducible component of $\supp (E_{(a,\tau([v]))})$.
By $\mathcal{D}E=a+\mathcal{D}E'$, we complete the proof of the claim.

\medskip

{\it Step 5.}\
We are going to take a closed subvariety $V\subset Z_o$ such that the restriction 
$$c|_{(A/B)\times V}:
(A/B)\times V\to \overline{c((A/B)\times Z_o)}$$ 
is dominant, and generically finite.
Here $\overline{c((A/B)\times Z_o)}\subset \mathrm{Hilb}$ is the Zariski closure in $\mathrm{Hilb}$.
Note that, the set $\tau(O_2)\subset Z_o$ is dense and constructible (cf. \cite[II, Ex. 3.19]{H}), hence contains non-empty Zariski open set.
Hence for generic $t\in Z_o$, by the assumptions (1), (2), we have $\mathrm{Stab^0}(X_t)=B$.
Hence the restriction $c|_{(A/B)\times \{ t\} }:(A/B)\times \{ t\} \to \mathrm{Hilb}$ is quasi-finite. 
We apply Lemma \ref{lem:20220614} to $(A/B)\times Z_o\to \overline{c((A/B)\times Z_o)}$ to get our $V\subset Z_o$.

Let $W\subset (A/B)\times V$ be a non-empty Zariski open subset such that $W\to \overline{c((A/B)\times Z_o)}$ is {\'e}tale.
We fix $a_0\in A/B$ such that 
\begin{equation}\label{eqn:202103301}
(\{ a_0\}\times V)\cap W\not=\emptyset.
\end{equation}

\medskip

{\it Step 6.}\
We construct a Zariski closed subsets $F_o$ and $\Sigma_0$ of $(A/B)\times S$.
(After $U\hookrightarrow (A/B)\times S$ is constructed, these would be irreducible components of $\overline{U}\cap p^{-1}(Z_o)$ and $\overline{U}\cap p^{-1}(\supp S^P)$, where $p:(A/B)\times S\to S$ is the second projection.) 
We consider the restriction
$$
c|_{(A/B)\times Z_o}:(A/B)\times Z_o\to \overline{c((A/B)\times Z_o)}\subset \mathrm{Hilb}.
$$
We take an irreducible component $F_o$ of a Zariski closed set $(c|_{(A/B)\times Z_o})^{-1}(\overline{c(\{a_0\}\times V)})$ such that $\{a_0\}\times V\subset F_o$.
Here $\overline{c(\{a_0\}\times V)}\subset \mathrm{Hilb}$ is the Zariski closure in $\mathrm{Hilb}$.
We have a map 
\begin{equation*}
p|_{F_o}:F_o\to Z_o
\end{equation*}
by the composition of the closed immersion $F_o\hookrightarrow (A/B)\times Z_o$ and the second projection $p|_{(A/B)\times Z_o}:(A/B)\times Z_o\to Z_o$.
Then by Lemma \ref{lem:202206142} (2), the map $p|_{F_o}:F_o\to Z_o$ is dominant and generically finite.

Next we construct $\Sigma_0\subset (A/B)\times S$ so that $F_o\subset \Sigma_0$.
Let $F$ be an irreducible component of $F_o\times_{Z_o}Z$ such that the natural maps $F\to F_o$ and $p'|_{F}:F\to Z$ are dominant, where $p':(A/B)\times S_{1,A/B}\to S_{1,A/B}$ is the second projection.
We remark that $F_o\times_{Z_o}Z$ is a Zariski closed set of $(A/B)\times Z$.
Hence $F\subset (A/B)\times Z$ is a Zariski closed subset.
Let $\Theta\subset (A/B)\times S^P$ be a Zariski closed subset defined by $\Theta=c^{-1}(\overline{c(\{a_0\}\times V)})$.
Then $\Theta\subset (A/B)\times S$ is a Zariski closed subset.
Let $(a,[v])\in \cap_{k\geq 2}(p'|_{F})^{-1}(O_k)\subset F$.
We take an irreducible component $E$ of $\supp E_{(a,\tau ([v]))}$ as in Claim 2 of step 4.
By $(a,\tau([v]))\in F_o$, we have $c((a,\tau([v])))\in \overline{c(\{ a_0\}\times V)}$.
Hence $E\subset \Theta$.
Hence there exists an irreducible component $\Theta'$ of $\Theta$ such that $E\subset \Theta'$.
Then $\Theta'\subset (A/B)\times S$ is a Zariski closed subset.
By \eqref{eqn:20210330}, we have $(a,[v])\in \mathcal{D}\Theta'$.
Hence, denoting by $\Theta_1,\ldots,\Theta_l$ all irreducible components of $\Theta$,
we have
$$
\cap_k(p'|_{F})^{-1}(O_k)
\subset \mathcal{D}\Theta_1\cup\cdots\cup\mathcal{D}\Theta_l.
$$
Since $\mathcal{D}\Theta_i\subset (A/B)\times S$ are closed subschemes, $\supp (F\cap \mathcal{D}\Theta_i)$ are Zariski closed subsets of $F$.
Note that $(p'|_{F})^{-1}(O_k)\subset F$ is a non-empty Zariski open set for each $k$.
Hence we may choose $\Sigma_0$ from $\Theta_1,\ldots,\Theta_l$ such that 
\begin{equation}\label{eqn:202108191}
F\subset \mathcal{D}\Sigma_0.
\end{equation}
Since the natural map $F\to F_o$ is surjective, we have
$$F_o\subset \Sigma_0.$$
In particular, $\{a_0\}\times V\subset \Sigma_0$.
Hence we may apply Lemma \ref{lem:202206142} (1) for $(A/B)\times S^P\to \overline{c((A/B)\times S^P)}$ to get that $p|_{\Sigma_0}:\Sigma_0\to S$ is unramified at generic $(a_0,t)\in \{a_0\}\times V$, where $p:(A/B)\times S\to S$ is the second projection.

\medskip

{\it Step 7.}\
Now we apply Lemma \ref{lem:20210325} to take $\Sigma \subset (A/B)\times S$ from $\Sigma_0\subset (A/B)\times S$ such that $\Sigma_0\subset \Sigma$ and $\Sigma \to S$ is generically finite and {\'e}tale at generic $(a_0,t)\in \{ a_0\}\times V\subset \Sigma_0$.
We take a Zariski open $U\subset \Sigma$ such that $q:U\to S$ is {\'e}tale, where $q$ is the restriction of $p$ onto $U$.
We may take $U$ so that $(\{a_0\}\times V)\cap U\not=\emptyset$, hence $F_o\cap U\not=\emptyset$.
Hence $F_o\cap U\subset F_o$ is Zariski dense.
By \eqref{eqn:202108191}, this shows that $F\cap PTU\subset F$ is Zariski dense, where $PTU=(\mathcal{D}\Sigma)|_U$ in $(A/B)\times S_{1,A/B}$.
Hence $p'(F\cap PTU)\subset Z$ is Zariski dense in $Z$.
Hence $p'(PTU)\cap Z\subset Z$ is Zariski dense in $Z$.
This shows that the immersion $U\hookrightarrow (A/B)\times S$ satisfies the property of Definition \ref{defn:20201122}.
Hence $Z$ is horizontally integrable.
This completes the proof of the lemma.
\hspace{\fill} $\square$

\section{Verification of the normality condition: Existence of horizontally integrable $Z$}\label{sec:9}

The purpose of this section is to prove Proposition \ref{cor:20210306} below.
To state this proposition, we introduce several terminologies.
We recall $\Pi(\mathcal{F})$ from Definition \ref{defn:20220601}.

\begin{defn}
Let $X\subset A$ be a closed subvariety.
Let $\overline{A}$ be an equivariant compactification.
Let $\mathcal{F}\subset \mathrm{Hol}(\mathbb D, X)$ be an infinite set of holomorphic maps.
We define $\Lambda_{X,\overline{A}}(\mathcal{F})$ to be the set of all semi-abelian subvarieties $B\subset A$ such that \begin{enumerate}
\item
$\mathcal{F}\to \mathrm{Sp}_B\overline{X}$, where $\overline{X}\subset \overline{A}$ is the compactification, and
\item
$B\not\in \Pi(\mathcal{F})$.
\end{enumerate}
\end{defn}

Let $\mathcal{F}\subset \mathrm{Hol}(\mathbb D,A)$ be an infinite set of holomorphic maps.
For $B\in \Pi(\mathcal{F})$, denoting by $\varpi_B:A\to A/B$ the quotient map, we set $$\mathcal{F}_{A/B}=( \varpi_B\circ f)_{f\in \mathcal{F}}.$$ 
This is an infinite indexed family in $\mathrm{Hol}(\mathbb D,A/B)$.
Then $\mathcal{F}_{A/B}$ contains at most finitely many constant maps.
We remove these constant maps from $\mathcal{F}_{A/B}$ to get an infinite subfamily $\mathcal{F}_{A/B}'$ of $\mathcal{F}_{A/B}$.
We consider the following assumption for $\mathcal{F}$.

\begin{ass}\label{rem:20200630}
For every $B\in\Pi(\mathcal{F})$, the infinite indexed family $\mathcal{F}_{A/B}'$ in $\mathrm{Hol}(\mathbb D,A/B)$ satisfies Assumption \ref{ass:202103061}.
\end{ass}

As we shall see later (cf. Lemma \ref{lem:20201227}), every infinite subset $\mathcal{F}\subset \mathrm{Hol}(\mathbb D,A)$ of non-constant holomorphic maps contains infinite subset which satisfies this assumption. 
Let $\mathcal{F}\subset \mathrm{Hol}(\mathbb D,A)$ be an infinite set of holomorphic maps which satisfies Assumption \ref{rem:20200630}.
Let $B\in \Pi(\mathcal{F})$.
Then there exists a unique $T_{k,A/B}\subset P_{k,A/B}$ as in Assumption \ref{ass:202103061}.
We define $Z_{k,A/B}\subset P_{k,A/B}$ by
\begin{equation}\label{eqn:20220623}
Z_{k,A/B}=\cup_{l\geq 0}\mathrm{Im}(T_{k+l,A/B}\hookrightarrow P_{k+l,A/B}\to P_{k,A/B}).
\end{equation}
Then by Lemma \ref{lem:20220219}, we have $\mathcal{F}_{P_{k,A/B}}\Rightarrow Z_{k,A/B}$, where we set $\mathcal{F}_{P_{k,A/B}}=(\mathcal{F}_{A/B}')_{P_{k,A/B}}$ (cf. \eqref{eqn:20230326}).

\begin{pro}\label{cor:20210306}
Let $X\subset A$ be a closed subvariety and let $\overline{A}$ be a smooth equivariant compactification, where $\overline{A}$ is projective.
Let $\mathcal{F}\subset \mathrm{Hol}(\mathbb D,X)$ be an infinite set of holomorphic maps which satisfies Assumption \ref{rem:20200630}.
Assume that $\Lambda_{X,\overline{A}}(\mathcal{F})=\emptyset$.
Then there exists a semi-abelian subvariety $B\subset A$ such that $\mathrm{Sp}_B\overline{X}\not= \emptyset$ with the following properties:
\begin{enumerate}
\item $B\in \Pi(\mathcal{F})$.
\item There exists $k$ such that $Z_{k+1,A/B}\subset P_{k+1,A/B}$ is horizontally integrable, where 
$$Z_{k+1,A/B}=\cup_{l\geq 1}\mathrm{Im}(T_{k+l,A/B}\hookrightarrow P_{k+l,A/B}\to P_{k+1,A/B}).$$
\item $\mathcal{F}\to \mathrm{Sp}_B\overline{X}$.
\end{enumerate}
\end{pro}

We note that the assumption $\Lambda_{X,\overline{A}}(\mathcal{F})=\emptyset$ in this proposition reads $X\subsetneqq A$.
Indeed if $X=A$, then we have $A\in \Lambda_{X,\overline{A}}(\mathcal{F})$, hence $\Lambda_{X,\overline{A}}(\mathcal{F})\not=\emptyset$.

\subsection{Auxiliary Lemmas 1}

In this subsection, we prove several lemmas related to the following definition.

\begin{defn}\label{defn:20220418}
Let $S$ be a variety and let $X$ be a scheme of finite type over $\mathbb C$.
Let $\psi:X\to S$ be a morphism and let $Z\subset S$ be a Zariski closed set.
Let $\mathcal{U}$ be the set of all Zariski open set $U\subset S$ such that $Z\cap U\not= \emptyset$.
For $U\in \mathcal{U}$, we consider the scheme theoretic closure $\overline{\psi^{-1}(U)}\subset X$.
We set $X[Z]=\cap_{U\in \mathcal{U}}\overline{\psi^{-1}(U)}$, which is a closed subcheme of $X$.
\end{defn}

\begin{rem}\label{rem:20220418}
Assume that $Z$ is irreducible.
Then there exists $U\in \mathcal{U}$ such that $X[Z]=\overline{\psi^{-1}(U)}$.
Indeed, by the Noetherian property, there exist $U_1,\ldots ,U_k\in\mathcal{U}$ such that $X[Z]=\overline{\psi^{-1}(U_1)}\cap\cdots\cap \overline{\psi^{-1}(U_k)}$.
Set $U=U_1\cap\cdots\cap U_k$.
Then since $Z$ is irreducible, we have $U\in\mathcal{U}$.
Then we have $X[Z]\subset \overline{\psi^{-1}(U)}\subset \overline{\psi^{-1}(U_1)}\cap\cdots\cap \overline{\psi^{-1}(U_k)}$.
Hence $X[Z]=\overline{\psi^{-1}(U)}$.
\end{rem}

\begin{lem}\label{lem:20220420}
Let $X\to S$ and $\varphi:S\to T$ be morphisms, where $S$ and $T$ are varieties and $X$ is a scheme of finite type over $\mathbb C$.
Let $\psi:X\to T$ be the composite of $X\to S\to T$.
Let $V\subset S$ and $W\subset T$ be irreducible Zariski closed subsets such that $W\subset \overline{\varphi(V)}$.
Then $X[V]\subset X[W]$.
Moreover assume that there exists a Zariski open set $U\subset T$ such that $U\cap W\not=\emptyset$ and that $\psi^{-1}(U)$ is integral.
Then $X[V]=X[W]$.
\end{lem}

{\it Proof.}\
By Remark \ref{rem:20220418}, we may take a Zariski open $U_1\subset T$ such that $U_1\cap W\not=\emptyset$ and 
\begin{equation}\label{eqn:20220618}
X[W]=\overline{\psi^{-1}(U_1)}.
\end{equation}
We have $\varphi^{-1}(U_1)\cap V\not=\emptyset$.
Hence by the definition of $X[V]$, we have
\begin{equation*}
X[V]\subset \overline{\psi^{-1}(U_1)}
= X[W].
\end{equation*}

Now we assume moreover that $\psi^{-1}(U)$ is integral for some Zariski open $U\subset T$ such that $U\cap W\not=\emptyset$.
We take $U_1\subset T$ as above.
We may assume $U_1\subset U$.
By $\varphi^{-1}(U_1)\cap V\not=\emptyset$, we may take $U_2\subset \varphi^{-1}(U_1)$ such that $\overline{p^{-1}(U_2)}=X[V]$, where $p:X\to S$ (cf. Remark \ref{rem:20220418}).
Since $\psi^{-1}(U_1)$ is integral, the inclusion $p^{-1}(U_2)\subset \psi^{-1}(U_1)$ is schematic dense.
Hence $\psi^{-1}(U_1)\subset X[V]$.
Thus by \eqref{eqn:20220618}, we get $X[W]\subset X[V]$.
This concludes the proof of the lemma.
\hspace{\fill} $\square$

\begin{lem}\label{lem:20220618}
Let $T$ be a variety and let $W\subset T$ be an irreducible Zariski closed subset.
Let $\psi:X\to T$ and $\psi':X'\to T$ be morphisms from schemes of finite type over $\mathbb C$.
Let $p:X\hookrightarrow X'$ be a closed immersion over $T$.
Assume that there exists a Zariski open set $U\subset T$ such that $U\cap W\not=\emptyset$ and the induced map $\psi^{-1}(U)\to(\psi')^{-1}(U)$ is an isomorphism.
Then $X[W]=X'[W]$.
\end{lem}

{\it Proof.}\
The immersion $p$ induces a closed immersion $X[W]\hookrightarrow X'[W]$.
Hence we prove the converse.
By Remark \ref{rem:20220418}, we may take a Zariski open $U\subset T$ such that  $U\cap W\not=\emptyset$ and $\overline{\psi^{-1}(U)}=X[W]$.
By replacing $U$ by a smaller Zariski open set, we may assume moreover that the induced map $\psi^{-1}(U)\to(\psi')^{-1}(U)$ is an isomorphism.
Hence the scheme theoretic closure $\overline{(\psi')^{-1}(U)}\subset X'$ factors $X[W]$, where $X[W]$ is a closed subscheme of $X'$ by $p:X\hookrightarrow X'$.
This shows $X'[W]\subset X[W]$.
\hspace{\fill} $\square$
\medskip

We recall Definition \ref{defn:20211007}.
Let $\varphi:S\to T$ be a morphism of varieties.
Let $V\subset S$ and $W\subset T$ be irreducible Zariski closed subsets such that $W\subset \overline{\varphi(V)}$.
Let $\psi:T'\to T$ be a $W$-admissible modification and set $S'=(S\times_TT')[W']$, where $W'\subset T'$ is the minimal transform.
Then the induced map $S'\to S$ is a $V$-admissible modification.
We denote by $V'\subset S'$ the minimal transform of $V\subset S$.
In this situation, we have the following lemma.

\begin{lem}\label{lem:20210909}
Let $X\to S$ and $\varphi:S\to T$ be morphisms, where $S$ and $T$ are varieties and $X$ is a scheme of finite type over $\mathbb C$.
Let $V\subset S$ and $W\subset T$ be irreducible Zariski closed subsets such that $W\subset \overline{\varphi(V)}$.
Let $T'\to T$ be a $W$-admissible modification and set $S'=(S\times_TT')[W']$.
Then $(X\times_{S}S')[V']\subset (X\times_TT')[W']$.
Assume moreover that $X$ is integral.
Then $(X\times_{S}S')[V']=(X\times_TT')[W']$.
\end{lem}

{\it Proof.}\
We first prove
\begin{equation}\label{eqn:20220418}
(X\times_{S}S')[W']= (X\times_TT')[W'].
\end{equation}
By the closed immersion $S'\hookrightarrow S\times_TT'$, we get a closed immersion $X\times_SS'\hookrightarrow X\times_TT'$.
Note that the following is the fiber product:
\begin{equation}\label{eqn:202206182}
\begin{CD}
X\times_SS' @>>> X\times_TT'@=X\times_S(S\times_TT')\\
@VVV @VVV \\
S'@>>> S\times_TT'
\end{CD}
\end{equation}
By Remark \ref{rem:20220418}, we may take a Zariski open $U_1\subset T'$ such that $U_1\cap W'\not=\emptyset$ and a Zariski open set $S\times_{T}U_1\subset S\times_TT'$ is schematic dense in $S'\subset S\times_TT'$.
In particular, we have an open immersion 
\begin{equation}\label{eqn:202206183}
S\times_{T}U_1\subset S'.
\end{equation}
Hence the map $S'\to S\times_TT'$ over $T'$ is an isomorphism over $U_1\subset T'$.
Hence the closed immersion $X\times_SS'\hookrightarrow X\times_TT'$ is an isomorphism over $U_1\subset T'$.
Hence by Lemma \ref{lem:20220618}, we get \eqref{eqn:20220418}.

Set $X'=X\times_{S}S'$.
Let $\psi:X'\to T'$ be the natural map.
Since \eqref{eqn:202206182} is a fiber product, the open immersion \eqref{eqn:202206183} yields 
\begin{equation}\label{eqn:20220421}
\psi^{-1}(U_1)=X\times_TU_1.
\end{equation}

Now let $\varphi':S'\to T'$ be the induced map.
We have $W'\subset \overline{\varphi'(V')}$.
We apply Lemma \ref{lem:20220420} to get $X'[V']\subset X'[W']$.
Hence combining with \eqref{eqn:20220418}, we get 
\begin{equation*}
(X\times_{S}S')[V']\subset (X\times_TT')[W'].
\end{equation*}

We assume moreover that $X$ is integral.
We take a Zariski open set $U_2\subset T'$ such that $T'\to T$ is an isomorphism over $U_2$ and $U_2\cap W'\not=\emptyset$.
We may assume $U_2\subset U_1$.
Then by \eqref{eqn:20220421}, we may consider $\psi^{-1}(U_2)$ as a Zariski open set of $X$.
Hence $\psi^{-1}(U_2)$ is integral.
Hence by Lemma \ref{lem:20220420}, we get $X'[V']=X'[W']$.
Hence by \eqref{eqn:20220418}, we conclude the proof of our lemma.
\hspace{\fill} $\square$

\begin{lem}\label{lem:20210903}
Let $Z\subset S$ be an irreducible Zariski closed set.
Let $S'\to S$ be a $Z$-admissible modification and let $S''\to S'$ be a $Z'$-admissible modification, where $Z'\subset S'$ is the minimal transform.
Let $X\to S$ be a morphism.
Then we have
$$
(X\times_SS'')[Z'']= ((X\times_SS')[Z']\times_{S'}S'')[Z''],
$$
where $Z''\subset S''$ is the minimal transform.
\end{lem}

{\it Proof.}\
The closed immersion $(X\times_SS')[Z']\hookrightarrow X\times_SS'$ induces a closed immersion 
$$
((X\times_SS')[Z']\times_{S'}S'')[Z'']\hookrightarrow (X\times_SS'')[Z''].
$$
We shall show that this is an isomorphism.
By Remark \ref{rem:20220418}, we may take a Zariski open set $U_1\subset S'$ such that $Z'\cap U_1\not=\emptyset$ and 
$$(X\times_SS')[Z']=\overline{X\times_SU_1}.$$
We denote by $\varphi:S''\to S'$ the natural map.
By Remark \ref{rem:20220418}, we may take a Zariski open set $U_2\subset \varphi^{-1}(U_1)$ such that $U_2\cap Z''\not=\emptyset$ and 
$$((X\times_SS')[Z']\times_{S'}S'')[Z'']=\overline{(X\times_SS')[Z']\times_{S'}U_2}.
$$
Hence, we have 
$$X\times_SU_2=(X\times_SU_1)\times_{S'}U_2 \subset
(X\times_SS')[Z']\times_{S'}U_2
\subset
 ((X\times_SS')[Z']\times_{S'}S'')[Z''].$$
Then by the definition of $(X\times_SS'')[Z'']$, we have $(X\times_SS'')[Z'']\subset \overline{X\times_SU_2}$.
Hence we get 
$$(X\times_SS'')[Z'']\subset ((X\times_SS')[Z']\times_{S'}S'')[Z''].$$
This completes the proof of our lemma.
\hspace{\fill} $\square$

\medskip

In the following two lemmas, we consider a closed subscheme $X\subset \mathbb P^N\times S$.

\begin{lem}\label{lem:20210905}
Let $S'\to S$ be a $Z$-admissible modification of varieties, where $Z\subset S$ is an irreducible Zariski closed set.
Let $X\to S$ be a projective morphism such that $X|_Z\to Z$ is flat, where $X|_Z=X\times_SZ$.
Let $Z'\subset S'$ be the minimal transform.
Set $X'=(X\times_SS')[Z']$ and $X'|_{Z'}=X'\times_{S'}Z'$.
Then $X'|_{Z'}=X|_Z\times_ZZ'$.
In particular, $X'|_{Z'}\to Z'$ is flat.
\end{lem}

{\it Proof.}\
The closed immersion $X'\hookrightarrow X\times_SS'$ induces a closed immersion 
\begin{equation}\label{eqn:20220420}
X'|_{Z'}\hookrightarrow (X\times_SS')|_{Z'}.
\end{equation}
By Remark \ref{rem:20220418}, we may take a Zariski open $U\subset S'$ such that $U\cap Z'\not=\emptyset$ and $\overline{X\times_SU}=X'$.
In particular, we have an open immersion
\begin{equation}\label{eqn:202204201}
(X\times_SS')|_{Z'\cap U}\subset X'|_{Z'}.
\end{equation}
Since the composite of $Z'\hookrightarrow S'\to S$ factors throw $Z\hookrightarrow S$, we have
\begin{equation}\label{eqn:202204202}
(X\times_SS')|_{Z'}=X|_Z\times_ZZ'.
\end{equation}
In particular, the morphism $(X\times_SS')|_{Z'}\to Z'$ is flat.
Hence by Lemma \ref{lem:20210417}, the inclusion $(X\times_SS')|_{Z'\cap U}\subset (X\times_SS')|_{Z'}$ is scheme theoretic dense.
Hence by \eqref{eqn:202204201}, we get $(X\times_SS')|_{Z'}\subset X'|_{Z'}$.
Thus by \eqref{eqn:20220420}, we get $X'|_{Z'}= (X\times_SS')|_{Z'}$.
By \eqref{eqn:202204202}, we get $X'|_{Z'}=X|_Z\times_ZZ'$.
\hspace{\fill} $\square$

\begin{lem}\label{lem:20210831}
Let $X\subset \mathbb P^N\times S$ be a closed subscheme where $S$ is a variety.
Let $Z\subset S$ be an irreducible Zariski closed set.
Then there exists a $Z$-admissible modification $S'\to S$ such that $(X\times_SS')[Z']|_{Z'}\to Z'$ is flat, where $Z'\subset S'$ is the minimal transform and $(X\times_SS')[Z']|_{Z'}=(X\times_SS')[Z']\times_{S'}Z'$.
\end{lem}

{\it Proof.}\
We apply Lemma \ref{lem:202110273} to get a $Z$-admissible modification $S'\to S$ and a Zariski open set $U\subset S'$ such that $U\cap Z'\not=\emptyset$ and that $X'|_{Z'}\to Z'$ is flat, where $X'\subset X\times_SS'$ is the scheme theoretic closure of $X\times_SU$.
Then we have $(X\times_SS')[Z']\subset X'$.
Hence $(X\times_SS')[Z'][Z']\subset X'[Z']$.
By Lemma \ref{lem:20210903}, applied to $S''=S'$, we get $(X\times_SS')[Z'][Z']=(X\times_SS')[Z']$.
Hence $(X\times_SS')[Z']\subset X'[Z']$.
On the other hand, by $X'\subset X\times_SS'$, we get $X'[Z']\subset (X\times_SS')[Z']$.
Thus $X'[Z']=(X\times_SS')[Z']$.
We apply Lemma \ref{lem:20210905} to $X'\to S'$ and $\mathrm{id}_{S'}:S'\to S'$.
The conclusion is $X'[Z']|_{Z'}=X'|_{Z'}$.
Hence $(X\times_SS')[Z']|_{Z'}=X'|_{Z'}$.
Thus $(X\times_SS')[Z']|_{Z'}\to Z'$ is flat.
\hspace{\fill} $\square$

\subsection{Auxiliary Lemmas 2}

\begin{lem}\label{lem:20210414}
Let $X\subset \overline{A}\times S$ be a closed subscheme, where $S$ is integral and $\overline{A}$ is a projective, equivariant compactification.
Assume that the induced map $X\to S$ is surjective.
Then there exist a non-empty Zariski open set $U\subset S$ and a semi-abelian subvariety $C\subset A$ such that $\mathrm{Stab^0}(X_y)=C$ for all $y\in U$.
\end{lem}

{\it Proof.}
Since $S$ is integral, by replacing $S$ by its non-empty Zariski open set, we may assume that $X\to S$ is flat.
For each $B\subset A$, we set $V_B=\{y\in S;B\subset \mathrm{Stab^0}(X_y)\}$.
Then $V_B\subset S$ is a Zariski closed set.
Indeed let $P$ be the Hilbert polynomial of $X_y$ for some (hence for all) $y\in S$.
Set $Y=\bigcap_{b\in B}(X+b)\subset X$, where the intersection is taken scheme theoretically. 
Then $y\in V_B$ if and only if the Hilbert polynomial of $Y_y$ is equal to $P$.
Hence $V_B$ is a Zariski closed set (cf. \cite[Lemma 3.1]{Y}).

Now let $\mathcal{X}\subset \overline{A}\times A\times S$ be the pull back of $X$ by the action $m:\overline{A}\times A\times S\to \overline{A}\times S$, where $(x,a,s)\mapsto(x+a,s)$.
Then $\mathcal{X}\to A\times S$ is flat.
Let $c:A\times S\to \mathrm{Hilb}$ be the classification map.
Let $\varphi:A\times S\to \mathrm{Hilb}\times S$ be the induced map such that $\varphi(a,s)=(c(a,s),s)$.
Then for each $y\in S$, we have $\supp (\varphi^{-1}(\varphi (0,y)))\subset A\times\{y\}=A$.
We have
$$
\supp (\varphi^{-1}(\varphi (0,y)))=\mathrm{Stab}(X_y).
$$
For each integer $d\geq 0$, let $E_d\subset S$ be the set of $y\in S$ such that $\dim_{(0,y)}\varphi^{-1}(\varphi (0,y))\geq d$.
Then $E_d\subset S$ is a Zariski closed set.
We take $d$ such that $E_d=S$ and $E_{d+1}\subsetneqq S$.
Then for each $y\in S-E_{d+1}$, we have $\dim \mathrm{Stab^0}(X_y)=d$.
Let $B_1,B_2,\ldots$ be the set of $d$-dimensional semi-abelian subvarieties of $A$.
Then $S-E_{d+1}\subset\cup V_{B_i}$.
Hence there exists $B_i$ such that $S=V_{B_i}$.
For $y\in S-E_{d+1}$, we have $B_i\subset \mathrm{Stab^0}(X_y)$.
Hence $\mathrm{Stab^0}(X_y)=B_i$.
We set $C=B_i$ and $U=S-E_{d+1}$ to conclude the proof.
\hspace{\fill} $\square$

\medskip

In the following lemma, we recall the definition $X_{\{0\},k}\subset \overline{A}\times S_{k,A}$ from \eqref{eqn:202109111}.
We have the isomorphism $(A\times S)_k\simeq A\times S_{k,A}$ as in \eqref{eqn:20211203}.

\begin{lem}\label{lem:202109031}
Let $X\subset \overline{A}\times S$ be a closed subscheme, where $S$ is a smooth algebraic variety and $\overline{A}$ is a smooth equivariant compactification.
Let $f\in \mathrm{Hol}(\mathbb D,A\times S)$ satisfies $f(\mathbb D)\subset \supp X$, where $f$ is non-constant.
Then $f_{[k]}(\mathbb D)\subset \supp X_{\{0\},k}$.
\end{lem}

{\it Proof.}\
We define $g:\mathbb D\times \mathbb D\to A\times A\times S$ by 
$$g(z,w)=(f_A(z),f_A(w)-f_A(z),f_S(w)).$$
Let $m:\overline{A}\times A\times S\to \overline{A}\times S$ be defined by $(x,a,s)\mapsto(x+a,s)$.
Then $m\circ g(z,w)=(f_A(w),f_S(w))\in \supp X$.
Hence $g(\mathbb D\times \mathbb D)\subset \supp \mathcal{X}$, where $\mathcal{X}\subset A\times A\times S$ is a closed subscheme obtained by the pull-back of $X\subset \overline{A}\times S$ by $m:\overline{A}\times A\times S\to \overline{A}\times S$.
By taking $k$-th derivative for $w$, we get $\partial_w^kg:\mathbb D\times \mathbb D\to A\times (A\times S)_k$.

\medskip

{\it Claim.}\
$\partial_w^kg(\mathbb D\times \mathbb D)\subset \supp \mathcal{P}_k\mathcal{X}$.

{\it Proof.}\
We prove this by the induction on $k$.
The case $k=0$ is obvious.
So we assume the case $k$ and prove the case for $k+1$.
By $\partial_w^kg:\mathbb D\times \mathbb D\to \mathcal{P}_k\mathcal{X}\subset A\times (A\times S)_k$, we get $\partial_w(\partial_w^kg):\mathbb D\times \mathbb D\to \mathcal{D}\mathcal{P}_k\mathcal{X}\subset (A\times (A\times S)_k)_1$.
Since $\partial_w(\partial_w^kg):\mathbb D\times \mathbb D\to A\times (A\times S)_{k+1}\subset (A\times (A\times S)_k)_1$, we get 
$$\partial_w(\partial_w^kg)(\mathbb D\times \mathbb D)\subset \mathcal{D}\mathcal{P}_k\mathcal{X}\cap (A\times (A\times S)_{k+1})=\mathcal{P}_{k+1}\mathcal{X}.
$$
This completes the induction step.
\hspace{\fill} $\square$

\medskip

By \eqref{eqn:202112035}, we have $\partial_w^kg=(f_A(z),f_A(w)-f_A(z),f_{S_{k,A}}(w))$.
Restricting this to the diagonal $\mathbb D\subset \mathbb D\times \mathbb D$, we get 
$$\partial_w^kg\circ \Delta(z)=(f_A(z),0_A,f_{S_{k,A}}(z))
\in A\times (A\times S)_k\vert_{A\times \{ 0\}\times S_{k,A}}.$$
Hence by $X_{\{0\},k}=\mathcal{P}_k\mathcal{X}\cap (\overline{A}\times \{ 0\}\times S_{k,A})$, the claim above implies $\partial_w^kg\circ \Delta(\mathbb D)\subset \supp X_{\{0\},k}$.
Now under the isomorphism 
$$\psi_k:A\times (A\times S)_k\vert_{A\times \{ 0\}\times S_{k,A}}\to (A\times S)_k,$$ 
we have $\psi_k\circ\partial_w^kg\circ\Delta=f_{[k]}$.
Hence, we have $f_{[k]}(\mathbb D)\subset \supp X_{\{0\},k}$.
\hspace{\fill} $\square$

\begin{lem}\label{lem:202206233}
Let $X\subset \overline{A}$ be a closed subvariety, where $\overline{A}$ is a smooth equivariant compactification.
Let $\mathcal{F}\subset \mathrm{Hol}(\mathbb D,A\cap X)$ be an infinite set of non-constant holomorphic maps.
Let $Z\subset P_{k,A}$ be a Zariski closed subset such that $\mathcal{F}_{P_{k,A}}\Rightarrow Z$.
Let $P_{k,A}'\to P_{k,A}$ be a $Z$-admissible modification and let $Z'\subset P_{k,A}'$ be the minimal transform.
Set $X_{\{0\},k}'=(X_{\{0\},k}\times_{P_{k,A}}P_{k,A}')[Z']$.
Then we have $f'_{[k]}(\mathbb D)\subset X_{\{0\},k}'$ for all $f\in \mathcal{F}$ with finite exception, where $f_{[k]}':\mathbb D\to A\times P_{k,A}'$ is the lift of $f_{[k]}:\mathbb D\to A\times P_{k,A}$.
\end{lem}

{\it Proof.}\
By Lemma \ref{lem:202105272}, $Z$ is irreducible.
By Remark \ref{rem:20220418}, we may take a Zariski open set $U\subset P_{k,A}'$ such that $U\cap Z'\not=\emptyset$ and $X_k'=\overline{\psi^{-1}(U)}$, where $\psi:X_k\times_{P_{k,A}}P_{k,A}'\to P_{k,A}'$ is the projection.
We may assume moreover that the map $P_{k,A}'\to P_{k,A}$ is isomorphic on $U$.
By Lemma \ref{lem:20220623}, we have $\mathcal{F}_{P_{k,A}'}\Rightarrow Z'$.
Hence we have $f_{P_{k,A}'}(\mathbb D)\cap U\not=\emptyset$ for all $f\in\mathcal{F}$ with finite exception.
We set $\Omega_f=f_{P_{k,A}'}^{-1}(U)$.
Then for all $f\in \mathcal{F}$ with finite exception, $\mathbb D-\Omega_f$ is discrete and $f_{P_{k,A}'}(\Omega_f)\subset U$.
By Lemma \ref{lem:202109031}, we have $f_{[k]}(\mathbb D)\subset X_{\{0\},k}$.
This shows $f_{[k]}'(\Omega_f)\subset \psi^{-1}(U)$, hence $f_{[k]}'(\mathbb D)\subset X_{\{0\},k}'$ for all $f\in \mathcal{F}$ with finite exception.
\hspace{\fill} $\square$

\subsection{Demailly jet spaces and quotient maps}

We recall Definition \ref{defn:20201225}.
Let $B\subset A$ be a semi-abelian subvariety.
For each $k\geq 1$, the quotient $A\to A/B$ canonically induces the map
\begin{equation}\label{eqn:202206221}
\mu_k:P_{k,A}\backslash E_{k,A,A/B}  \to P_{k,A/B}
\end{equation}
inductively as follows.
When $k=1$, we have $P_{1,A}=P(\mathrm{Lie}(A))$, $E_{1,A,A/B}=P(\mathrm{Lie}(B))$ and $P_{1,A/B}=P(\mathrm{Lie}(A/B))$.
The map $\mu_1$ is defined by the projectivization of the quotient map $\mathrm{Lie}(A)\to\mathrm{Lie}(A/B)$.
Suppose we have $\mu_{k-1}$.
Then by $(\mu_{k-1})_{*}:T(P_{k-1,A}\backslash E_{k-1,A,A/B})\to TP_{k-1,A/B}$ and the quotient map $\mathrm{Lie}(A)\to\mathrm{Lie}(A/B)$, we get
\begin{equation}\label{eqn:202206211}
T(P_{k-1,A}\backslash E_{k-1,A,A/B} )\times\mathrm{Lie}(A)\to TP_{k-1,A/B}\times\mathrm{Lie}(A/B).
\end{equation}
The kernel of this map is contained in 
$$T(P_{k-1,A}\backslash E_{k-1,A,A/B} )\times\mathrm{Lie}(B)\subset T(P_{k-1,A}\backslash E_{k-1,A,A/B} )\times\mathrm{Lie}(A).$$
Let $\tau:P_{k,A}\to P_{k-1,A}$ be the projection.
Then the restriction of the projectivization of \eqref{eqn:202206211} onto $P_{k,A}$ yields the map
$$
(P_{k,A}\backslash \tau^{-1}(E_{k-1,A,A/B}))\backslash E_{k,A,A/B}\to P_{k,A/B}.
$$
Hence by Lemma \ref{lem:202206211}, we get $\mu_{k}:P_{k,A}\backslash E_{k,A,A/B} \to P_{k,A/B}$, which is the map \eqref{eqn:202206221} for $k$.

Now let $f:\mathbb D\to A$ be a holomorphic map.
Let $f_{A/B}:\mathbb D\to A/B$ be the composition of $f$ and the quotient map $A\to A/B$.
Assume that $f_{A/B}$ is non-constant.
Then we get $f_{P_{k,A}}:\mathbb D\to P_{k,A}$ and $(f_{A/B})_{P_{k,A/B}}:\mathbb D\to P_{k,A/B}$, where $f_{P_{k,A}}(\mathbb D)\not\subset E_{k,A,A/B}$.
Then by the construction above, we have
$$
(f_{A/B})_{P_{k,A/B}}=\mu_k\circ f_{P_{k,A}}.
$$

\begin{lem}\label{lem:202206212}
Let $C\subset A$ be a semi-abelian subvariety such that $B\subset C\subset A$, then $\mu_k^{-1}(E_{k,A/B,A/C})\subset E_{k,A,A/C}\cap (P_{k,A}\backslash E_{k,A,A/B})$.
\end{lem}

{\it Proof.}\
Let $(x,[v])\in P_{k,A}\backslash E_{k,A,A/C}$, where $x\in P_{k-1,A}$ and $v\in V_{k-1}^{\dagger}\backslash\{0\}\subset TP_{k-1,A}\times\mathrm{Lie}(A)$ (cf. \eqref{eqn:202112033}).
Then $v\not\in TP_{k-1,A}\times\mathrm{Lie}(C)$.
Then the image of $v$ under the map \eqref{eqn:202206211} is not contained in $TP_{k-1,A/B}\times\mathrm{Lie}(C/B)$.
Hence $\mu_k((x,[v]))\not\in E_{k,A/B,A/C}$.
This shows $\mu_k^{-1}(E_{k,A/B,A/C})\subset E_{k,A,A/C}$ in $P_{k,A}\backslash E_{k,A,A/B}$.
\hspace{\fill} $\square$

\medskip

We recall $Z_{k,A/B}\subset P_{k,A/B}$ from \eqref{eqn:20220623}.

\begin{lem}\label{lem:202206234}
Let $\mathcal{F}\subset \mathrm{Hol}(\mathbb D,A)$ be an infinite set of holomorphic maps which satisfies Assumption \ref{rem:20200630}.
Let $B,C\in\Pi(\mathcal{F})$ such that $B\subset C$.
Let $P_{k,A/B}'\to P_{k,A/B}$ be a $Z_{k,A/B}$-admissible modification such that the rational map $P_{k,A/B}\dashrightarrow P_{k,A/C}$ induces a morphism $\mu_k:P_{k,A/B}'\to P_{k,A/C}$.
Then $Z_{k,A/C}\subset \mu_k(Z_{k,A/B}')$, where $Z_{k,A/B}'\subset P_{k,A/B}'$ is the minimal transform. 
\end{lem}

{\it Proof.}\
We take $l\in \mathbb Z_{\geq 0}$ such that $T_{k+l,A/C}\to Z_{k,A/C}$ and $T_{k+l,A/B}\to Z_{k,A/B}$ are surjective maps.
The rational map $P_{k+l,A/B}\dashrightarrow P_{k+l,A/C}$ is holomorphic outside $E_{k+l,A/B,A/C}\subset P_{k+l,A/B}$.
By $C/B\in \Pi(\mathcal{F}_{A/B})$, we have $T_{k+l,A/B}\not\subset E_{k+l,A/B,A/C}$.
Hence there exists a $T_{k+l,A/B}$-admissible modification $P_{k+l,A/B}'\to P_{k+l,A/B}$ such that the rational map $P_{k+l,A/B}\dashrightarrow P_{k+l,A/C}$ induces a regular map 
$$\mu_{k+l}:P_{k+l,A/B}'\to P_{k+l,A/C}.
$$
By $\mathcal{F}_{P_{k+l,A/B}}\Rightarrow T_{k+l,A/B}$, we have $\mathcal{F}_{P_{k+l,A/B}'}\Rightarrow T_{k+l,A/B}'$ (cf. Lemma \ref{lem:20220623}).
Hence by Lemma \ref{lem:20220219}, we get
$$\mathcal{F}_{P_{k+l,A/C}}\Rightarrow \mu_{k+l}(T_{k+l,A/B}').$$
By Lemma \ref{lem:202105271}, we have either $T_{k+l,A/C}\subset  \mu_{k+l}(T_{k+l,A/B}')$ or $\mu_{k+l}(T_{k+l,A/B}')\subsetneqq T_{k+l,A/C}$.
To show $T_{k+l,A/C}\subset  \mu_{k+l}(T_{k+l,A/B}')$, we assume contrary $\mu_{k+l}(T_{k+l,A/B}')\subsetneqq T_{k+l,A/C}$.
By Lemma \ref{lem:202105272}, $\mu_{k+l}(T_{k+l,A/B}')$ is irreducible.
Hence by the definition of $T_{k+l,A/C}$, there exists $C'/C\in\Pi(\mathcal{F}_{A/C})$ such that $\mu_{k+l}(T_{k+l,A/B}')\subset E_{k+l,A/C,A/C'}$.
Hence $T_{k+l,A/B}'\subset \mu_{k+l}^{-1}(E_{k+l,A/C,A/C'})$.
Since $T_{k+l,A/B}$ is irreducible, this implies $T_{k+l,A/B}\subset E_{k+l,A/B,A/C'}$ (cf. Lemma \ref{lem:202206212}).
By $C'/B\in \Pi(\mathcal{F}_{A/B})$, this contradicts to the definition of $T_{k+l,A/B}$.
Hence $T_{k+l,A/C}\subset  \mu_{k+l}(T_{k+l,A/B}')$.

Now $Z_{k,A/C}\subset P_{k,A/C}$ is contained in the image of $T_{k+l,A/B}'\subset P_{k+l,A/B}'$ under the composition of $\mu_{k+l}:P_{k+l,A/B}'\to P_{k+l,A/C}$ and $P_{k+l,A/C}\to P_{k,A/C}$.
Thus we get $Z_{k,A/C}\subset \mu_k(Z_{k,A/B}')$.
\hspace{\fill} $\square$

\medskip

Let $B$ and $C$ be semi-abelian subvarieties of $A$ such that $B\subset C$.
Let $S\subset \overline{A}$ be a Zariski closed set which is $C$-invariant.
Then we get $S_{B,k}\subset \overline{A}\times P_{k,A/B}$ and $S_{C,k}\subset \overline{A}\times P_{k,A/C}$.

\begin{lem}\label{lem:202206222}
Let $\varphi:\overline{A}\times (P_{k,A/B}\backslash E_{k,A/B,A/C})\to \overline{A}\times P_{k,A/C}$ be the map induced from the regular map $P_{k,A/B}\backslash E_{k,A/B,A/C}\to P_{k,A/C}$.
Then we have
$$\varphi^{*}(S_{C,k})=S_{B,k}|_{P_{k,A/B}\backslash E_{k,A/B,A/C}}.$$
\end{lem}

{\it Proof.}\
Let $\mathcal{S}\subset\overline{A}\times A$ be the pull-back of $S\subset \overline{A}$ by the action $m:\overline{A}\times A\to \overline{A}$ so that $\mathcal{S}_a=S-a$, where $\mathcal{S}_a\subset \overline{A}$ is the fiber of $\mathcal{S}\to A$ over $a\in A$.
Let $\mathcal{S}_B\subset \overline{A}\times(A/B)$ be defined so that the pull-back of $\mathcal{S}_B$ by $\overline{A}\times A\to \overline{A}\times (A/B)$ is equal to $\mathcal{S}$.
We define $\mathcal{S}_C\subset \overline{A}\times (A/C)$ similarly.
Then the pull-back of $\mathcal{S}_C$ by $\overline{A}\times (A/B)\to \overline{A}\times (A/C)$ is equal to $\mathcal{S}_B$.
Let $\phi:\overline{A}\times ((A/B)\times P_{k,A/B})\dashrightarrow \overline{A}\times ((A/C)\times P_{k,A/C})$ be the induced rational map, which is regular on $\overline{A}\times ((A/B)\times (P_{k,A/B}\backslash E_{k,A/B,A/C}))$.
Then we have 
$$
\phi^*(\mathcal{P}_k\mathcal{S}_C)=(\mathcal{P}_k\mathcal{S}_B)\vert_{(A/B)\times (P_{k,A/B}\backslash E_{k,A/B,A/C})}.
$$
Hence by \eqref{eqn:202109111}, we get $\varphi^{*}(S_{C,k})=S_{B,k}|_{P_{k,A/B}\backslash E_{k,A/B,A/C}}$.
\hspace{\fill} $\square$

\subsection{Main lemma for the proof of Proposition \ref{cor:20210306}}

Let $A$ be a semi-abelian variety and let $B\subset A$ be a semi-abelian subvariety.
Given a Zariski closed set $V\subset\overline{A}$, we set $Y=\mathrm{Sp}_BV$.
Then the Zariski closed set $Y\subset \overline{A}$ is $B$-invariant.
Hence we get the closed subscheme $Y_{B,k}\subset \overline{A}\times P_{k,A/B}$.
For each $y\in P_{k,A/B}$, the fiber of $Y_{B,k}\to P_{k,A/B}$ over $y$ is denoted by $(Y_{B,k})_y$, which is a closed subscheme of $\overline{A}$.

Let $X\subset A$ be a closed subvariety and let $\overline{A}$ be a smooth equivariant compactification.
Let $\overline{X}\subset \overline{A}$ be the Zariski closure.
In this subsection, we write $\overline{X}_k=(\overline{X})_{\{ 0\},k}$ for short.
Then $\overline{X}_k\subset \overline{A}\times P_{k,A}$.
Although this $\overline{X}_k$ is not the same as the Demailly jet space of $\overline{X}$ discussed in Section \ref{sec:dem}, no confusion will occur.

We recall $Z_{k,A/B}\subset P_{k,A/B}$ from \eqref{eqn:20220623}.
We set $Z_k=Z_{k,A}\subset P_{k,A}$ for short.

\begin{lem}\label{pro:410}
Let $A$ be a non-trivial semi-abelian variety and let $X\subset A$ be a closed subvariety.
Let $\overline{A}$ be a smooth equivariant compactification, which is projective.
Let $\mathcal{F}\subset \mathrm{Hol}(\mathbb D,X)$ be an infinite set of holomorphic maps which satisfies Assumption \ref{rem:20200630}.
Assume that $\Lambda_{X,\overline{A}}(\mathcal{F})=\emptyset$.
Then there exists $B\subset A$ such that $\mathrm{Sp}_B\overline{X}\not= \emptyset$ with the following properties:
\begin{enumerate}
\item $B\in \Pi(\mathcal{F})$, in particular $B\subsetneqq A$.
\item $Z_{k,A/B}\subset p_{k,B}(Y_{B,k})$ for sufficiently large $k$, where $Y=\mathrm{Sp}_B\overline{X}$ and $p_{k,B}:\overline{A}\times P_{k,A/B}\to P_{k,A/B}$ is the second projection.
\item $\mathrm{Stab^0}((Y_{B,k})_y)=B$ for generic $y\in Z_{k,A/B}$ and sufficiently large $k$.
\item
For sufficiently large $k\geq 0$, there exist a $Z_k$-admissible modification $\hat{P}_{k,A}\to P_{k,A}$ with a regular map $\sigma :\hat{P}_{k,A}\to P_{k,A/B}$ and a Zariski closed set $V_k\subset \hat{P}_{k,A}$ such that:
\begin{enumerate}
\item
$Z_{k,A/B}\subset \sigma(V_k)$ and $V_k\subset \hat{Z}_k$, where $\hat{Z}_k\subset \hat{P}_{k,A}$ is the minimal transform.
\item
Set $\hat{X}_k=(\overline{X}_k\times_{P_{k,A}}\hat{P}_{k,A})[\hat{Z}_k]\subset \overline{A}\times \hat{P}_{k,A}$.
Then $\hat{X}_k|_{\hat{Z}_k}\to \hat{Z}_{k}$ is flat.
\item 
Let $P_{k,A/B}'\to P_{k,A/B}$ be a $Z_{k,A/B}$-admissible modification with the minimal transform $Z_{k,A/B}'\subset P_{k,A/B}'$, let
$\hat{P}_{k,A}'=(\hat{P}_{k,A}\times_{P_{k,A/B}}P_{k,A/B}')[Z_{k,A/B}']$ and $V'_k=(V_k\times_{P_{k,A/B}}P_{k,A/B}')[Z_{k,A/B}']
\subset \hat{P}_{k,A}'$.
Then $\mathcal{F}_{\hat{P}_{k,A}'}\to V'_k$.
\item
The image of $\supp \hat{X}_k|_{V_k} \subset \overline{A}\times \hat{P}_{k,A}$ under the map $\overline{A}\times \hat{P}_{k,A}\to \overline{A}\times P_{k,A/B}$ is contained in $\supp Y_{B,k}\subset \overline{A}\times P_{k,A/B}$.
\end{enumerate}
\item
If $k$ is sufficiently large, then for generic $y\in Z_{k,A/B}$, the natural map $(Y_{B,k})_y\to (Y_{B,k-1})_{y_0}$ is an isomorphism, where $y_0\in P_{k-1,A/B}$ is the image of $y$ under the map $P_{k,A/B}\to P_{k-1,A/B}$.
\end{enumerate}
\end{lem}

In the proof of this lemma, we use the following notation:
Let $V\subset \overline{A}\times S$ be a Zariski closed set.
Let $B\subset A$ be a semi-abelian variety.
We set
$$
\mathrm{Sp}_BV=\bigcap_{b\in B}(V+b)\subset V.
$$
Then $\mathrm{Sp}_BV\subset \overline{A}\times S$ is a Zariski closed subset.

\medskip

{\it Proof of Lemma \ref{pro:410}.}\
Let $\mathcal{B}$ be the set of all semi-abelian subvarieties $B\subset A$ such that $B$ satisfies the three properties: (1), (4) and $\mathrm{Sp}_B\overline{X}\not= \emptyset$.
To show that $\mathcal{B}$ is non-empty, we shall prove $\{ 0\}\in \mathcal{B}$.
We note that $\{0\}\in\Pi(\mathcal{F})$.
Indeed, by $\mathrm{Sp}_{\{0\}}\overline{X}=\overline{X}$, we have $\mathcal{F}\to \mathrm{Sp}_{\{0\}}\overline{X}$.
Hence if $\{ 0\}\not\in \Pi(\mathcal{F})$, then $\{0\}\in\Lambda_{X,\overline{A}}(\mathcal{F})$.
This contradicts to $\Lambda_{X,\overline{A}}(\mathcal{F})=\emptyset$.
Hence $\{0\}\in\Pi(\mathcal{F})$.
To prove that $\{0\}$ satisfies (4), we take a $Z_k$-admissible modification $\sigma:\hat{P}_{k,A}\to P_{k,A}$ such that $\hat{X}_k|_{\hat{Z}_k}\to \hat{Z}_k$ is flat, where $\hat{X}_k$ and $\hat{Z}_k$ are defined as in the statement of the lemma.
The existence of such modification follows from Lemma \ref{lem:20210831}.
We set $V_k=\hat{Z}_k$.
Then (4) is satisfied.
Thus $\{ 0\}\in \mathcal{B}$.
In particular $\mathcal{B}$ is non-empty.

We remark that if $B\in \mathcal{B}$, then (2) is satisfied for sufficiently large $k$ satisfying (4).
Indeed, by $\mathcal{F}_{P_{k,A}}\Rightarrow Z_{k}$, we may apply Lemma \ref{lem:202206233} to get $\hat{Z}_k\subset p_k(\hat{X}_k)$, where we continue to write the induced map $p_k:\overline{A}\times \hat{P}_{k,A}\to\hat{P}_{k,A}$.
Hence by $V_k\subset \hat{Z}_k$ (cf. (4a)), we have $V_k= p_k(\hat{X}_k|_{V_k})$.
By (4d), we have $\sigma(V_k)\subset p_{k,B}(Y_{B,k})$, where $p_{k,B}:\overline{A}\times P_{k,A/B}\to P_{k,A/B}$ is the second projection.
By $Z_{k,A/B}\subset \sigma(V_k)$ (cf. (4a)), we get $Z_{k,A/B}\subset p_{k,B}(Y_{B,k})$.
Hence (2) is true for $B\in\mathcal{B}$.

\medskip

{\bf Claim 1.}\
If $B\in\mathcal{B}$, then the assertion (5) is satisfied.

\medskip

{\it Proof.}\
Since $Z_{k,A/B}$ is integral, the generic flatness yields that there exists a non-empty Zariski open set $U_k\subset  Z_{k,A/B}$ such that $Y_{B,k}|_{Z_{k,A/B}}\to Z_{k,A/B}$ is flat over $U_k$.
We may assume that the image of $U_k$ under $Z_{k,A/B}\to Z_{k-1,A/B}$ is contained in $U_{k-1}$.
For each $k$, note that the Hilbert polynomials of the fibers $(Y_{B,k})_y$ are all the same for $y\in U_k$.
We denote this polynomial by $H_k$.
For $y\in U_k$, we have $(Y_{B,k})_y\subset (Y_{B,k-1})_{y_0}\subset \overline{A}$, where $y_0\in P_{k-1,A/B}$ is the image of $y$ under the map $P_{k,A/B}\to P_{k-1,A/B}$.
Hence $H_k\leq H_{k-1}$.
Thus we get $H_k\geq H_{k+1}\geq H_{k+2}\geq \cdots$.
By \cite[Lemma 8.2]{Y}, there exists $k_0$ such that $H_{k_0}=H_{k_0+1}=\cdots$.
Hence if $k\geq k_0+1$, we have $(Y_{B,k})_y= (Y_{B,k-1})_{y_0}$.
This completes the proof.
\hspace{\fill} $\square$

\medskip

In the following, we shall prove that a maximal element in $\mathcal{B}$ satisfies (3).
We take $B\in\mathcal{B}$.
We consider $Y_{B,k}|_{Z_{k,A/B}}\to Z_{k,A/B}$.
Then by the property (2), this map is surjective.
Set $C_k=\mathrm{Stab}^0(Y_{B,k})_y\subset A$ for generic $y\in Z_{k,A/B}$
(cf. Lemma \ref{lem:20210414}).
Note that $(Y_{B,k})_y\subset \overline{A}$ is $B$-invariant.
Hence
$$
B\subset C_k.
$$
By Claim 1 above, there exists $C$ such that $C_k=C$ for all sufficiently large $k$.

We shall show $C\in\mathcal{B}$.
By the construction, we have $\mathrm{Sp}_CY\not=\emptyset$.
Hence by $\mathrm{Sp}_CY\subset\mathrm{Sp}_C\overline{X}$, we have $\mathrm{Sp}_C\overline{X}\not=\emptyset$.
By $B\in \mathcal{B}$, we may take $\hat{P}_{k,A}\to P_{k,A}$, $\sigma:\hat{P}_{k,A}\to P_{k,A/B}$ and $V_k\subset \hat{P}_{k,A}$, which are described in (4).
Here and what follows, we assume that $k$ is sufficiently large satisfying (4) for $B$.
Using Lemma \ref{lem:20210831}, we take a $Z_{k,A/B}$-admissible modification $P_{k,A/B}'\to P_{k,A/B}$ such that
\begin{itemize}
\item $V'_k|_{Z_{k,A/B}'}\to Z_{k,A/B}'$ is flat,  where $Z_{k,A/B}'\subset P_{k,A/B}'$ is the minimal transform and $V_k'\subset \hat{P}_{k,A}'$ is defined by $V'_k=(V_k\times_{P_{k,A/B}}P_{k,A/B}')[Z_{k,A/B}']$ as in the statement of (4c).
\end{itemize}

\begin{equation*}
\begin{CD}
\hat{P}_{k,A}'@>\sigma'>> P_{k,A/B}'\\
@VVV @VVV   \\
\hat{P}_{k,A} @>>\sigma> P_{k,A/B}
\end{CD}
\end{equation*}

Note that by $\{0\}\in\Pi(\mathcal{F})$, Lemma \ref{lem:202206234} yields that $Z_{k,A/B}\subset \sigma(\hat{Z}_k)$.
Hence $\hat{P}_{k,A}'\to \hat{P}_{k,A}$ is a $\hat{Z}_k$-admissible modification.
Hence we may define the minimal transform $\hat{Z}_k'\subset \hat{P}_{k,A}'$.
Set $\hat{X}_k'=(\hat{X}_k\times_{\hat{P}_{k,A}}\hat{P}_{k,A}')[\hat{Z}_k']\subset \overline{A}\times \hat{P}_{k,A}'$.
We set 
$$W_k=\supp V'_k|_{Z_{k,A/B}'}\subset \hat{P}_{k,A}'.
$$
We claim that

\medskip

{\bf Claim 2.}\
The image of $\supp \hat{X}_k'|_{W_k} \subset \overline{A}\times \hat{P}_{k,A}'$ under the map $\overline{A}\times \hat{P}_{k,A}'\to \overline{A}\times P_{k,A/B}$ is contained in $\mathrm{Sp}_{C_k}(\supp Y_{B,k})$.

\medskip

{\it Proof.}\
By the definition of $C_k$, there exists a dense Zariski open $U\subset Z_{k,A/B}$ such that $\supp Y_{B,k}|_U\subset \mathrm{Sp}_{C_k}(\supp Y_{B,k})$.
By shrinking $U$, we may asuume that $U\subset Z_{k,A/B}'$.
Since $V_k'|_{Z_{k,A/B}'}\to Z_{k,A/B}'$ is flat, Lemma \ref{lem:20210417} yields that $W_k|_U\subset W_k$ is dense.
Since $\hat{X}_k|_{\hat{Z}_k}\to \hat{Z}_{k}$ is flat, Lemma \ref{lem:20210905} yields that $\hat{X}_k'|_{\hat{Z}_k'}\to \hat{Z}_k'$ is flat.
By $W_k\subset V_k'\subset \hat{Z}_k'$, we get that $\hat{X}_k'|_{W_k}\to W_k$ is flat.
Hence Lemma \ref{lem:20210417} yields that $\supp \hat{X}_k'|_{(W_k|_U)}\subset \supp \hat{X}_k'|_{W_k}$ is dense.

Now by (4d), the image of $\supp \hat{X}_k'|_{(W_k|_U)}$ under $\overline{A}\times \hat{P}_{k,A}'\to \overline{A}\times P_{k,A/B}$ is contained in $\supp Y_{B,k}|_U$, hence $\mathrm{Sp}_{C_k}(\supp Y_{B,k})$.
Thus the image of $\supp \hat{X}_k'|_{W_k}$ is contained in $\mathrm{Sp}_{C_k}(\supp Y_{B,k})$.
\hspace{\fill} $\square$

\medskip

{\bf Claim 3.}\ 
Let $P_{k,A/B}''\to P_{k,A/B}'$ be a $Z_{k,A/B}'$-admissible modification. 
Let $Z_{k,A/B}''\subset P_{k,A/B}''$ be the minimal transform.
Set $\hat{P}_{k,A}''=(\hat{P}_{k,A}'\times_{P_{k,A/B}'}P_{k,A/B}'')[Z_{k,A/B}'']$ and $W_k^+=(W_k\times_{P_{k,A/B}'}P_{k,A/B}'')[Z_{k,A/B}'']\subset \hat{P}_{k,A}''$.
Then $\mathcal{F}_{\hat{P}_{k,A}''}\to W_k^+$.

\begin{equation*}
\begin{CD}
\hat{P}_{k,A}''@>\sigma''>> P_{k,A/B}'' \\
@VVV @VVV \\
\hat{P}_{k,A}'@>\sigma'>> P_{k,A/B}'
\end{CD}
\end{equation*}

{\it Proof.}\
Set $V_k''=(V_k'\times_{P_{k,A/B}'}P_{k,A/B}'')[Z_{k,A/B}'']\subset \hat{P}_{k,A}''$.
Since $V_k'|_{Z_{k,A/B}'}\to Z_{k,A/B}'$ is flat, Lemma \ref{lem:20210905} yields that $V_k''|_{Z_{k,A/B}''}=V_k'|_{Z_{k,A/B}'}\times_{Z_{k,A/B}'}Z_{k,A/B}''$.
Similarly, applying Lemma \ref{lem:20210905} to $V_{k}'|_{Z_{k,A/B}'}\to P_{k,A/B}'$, we have 
$$((V_{k}'|_{Z_{k,A/B}'})\times_{P_{k,A/B}'}P_{k,A/B}'')[Z_{k,A/B}'']=V_k'|_{Z_{k,A/B}'}\times_{Z_{k,A/B}'}Z_{k,A/B}''.$$
Note that $W_k^+=\supp (((V_{k}'|_{Z_{k,A/B}'})\times_{P_{k,A/B}'}P_{k,A/B}'')[Z_{k,A/B}''])$.
Hence we get 
$$\supp V_k''|_{Z_{k,A/B}''}= W_k^+.$$

Now by Lemma \ref{lem:20210903}, we have $V_k''=(V_k\times_{P_{k,A/B}}P_{k,A/B}'')[Z_{k,A/B}'']$ and $\hat{P}_{k,A}''=(\hat{P}_{k,A}\times_{P_{k,A/B}}P_{k,A/B}'')[Z_{k,A/B}'']$.
Hence by (4c), we have $\mathcal{F}_{\hat{P}_{k,A}''}\to V_k''$.
By Assumption \ref{rem:20200630}, we have $\mathcal{F}_{P_{k,A/B}''}\to Z_{k,A/B}''$.
Thus by Lemmas \ref{lem:20220626} and \ref{lem:202206260},
we get $\mathcal{F}_{\hat{P}_{k,A}''}\to W_k^+$.
This conclude the proof of the claim.
\hspace{\fill} $\square$

\medskip

We are going to prove $C\in \mathcal{B}$.
We first prove $C\in\Pi(\mathcal{F})$.
By Lemma \ref{lem:20210903}, we have $\hat{X}_k'=(\overline{X}_k\times _{P_{k,A}}\hat{P}_{k,A}')[\hat{Z}_k']$.
Hence by Lemma \ref{lem:202206233}, we have
\begin{equation}\label{eqn:202204212}
f_{[k]}(\mathbb D)\subset \hat{X}_k'
\end{equation}
for all $f\in \mathcal{F}$ with finite exception.
We have $\mathcal{F}_{\hat{P}_{k,A}'}\to W_k$.
Hence by Lemma \ref{lem:202206260}, $\{ f_{[k]}\}_{f\in\mathcal{F}}\to \hat{X}_k'|_{W_k}$.
By claim 2 above, we have $\mathcal{F}\to \mathrm{Sp}_{C_k}Y$, hence $\mathcal{F}\to \mathrm{Sp}_{C_k}\overline{X}$.
Hence by $\Lambda_{X,\overline{A}}(\mathcal{F})=\emptyset$, we get $C_k\in\Pi(\mathcal{F})$.
Thus $C\in\Pi(\mathcal{F})$.
In particular, we have $C\not= A$.

By $C\in\Pi(\mathcal{F})$, the definition of $T_{k,A/B}$ yields that $T_{k,A/B}\not\subset E_{k,A/B,A/C}$.
Hence 
\begin{equation}\label{eqn:202109112}
Z_{k,A/B}\not\subset E_{k,A/B,A/C}.
\end{equation}
Note that the rational map $P_{k,A/B}\dashrightarrow P_{k,A/C}$ is regular outside $E_{k,A/B,A/C}$.
Thus we may moreover assume for the $Z_{k,A/B}$-admissible map $P_{k,A/B}'\to P_{k,A/B}$ that
\begin{itemize}
\item $\mu:P_{k,A/B}'\to P_{k,A/C}$ exists.
\end{itemize}
We get the following:
\begin{equation*}
\begin{CD}
\hat{P}_{k,A}'@>\sigma'>> P_{k,A/B}'@>\mu>> P_{k,A/C} \\
@VVV @VVV   \\
\hat{P}_{k,A} @>>\sigma> P_{k,A/B}
\end{CD}
\end{equation*}
By Lemma \ref{lem:202206234}, we have 
\begin{equation}\label{eqn:20220425}
Z_{k,A/C}\subset \mu(Z_{k,A/B}').
\end{equation}

Now we shall show that $\tau:\hat{P}_{k,A}'\to P_{k,A/C}$ and $W_k\subset \hat{P}_{k,A}'$ satisfy the condition (4).
By Lemma \ref{lem:20210903}, we have $\hat{X}_k'=(X_k\times _{P_{k,A}}\hat{P}_{k,A}')[\hat{Z}_k']$.
Note that $Z_{k,A/C}\subset \tau(W_k)$ follows from $\sigma' (W_k)=Z_{k,A/B}'$ (cf. \eqref{eqn:20220425}).
Also by $V_k\subset \hat{Z}_k$, we have $W_k\subset V_k'\subset \hat{Z}_k'$.
These show (4a).
Since $\hat{X}_k|_{\hat{Z}_k}\to \hat{Z}_{k}$ is flat, Lemma \ref{lem:20210905} yields that $\hat{X}_k'|_{\hat{Z}_k'}\to \hat{Z}_{k}'$ is flat.
This shows (4b).

Next we prove (4c).
Let $P_{k,A/C}'\to P_{k,A/C}$ be a $Z_{k,A/C}$-admissible modification.
We consider the following
\begin{equation*}
\begin{CD}
\hat{P}_{k,A}''@>\sigma''>> P_{k,A/B}'' @>>> P_{k,A/C}'\\
@V\psi VV @VVV @VVV \\
\hat{P}_{k,A}'@>\sigma'>> P_{k,A/B}'@>\mu>> P_{k,A/C} 
\end{CD}
\end{equation*}
Here $P_{k,A/B}''=(P_{k,A/B}'\times_{P_{k,A/C}}P_{k,A/C}')[Z_{k,A/C}']$.
Then $P_{k,A/B}''\to P_{k,A/B}'$ is a $Z_{k,A/B}'$-admissible modification (cf. \eqref{eqn:20220425}).
Since $\hat{P}_{k,A}'$ is integral, we may apply Lemma \ref{lem:20210909} to get $\hat{P}_{k,A}''=(\hat{P}_{k,A}'\times_{P_{k,A/C}}P_{k,A/C}')[Z_{k,A/C}']$ (cf. \eqref{eqn:20220425}).
Set $W_k'=(W_k\times_{P_{k,A/C}}P_{k,A/C}')[Z_{k,A/C}']$.
Then by Lemma \ref{lem:20210909}, we get $W_k^+\subset W_k'$ (cf. \eqref{eqn:20220425}).
Thus by Claim 3, we get $\mathcal{F}_{\hat{P}_{k,A}''}\to W_k'$.
This shows (4c).

Finally to check (4d), we prove that the image $\hat{X}_k'|_{W_k} \subset \overline{A}\times \hat{P}_{k,A}'\to \overline{A}\times P_{k,A/C}$ is contained in $S_{C,k}\subset \overline{A}\times P_{k,A/C}$, where $S=\mathrm{Sp}_C\overline{X}$.
We have $\mathrm{Sp}_CY=S$.
Hence we get $\mathrm{Sp}_C(Y_{B,k})=S_{B,k}$.
Hence Lemma \ref{lem:202206222} yields
\begin{equation}\label{eqn:20210831}
\mathrm{Sp}_C(Y_{B,k})|_{(P_{k,A/B}\backslash E_{k,A/B,A/C})}\subset \varphi^{-1}(S_{C,k}),
\end{equation}
where $\varphi:\overline{A}\times (P_{k,A/B}\backslash E_{A/B,A/C,k})\to \overline{A}\times P_{k,A/C}$.
We take $U\subset P_{k,A/B}$ such that $U\cap Z_{k,A/B}\not=\emptyset$ and $P_{k,A/B}'\to P_{k,A/B}$ is an isomorphism over $U$.
We consider as $U\subset P_{k,A/B}'$.
By \eqref{eqn:202109112}, we may assume
$$U \subset P_{k,A/B}-E_{k,A/B,A/C}.$$
Since $V_k'|_{Z_{k,A/B}'}\to Z_{k,A/B}'$ is flat, Lemma \ref{lem:20210417} yields that $W_k|_{(Z_{k,A/B}'\cap U)}\subset W_k$ is dense.
By Claim 2 and \eqref{eqn:20210831}, the image of $\hat{X}_k'|_{W_k|_{(Z_{k,A/B}'\cap U)}}$ under $\overline{A}\times \hat{P}_{k,A}'\to \overline{A}\times P_{k,A/C}$ is contained in $S_{C,k}$, provided $k$ is sufficiently large so that $C_k=C$.
Since $\hat{X}_k'|_{W_k}\to W_k$ is flat, Lemma \ref{lem:20210417} yields that the inclusion $\hat{X}_k'|_{W_k|_{(Z_{k,A/B}'\cap U)}}\subset \hat{X}_k'|_{W_k}$ is dense.
Hence the image $\supp \hat{X}_k'|_{W_k}\subset \overline{A}\times \hat{P}_{k,A}'\to \overline{A}\times P_{k,A/C}$ is contained in $S_{C,k}\subset \overline{A}\times P_{k,A/C}$.
Thus we have proved $C\in \mathcal{B}$.

Now we finish the proof of the lemma.
We take maximal $B\subset A$ such that $B\in\mathcal{B}$.
Set $C=\mathrm{Stab}^0(Y_{B,k})_y\subset A$ for generic $y\in Z_{k,A/B}$ and sufficiently large $k$.
Then $B\subset C$.
By $C\in \mathcal{B}$ and the maximal property of $B\in \mathcal{B}$, we have $C=B$.
This conclude the proof of Lemma \ref{pro:410}.
\hspace{\fill} $\square$

\subsection{Proof of Proposition \ref{cor:20210306}}
Since $\mathcal{F}$ satisfies Assumption \ref{rem:20200630} and $\Lambda_{X,\overline{A}}(\mathcal{F})=\emptyset$, we may take $B$ as in Lemma \ref{pro:410}.
Then $B\in \Pi(\mathcal{F})$.
In particular, $B\not= A$.

Next we prove the property (2) of Proposition \ref{cor:20210306}.
We take sufficiently large $k_0$ such that the properties of Lemma \ref{pro:410} (2), (3), (4) and (5) are true for $k\geq k_0$.
For each $k\in \mathbb Z_{k\geq 1}$, let $Z_{k,A/B}\subset P_{k,A/B}$ be defined by \eqref{eqn:20220623}.
Then $Z_{k,A/B}$ is an irreducible Zariski closed set such that $T_{k,A/B}\subset Z_{k,A/B}$.
Set $Y=\mathrm{Sp}_B\overline{X}$.
If $k\geq k_0$, then for generic $y\in Z_{k,A/B}$, the natural map $(Y_{B,k})_y\to (Y_{B,k-1})_{y_0}$ is an isomorphism, where $y_0\in P_{k-1,A/B}$ is the image of $y$ under the map $P_{k,A/B}\to P_{k-1,A/B}$
(cf. Lemma \ref{pro:410} (5)).
We set $S=P_{k_0,A/B}$.
Then $P_{k_0+l,A/B}\subset S_{l,A/B}$ (cf. Remark \ref{rem:20230312}).
We consider as $Z_{k_0+l,A/B}\subset S_{l,A/B}$. 
We set $Z=Z_{k_0+1,A/B}\subset S_{1,A/B}$.
By Lemma \ref{lem:alglem}, $Z$ is horizontally integrable.
Indeed, the assumption (1) of Lemma \ref{lem:alglem} follows from Lemma \ref{pro:410} (3).
The assumption (2) of Lemma \ref{lem:alglem} directly follows from Lemma \ref{pro:410} (5).
By $B\in \Pi(\mathcal{F})$, the definition of $T_{k_0+l,A/B}$ yields that $T_{k_0+l,A/B}\not\subset E_{k_0+l,A/B,A/B}$ as subsets of $P_{k_0+l,A/B}$.
By $P_{k_0+l,A/B}^{\mathrm{sing}}= E_{k_0+l,A/B,A/B}$ (cf. \eqref{eqn:202206235}), we get $T_{k_0+l,A/B}\not\subset P_{k_0+l,A/B}^{\mathrm{sing}}$, so $Z_{k_0+l,A/B}\not\subset P_{k_0+l,A/B}^{\mathrm{sing}}$.
Note that $P_{k_0+l,A/B}\subset S_{l,A/B}\subset P(TS_{l-1,A/B}\times \mathrm{Lie}(A/B))$.
Then we have
\begin{equation*}
\begin{split}
S_{l,A/B}^{\mathrm{sing}}\cap P_{k_0+l,A/B}
&=P(T_{S_{l-1,A/B}/S}\times\{0\})\cap P_{k_0+l,A/B}
\\
&\subset P(TS_{l-1,A/B}\times\{0\})\cap P_{k_0+l,A/B}= P_{k_0+l,A/B}^{\mathrm{sing}}.
\end{split}
\end{equation*}
Thus $Z_{k_0+l,A/B}\not\subset S_{l,A/B}^{\mathrm{sing}}$.
Hence for generic $y\in Z$, there exists $y_l\in Z_{k_0+l,A/B}-S_{l,A/B}^{\mathrm{sing}}$ such that $(Y_{B,k_0+1})_y=(Y_{B,k_0+l})_{y_l}$.
We note that
$$Y_{B,k_0+l}=(Y_{B,k_0})_{B,l}\cap (\overline{A}\times P_{k_0+l,A/B})$$ 
for $l\geq 1$, which follows from the definitions \eqref{eqn:202206236} and \eqref{eqn:202109111}.
Hence the assumption (3) of Lemma \ref{lem:alglem} is satisfied.
Thus $Z$ is horizontally integrable.

Now we prove $\mathcal{F}\to \mathrm{Sp}_B\overline{X}$.
We take $\hat{P}_{k,A}\to P_{k,A}$ and $\hat{X}_k\subset \overline{A}\times \hat{P}_{k,A}$ as in Lemma \ref{pro:410} (4), where we fix $k\geq k_0$.
By Lemma \ref{lem:202206233}, we have $f_{[k]}(\mathbb D)\subset \hat{X}_k$ for all $f\in \mathcal{F}$ with finite exception.
Hence by $\mathcal{F}_{\hat{P}_{k,A}}\to V_k$ (cf. Lemma \ref{pro:410} (4c)), we have $\{ f_{[k]}\}_{f\in \mathcal{F}}\to \hat{X}_k|_{V_k}$ (cf. Lemma \ref{lem:202206260}).
Hence by Lemma \ref{pro:410} (4d), we have $\mathcal{F}\to \mathrm{Sp}_B\overline{X}$.
\hspace{\fill} $\square$

\section{Proof of Theorem \ref{thm:mstrong}}\label{sec:10}

The aim of this section is to prove Theorem \ref{thm:mstrong}.

\begin{lem}\label{lem:202205233}
Let $\mathcal{F}=(f_i)_{i\in I}$ be an infinite indexed family of holomorphic maps in $\mathrm{Hol}(\mathbb D,A)$.
Then replacing $\mathcal{F}$ by its infinite subfamily, we have $\Pi(\mathcal{G})=\Pi(\mathcal{F})$ for all infinite subfamily $\mathcal{G}$ of $\mathcal{F}$.
\end{lem}

{\it Proof.}\
There are only countably many semi-abelian subvarieties $B\subset A$ (cf. \cite[Cor. 5.1.9]{NW}).
So we enumerate them as $B_1,B_2,\ldots$.
Note that this is possibly finite.
If $I$ contains an infinite subset $I'\subset I$ such that $(|(\varpi_{B_1}\circ f_i)'|_{\omega_{A/B_1}})_{i\in I'}$ converges uniformly on compact subsets of $\mathbb D$ to $0$, then we set $I_1=I'$.
If there is no such $I'$, then we set $I_1=I$.
If $I_1$ contains an infinite subset $I'\subset I_1$ such that $\{ |(\varpi_{B_2}\circ f_i)'|_{\omega_{A/B_2}}\}_{i\in I'}$ converges uniformly on compact subsets of $\mathbb D$ to $0$, then we set $I_2=I'$.
If there is no such $I'$, then we set $I_2=I_1$.
We continue this process to get (possibly finite) decreasing sequence of infinite sets $I\supset I_1\supset I_2\supset \ldots$.
If this sequence is finite $I_1\supset I_2\supset \cdots\supset I_L$, then we continue infinitely by letting $I_L=I_{L+1}=I_{L+2}=\cdots$.
We define an infinite sequence $i_1,i_2,i_3,\ldots$ of distinct elements in $I$ so that for each $l\geq 1 $, we have 
\begin{equation}\label{eqn:202205233}
i_k\in I_l
\end{equation}
for all $k\geq l$.
This sequence is constructed inductively as follows.
We take $i_1\in I_1$.
Suppose distinct elements $i_1,i_2,\ldots,i_n$ are chosen so that \eqref{eqn:202205233} holds for all $k,l\in \mathbb N$ with $1\leq l\leq k\leq n$.
Then we choose $i_{n+1}\in I_{n+1}-\{ i_1,\ldots,i_n\}$.
Thus we have constructed the sequence $i_1,i_2,\ldots$ with the desired property \eqref{eqn:202205233}.

Now we set $J=\{ i_1,i_2,i_3,\ldots\}$.
Then for all semi-abelian subvariety $B$, we have either 
\begin{itemize}
\item
$(|(\varpi_{B}\circ f_i)'|_{\omega_{A/B}})_{i\in J}$ converges uniformly on compact subsets of $\mathbb D$ to $0$, or 
\item
no infinite subfamily of $(|(\varpi_{B}\circ f_i)'|_{\omega_{A/B}})_{i\in J}$ converges uniformly on compact subsets of $\mathbb D$ to $0$.
\end{itemize}
Hence replacing $\mathcal{F}$ by $(f_j)_{j\in J}$, we get $\Pi(\mathcal{G})=\Pi(\mathcal{F})$ for all infinite subfamily $\mathcal{G}$ of $\mathcal{F}$.
\hspace{\fill} $\square$

\begin{lem}\label{lem:20220523}
Let $\mathcal{F}=(f_i)_{i\in I}$ be an infinite indexed family of non-constant holomorphic maps in $\mathrm{Hol}(\mathbb D,A)$.
Then there exists an infinite subfamily $\mathcal{G}$ of $\mathcal{F}$ such that all infinite subfamilies of $\mathcal{G}$ satisfy Assumption \ref{ass:202103061}. 
\end{lem}

{\it Proof.}\
By replacing $\mathcal{F}$ by its infinite subfamily, we have
\begin{equation}\label{eqn:202205231}
\Pi(\mathcal{F}')=\Pi(\mathcal{F})
\end{equation}
for all infinite subfamily $\mathcal{F}'$ of $\mathcal{F}$ (cf. Lemma \ref{lem:202205233}).
Let $k\geq 1$ and $B\in\Pi(\mathcal{F})$.
Then only finitely many $i\in I$ satisfies $(f_i)_{P_{k,A}}(\mathbb D)\subset E_{k,A,A/B}$, for if $(f_i)_{P_{k,A}}(\mathbb D)\subset E_{k,A,A/B}$, then $(\varpi_B\circ f_i)'=0$.
Hence, given infinite subfamily $\mathcal{F}'$ of $\mathcal{F}$, we may apply Lemma \ref{lem:202105284} to get an infinite subfamily $\mathcal{H}$ of $\mathcal{F}'$ such that $\mathrm{LIM}(\mathcal{H}_{P_{k,A}},\{E_{k,A,A/B}\}_{B\in\Pi(\mathcal{F})})$ exists.

We apply the argument above for $k=1$ to get an infinite subfamily $\mathcal{F}_1=(f_i)_{i\in I_1}$, where $I_1\subset I$,  such that $\mathrm{LIM}((\mathcal{F}_1)_{P_{1,A}},\{E_{1,A,A/B}\}_{B\in\Pi(\mathcal{F})})$ exists.
Again we apply the argument above for $k=2$ to get an infinite subfamily $\mathcal{F}_2=(f_i)_{i\in I_2}$, where $I_2\subset I_1$,  such that $\mathrm{LIM}((\mathcal{F}_2)_{P_{2,A}},\{E_{2,A,A/B}\}_{B\in\Pi(\mathcal{F})})$ exists.
Continue this process to get a decreasing sequence of infinite sets $I\supset I_1\supset I_2\supset \ldots$.
We define a countable infinite subset $J=\{ i_1,i_2,i_3,\ldots\}$ of $I$ such that  \eqref{eqn:202205233} holds for all $l\geq 1$ and $k\geq l$ (cf. the proof of Lemma \ref{lem:202205233}).
Set $\mathcal{G}=(f_i)_{i\in J}$.
Then $\mathrm{LIM}(\mathcal{G}_{P_{k,A}},\{E_{k,A,A/B}\}_{B\in\Pi(\mathcal{F})})$ exists for all $k\geq  1$ (cf. Remark \ref{rem:20210805} (3)).
Let $\mathcal{H}$ be an infinite subfamily of $\mathcal{G}$.
Then $\mathrm{LIM}(\mathcal{H}_{P_{k,A}},\{E_{k,A,A/B}\}_{B\in\Pi(\mathcal{F})})$ exists for all $k\geq  1$ (cf. Remark \ref{rem:20210805} (3)).
Hence, by \eqref{eqn:202205231}, $\mathrm{LIM}(\mathcal{H}_{P_{k,A}},\{E_{k,A,A/B}\}_{B\in\Pi(\mathcal{H})})$ exists for all $k\geq  1$.
Thus $\mathcal{H}$ satisfies Assumption \ref{ass:202103061}. 
The proof is completed.
\hspace{\fill} $\square$

\begin{lem}\label{lem:20201227}
Let $\mathcal{F}\subset \mathrm{Hol}(\mathbb D,A)$ be an infinite set of non-constant holomorphic maps.
By replacing $\mathcal{F}$ by its infinite subset, we may assume Assumption \ref{rem:20200630}.
\end{lem}

{\it Proof.}
By replacing $\mathcal{F}$ by its infinite subset, we may assume \eqref{eqn:202205231} for all infinite subset $\mathcal{F}'\subset \mathcal{F}$ (cf. Lemma \ref{lem:202205233}).
We enumerate the elements of $\Pi(\mathcal{F})$ as $B_1,B_2,\ldots$.
Note that this is possibly finite.
We apply Lemma \ref{lem:20220523} to get an infinite subset $\mathcal{F}_1\subset \mathcal{F}$ such that all infinite subfamilies of $(\mathcal{F}_1)_{A/B_1}$ satisfy Assumption \ref{ass:202103061}.
Again we apply Lemma \ref{lem:20220523} to get an infinite subset $\mathcal{F}_2\subset \mathcal{F}_1$ such that all infinite subfamilies of $(\mathcal{F}_2)_{A/B_2}$ satisfy Assumption \ref{ass:202103061}.
We continue this process to get (possibly finite) descending sequence of infinite sets $\mathcal{F}\supset \mathcal{F}_1\supset \mathcal{F}_2\supset \ldots$.
If this sequence is finite $\mathcal{F}_1\supset \mathcal{F}_2\supset \cdots\supset \mathcal{F}_L$, then we continue infinitely by letting $\mathcal{F}_L=\mathcal{F}_{L+1}=\mathcal{F}_{L+2}=\cdots$.
We define a countably infinite subset $\mathcal{G}=\{ f_1,f_2,f_3,\ldots\}$ such that $f_k\in \mathcal{F}_l$ for all $l\geq 1$ and $k\geq l$ (cf. the proof of Lemma \ref{lem:202205233}).
Then, for all $B\in\Pi(\mathcal{F})$, $\mathcal{G}_{A/B}$ satisfies Assumption \ref{ass:202103061}. 
By \eqref{eqn:202205231}, we have $\Pi(\mathcal{G})=\Pi(\mathcal{F})$.
Hence for all $B\in\Pi(\mathcal{G})$, $\mathcal{G}_{A/B}$ satisfies Assumption \ref{ass:202103061}. 
Hence $\mathcal{G}$ satisfies Assumption \ref{rem:20200630}.
\hspace{\fill} $\square$

\begin{lem}\label{thm:m}
Let $X\subsetneqq A$.
Let $\overline{A}$ be a smooth equivariant compactification and let $\overline{X}\subset \overline{A}$ be the compactification.
Let $\mathcal{F}\subset \mathrm{Hol}(\mathbb D,X)$ be an infinite set of holomorphic maps.
Then there exist a semi-abelian subvariety $B\subset A$ and an infinite subset $\mathcal{G}\subset \mathcal{F}$ with the following two properties:
\begin{enumerate}
\item Let $\overline{A/B}$ be a smooth equivariant compactification.
Then there exists an infinite subset $\mathcal{H}\subset \mathcal{G}$ such that $\{ \varpi\circ f\}_{f\in \mathcal{H}}$ converges uniformly on compact subsets of $\mathbb D$ to $g:\mathbb D\to \overline{A/B}$, where $\varpi:A\to A/B$ is the quotient map.
\item $\mathcal{G}\to \mathrm{Sp}_B\overline{X}$.
\end{enumerate}
\end{lem}

{\it Proof.}\
We may assume that $\mathcal{F}$ contains only finitely many constant maps, for otherwise our lemma is trivial by setting $B=\{ 0\}$ and $\mathcal{G}$ is the set of constant mappings in $\mathcal{F}$.
We remove those constant maps from $\mathcal{F}$ and assume that all elements of $\mathcal{F}$ are non-constant. 

We first reduce to the case $\overline{A}$ is projective.
When this is not satisfied, we may take an equivariant modification $\hat{A}\to \overline{A}$ such that $\hat{A}$ is smooth and projective (cf. Lemma \ref{lem:202206031}).
Let $\hat{X}\subset \hat{A}$ be the compactification.
Then the image of $\mathrm{Sp}_B\hat{X}\subset \hat{X}$ under the natural map $\hat{X}\to \overline{X}$ is contained in $\mathrm{Sp}_B\overline{X}$.
Hence $\mathcal{G}\to \mathrm{Sp}_B\hat{X}$ implies $\mathcal{G}\to \mathrm{Sp}_B\overline{X}$.
Hence we have reduced to the case that $\overline{A}$ is projective.

By Lemma \ref{lem:20201227}, we may take an infinite subset $\mathcal{G}\subset \mathcal{F}$ such that $\mathcal{G}$ satisfies Assumption \ref{rem:20200630}.
We prove the lemma in two cases.

If $\Lambda_{X,\overline{A}}(\mathcal{G})\not= \emptyset$, then we take $B\in \Lambda_{X,\overline{A}}(\mathcal{G})$.
Then $\mathcal{G}\to \mathrm{Sp}_B\overline{X}$ and there exists an infinite subset $\mathcal{H}\subset \mathcal{G}$ such that $\{ |(\varpi\circ f)'|_{\omega_{A/B}}\}_{f\in \mathcal{H}}$ converges uniformly on compact subsets of $\mathbb D$ to $0$. 
Hence we conclude that $\{ \varpi\circ f\}_{f\in \mathcal{H}}$ converges to some constant map.
Thus our lemma is valid in this case.

If $\Lambda_{X,\overline{A}}(\mathcal{G})=\emptyset$, then we take $B$ as in Proposition \ref{cor:20210306}.
We have $\mathcal{G}\to \mathrm{Sp}_B\overline{X}$.
By Proposition \ref{cor:20210306} (1), we have $B\in\Pi(\mathcal{G})$.
Hence an infinite indexed family $\mathcal{G}_{A/B}$ in $\mathrm{Hol}(\mathbb D,A/B)$ contains only finitely many constant maps and satisfies Assumption \ref{ass:202103061}.
By Proposition \ref{cor:20210306} (2), we get $k$ such that $Z_{k+1,A/B}\subset P_{k+1,A/B}$ is horizontally integrable, where $Z_{k+1,A/B}$ is defined in \eqref{eqn:20220623}.
Let $\overline{A/B}$ be a smooth equivariant compactification. 
By Corollary \ref{cor:20220424}, we may choose an infinite subset $\mathcal{H}\subset \mathcal{G}$ such that $\{ \varpi\circ f\}_{f\in\mathcal{H}}$ converges uniformly on compact subsets of $\mathbb D$ to $g:\mathbb D\to \overline{A/B}$.
The proof is completed.
\hspace{\fill} $\square$

\medskip

Let $B\subset A$ be a semi-abelian subvariety.
Let $\hat{A}$ and $\overline{A/B}$ be equivariant compactifications with an equivariant map $p:\hat{A}\to \overline{A/B}$.
We consider the following assumption:

\begin{ass}\label{ass:20220626}
For every $x\in \hat{A}$, if a semi-abelian subvariety $C\subset A$ satisfies $C\cdot p(x)=p(x)$, then $C\cdot x\subset B\cdot x$ as subsets of $\hat{A}$.
\end{ass}

By choosing an equivariant blow-up $\varphi:\hat{A}\to \bar{A}$ and a particular equivariant compactification $\overline{A/B}$, we may assume this assumption.
We may assume that $\overline{A/B}$ is smooth.
This is a consequence of Lemma \ref{lem:20210513}.

\begin{lem}\label{lem:202206263}
Let $X\subsetneqq A$ be a proper closed subvariety.
Let $B\subset A$ be a semi-abelian subvariety.
Let $\overline{A/B}$ be a smooth equivariant compactification and let $\hat{A}$ be an equivariant compactification such that an equivariant map $p:\hat{A}\to \overline{A/B}$ with Assumption \ref{ass:20220626} exists.
Let $\mathcal{F}\subset \mathrm{Hol}(\mathbb D,X)$ be an infinite set of holomorphic maps such that $\mathcal{F}\to \mathrm{Sp}_B\hat{X}$, where $\hat{X}\subset \hat{A}$ is the Zariski closure of $X\subset A$.
Assume that $(p\circ f)_{f\in\mathcal{F}}$ converges uniformly on compact subsets of $\mathbb D$ to $g:\mathbb D\to \overline{A/B}$ such that $g(\mathbb D)\subset \partial (A/B)$.
Then there exists a semi-abelian subvariety $B'\subset A$ with $B\subset B'$ such that the following two properties hold:
\begin{enumerate}
\item $(p'\circ f)_{f\in\mathcal{F}}$ converges uniformly on compact subsets of $\mathbb D$ to $h:\mathbb D\to A/B'$, where $p':A\to A/B'$ is the quotient map.
\item $\mathcal{F}\to \mathrm{Sp}_{B'}\hat{X}$.
\end{enumerate}
\end{lem}

{\it Proof.}\
Let $Y\subset \hat{X}$ be the smallest Zariski closed set such that $\mathcal{F}\to Y$ (cf. Lemma \ref{lem:202206261}).
We have 
\begin{equation}\label{eqn:202206271}
Y\subset \mathrm{Sp}_B(\hat{X}).
\end{equation}
Note that $p(Y)$ is the smallest Zariski closed set such that $(p\circ f)_{f\in\mathcal{F}}\to p(Y)$.
Hence $p(Y)$ is the Zariski closure of $g(\mathbb D)$ (cf. Remark \ref{rem:20220626}).
By $g(\mathbb D)\subset \partial (A/B)$, we have $p(Y)\subset \partial (A/B)$.
Let $V\subset \overline{A/B}$ be the Zariski closure of $p(Y)+(A/B)\subset \overline{A/B}$.
Then since $p(Y)$ is irreducible, $V$ is irreducible.
Note that $V$ is $A/B$-invariant.
Hence by Lemma \ref{lem:vect}, there exists an $A/B$-invariant Zariski open set $W\subset \overline{A/B}$ such that $V\subset W$ and that an equivariant morphism $\tau: W\to V$ exists. 
Let $I\subset A/B$ be the isotropy group for $V$.
Let $O\subset V$ be the unique dense $A/B$-orbit.
Then $O$ is isomorphic to $(A/B)/I$ and $O\cap p(Y)\not=\emptyset$.
Hence 
\begin{equation}\label{eqn:29229627}
g(\mathbb D)\not\subset V\backslash O.
\end{equation}
Since $(p\circ f)_{f\in\mathcal{F}}$ converges uniformly on compact subsets on $\mathbb D$ to $g$, the indexed family $(\tau\circ p\circ f)_{f\in\mathcal{F}}$ converges uniformly on compact subsets on $\mathbb D$ to $\tau\circ g=g$.
We have $\tau\circ p\circ f(\mathbb D)\subset O$ for all $f\in\mathcal{F}$.
Hence by \eqref{eqn:29229627}, we have $g(\mathbb D)\subset O$.
Let $B'$ be $B\subset B'\subset A$ such that $B'/B=I$.
Then under the isomorphism $O\simeq A/B'$, we have $p'\circ f=\tau\circ p\circ f$.
Hence we have proved (1).

In the following, we shall prove $\mathcal{F}\to \mathrm{Sp}_{B'}(\hat{X})$.
Note that the action of $B'$ on $p(Y)$ is trivial.
Hence by Assumption \ref{ass:20220626}, we have $B'+Y\subset B+Y$.
Since $\mathrm{Sp}_B(\hat{X})$ is $B$-invariant, we have $B+Y\subset \mathrm{Sp}_B\hat{X}$ (cf. \eqref{eqn:202206271}).
Hence we get
$$B'+Y\subset \mathrm{Sp}_B\hat{X}\subset \hat{X}.$$
Hence $Y\subset \mathrm{Sp}_{B'}\hat{X}$, so $\mathcal{F}\to \mathrm{Sp}_{B'}\hat{X}$.
\hspace{\fill} $\square$

\begin{lem}\label{lem:20200910}
Let $X\subsetneqq A$ be a proper closed subvariety.
Then there exists a smooth equivariant compactification $\overline{A}$ such that $\mathrm{Sp}_A(\overline{X})= \emptyset$. 
\end{lem}

{\it Proof.}\
Let $\tilde{A}$ be an equivariant compactification and let $\tilde{X}\subset \tilde{A}$ be the closure.
Set $Z=\bigcap_{a\in A}(a+\tilde{X})$ as scheme.
By $X\not=A$, we have $\supp Z\subset \partial A$.
Let $m:A\times \tilde{A}\to\tilde{A}$ be the $A$-action on the equivariant compactification $\tilde{A}$.
We first prove
\begin{equation}\label{eqn:202206246}
m^*Z= A\times Z
\end{equation}
Indeed we have $a+Z=Z$ for all $a\in A$.
Hence by Lemma \ref{lem:20211030}, we have $A\times Z\subset m^*Z$.
To prove the converse $m^*Z\subset A\times Z$, we consider an isomorphism $\mu:A\times\tilde{A}\to A\times \tilde{A}$ defined by $(a,x)\mapsto (a,m(-a,x))$.
By Lemma \ref{lem:20211030}, we also have $A\times Z\subset \mu^*(A\times Z)$.
Note that $m\circ \mu:A\times \tilde{A}\to\tilde{A}$ is the second projection.
Hence $\mu^*m^*Z=A\times Z$.
Thus $\mu^*m^*Z\subset \mu^*(A\times Z)$, hence $m^*Z\subset A\times Z$.
Thus we get \eqref{eqn:202206246}.

We claim that the blow-up $\hat{A}=\mathrm{Bl}_Z\tilde{A}$ is an equivariant compactification.
Indeed the map $m:A\times \tilde{A}\to\tilde{A}$ induces $m':A\times \mathrm{Bl}_Z\tilde{A}\to \tilde{A}$.
Note that $A\times \mathrm{Bl}_Z\tilde{A}=\mathrm{Bl}_{A\times Z}(A\times\tilde{A})$.
Hence by \eqref{eqn:202206246}, $(m')^*Z$ is a Cartier divisor.
Hence $m'$ factors as $A\times \mathrm{Bl}_Z\tilde{A}\to \mathrm{Bl}_Z\tilde{A}$.
Hence we have proved that $\mathrm{Bl}_Z\tilde{A}$ is an equivariant compactification.

Now we may take finite points $a_1,\ldots,a_l\in A$ such that $Z=\bigcap_{i=1}^l(a_i+\tilde{X})$ as scheme.
Then we have 
\begin{equation}\label{eqn:202206247}
\bigcap_{i=1}^l(a_i+\hat{X})=\emptyset,
\end{equation}
where $\hat{X}\subset \hat{A}$ is the Zariski closure of $X$.
We prove this.
Let $\varphi:\hat{A}\to\tilde{A}$ be the morphism.
Since $\varphi^*Z\subset \hat{A}$ is a Cartier divisor, $\mathcal{I}_{\varphi^*Z}\subset \mathcal{O}_{\hat{A}}$ is an invertible sheaf.
By $\mathcal{I}_{\varphi^*(a_i+\tilde{X})}\subset \mathcal{I}_{\varphi^*Z}\subset\mathcal{O}_{\hat{A}}$, we take $Y_i\subset \hat{A}$ such that $\mathcal{I}_{Y_i}=\mathcal{I}_{\varphi^*(a_i+\tilde{X})}\otimes (\mathcal{I}_{\varphi^*Z})^{-1}\subset \mathcal{O}_{\hat{A}}$.
Then by $\mathcal{I}_{\varphi^*(a_1+\tilde{X})}+\cdots+\mathcal{I}_{\varphi^*(a_l+\tilde{X})}=\mathcal{I}_{\varphi^*Z}$, we have $\mathcal{I}_{\varphi^*Z}\cdot (\mathcal{I}_{Y_1}+\cdots+\mathcal{I}_{Y_l})=\mathcal{I}_{\varphi^*Z}$ so that $\mathcal{I}_{Y_1}+\cdots+\mathcal{I}_{Y_l}=\mathcal{O}_{\hat{A}}$.
Hence $Y_1\cap\cdots\cap Y_l=\emptyset$.
Note that $a_i+\hat{X}\subset Y_i$.
Hence we get \eqref{eqn:202206247}.
By $\mathrm{Sp}_A(\hat{X})\subset \bigcap_{i=1}^l(a_i+\hat{X})$, we have $\mathrm{Sp}_A(\hat{X})= \emptyset$. 
By Lemma \ref{lem:20211030111}, we may take a smooth equivariant modification $\overline{A}\to\hat{A}$. 
Then $\overline{A}$ satisfies our assertion.
(See also \cite[Prop. 5.6.7]{NW}.)
\hspace{\fill} $\square$

\medskip

Let $X\subset A$ be a closed subvariety and $\mathcal{F}\subset \mathrm{Hol}(\mathbb D,X)$ be an infinite set.
Given an equivariant compactification $\bar{A}$, we set
$$
Q(\bar{A},\mathcal{F})=\{ B\subset A\ ;\ \mathcal{F}\to \mathrm{Sp}_B\bar{X}\},
$$
where $\bar{X}\subset \bar{A}$ is the Zariski closure.
If $\mathcal{G}\subset \mathcal{F}$ is an infinite subset, then we have $Q(\bar{A},\mathcal{F})\subset Q(\bar{A},\mathcal{G})$.
Indeed, if $B\in Q(\bar{A},\mathcal{F})$, then $\mathcal{F}\to \mathrm{Sp}_B\bar{X}$.
Hence $\mathcal{G}\to \mathrm{Sp}_B\bar{X}$, which shows $B\in Q(\bar{A},\mathcal{G})$.

\begin{lem}\label{lem:202206254}
Replacing $\mathcal{F}$ by its infinite subset, we have $Q(\bar{A},\mathcal{F})=Q(\bar{A},\mathcal{G})$ for all infinite subset $\mathcal{G}\subset \mathcal{F}$ and all smooth equivariant compactification $\bar{A}$.
\end{lem}

{\it Proof.}\
The pairs $(\hat{A},B)$ of smooth equivariant compactifications $\hat{A}$ and semi-abelian subvarieties $B\subset A$ are countable (cf. Lemma \ref{lem:20220728} and \cite[Cor. 5.1.9]{NW}).
So we enumerate them as $(\hat{A}_1,B_1),(\hat{A}_2,B_2),\ldots$.
In the following, we construct a descending sequence $\mathcal{F}\supset \mathcal{F}_1\supset \mathcal{F}_2\supset \cdots$ of infinite subsets such that for all $k\geq 1$, one of the followings occurs:
\begin{enumerate}
\item 
$\mathcal{F}_k\to \mathrm{Sp}_{B_k}\hat{X}_k$, where $\hat{X}_k\subset \hat{A}_k$ is the Zariski closure of $X$ in $\hat{A}_k$, or
\item
$\mathcal{F}'\not\to \mathrm{Sp}_{B_k}\hat{X}_k$ for all infinite subset $\mathcal{F}'\subset \mathcal{F}_k$.
\end{enumerate}
The construction is the following.
If there exists an infinite subset $\mathcal{F}'\subset \mathcal{F}$ such that $\mathcal{F}'\to \mathrm{Sp}_{B_1}\hat{X}_1$, then we set $\mathcal{F}_1=\mathcal{F}'$, otherwise we set $\mathcal{F}_1=\mathcal{F}$.
If there exists $\mathcal{F}'\subset \mathcal{F}_1$ such that $\mathcal{F}'\to \mathrm{Sp}_{B_2}\hat{X}_2$, then we set $\mathcal{F}_2=\mathcal{F}'$, otherwise $\mathcal{F}_2=\mathcal{F}_1$.
In this way, we get $\mathcal{F}\supset \mathcal{F}_1\supset \mathcal{F}_2\supset \cdots$.

We define an infinite sequence $f_1,f_2,f_3,\ldots$ of distinct elements in $\mathcal{F}$ so that for each $l\geq 1$, we have 
\begin{equation}\label{eqn:202206254}
f_k\in \mathcal{F}_l
\end{equation}
for all $k\geq l$.
This sequence is constructed inductively as follows.
We take $f_1\in \mathcal{F}_1$.
Suppose distinct elements $f_1,f_2,\ldots,f_n$ are chosen so that \eqref{eqn:202206254} holds for all $k,l\in \mathbb N$ with $1\leq l\leq k\leq n$.
Then we choose $f_{n+1}\in \mathcal{F}_{n+1}-\{ f_1,\ldots,f_n\}$.
Thus we have constructed the sequence $f_1,f_2,\ldots$ with the desired property \eqref{eqn:202206254}.

Now we set $\mathcal{F}_o=\{ f_1,f_2,\ldots\}$.
Then $\mathcal{F}_o\subset \mathcal{F}$ is an infinite subset.
Let $\mathcal{G}\subset \mathcal{F}_o$ be an infinite subset and $\bar{A}$ be a smooth equivariant compactification.
We claim $Q(\bar{A},\mathcal{F}_o)=Q(\bar{A},\mathcal{G})$.
To show this, it is enough to show $Q(\bar{A},\mathcal{G})\subset Q(\bar{A},\mathcal{F}_o)$.
Let $B\in Q(\bar{A},\mathcal{G})$.
We may take $k$ such that $(\bar{A},B)=(\hat{A}_k,B_k)$.
Since $\mathcal{G}\cap\mathcal{F}_k$ is infinite, $\mathcal{G}\to \mathrm{Sp}_{B_k}\hat{X}_k$ yields that $\mathcal{F}_k\to \mathrm{Sp}_{B_k}\hat{X}_k$.
Hence $\mathcal{F}_o\to \mathrm{Sp}_{B_k}\hat{X}_k$, so that $B\in Q(\bar{A},\mathcal{F}_o)$.
Hence $Q(\bar{A},\mathcal{G})\subset Q(\bar{A},\mathcal{F}_o)$.
We replace $\mathcal{F}$ by $\mathcal{F}_o$ to conclude the proof.
\hspace{\fill} $\square$

\medskip

In the following, we write $Q(\bar{A})=Q(\bar{A},\mathcal{F})$ if no confusion may occur.

\begin{lem}\label{lem:20210515}
$Q(\bar{A})$ contains finite number of maximal element,
i.e.,
there exist $B_1,\ldots,B_l\in Q(\bar{A})$ such that for all $B\in Q(\bar{A})$, there exists $B_i$ such that $B\subset B_i$.
\end{lem}

{\it Proof.}\
By Lemma \ref{lem:202206261}, there exists a Zariski closed set $Y\subset \bar{X}$ such that for Zariski closed subsets $V\subset \bar{X}$, we have $\mathcal{F}\to V$ if and only if $Y\subset V$.
If $B\in Q(\bar{A})$, then $Y\subset \mathrm{Sp}_B\bar{X}$.
Since $\mathrm{Sp}_B\bar{X}$ is $B$-invariant, we have $B+Y\subset \mathrm{Sp}_B\bar{X}$, hence $B+Y\subset \bar{X}$.
Conversely, if $B\subset A$ is a semi-abelian subvariety such that $B+Y\subset \bar{X}$, then we have $Y\subset \mathrm{Sp}_B\bar{X}$.
Hence we have $\mathcal{F}\to \mathrm{Sp}_B\bar{X}$, so $B\in Q(\bar{A})$.
Thus $B\in Q(\bar{A})$ if and only if $B+Y\subset \bar{X}$.

Let $\varphi:A\times Y\to \bar{A}$ be the restriction of the $A$-action $A\times \bar{A}\to \bar{A}$.
Then $\varphi^{-1}(\bar{X})\subset A\times Y$ is Zariski closed.
We define $Z\subset A$ to be the set of $a\in A$ such that $\{ a\}\times Y\subset \varphi^{-1}(\bar{X})$.
Then $Z\subset A$ is Zariski closed.
We have 
\begin{equation}\label{eqn:202206242}
B\in Q(\bar{A})\iff B\subset Z.
\end{equation}

We set $\mathcal{U}=\mathrm{Lie}(A)-\{ 0\}$.
Let $\phi:\mathbb C\times \mathcal{U}\to A$ be an analytic map defined by $\phi(z,u)=\mathrm{exp}(zu)$.
For $t\in \mathbb C$, we define $\iota_t:\mathcal{U}\to \mathbb C\times\mathcal{U}$ by $\iota_t(u)=(t,u)$.
We set
$$
\mathcal{Z}=\bigcap_{t\in\mathbb C}(\phi\circ\iota_t)^{-1}(Z).
$$
Then $\mathcal{Z}\subset \mathcal{U}$ is an analytic subset such that 
\begin{equation}\label{eqn:20220624}
u\in \mathcal{Z} \iff \phi_u(\mathbb C)\subset Z,
\end{equation}
where $\phi_u:\mathbb C\to A$ is a one parameter group defined by $\phi_{u}(z)=\exp (zu)$.
For each $b\in\mathbb C^*$, we have $\phi_{bu}(z)=\phi_u(bz)$.
Hence $\mathcal{Z}\subset \mathcal{U}$ is invariant under the natural $\mathbb C^*$-action on $\mathcal{U}$.
Thus $\mathcal{Z}/\mathbb C^*$ is an analytic subset of $\mathcal{U}/\mathbb C^*=\mathbb P(\mathrm{Lie}(A))$.
Indeed let $\Omega\subset \mathbb P(\mathrm{Lie}(A))$ be an open set with a local analytic section $\sigma:\Omega \to \mathcal{U}$ of the natural projection $\mathcal{U}\to \mathbb P(\mathrm{Lie}(A))$.
Then we have $\sigma^{-1}(\mathcal{Z})=(\mathcal{Z}/\mathbb C^*)\cap \Omega$.
Since $\sigma^{-1}(\mathcal{Z})\subset \Omega$ is an analytic subset, $\mathcal{Z}/\mathbb C^*$ is an analytic subset of $\mathbb P(\mathrm{Lie}(A))$.
Hence by Chow's lemma, $\mathcal{Z}/\mathbb C^*$ is an algebraic subset of $\mathbb P(\mathrm{Lie}(A))$.

For $[u]\in \mathbb P(\mathrm{Lie}(A))$, where $u\in \mathrm{Lie}(A)-\{0\}$, let $B_{[u]}$ be the Zariski closure of $\phi_u(\mathbb C)\subset A$.
Then $B_{[u]}\subset A$ is a semi-abelian subvariety.
Then by \eqref{eqn:20220624}, we have 
\begin{equation}\label{eqn:202206244}
[u]\in \mathcal{Z}/\mathbb C^* \iff B_{[u]}\subset Z.
\end{equation} 
Hence $[u]\in \mathcal{Z}/\mathbb C^*$ if and only if $\mathbb P(\mathrm{Lie}(B_{[u]}))\subset \mathcal{Z}/\mathbb C^*$.
Thus $\mathcal{Z}/\mathbb C^*=\cup_{[u]\in \mathcal{Z}/\mathbb C^*}\mathbb P(\mathrm{Lie}(B_{[u]}))$.
Note that there are only countably many semi-abelian subvarieties of $A$ (cf. \cite[Cor. 5.1.9]{NW}).
Hence there exist $[u_1],\ldots,[u_l]\in \mathcal{Z}/\mathbb C^*$ such that $\mathcal{Z}/\mathbb C^*=\mathbb P(\mathrm{Lie}(B_{[u_1]}))\cup \cdots\cup\mathbb P(\mathrm{Lie}(B_{[u_l]}))$.

Now if $B\in Q(\bar{A})$, then $B\subset Z$ (cf. \eqref{eqn:202206242}).
Hence we have
$\mathbb P(\mathrm{Lie}(B))\subset \mathcal{Z}/\mathbb C^*$ (cf. \eqref{eqn:202206244}).
Hence there exists $B_{[u_i]}$ such that $B\subset B_{[u_i]}$.
This conclude the proof.
\hspace{\fill} $\square$

\medskip

Let $\Sigma$ be the set of all semi-abelian subvarieties of $A$.
Let $\mathcal{Q}$ be the set of all subset $Q\subset \Sigma$ such that
\begin{enumerate}
\item $Q$ contains finite number of maximal elements, and
\item If $B,B'\in \Sigma$ such that $B\in Q$ and $B'\subset B$, then $B'\in Q$.
\end{enumerate}
Then $\mathcal{Q}$ satisfies the following

\begin{lem}\label{lem:202206251}
Let $Q_1\supset Q_2\supset Q_3\supset\cdots$ be a descending sequence in $\mathcal{Q}$.
Then there exists $k\geq 1$ such that $Q_k=Q_{k+1}=Q_{k+2}=\cdots$.
\end{lem}

{\it Proof.}\
For each $Q\in\mathcal{Q}$, we denote by $Q_{\mathrm{max}}\subset Q$ the set of the maximal elements in $Q$.
Then $Q_{\mathrm{max}}\subset Q$ is a finite subset.
We set
$$
P_{Q}=\bigcup_{B\in Q_{\mathrm{max}}}B.
$$
Then $P_Q\subset A$ is a Zariski closed set.
We claim that $Q\subset Q'$ if and only if $P_{Q}\subset P_{Q'}$.
Indeed, assume $Q\subset Q'$.
Then for $B\in Q_{\mathrm{max}}$, we have $B\in Q'$, hence there exists $B'\in Q'_{\mathrm{max}}$ such that $B\subset B'$.
This shows $B\subset P_{Q'}$, hence $P_{Q}\subset P_{Q'}$.
Conversely, suppose $P_{Q}\subset P_{Q'}$.
Let $B\in Q$.
Then $B\subset P_{Q}\subset P_{Q'}$.
Hence there exists $B'\in Q'_{\mathrm{max}}$ such that $B\subset B'$.
Thus by the property (2) for $\mathcal{Q}$, we have $B\in Q'$.
Hence $Q\subset Q'$.
Thus we have proved that $Q\subset Q'$ if and only if $P_{Q}\subset P_{Q'}$.

Now let $Q_1\supset Q_2\supset Q_3\supset\cdots$ be a descending sequence in $\mathcal{Q}$.
Then we have $P_{Q_1}\supset P_{Q_2}\supset P_{Q_3}\supset \cdots$.
Thus by the Noetherian property, there exists $k$ such that $P_{Q_k}=P_{Q_{k+1}}=P_{Q_{k+2}}=\cdots$.
Then we have $Q_k=Q_{k+1}=Q_{k+2}=\cdots$.
The proof is completed.
\hspace{\fill} $\square$

\medskip

For each equivariant compactification $\bar{A}$, we have $Q(\bar{A})\in\mathcal{Q}$.
Indeed, by Lemma \ref{lem:20210515}, $Q(\bar{A})$ contains finite number of maximal element.
If $B\in Q(\bar{A})$ and $B'\subset B$, then by $\mathrm{Sp}_B\bar{X}\subset \mathrm{Sp}_{B'}\bar{X}$, we have $\mathcal{F}\to \mathrm{Sp}_{B'}\bar{X}$, hence $B'\in Q(\bar{A})$.
Thus $Q(\bar{A})\in\mathcal{Q}$.
For an equivariant compactification $\hat{A}\to\bar{A}$, we have $Q(\hat{A})\subset Q(\bar{A})$.
Indeed the image of $\mathrm{Sp}_B(\hat{X})\subset \hat{A}$ under the map $\hat{A}\to\bar{A}$ is contained in $\mathrm{Sp}_B(\bar{X})\subset \bar{A}$, which implies that if $B\in Q(\hat{A})$, then $B\in Q(\bar{A})$.

\begin{lem}\label{lem:202206256}
There exists an equivariant compactification $\bar{A}$ such that for  every equivariant compactification $\hat{A}\to\bar{A}$, we have $Q(\bar{A})=Q(\hat{A})$.
Moreover we may take $\bar{A}$ to be smooth.
\end{lem}

{\it Proof.}\
Assume contrary to get a sequence $\bar{A}_1\leftarrow \bar{A}_2\leftarrow\cdots$ such that $Q(\bar{A}_1)\supsetneqq Q(\bar{A}_2)\supsetneqq\cdots$.
But this does not occur from Lemma \ref{lem:202206251}.
Hence there exists $\bar{A}$ of desired property.
If our $\bar{A}$ is not smooth, we may replace $\bar{A}$ by its smooth equivariant blow-up.
This exists by Lemma \ref{lem:20211030111}.
\hspace{\fill} $\square$

\medskip

We reduce Theorem \ref{thm:mstrong} to the following equivalent statement.

\begin{thm}\label{thm:mstrong2}
Let $A$ be a semi-abelian variety.
Let $X\subsetneqq A$ be a closed subvariety.
Let $\mathcal{F}\subset \mathrm{Hol}(\mathbb D,X)$ be an infinite set of holomorphic maps.
Then there exist a proper semi-abelian subvariety $B\subsetneqq A$ and an infinite subset $\mathcal{G}\subset \mathcal{F}$ with the following two properties:
\begin{enumerate}
\item $(\varpi\circ f)_{f\in \mathcal{G}}$ converges uniformly on compact subsets of $\mathbb D$ to $g:\mathbb D\to A/B$, where $\varpi:A\to A/B$ is the quotient map.
\item Let $\overline{A}$ be an equivariant compactification and let $\overline{X}\subset \overline{A}$ be the compactification.
Then $\mathcal{G}\to \mathrm{Sp}_B\overline{X}$.
\end{enumerate}
\end{thm}

{\it Proof.}\
By Lemma \ref{lem:202206254}, replacing $\mathcal{F}$ by its infinite subset, we may assume that 
\begin{equation}\label{eqn:202206255}
Q(\tilde{A},\mathcal{F})=Q(\tilde{A},\mathcal{G})
\end{equation}
for every smooth equivariant compactification $A\subset \tilde{A}$ and every infinite subset $\mathcal{G}\subset \mathcal{F}$.
By Lemma \ref{lem:202206256}, we may take a smooth equivariant compactification $\bar{A}$ such that
\begin{equation}\label{eqn:202206256}
Q(\bar{A},\mathcal{F})=Q(\hat{A},\mathcal{F})
\end{equation}
for all $\hat{A}\to \bar{A}$.

We claim that for every equivariant compactification $\tilde{A}$ and infinite subset $\mathcal{G}\subset \mathcal{F}$, we have
\begin{equation}\label{eqn:202206272}
Q(\bar{A},\mathcal{G})\subset Q(\tilde{A},\mathcal{G}).
\end{equation}
To prove this, we apply Lemmas \ref{lem:20211030111} and \ref{lem:202109231} to get a smooth equivariant blow-up $\hat{A}\to \bar{A}$ such that $\hat{A}\to \tilde{A}$ exists.
Then by \eqref{eqn:202206255} and \eqref{eqn:202206256}, we get
$$
Q(\bar{A},\mathcal{G})=
Q(\bar{A},\mathcal{F})=
Q(\hat{A},\mathcal{F})=
Q(\hat{A},\mathcal{G}).
$$
On the other hand, the existence of the map $\hat{A}\to \tilde{A}$ yields $Q(\hat{A},\mathcal{G})\subset Q(\tilde{A},\mathcal{G})$.
Hence we get \eqref{eqn:202206272}.

Now we apply Lemma \ref{thm:m} for $\mathcal{F}$ and $\bar{A}$ to get $B_0\subset A$ and $\mathcal{G}\subset \mathcal{F}$.
Then we have $\mathcal{G}\to \mathrm{Sp}_{B_0}\bar{X}$, hence $B_0\in Q(\bar{A},\mathcal{G})$.
We choose an equivariant blow-up $\varphi:\hat{A}\to \bar{A}$ and an equivariant compactification $\overline{A/B_0}$ such that $p:\hat{A}\to \overline{A/B_0}$ exists and satisfies Assumption \ref{ass:20220626} (cf. Lemma \ref{lem:20210513}).
We may assume that $\overline{A/B_0}$ is smooth.
By \eqref{eqn:202206272}, we have $B_0\in Q(\hat{A},\mathcal{G})$.
Thus $\mathcal{G}\to \mathrm{Sp}_{B_0}(\hat{X})$.
By Lemma \ref{thm:m}, we get an infinite subset $\mathcal{H}\subset \mathcal{G}$ with a limit $g:\mathbb D\to \overline{A/B_0}$ of $(p\circ f)_{f\in \mathcal{H}}$.

We first consider the case $g(\mathbb D)\subset A/B_0$.
We take arbitrary $\tilde{A}$.
Then by \eqref{eqn:202206272}, we have $B_0\in Q(\tilde{A},\mathcal{G})$.
Hence $\mathcal{H}\to \mathrm{Sp}_{B_0}\tilde{X}$.
Note that $B_0$ satisfies $\mathrm{Sp}_{B_0}(\tilde{X})\not=\emptyset$ for all $\tilde{A}$.
Hence by Lemma \ref{lem:20200910}, we have $B_0\subsetneqq A$.
We replace $\mathcal{G}$ by $\mathcal{H}$ and set $B=B_0$ to conclude the proof in the case $g(\mathbb D)\subset A/B_0$.

Next we assume $g(\mathbb D)\not\subset A/B_0$.
Then we have $g(\mathbb D)\subset \partial (A/B_0)$.
Then by Lemma \ref{lem:202206263}, we get $B$ with $B_0\subset B\subset A$ such that $\{ \varpi\circ f\}_{f\in \mathcal{H}}$ converges uniformly on compact subsets of $\mathbb D$ to $g:\mathbb D\to A/B$ and that $\mathcal{H}\to \mathrm{Sp}_{B}(\hat{X})$.
Hence the existence of $\hat{A}\to\bar{A}$ shows that $\mathcal{H}\to \mathrm{Sp}_{B}\bar{X}$.
We take arbitrary $\tilde{A}$.
Then by \eqref{eqn:202206272}, we have $B\in Q(\tilde{A},\mathcal{H})$.
Hence $\mathcal{H}\to \mathrm{Sp}_{B}\tilde{X}$.
Our $B$ satisfies $\mathrm{Sp}_{B}(\tilde{X})\not=\emptyset$ for all $\tilde{A}$.
Hence by Lemma \ref{lem:20200910}, we have $B\subsetneqq A$.
We replace $\mathcal{G}$ by $\mathcal{H}$ to conclude the proof.
\hspace{\fill} $\square$

\medskip

{\it Proof of Theorem \ref{thm:mstrong}.}\
Let $(f_n)_{n\in\mathbb N}$ be a sequence in $\mathrm{Hol}(\mathbb D,X)$.
We may assume that the set $\mathcal{F}=\{ f_n;\ n\in\mathbb N\}\subset \mathrm{Hol}(\mathbb D,X)$ is infinite, for otherwise, we may take a subsequence consists of the same maps, hence the theorem is trivially valid for $B=\{ 0\}$.
Then we apply Theorem \ref{thm:mstrong2} to get $B$ and $\mathcal{G}\subset \mathcal{F}$.
To conclude the proof, we take a subsequence $(f_{n_k})_{k\in\mathbb N}$ of distinct elements in $\mathcal{G}$.
\hspace{\fill} $\square$

\section{Proof of Theorem \ref{thm:pkob}}\label{sec:11}

In this section, we prove Theorem \ref{thm:pkob}.
We recall the Zariski closed subset $S\subset \overline{X}$ from the statement of Theorem \ref{thm:pkob}.

\begin{lem}\label{lem:202105151}
Let $X\subsetneqq A$ be a closed algebraic subvariety.
Let $\overline{A}$ be a smooth equivariant compactification and let $\overline{X}\subset \overline{A}$ be the compactification.
Let $(f_i)_{i\in I}$ be an infinite indexed family of holomorphic maps in $\mathrm{Hol}(\mathbb D,X)$ such that $(f_i)_{i\in I}\not\to S$.
Then there exists an infinite subset $J\subset I$ such that $(f_{i})_{i\in J}$ converges uniformly on compact subsets of $\mathbb D$ to some $h:\mathbb D\to \overline{X}$. 
\end{lem}

In the following proof, we recall $C_H^1$ from \eqref{20220604}.

{\it Proof.}\
Set $\mathcal{F}=\{ f_i;\ i\in I\}\subset \mathrm{Hol}(\mathbb D,X)$.
Then we may assume that the map $I\to\mathcal{F}$ defined by $i\mapsto f_i$ is a finite-to-one mapping, for otherwise we may choose an infinite indexed subfamily which consists of the same elements in $\mathrm{Hol}(\mathbb D,X)$.
In particular, $\mathcal{F}$ is infinite.
By $(f_i)_{i\in I}\not\to S$, we have $\mathcal{F}\not\to S$.
We shall show the existence of an infinite subset $\mathcal{G}\subset \mathcal{F}$ which converges uniformly on compact subsets of $\mathbb D$ to some $h:\mathbb D\to \overline{X}$. 

By replacing $\mathcal{F}$ by its infinite subset, we may assume that $\mathcal{F}'\not\to S$ for all infinite subset $\mathcal{F}'\subset \mathcal{F}$ (cf. Lemma \ref{lem:20210803}).
We take $B\subset A$ and $\mathcal{G}\subset \mathcal{F}$ as in Theorem \ref{thm:mstrong2}.
If $B=\{ 0\}$, then our lemma is valid.
Hence we assume $\dim B\geq 1$.
Let $Z\subset \partial X$ be the locus where $B$ acts trivially.
Then $\mathrm{Sp}_B\overline{X}\subset S\cup Z$, hence $\mathcal{G}\to S\cup Z$.
Let $\mathcal{Z}$ be the set of all Zariski closed sets $Z'\subset Z$ such that $\mathcal{G}\to S\cup Z'$.
By $\mathcal{G}\to S\cup Z$, we have $Z\in \mathcal{Z}$, hence $\mathcal{Z}\not=\emptyset$. 
By the Noetherian property, we may take a minimum element $Z'\in\mathcal{Z}$.
By $\mathcal{G}\not\to S$, we have $Z'\not=\emptyset$.
Let $T$ be an irreducible component of $Z'$.
Set $S_T=S\cup \overline{(Z'\backslash T)}$.
By $\overline{(Z'\backslash T)}\subsetneqq Z'$, we have $\overline{(Z'\backslash T)}\not\in\mathcal{Z}$.
Hence
\begin{equation}\label{eqn:202110041}
\mathcal{G}\not\to S_T.
\end{equation}

Let $C\subset A$ be the isotropy group of the generic points of $T$.
Then $B\subset C$.
Let $V\subset \partial A$ be the Zariski closure of $A+T\subset \partial A$.
Then $C$ is the isotropy group of $V$.
Since $T$ is irreducible, $V$ is irreducible.
Let 
$$p:W\to V$$ 
be the vector bundle described in Lemma \ref{lem:vect} so that $V=\overline{A/C}$, where $V\subset W$ is considered as 0-section.
Note that
$$
S\cup Z'\subset S_T\cup T\subset  S_T\cup V.
$$
By $\mathcal{G}\to S\cup Z'$, we have
\begin{equation}\label{eqn:20211004}
\mathcal{G}\to S_T\cup V.
\end{equation}
From the limit map $g_0:\mathbb D\to A/B$, we get $g:\mathbb D\to A/C\subset V$, which is the limit of $(p\circ f)_{f\in \mathcal{G}}$.

We are going to prove $\mathcal{G}\to V$.
We first show $g(\mathbb D)\not\subset S_T\cap V$.
To prove this, we assume contrary that $g(\mathbb D)\subset S_T\cap V$.
Let $W\subset \overline{W}$ be a compactification such that $\varphi:\overline{W}\to \bar{A}$ and $\bar{p}:\overline{W}\to V$ exists.
By \eqref{eqn:20211004}, we have $\mathcal{G}\to \varphi^{-1}(S_T)\cup V$.
Since we are assuming $g(\mathbb D)\subset S_T\cap V$, we have $(p\circ f)_{f\in \mathcal{G}}\to S_T\cap V$ (cf. Remark \ref{rem:20220626}).
Hence by Lemma \ref{lem:202206260}, we have $\mathcal{G}\to \bar{p}^{-1}(S_T\cap V)$.
Hence $\mathcal{G}\to (\varphi^{-1}(S_T)\cup V)\cap \bar{p}^{-1}(S_T\cap V)$ (cf. Lemma \ref{lem:20220626}).
By 
$$
V\cap \bar{p}^{-1}(S_T\cap V)=S_T\cap V,
$$
we have
\begin{equation*}
\begin{split}
(\varphi^{-1}(S_T)\cup V)\cap \bar{p}^{-1}(S_T\cap V)
&=(\varphi^{-1}(S_T)\cap \bar{p}^{-1}(S_T\cap V))\cup (V\cap \bar{p}^{-1}(S_T\cap V))\\
&\subset 
\varphi^{-1}(S_T).
\end{split}
\end{equation*}
Hence $\mathcal{G}\to \varphi^{-1}(S_T)$, so $\mathcal{G}\to S_T$.
This contradicts to \eqref{eqn:202110041}.
Hence $g(\mathbb D)\not\subset S_T\cap V$.

Now we replace $\mathcal{G}$ by its infinite subset so that $\mathcal{G}'\not\to S_T$ for all infinite subset $\mathcal{G}'\subset \mathcal{G}$ (cf. Lemma \ref{lem:20210803}).
We take an open neighbourhood $S_T\subset U_0$ and positive constants $s_0\in (0,1)$ and $\gamma_0>0$ such that for all $f\in \mathcal{G}\backslash \mathcal{E}_1$, where $\mathcal{E}_1\subset \mathcal{G}$ is a finite subset, we have
\begin{equation}\label{eqn:202110043}
C_H^1(\mathbb D(s_0)\backslash f^{-1}(U_0))\geq  \gamma_0.
\end{equation}
To prove $\mathcal{G}\to V$, we take an open neighbourhood $V\subset U_1$ and positive constants $s\in (s_0,1)$ and $\gamma\in (0,\gamma_0)$ arbitrary.
We denote by $K\subset \mathbb D$ the finite union of small open discs centered at the points of the finite set $g^{-1}(S_T\cap V)\cap \overline{\mathbb D(\frac{s+1}{2})}$ so that $g^{-1}(S_T\cap V)\subset K$ and
\begin{equation}\label{eqn:202110045}
C_H^1(K)<\gamma/4.
\end{equation}
We take an open neighbourhood $U_3\subset W$ of $S_T\cap V$ such that $U_3\subset p^{-1}(U_3\cap V)$ and $g(\overline{\mathbb D(\frac{s+1}{2})}\backslash K)\subset V\backslash \overline{U_3\cap V}$.
Since $(p\circ f)_{f\in \mathcal{G}}$ converges uniformly on $\overline{\mathbb D(\frac{s+1}{2})}$ to $g$, we have $p\circ f(\overline{\mathbb D(\frac{s+1}{2})}\backslash K)\subset V\backslash \overline{U_3\cap V}$ for all $f\in \mathcal{G}\backslash \mathcal{E}_2$, where $\mathcal{E}_2\subset \mathcal{G}$ is a finite subset.
Hence by $U_3\subset p^{-1}(U_3\cap V)$, we have $f(\overline{\mathbb D(\frac{s+1}{2})}\backslash K)\subset \overline{A}\backslash \overline{U_3}$ for all $f\in \mathcal{G}\backslash \mathcal{E}_2$.
By replacing $U_0$ and $U_1$ by smaller open neighbourhoods of $S_T$ and $V$, respectively, we may assume that $U_0\cap U_1\subset U_3$.
Hence we have $f(\mathbb D((s+1)/2)\backslash K)\subset \overline{A}\backslash (\overline{U_0\cap U_1})$, so $\mathbb D((s+1)/2)\cap f^{-1}(\overline{U_0\cap U_1})\subset K$.
Hence by \eqref{eqn:202110045}, we have
$$
C_H^1(\mathbb D((s+1)/2)\cap f^{-1}(\overline{U_0\cap U_1}))<\gamma/4
$$
for all $f\in \mathcal{G}\backslash \mathcal{E}_2$.
Now by \eqref{eqn:20211004}, for all $f\in \mathcal{G}\backslash\mathcal{E}_3$, where $\mathcal{E}_3\subset \mathcal{G}$ is a finite subset, we have
\begin{equation*}
C_H^1(\mathbb D((s+1)/2)\backslash f^{-1}(U_0\cup U_1))< \gamma/4.
\end{equation*}
Then for all $f\in\mathcal{G}\backslash (\mathcal{E}_2\cup\mathcal{E}_3)$, we have 
\begin{equation*}
C_H^1(\mathbb D((s+1)/2)\backslash f^{-1}(U_0'\cup U_1'))<\gamma/2,
\end{equation*}
where $U_0'=U_0\backslash \overline{U_0\cap U_1}$ and $U_1'=U_1\backslash \overline{U_0\cap U_1}$.
By Lemma \ref{lem:20220529}, there exists an open set $\Omega_f\subset \mathbb D(s)\cap f^{-1}(U_0'\cup U_1')$ such that $\Omega_f$ is connected and 
$$C_H^1(\mathbb D(s)\backslash\Omega_f)<\gamma.$$
Note that $U_0'\cap U_1'=\emptyset$.
Since $f(\Omega_f)$ is connected, we have either $f(\Omega_f)\subset U_0'$ or $f(\Omega_f)\subset U_1'$.
By $\gamma_0>\gamma$ and \eqref{eqn:202110043}, we have $f(\Omega_f)\not\subset U_0$ for all $f\in\mathcal{G}\backslash(\mathcal{E}_1\cup\mathcal{E}_2\cup\mathcal{E}_3)$.
Hence $f(\Omega_f)\subset U_1'\subset U_1$.
Thus $\mathcal{G}\to V$.

Now by Lemma \ref{lem:20210303}, $\mathcal{G}$ converges uniformly on compact subsets of $\mathbb D$ to $g:\mathbb D\to V\cap \bar{X}$.
We take an infinite subset $J\subset I$ such that $f_i\in\mathcal{G}$ for all $i\in J$.
Then $(f_i)_{i\in J}$ converges uniformly on compact subsets of $\mathbb D$ to $g$.
This completes the proof of the lemma.
\hspace{\fill} $\square$

\medskip

{\it Proof of Theorem \ref{thm:pkob}.}\
If $X=A$, then $S=A$.
Hence our theorem is trivial.
In the following, we assume $X\subsetneqq A$.

Let $(f_n)_{n\in\mathbb N}$ be a sequence of holomorphic maps in $\mathrm{Hol}(\mathbb D,X)$.
We assume that the condition (2) of the definition of `tautly imbedded modulo' is not satisfied.
Then there exist compact sets $K\subset \mathbb D$ and $L\subset \overline{X}-S$ such that, by replacing $(f_n)_{n\in\mathbb N}$ by its subsequence, we have $f_n(K)\cap L\not=\emptyset$ for all $n\in\mathbb N$.

We shall show $(f_n)_{n\in\mathbb N}\not\to S$.
To prove this, we assume, contrary, that $(f_n)_{n\in\mathbb N}\to S$.
We take open neighbourhood $U\subset \bar{X}$ of $S$ such that $\bar{U}\cap L=\emptyset$.
For each $n\in \mathbb N$, there exists $a_n\in K$ such that $f_n(a_n)\in L$, which follows from $f_n(K)\cap L\not=\emptyset$.
We choose an automorphism $\varphi_n:\mathbb D\to\mathbb D$ so that $\varphi_n(0)=a_n$.
Then $f_n\circ\varphi_n(0)\in L$.
In particular, $f_n\circ\varphi_n(0)\not\in \bar{U}$.
We define $\delta_n>0$ by
\begin{equation*}
\delta_n=\sup\{ t\in (0,1);\ f_n\circ\varphi_n(\overline{\mathbb D(t)})\subset \bar{X}- \bar{U}\}.
\end{equation*}
Then we have
\begin{equation}\label{eqn:202205211}
f_n\circ\varphi_n(\mathbb D(\delta_n))\cap U=\emptyset
\end{equation}
and
\begin{equation}\label{eqn:202205212}
f_n\circ\varphi_n(\partial \mathbb D(\delta_n))\cap \bar{U}\not=\emptyset.
\end{equation}
Since we are assuming $(f_n)_{n\in\mathbb N}\to S$, we have
\begin{equation}\label{eqn:20220521}
\lim_{n\to\infty}\delta_n=0.
\end{equation}
Set $r_n=\frac{1}{2\delta_n}$.
We define $g_n:\mathbb D(r_n)\to X$ by $g_n(z)=f_n\circ\varphi_n(2\delta_nz)$.
Here we continue to write $\mathbb D(r)=\{z\in \mathbb C; |z|<r\}$ for $r\in\mathbb R_{>0}$ without assuming $r<1$.

We define a descending sequence of infinite subsets $\mathbb N\supset I_2\supset I_2\supset \cdots$ inductively as follows.
Set $I_0=\mathbb N$.
Suppose $I_{l-1}$ is defined.
By \eqref{eqn:20220521}, we may take an infinite subset $I_{l}'\subset I_{l-1}$ such that $\delta_i<\frac{1}{2l}$ for all $i\in I_l'$.
Hence $r_i>l$ for all $i\in I_l'$.
Hence for each $i\in I_l'$, we have the restriction $g_i|_{\mathbb D(l)}:\mathbb D(l)\to X$.
We get an infinite indexed family $(g_i|_{\mathbb D(l)})_{i\in I_l'}$ in $\mathrm{Hol}(\mathbb D(l),X)$.
By \eqref{eqn:202205211}, we have $g_i(\mathbb D(1/2))\cap U=\emptyset$ for all $i\in I_l'$.
Hence $(g_i|_{\mathbb D(l)})_{i\in I_l'}\not\to S$, where we consider $\mathbb D(l)$ as '$\mathbb D$'.
We apply Lemma \ref{lem:202105151}.
Then we get an infinite subset $I_l\subset I_l'$ such that $(g_i|_{\mathbb D(l)})_{i\in I_l}$ converges uniformly on compact subsets of $\mathbb D(l)$ to $h_l:\mathbb D(l)\to \bar{X}$.
Since $h_{l-1}$ is the limit of $(g_i|_{\mathbb D(l-1)})_{i\in I_l}$, which is a subfamily of $(g_i|_{\mathbb D(l-1)})_{i\in I_{l-1}}$, we have $h_l|_{\mathbb D(l-1)}=h_{l-1}$.

Now we get a holomorphic map $h:\mathbb C\to \bar{X}$ such that $h|_{\mathbb D(l)}=h_l$ for all $l\in\mathbb N$.
We shall show that $h$ is non-constant.
Note that $h_1:\mathbb D\to \bar{X}$ is the limit of $(g_i|_{\mathbb D})_{i\in I_1}$.
By \eqref{eqn:202205212}, we have $g_i(\partial\mathbb D(1/2))\cap \bar{U}\not=\emptyset$ for all $i\in I_1$.
Thus $h_1(\partial \mathbb D(1/2))\cap \overline{U}\not=\emptyset$.
On the other hand, by $g_i(0)\in L$ for all $i\in I_1$, we have $h_1(0)\in L$.
By $L\cap \bar{U}=\emptyset$, $h_1$ is non-constant.
Hence $h$ is non-constant.

For each irreducible component $D\subset \partial A$, we have either $h(\mathbb C)\subset D$ or $h(\mathbb C)\cap D=\emptyset$.
Thus there exists an $A$-orbit $O\subset \bar{A}$ such that $h(\mathbb C)\subset O$.
Note that $O$ is isomorphic to a semi-abelian variety.
Hence the Zariski closure of $h(\mathbb C)$ is a translate of a semi-abelian subvariety of $O$ (cf. Theorem \ref{thm:lbo}, or \cite[Thm 3.9.19]{Kbook}).
Hence $h(\mathbb C)\subset S$.
This contradicts to $h(0)=h_1(0)\in L$.
Hence $(f_n)_{n\in\mathbb N}\not\to S$.

Now by Lemma \ref{lem:202105151}, there exists a subsequence $(f_{n_k})_{k\in\mathbb N}$ which converges uniformly on compact subsets of $\mathbb D$ to $f:\mathbb D\to \overline{X}$.
Thus we have proved that $X$ is tautly imbedded modulo $S$ in $\bar{X}$.
\hspace{\fill} $\square$

\section{Proof of Theorem \ref{thm:20220429}}
\label{sec:20220429}

When $A$ is compact, we may eliminate the term `$\gamma$-almost' from the statement (2) in Theorem \ref{thm:mstrong} to get the following.

\begin{cor}\label{cor:20220729}
Let $A$ be an abelian variety.
Let $X\subsetneqq A$ be a proper closed algebraic subvariety.
Let $(f_n)_{n\in \mathbb N}$  be a sequence of holomorphic maps in $\mathrm{Hol}(\mathbb D,X)$.
Then there exist a proper abelian subvariety $B\subsetneqq A$ and a subsequence $(f_{n_k})_{k\in\mathbb N}$ with the following two properties:
\begin{enumerate}
\item $(\varpi\circ f_{n_k})_{k\in \mathbb N}$ converges uniformly on compact subsets of $\mathbb D$ to a holomorphic map $g:\mathbb D\to A/B$, where $\varpi:A\to A/B$ is the quotient map.
\item For every $0<s<1$ and open neighbourhood $U\subset A$ of $\mathrm{Sp}_BX$, there exists $k_0\in \mathbb N$ such that, for all $k\geq k_0$, we have $f_{n_k}(z)\in U$ for all $z\in \mathbb D(s)$.
\end{enumerate}
\end{cor}

{\it Proof.}\
We apply Theorem \ref{thm:mstrong} to get $B\subsetneqq A$ and $(f_{n_k})_{k\in\mathbb N}$.
The assertion (1) follows from that of Theorem \ref{thm:mstrong}.
To prove the assertion (2), we set $Z=\varpi (\mathrm{Sp}_BX)$.
Then since each fiber of $\varpi:A\to A/B$ consists of one $B$-orbit, we have
$$\varpi^{-1}(Z)=\mathrm{Sp}_BX.$$
Let $U\subset A$ be an open neighbourhood of $\mathrm{Sp}_BX$.
Then since $A\backslash U$ is compact, $\varpi(A\backslash U)\subset A/B$ is a closed set.
We have $Z\cap \varpi(A\backslash U)=\emptyset$.
Note that $Z\subset A/B$ is a Zariski closed set.
Hence there exists an open neighbourhood $W\subset A/B$ of $Z$ such that $W\cap \varpi(A\backslash U)=\emptyset$.
Then $\varpi^{-1}(W)\subset U$.

Now let $s\in (0,1)$.
By $( f_{n_k})_{k\in\mathbb N}\to \mathrm{Sp}_BX$, we have $(\varpi\circ f_{n_k})_{k\in\mathbb N}\to Z$.
Hence $g(\mathbb D)\subset Z$ (cf. Remark \ref{rem:20220626}).
Hence there exists $k_0\in\mathbb N$ such that $\varpi\circ f_{n_k}(\mathbb D(s))\subset W$ for all $k\geq k_0$.
Hence we have $f_{n_k}(\mathbb D(s))\subset U$ for all $k\geq k_0$.
\hspace{\fill} $\square$

\begin{rem}
When $A$ is compact, Theorem \ref{thm:pkob} is easily derived from Corollary \ref{cor:20220729} as follows.
Let $(f_n)_{n\in\mathbb N}$ be a sequence of holomorphic maps in $\mathrm{Hol}(\mathbb D,X)$.
We assume that the condition (2) of the definition of `tautly imbedded modulo' is not satisfied.
Then there exist compact sets $K\subset \mathbb D$ and $L\subset X-S$ such that, by replacing $(f_n)_{n\in\mathbb N}$ by its subsequence, we have $f_n(K)\cap L\not=\emptyset$ for all $n\in\mathbb N$.
We apply Corollary \ref{cor:20220729} to get $B\subset A$ and $(f_{n_k})_{k\in\mathbb N}$.
To prove $\dim B=0$, we assume contrary $\dim B>0$.
Then $\mathrm{Sp}_BX\subset S$.
Hence the assertion (2) of Corollary \ref{cor:20220729} contradicts to the condition $f_n(K)\cap L\not=\emptyset$ for all $n\in\mathbb N$.
Hence $\dim B=0$.
Then by the assertion (1) of Corollary \ref{cor:20220729}, the sequence $(f_{n_k})_{k\in\mathbb N}$ converges uniformly on compact subsets of $\mathbb D$.
Thus Theorem \ref{thm:pkob} for the compact case is reproved.
\end{rem}

{\it We prove Theorem \ref{thm:20220429}.}
The case $X=A$ is trivial, so we assume $X\subsetneqq A$. 
Let $v\in \breve{T}_xX$ satisfies $F_X(v)=0$, where $x\in X$.
We consider $v\in TA$ by the natural inclusion $\breve{T}X\subset TA$.
There exists a sequence $(f_n)_{n\in\mathbb N}$ in $\mathrm{Hol}(\mathbb D,X)$ such that $f_n'(0)=nv$.
Then $f_n(0)=x$.
We apply Corollary \ref{cor:20220729} to get $B\subsetneqq A$ and $( f_{n_k})_{k\in\mathbb N}$.
We first show $\varpi_*(v)=0$, where $\varpi:A\to A/B$ is the quotient map.
Note that $\varpi\circ f_{n_k}:\mathbb D\to A/B$ converges uniformly on compact subsets of $\mathbb D$ to $g:\mathbb D\to A/B$.
Hence $(\varpi\circ f_{n_k})'(0)=n_k\varpi_*(v)$ converges to $g'(0)$.
Thus we have $\varpi_*(v)=0$.

Now by the assertion (2) of Corollary \ref{cor:20220729}, we have 
$$x\in \mathrm{Sp}_BX$$
(thanks to the absence of the term `$\gamma$-almost' in Corollary \ref{cor:20220729}).
Hence $x+B\subset X$.
Since $\varpi_*(v)=0$, we have $v\in T(x+B)$.
Hence there exists $f:\mathbb C\to (x+B)\subset X$ such that $f'(0)=v$.
The proof is completed.
\hspace{\fill} $\square$

\medskip

The following example shows that the condition $F_X(v)=0$ does not necessarily imply the existence of $f:\mathbb C\to X$ with $f'(0)=v$.

\begin{examp}\label{examp:20220518}
Let $X$ be a smooth surface which is Kobayashi hyperbolic.
Let $\tilde{X}\to X$ be a blowing-up along one point $p\in X$. 
Let $E\subset \tilde{X}$ be the exceptional divisor.
Let $a_1,a_2,a_3\in E$ be three distinct points. 
Then $\tilde{X}-Z$ is Brody hyperbolic, where $Z=\{a_1,a_2,a_3\}$.
Namely there is no non-constant holomorphic map $f:\mathbb C\to \tilde{X}-Z$.
On the other hand, for $q\in E-Z\subset \tilde{X}-Z$ and $v\in T_qE$, we have $F_{\tilde{X}-Z}(v)=0$.
To show this, we take an open neighbourhood $W$ of $E$ which is biholomorphic to
$$
\{ (x,y,[s,t])\in \mathbb C^2\times \mathbb P^1;\ |x|<1,|y|<1, xt=ys\}.
$$
We may assume that $q\in E$ corresponds to $(0,0,[0,1])$.
For $n\in\mathbb N$, we define $f_n\in\mathrm{Hol}(\mathbb D,W)$ by $f_n(z)=(z^3,\frac{1}{n}z^2,[nz,1])$.
Set $v=f_1'(0)$.
Then $v\in T_qE$ and $v\not= 0$.
We have $f_n'(0)=nv$.
We have $f_n(\mathbb D-\{ 0\})\subset W-E$ and $f_n(0)=q\in \tilde{X}-Z$.
Hence $f_n(\mathbb D)\subset \tilde{X}-Z$.
Thus $F_{\tilde{X}-Z}(v)=0$.
\end{examp}

\begin{rem}
We do not know whether the similar statement for Theorem \ref{thm:20220429} holds for non-compact semi-abelian varieties. 
\end{rem}

\section{Proof of Theorem \ref{cor:bc}}\label{sec:12}
We have an equivariant compactification $(\mathbb G_m)^{p-1}\subset \mathbb P^{p-1}$ such that the inclusion is defined by $(g_1,\ldots,g_{p-1})\mapsto [g_1:\ldots:g_{p-1}:1]$.
The action $(\mathbb G_m)^{p-1}\times \mathbb P^{p-1}\to\mathbb P^{p-1}$ is defined by 
\begin{equation}\label{eqn:20230331}
(g_1,\ldots,g_{p-1})\cdot [x_1:\ldots:x_{p-1}:x_p]=[g_1x_1:\ldots:g_{p-1}x_{p-1}:x_p].
\end{equation}
Here $x_1,\ldots,x_p$ are the homogeneous coordinates of $\mathbb P^{p-1}$.
Let $V\subset (\mathbb G_m)^{p-1}$ be a closed subvariety such that $\overline{V}\subset \mathbb P^{p-1}$ is defined by $x_1+x_2+\cdots +x_p=0$.

Let $\mathcal{I}=\{ I_1,\ldots,I_l\}$ be a disjoint partition $I_1\sqcup \cdots\sqcup I_l=\{1,\ldots,p\}$.
We define a subtorus $G_{\mathcal{I}}\subset (\mathbb G_m)^{p-1}$ as follows. 
For $i\in \{1,\ldots,p\}$, we define $\tau(i)\in \{1,\ldots,l\}$ by $i\in I_{\tau(i)}$.
We define a linear subspace $L(\mathcal{I})\subset \mathbb C^p$ by the equations $x_i=x_j$ for all $i,j\in \{1,\ldots,p\}$ such that $\tau(i)=\tau(j)$.
We consider the immersion $(\mathbb G_m)^{p-1}\subset \mathbb C^p$ by $(g_1,\ldots,g_{p-1})\mapsto (g_1,\ldots,g_{p-1},1)$.
We set
$$
G_{\mathcal{I}}=(\mathbb G_m)^{p-1}\cap L(\mathcal{I}).
$$
Then $G_{\mathcal{I}}\subset (\mathbb G_m)^{p-1}$ is a subtorus.

\begin{lem}\label{lem:202206302}
We have
\begin{equation*}
\mathrm{Sp}_{G_{\mathcal{I}}}\overline{V}=\left\{ [x_1:\ldots:x_p]\in \mathbb P^{p-1};\ 
\sum_{i\in I_k}x_i=0\  \text{for all}\  1\leq k\leq l\right\}.
\end{equation*}
\end{lem}

{\it Proof.}\
Note that $\mathrm{Sp}_{G_{\mathcal{I}}}\overline{V}$ is defined by the simultaneous equations
\begin{equation}\label{eqn:20220630}
g_1x_1+\cdots +g_{p-1}x_{p-1}+x_p=0, \quad \text{for all $(g_1,\ldots,g_{p-1})\in G_{\mathcal{I}}$}.
\end{equation}
By changing the indices of $I_1,\ldots,I_l$, we may assume $p\in I_l$.
We have an isomorphism $(\mathbb G_m)^{l-1}\to G_{\mathcal{I}}$ by $(a_1,\ldots,a_{l-1})\mapsto (a_{\tau(1)},\ldots,a_{\tau(p-1)})$, where we set $a_l=1$.
For each $k\in \{1,\ldots,l\}$, we choose $\iota(k)\in \{1,\ldots,p\}$ such that $\iota(k)\in I_{k}$.
Then $\iota$ is a section of $\tau:\{1,2,\ldots,p\}\to \{ 1,\ldots,l\}$.
For $(g_1,\ldots,g_{p-1})\in G_{\mathcal{I}}$, we have
$$
g_1x_1+\cdots +g_{p-1}x_{p-1}+x_p=\sum_{k=1}^{l-1}\left(g_{\iota(k)}\sum_{i\in I_k}x_i\right)+\sum_{i\in I_l}x_i.
$$
Hence \eqref{eqn:20220630} is equivalent to $\sum_{i\in I_k}x_i=0$ for all $k\in \{1,\ldots,l\}$.
The proof is completed.
\hspace{\fill} $\square$

\begin{lem}\label{lem:20220630}
Let $G\subset (\mathbb G_m)^{p-1}$ be a subtorus.
Then there exists a disjoint partition $\mathcal{I}=\{ I_1,\ldots,I_l\}$ of $\{1,2,\ldots,p\}$ such that
\begin{enumerate}
\item $G\subset G_{\mathcal{I}}$, and
\item $\mathrm{Sp}_{G_\mathcal{I}}\overline{V}=\mathrm{Sp}_{G}\overline{V}$.
\end{enumerate}
\end{lem}

{\it Proof.}\
By the inclusion $(\mathbb G_m)^{p-1}\subset \mathbb C^p$, we have $G\subset\mathbb C^p$.
Let $\chi_1,\ldots,\chi_p\in \mathrm{Hom}(G,\mathbb G_m)$ be the composite of $G\subset \mathbb C^{p}$ and the $i$-th projections $ \mathbb C^p\to \mathbb C$, where $\chi_p\equiv 1$.
Then $\chi_1,\ldots\chi_p$ are group homomorphisms.
We define an equivalence relation $\sim$ on $\{1,\ldots,p\}$ such that $i\sim j$ if and only if $\chi_i=\chi_j$.
Thus we get a disjoint partition $\mathcal{I}=\{ I_1,\ldots,I_l\}$ of $\{1,\ldots,p\}$ such that each $I_k$ is equivalence class of the equivalence relation.

Let $L\subset \mathbb C^p$ be the linear subspace spanned by $G\subset \mathbb C^p$.
We claim
\begin{equation}\label{eqn:202206301}
L=L(\mathcal{I}).
\end{equation}

We prove this.
By the definition of $\mathcal{I}$, we have $\chi_i=\chi_j$ if and only if $\tau(i)=\tau(j)$.
Hence we have $G\subset L(\mathcal{I})$.
Hence $L\subset L(\mathcal{I})$.
Let $\iota:\{ 1,\ldots,l\}\to \{1,\ldots,p\}$ be a section of $\tau:\{1,\ldots,p\}\to \{ 1,\ldots,l\}$.
Note that $x_{\iota(1)}|_{L(\mathcal{I})},\ldots,x_{\iota(l)}|_{L(\mathcal{I})}$ form a basis of  the dual space of $L(\mathcal{I})$, where $x_1,\ldots,x_p$ are the coordinate functions of $\mathbb C^p$.
On the other hand, $\{ \chi_{\iota(k)}\}_{k\in \{ 1,\ldots,l\}}\subset \mathrm{Hom}(G,\mathbb G_m)$ is linearly independent (cf. \cite[Lemma 8.1]{Borellinear}).
Hence $x_{\iota(1)}|_{L},\ldots,x_{\iota(l)}|_{L}$ are linearly independent on the dual space of $L$.
Hence $\mathrm{dim}L\geq l=\mathrm{dim}L(\mathcal{I}_G)$.
Thus $L=L(\mathcal{I})$.
Thus we get \eqref{eqn:202206301}.

Now we have $G\subset G_{\mathcal{I}}$, which follows from $G\subset L(\mathcal{I})$.
Hence we have $\mathrm{Sp}_{G_{\mathcal{I}}}\overline{V}\subset \mathrm{Sp}_G\overline{V}$.
It remains to prove $\mathrm{Sp}_G\overline{V}\subset \mathrm{Sp}_{G_{\mathcal{I}}}\overline{V}$.
Note that $\mathrm{Sp}_G\overline{V}$ is defined by simultaneous equations
\begin{equation}\label{eqn:202206302}
b_1x_1+\cdots +b_{p-1}x_{p-1}+x_p=0, \quad \text{for all $(b_1,\ldots,b_{p-1})\in G\subset (\mathbb G_m)^{p-1}$}.
\end{equation}
For $g=(g_1,\ldots,g_{p-1},1)\in G_{\mathcal{I}}\subset \mathbb C^p$,  we have $g\in L(\mathcal{I})$.
Thus by \eqref{eqn:202206301}, we have $g\in L$.
Hence there exist $(b_{1,1},\ldots ,b_{1,p-1},1),\cdots,(b_{s,1},\ldots,b_{s,p-1},1)\in G\subset \mathbb C^p$ and $\alpha_1,\ldots,\alpha_s\in \mathbb C$ such that 
$$
(g_1,\ldots,g_{p-1},1)=\alpha_1(b_{1,1},\ldots ,b_{1,p-1},1)+\cdots+\alpha_s(b_{s,1},\ldots,b_{s,p-1},1).
$$
Hence the solutions of the simultaneous equations \eqref{eqn:202206302} satisfy the equation 
$$g_1x_1+\cdots+g_{p-1}x_{p-1}+x_p=0.$$
Thus $\mathrm{Sp}_G\overline{V}\subset \mathrm{Sp}_{G_{\mathcal{I}}}\overline{V}$.
\hspace{\fill} $\square$

\begin{lem}\label{lem:202206301}
Suppose $G=(\mathbb G_m)^{p-1}$.
Then we have
$\mathrm{Sp}_{G}\overline{V}=\emptyset$.
\end{lem}

{\it Proof.}\
Note that $\mathrm{Sp}_G\bar{V}$ is defined by simultaneous equations
$$
g_1x_1+\cdots +g_{p-1}x_{p-1}+x_p=0, \quad \text{for all $(g_1,\ldots,g_{p-1})\in  (\mathbb G_m)^{p-1}$},
$$
which has only trivial solution.
Hence $\mathrm{Sp}_G\overline{V}=\emptyset$.
\hspace{\fill} $\square$

\medskip

{\it Proof of Theorem \ref{cor:bc}.}\
We consider $\mathbb P^{p-1}$ as an equivariant compactification of $(\mathbb G_m)^{p-1}$, where the action is defined by \eqref{eqn:20230331}.
Let the closed subvariety $V\subset (\mathbb G_m)^{p-1}$ be defined so that $\overline{V}\subset \mathbb P^{p-1}$ becomes
$$
x_1+x_2+\cdots +x_p=0,
$$
where $x_1,\ldots,x_p$ are homogeneous coordinates of $\mathbb P^{p-1}$.
For each $f=(f_1,\ldots,f_p)\in\mathcal{F}$, we consider a holomorphic map $\hat{f}:\mathbb D\to \mathbb P^{p-1}$ defined by $\hat{f}(z)=[f_1(z):\ldots:f_p(z)]$.
Then we have $\hat{f}(\mathbb D)\subset V$, hence $\hat{f}\in\mathrm{Hol}(\mathbb D,V)$. 
Hence $(\hat{f})_{f\in\mathcal{F}}$ is an infinite indexed family in $\mathrm{Hol}(\mathbb D,V)$. 
We replace $\mathcal{F}$ by its countably infinite subset.
Then we may consider $(\hat{f})_{f\in\mathcal{F}}$ as a sequence in $\mathrm{Hol}(\mathbb D,V)$. 
We continue to write this sequence $\mathcal{F}$.

We apply Theorem \ref{thm:mstrong} to get $B\subset (\mathbb G_m)^{p-1}$ and an infinite subset $\mathcal{G}\subset\mathcal{F}$.
Then by $(\hat{f})_{f\in\mathcal{G}}\to \mathrm{Sp}_B\bar{V}$, we have $\mathrm{Sp}_B\bar{V}\not=\emptyset$.
We apply Lemma \ref{lem:20220630} to get a disjoint partition $\mathcal{I}=\{ I_1,\ldots,I_l\}$ of $\{ 1,\ldots,p\}$ such that $B\subset G_{\mathcal{I}}$ and $\mathrm{Sp}_{G_{\mathcal{I}}}\overline{V}=\mathrm{Sp}_{B}\overline{V}$.
Then we have $\mathrm{Sp}_{G_{\mathcal{I}}}\overline{V}\not=\emptyset$.
Hence by Lemma \ref{lem:202206301}, we have $G_{\mathcal{I}}\not=(\mathbb G_m)^{p-1}$.
Since the limit function satisfies $\mathbb D\to (\mathbb G_m)^{p-1}/B$, the sequence $(\hat{f})_{f\in\mathcal{G}}$ converges under the quotient $(\mathbb G_m)^{p-1}/G_{\mathcal{I}}$.
Hence by replacing $B$ by $G_{\mathcal{I}}$, we may assume that $B=G_{\mathcal{I}}$.

Let $I_k$, $1\leq k\leq l$.
By Lemma \ref{lem:202206302}, we have $\mathrm{Sp}_B\overline{V}\subset \{ \sum_{j\in I_k}x_j=0\}$ as Zariski closed subsets of $\mathbb P^{p-1}$.
Hence we have
\begin{equation}\label{eqn:20200917}
(\hat{f})_{f\in\mathcal{G}}\to \left\{ \sum_{j\in I_k}x_j=0\right\}.
\end{equation}
Let $\psi_k:\mathbb P^{p-1}\dashrightarrow \mathbb P^{|I_k|-1}$ be defined by $[x_1:\ldots :x_p]\mapsto [x_j]_{j\in I_k}$.
Here we note that $\mathbb P^0=\mathrm{pt}$, when $|I_k|=1$.
Then $\psi_k$ is regular on $(\mathbb G_m)^{p-1}\subset \mathbb P^{p-1}$.
We have $(\mathbb G_m)^{|I_k|-1}\subset \mathbb P^{|I_k|-1}$.
Note that $\psi_k$ induces a group homomorphism
$$
\psi_k|_{(\mathbb G_m)^{p-1}}:(\mathbb G_m)^{p-1}\to (\mathbb G_m)^{|I_k|-1}.
$$
This $\psi_k|_{(\mathbb G_m)^{p-1}}$ is invariant under the action of $B$ on $(\mathbb G_m)^{p-1}$, hence factors the quotient map $(\mathbb G_m)^{p-1}\to (\mathbb G_m)^{p-1}/B$.
Hence the sequence $(\psi_k\circ \hat{f})_{f\in\mathcal{G}}$ converges uniformly on compact subsets of $\mathbb D$ to $g_k:\mathbb D\to (\mathbb G_m)^{|I_k|-1}$.
Set $H_k=\{ \sum_{j\in I_k}x_j=0\}\subset \mathbb P^{|I_k|-1}$.
We note that $H_k=\emptyset$ if $|I_k|=1$.

\medskip

{\it Claim.}\
Suppose $g_k(\mathbb D)\not\subset H_k$.
Let $\varepsilon\in (0,1)$, $s\in (0,1)$ and $\gamma>0$.
Then there exists a finite subset $\mathcal{E}\subset \mathcal{G}$ such that for all $f\in\mathcal{G}\backslash \mathcal{E}$, we have
\begin{equation}\label{eqn:202110051}
\frac{\sqrt{\sum_{j\in I_k}|f_j(z)|^2}}{\sqrt{\sum_{1\leq i\leq p}|f_i(z)|^2}}<\frac{\varepsilon}{2p}
\end{equation}
for $\frac{\gamma}{2p}$-almost all $z\in\mathbb D(s)$.

\medskip

We prove this.
We define $E_k\subset \mathbb P^{p-1}$ by $E_k=\cap_{i\in I_k}\{ x_i=0\}$.
Then $\psi_k$ induces a regular map $\varphi_k:\mathrm{Bl}_{E_k}\mathbb P^{p-1}\to \mathbb P^{|I_k|-1}$.
Let $E_k'\subset \mathrm{Bl}_{E_k}\mathbb P^{p-1}$ be the exceptional divisor.
Let $\mu:\mathbb P^{p-1}\to\mathbb R_{\geq 0}$ be defined by
$$
\mu([x_1:\ldots:x_p])=\frac{\sqrt{\sum_{j\in I_k}|x_j|^2}}{\sqrt{\sum_{1\leq i\leq p}|x_i|^2}}.
$$
Set $U=\{ \mu<\varepsilon/2p\}\subset \mathbb P^{p-1}$.
Then $U$ is an open neighbourhood of $E_k$.
Let $U'\subset  \mathrm{Bl}_{E_k}\mathbb P^{p-1}$ be the inverse image of $U$ under $\mathrm{Bl}_{E_k}\mathbb P^{p-1}\to \mathbb P^{p-1}$.
Then $U'$ is an open neighbourhood of $E_k'$.
By $g_k(\mathbb D)\not\subset H_k$, we may take an open neighbourhood $W_0\subset \mathbb P^{|I_k|-1}$ of $H_k$ such that 
$$C_H^1(\mathbb D(s)\cap g_k^{-1}(W_0))<\gamma/4p.$$
We take an open neighbourhood $W\Subset W_0$ of $H_k$.
Then since $(\varphi_k\circ \hat{f})_{f\in\mathcal{G}}$ converges uniformly on compact subsets of $\mathbb D$ to $g_k$, we have 
$$
\mathbb D(s)\cap \hat{f}^{-1}(\varphi_k^{-1}(W))\subset \mathbb D(s)\cap g_k^{-1}(W_0)
$$
for all but finitely many $f\in\mathcal{G}$.
Hence we have 
$$C_H^1(\mathbb D(s)\cap \hat{f}^{-1}(\varphi_k^{-1}(W)))<\gamma/4p$$ 
for all but finitely many $f\in\mathcal{G}$.
By \eqref{eqn:20200917}, we have $(\hat{f})_{f\in\mathcal{G}}\to E_k'\cup \varphi_k^{*}H_k$.
Hence, we have 
$$C_H^1(\mathbb D(s)\backslash \hat{f}^{-1}(U'\cup \varphi_k^{-1}(W)))<\gamma/4p$$ for all but finitely many $f\in\mathcal{G}$.
Hence we have 
$$C_H^1(\mathbb D(s)\backslash \hat{f}^{-1}(U))<\gamma/2p$$ 
for all but finitely many $f\in\mathcal{G}$.
This proves our claim.

We define $\Lambda\subset \{1,\ldots,l\}$ to be the set of $k\in \{1,\ldots,l\}$ such that $g_k(\mathbb D)\not\subset H_k$.
If $|I_k|=1$, then $k\in \Lambda$.
We note that \eqref{eqn:202110051} shows that $\Lambda\not=\{1,\ldots,l\}$.
By changing the indexes, we may assume that $\{1,\ldots,n\}=\{1,\ldots,l\}\backslash\Lambda$.
Then we have $n\geq 1$.

We show $I_1,\ldots,I_n$ satisfy the assertions (2) and (3) of Theorem \ref{cor:bc}.
We take $k\in \{1,\ldots,n\}$.
Then $|I_k|\geq 2$.
For $i,j\in I_k$, we define $\tau_{ij}:(\mathbb G_m)^{p-1}\to\mathbb G_m$ by $(b_1,\ldots,b_{p-1})\mapsto b_i/b_j$, where we set $b_p=1$.
Then $\tau_{ij}$ factors $\psi_k|_{(\mathbb G_m)^{p-1}}:(\mathbb G_m)^{p-1}\to (\mathbb G_m)^{|I_k|-1}$.
Hence $(\tau_{ij}\circ\hat{f})_{f\in\mathcal{G}}$ converges uniformly on compact subsets of $\mathbb D$ to the composite of $g_k:\mathbb D\to (\mathbb G_m)^{|I_k|-1}$ and $\kappa_{ij}:(\mathbb G_m)^{|I_k|-1}\to\mathbb G_m$.
Note that $\tau_{ij}\circ\hat{f}=f_i/f_j$.
Hence $(f_i/f_j)_{f\in\mathcal{G}}$ converges uniformly on compact subsets of $\mathbb D$ to $\kappa_{ij}\circ g_k:\mathbb D\to\mathbb G_m$.
By $g_k(\mathbb D)\subset H_k$, we have
$$
\sum_{i\in I_k}\kappa_{ij}\circ g_k=0.$$
Hence $(\sum_{i\in I_k}f_i/f_j)_{f\in\mathcal{G}}$ converges uniformly on compact subsets of $\mathbb D$ to $0$.
Thus the assertion (2) of Theorem \ref{cor:bc} is true.

We prove Theorem \ref{cor:bc} (3).
Set $I=I_1\sqcup\cdots\sqcup I_n$.
We have
$$
\frac{\sqrt{\sum_{i\in I}|f_i|^2}}{\sqrt{\sum_{1\leq i\leq p}|f_i|^2}}
\geq
\frac{\sqrt{\sum_{1\leq i\leq p}|f_i|^2}-\sum_{j\not\in I}|f_j|}{\sqrt{\sum_{1\leq i\leq p}|f_i|^2}}
=1-\sum_{j\not\in I}\frac{|f_j|}{\sqrt{\sum_{1\leq i\leq p}|f_i|^2}}.
$$
Hence by \eqref{eqn:202110051}, we have
$$
\frac{\sqrt{\sum_{1\leq i\leq p}|f_i(z)|^2}}{\sqrt{\sum_{i\in I}|f_i(z)|^2}}<1+\varepsilon<2
$$
for $\frac{\gamma}{2}$-almost all $z\in\mathbb D(s)$.
Combining this with \eqref{eqn:202110051}, we get the assertion (3) of Theorem \ref{cor:bc}.
The proof of Theorem \ref{cor:bc} is completed.
\hspace{\fill} $\square$

\appendix

\section{Semi-abelian varieties} \label{sec:a}
In this appendix, we describe on semi-abelian varieties.
We only treat the definitions and properties which are needed in this paper.
There are several good references on semi-abelian varieties including \cite[Sec. 5.4]{Brion}, \cite{Iitaka}, \cite[Chap 5]{NW}, \cite[Chap VI]{Serre}.
All algebraic groups are defined over $\mathbb C$.

A semi-abelian variety $A$ is an algebraic group with a (unique) expression 
$$0\to T\to A\to A_0\to 0,$$
where $A_0$ is an abelian variety and $T\simeq \mathbb G_m^l$ is an algebraic torus.
Then $A$ is smooth, connected and commutative (cf. \cite[Rem 5.4.2 (ii)]{Brion}).
By \cite[Thm 2.7.2]{Brion}, the map $A\to A_0$ is a $T$-torsor, i.e., we have the following Cartesian diagram
\begin{equation*}
\begin{CD}
A@<\varphi<< T\times A\\
@VVV    @VVpV\\
A_0@<<< A
\end{CD}
\end{equation*}
Here $\varphi(a,t)=a+t$ and $p$ is the second projection.

\begin{rem}\label{rem:20220613}
Let $A$ be a semi-abelian variety and let $B\subset A$ be a connected algebraic subgroup.
By \cite[Cor 5.4.6]{Brion}, $B$ is a semi-abelian variety.
We have a quotient $q:A\to A/B$, which is a $B$-torsor (cf. \cite[Thm 2.7.2]{Brion}).
In particular, $q$ is faithfully flat and quasi-compact.
By \cite[Cor 5.4.6]{Brion}, $A/B$ is a semi-abelian variety.
\end{rem}

By an equivariant compactification $\overline{A}$ of $A$, we mean that (1) $\overline{A}$ is compact, (2) an open immersion $A\subset \overline{A}$ exists, and (3) the group morphism $A\times A\to A$ extends to $A\times \overline{A}\to \overline{A}$.

\subsection{Construction of an equivariant compactification}
The main purpose of this subsection is to introduce an equivariant compactification $\overline{A}$ of $A$ constructed from an equivariant compactification $\overline{T}$ of $T$.
This compactification is described in \cite[p. 1414]{kuh}, \cite[Lemma 2.2]{Vojta2}.

Let $V$ be an algebraic variety which admits a $T$-action.
Then we have a $T$-action on $V\times A$ defined by $(x,a)\mapsto (x-t,a+t)$ for $t\in T$.

\begin{lem}
The categorical quotient $V\times A\to (V\times A)/T$ exists.
Namely every $T$-invariant morphism $V\times A\to W$ factors uniquely through $V\times A\to (V\times A)/T$.
\end{lem}

{\it Proof.}\
Let $\{ U_i\}$ be a Zariski open covering of $A_0$ with $T$-equivariant trivialization $\phi_i:A|_{U_i}\to T\times U_i$ (cf. \cite[p. 169]{Serre}, \cite[Lem 2.2]{Vojta}, \cite[Prop 16.55]{Milne}).
Then for each $i,j$, we get an isomorphism
$$
\phi_j\circ \phi_i^{-1}|_{T\times(U_i\cap U_j)}:T\times (U_i\cap U_j)\to T\times (U_i\cap U_j).
$$
Let $s_{ij}:U_i\cap U_j\to T$ be defined by $\phi_j\circ \phi_i^{-1}|_{T\times (U_i\cap U_j)}(0_{T},u)=(s_{ij}(u),u)$.
Note that $A$ is reconstructed from a gluing of trivial $T$-torsors $T\times U_i$ by the \v{C}ech cocycle $\{ s_{ij}\}\in \check{\mathrm{H}}(\{ U_i\},T)$.
Now we glue $V\times U_i$ by the same \v{C}ech cocycle $\{ s_{ij}\}$ to get $(V\times A)/T$.
The $T$-invariant morphism $V\times A\to (V\times A)/T$ is described as follows.
For each $i$, we define a map 
$$\mu_i:V\times T\times U_i\to V\times U_i$$ 
by $\mu_i(x,t,u)=(x+t,u)$.
Then $\mu_i$ is $T$-invariant under the $T$-action on $V\times T\times U_i$ defined by $(x,t,u)\mapsto (x-\tau,t+\tau,u)$, where $\tau\in T$.
We note that $\mu_i$ is a $T$-torsor with respect to this $T$-action.
The space $V\times A$ is described by a gluing of spaces $V\times T\times U_i$ by the isomorphisms $V\times T\times (U_i\cap U_j)\to V\times T\times (U_i\cap U_j)$ defined by $(x,t,u)\mapsto (x,t+s_{ij}(u),u)$.
Then we may glue $\mu_i$ to get a $T$-torsor 
$$\mu:V\times A\to (V\times A)/T.$$ 
In particular, $\mu$ is a categorical quotient (cf. \cite[Prop 2.6.4]{Brion}).
We remark that the space $(V\times A)/T$ does not depend on the choices of $\{ U_i\}$ and $\phi_i$.
\hspace{\fill} $\square$

\begin{rem}\label{rem:20220725}
Suppose $T$-equivariant map $f:V'\to V$ exists.
Then by the above construction, $f$ induces $f':(V'\times A)/T\to (V\times A)/T$.
If $f$ is an open (resp. closed) immersion, then $f'$ is an open (resp. closed) immersion.
If $V$ is smooth, then $(V\times A)/T$ is smooth.
Indeed $(V\times A)/T$ is constructed by a gluing of the spaces $V\times U_i$, which are smooth.
\end{rem}

\begin{lem}
Let $\overline{T}$ be an equivariant compactification of $T$.
Then $(\overline{T}\times A)/T$ is an equivariant compactification of $A$.
\end{lem}

{\it Proof.}\
We have $A=(T\times A)/T$.
Hence the open immersion $T\subset \overline{T}$ induces an open immersion $A\subset (\overline{T}\times A)/T$.
The $A$-action on $(\overline{T}\times A)/T$ is described as follows.
Let 
$$h:A\times (\overline{T}\times A)\to A\times ((\overline{T}\times A)/T)$$
be defined by $h(a,t,a')=(a,\mu(t,a'))$, where $\mu:\overline{T}\times A\to(\overline{T}\times A)/T$ is the quotient map.
Then $h$ is a $T$-torsor, hence a categorical quotient under the $T$-action on $A\times (\overline{T}\times A)$ defined by
\begin{equation}\label{eqn:202207253}
(a,t,a')\mapsto (a,t-\tau,a'+\tau),
\end{equation}
where $\tau\in T$.
Let 
\begin{equation}\label{eqn:20220724}
\varphi:A\times  (\overline{T}\times A)\to (\overline{T}\times A)/T
\end{equation}
be defined by $\varphi(a,t,a')=\mu(t,a+a')$.
Then $\varphi$ is invariant under the $T$-action \eqref{eqn:202207253}.
Hence $\varphi$ is the composite of $h$ and a unique map 
\begin{equation}\label{eqn:202207241}
\psi:A\times ((\overline{T}\times A)/T)\to (\overline{T}\times A)/T.
\end{equation}
This is our $A$-action.
Note that the map $A\to A_0$ extends to $(\overline{T}\times A)/T\to A_0$, which is proper.
Hence $(\overline{T}\times A)/T$ is compact.
\hspace{\fill} $\square$

\begin{lem}
Let $\overline{T}$ be an equivariant compactification of $T$.
Let $Z\subset \overline{T}$ be an irreducible locally closed set which is $T$-invariant.
Then $(Z\times A)/T\subset (\overline{T}\times A)/T$ is an $A$-invariant, irreducible locally closed set.
Moreover, every $A$-invariant, irreducible locally closed set $V\subset (\overline{T}\times A)/T$ is obtained in this way.
\end{lem}

{\it Proof.}\
By the construction, $(Z\times A)/T\subset (\overline{T}\times A)/T$ is a locally closed set (cf. Remark \ref{rem:20220725}).
We have $\varphi(A\times Z\times A)=(Z\times A)/T$, where $\varphi$ is the same as \eqref{eqn:20220724}.
Hence $\psi(A\times ((Z\times A)/T))=(Z\times A)/T$, where $\psi$ is the same as \eqref{eqn:202207241}.
Hence $(Z\times A)/T\subset (\overline{T}\times A)/T$ is $A$-invariant.

Let $V\subset (\overline{T}\times A)/T$ be an $A$-invariant, irreducible locally closed set.
Set $Z=V\cap \overline{T}$, where $\overline{T}$ is identified with the fibers of $(\overline{T}\times A)/T\to A_0$.
Then we have $V=(Z\times A)/T$.
\hspace{\fill} $\square$

\begin{lem}\label{lem:20200906}
Let $\overline{A}$ be a smooth equivariant compactification of $A$.
Then there exists a smooth equivariant compactification $T\subset \overline{T}$ so that $\overline{A}=(\overline{T}\times A)/T$.
\end{lem}

{\it Proof.}
There exists a rational map $\overline{A}\dashrightarrow A_0$, but every rational map from smooth variety to an abelian variety is regular.
So we have a regular map $p:\overline{A}\to A_0$.
Note that this map is $A$-equivariant.
Indeed denoting by $m:A\times \overline{A}\to \overline{A}$ and $m':A\times A_0\to A_0$ the actions, we have two maps $p\circ m:A\times \overline{A}\to A_0$ and $m'\circ (\mathrm{id}_A,p):A\times \overline{A}\to A_0$.
These coincide on $A\times A\subset A\times\overline{A}$.
Hence $p\circ m=m'\circ (\mathrm{id}_A,p)$, which shows that $p:\overline{A}\to A_0$ is $A$-equivariant.
Since $\overline{A}$ is smooth, there exists a non-empty Zariski open set $U\subset A_0$ such that $p^{-1}(U)\to U$ is smooth (\cite[III, Cor 10.7]{H}).
By translation, $p:\overline{A}\to A_0$ is smooth.

Now set $\overline{T}=p^{-1}(0_{A_0})$.
Then $\overline{T}$ is a smooth equivariant compactification of $T$.
By restricting $m:A\times \overline{A}\to \overline{A}$ to $A\times \overline{T}\subset A\times\overline{A}$, we get $A\times \overline{T}\to\overline{A}$.  
This map is $T$-invariant.
Hence we get $\psi:(\overline{T}\times A)/T\to \overline{A}$.
Note that the restriction $\psi|_A:A\to \overline{A}$ is the open immersion.
Hence $\psi$ is a map over $A_0$.
Note that $\psi$ is injective on the fiber $\overline{T}\subset (\overline{T}\times A)/T$ over $0_{A_0}\in A_0$.
Hence $\psi$ is injective.
Since $\overline{A}$ is smooth, Zariski's main theorem yields that $\psi$ is an open immersion.
Since $(\overline{T}\times A)/T$ is complete, $\psi$ is an isomorphism.
Thus $\overline{A}$ is obtained by $(\overline{T}\times A)/T$.
\hspace{\fill} $\square$

\subsection{Fans and torus embeddings}\label{sec:20230331}

We describe on torus embeddings associated to fans.
We basically follow the notations described in \cite{Oda}.
See also \cite{Cox}, \cite{Fulton}, \cite{Toroidal}.

Let $N=\mathbb Z^r$.
Set $N_{\mathbb R}=N\otimes \mathbb R$, $M=\mathrm{Hom}(N,\mathbb Z)$ and $M_{\mathbb R}=M\otimes \mathbb R$.
A subset $C\subset N_{\mathbb R}$ is a {\it convex rational polyhedral cone} if there exist $n_1,\ldots,n_s\in N$ such that $C=\mathbb R_{\geq 0}n_1+\cdots+\mathbb R_{\geq 0}n_s$.
We set 
$$C^{\vee}=\{ u\in M_{\mathbb R}; u(x)\geq0 \ \text{for all $x\in C$}\}.$$
Then $C^{\vee}\subset M_{\mathbb R}$ is a convex rational polyhedral cone (cf. \cite[p. 2]{Oda}).
A subset $F\subset C$ is called a {\it face} and is denoted $F\prec C$ if there exists $u\in C^{\vee}$ such that $F=C\cap\{u=0\}$.
Then $F$ is a convex rational polyhedral cone (cf. \cite[p. 2]{Oda}).
We call $C$ {\it strongly convex} if $C\cap(-C)=\{ 0\}$.
A {\it fan} in $N$ is a nonempty collection $\Delta$ of strongly convex rational polyhedral cones in $N_{\mathbb R}$ satisfying the following conditions: 
\begin{enumerate}
\item
Every face of any $\sigma\in\Delta$ is contained in $\Delta$.
\item
For any $\sigma,\sigma'\in \Delta$, the intersection $\sigma\cap\sigma'$ is a face of both $\sigma$ and $\sigma'$.
\end{enumerate}
We call $\Delta$ {\it complete} if $\Delta$ is finite and $\cup_{\sigma\in\Delta}\sigma=N_{\mathbb R}$.
We call $\Delta$ {\it non-singular} if for all $\sigma\in\Delta$, there exists a $\mathbb Z$-basis $n_1,\ldots,n_r$ of $N$ such that $\sigma=\mathbb R_{\geq0}n_1+\cdots+\mathbb R_{\geq 0}n_s$, where $s\leq r$.

Let $N'=\mathbb Z^{r'}$ and let $\Delta'$ be a fan in $N'_{\mathbb R}$.
A {\it map of fans} $\varphi:(N',\Delta')\to (N,\Delta)$ is a $\mathbb Z$-linear map $\varphi:N'\to N$ such that for each $\sigma'\in\Delta'$, there exists $\sigma\in \Delta$ such that $\varphi(\sigma')\subset \sigma$.
Here we continue to write the scalar extension by $\varphi:N'_{\mathbb R}\to N_{\mathbb R}$.

Let $T$ be an algebraic torus.
Let $M=\mathrm{Hom}(T,\mathbb G_m)$ and $N=\mathrm{Hom}(\mathbb G_m,T)$.
By definition (cf. \cite[Thm 1.4]{Oda}), a {\it torus embedding} $T\subset T\mathrm{emb}(\Delta)$ is defined from a fun $(N,\Delta)$ as follows. 
For each $\sigma\in \Delta$, we construct the affine scheme $\mathrm{Spec}\ (\mathbb C[M\cap \sigma^{\vee}])$.
Then for $\sigma'\in\Delta$ with $\sigma'\prec\sigma$, we get the open immersion $U_{\sigma'}\subset U_{\sigma}$.
By this, we may glue $\{ U_{\sigma}\}_{\sigma\in \Delta}$ to get $T\mathrm{emb}(\Delta)$.
Then $T\mathrm{emb}(\Delta)$ is normal and contains $T=U_{\{0\}}$ as a dense Zariski open set.
Moreover $T$ acts on $T\mathrm{emb}(\Delta)$.
A torus embedding $T\mathrm{emb}(\Delta)$ is complete if and only if $\Delta$ is complete, and smooth if and only if $\Delta$ is non-singular.
Conversely, by Sumihiro's theorem, every normal, equivariant compactification $\overline{T}$ of $T$ is a torus embedding associated to some complete fan (\cite[Thm 1.5]{Oda}).
Let $\psi:T'\to T$ be a homomorphism of algebraic tori.
Then we have $\psi_*:N'\to N$, where $N'=\mathrm{Hom}(\mathbb G_m,T')$.
Let $\Delta'$ be a fan in $N'$.
Suppose $\psi_*:(N',\Delta')\to (N,\Delta)$ is a map of fans.
Then $\psi$ extends to an equivariant map $T'\mathrm{emb}(\Delta')\to T\mathrm{emb}(\Delta)$ (cf. \cite[Thm 1.13]{Oda}).

Let $\Delta$ be a complete fan in $N$.
Set $\bar{T}=T\mathrm{emb}\Delta$.
For each $\sigma\in\Delta$, there exists a unique $T$-orbit $\mathrm{orb}(\sigma)\subset U_{\sigma}$ which is Zariski closed in $U_{\sigma}$.
Then by \cite[Prop 1.6 (iv)]{Oda}, we have
\begin{equation}\label{eqn:202207255}
\coprod_{\tau \prec\sigma}\mathrm{orb}(\tau)=U_{\sigma}.
\end{equation}
Set $\sigma^{\perp}=\{ u\in M_{\mathbb R}; u(x)=0 \ \text{for all $x\in \sigma$}\}$.
Then $\mathrm{orb}(\sigma)$ is identified with $\mathrm{Spec}\mathbb C[M\cap \sigma^{\perp}]$, where the closed immersion $\mathrm{orb}(\sigma)\subset U_{\sigma}$ is described by the surjection $\mathbb C[M\cap \sigma^{\vee}]\to \mathbb C[M\cap \sigma^{\perp}]$ defined by $\chi^m\mapsto 0$ for all $m\in (M\cap \sigma^{\vee})\backslash (M\cap \sigma^{\perp})$ (cf. \cite[Prop 1.6 (iv)]{Oda}).
Here for the notation $\chi^m$ we refer the readers to \cite[p. 15]{Fulton}.
We note that $\mathrm{Spec}\mathbb C[M\cap \sigma^{\perp}]$ is a quotient torus defined by the inclusion $\mathbb C[M\cap \sigma^{\perp}]\subset \mathbb C[M]$.
The inclusion $\mathbb C[M\cap \sigma^{\perp}]\subset \mathbb C[M\cap \sigma^{\vee}]$ yields the map
\begin{equation}\label{eqn:20220718}
U_{\sigma}\to \mathrm{orb}(\sigma)
\end{equation}
whose restriction to $T\subset U_{\sigma}$ is the quotient map $T\to \mathrm{orb}(\sigma)$.
The closed immersion $\mathrm{orb}(\sigma)\hookrightarrow U_{\sigma}$ is a section of the map \eqref{eqn:20220718}.

Let $I_{\sigma}\subset T$ be the kernel of the quotient map $T\to\mathrm{orb}(\sigma)$.
Then $I_{\sigma}$ is the isotropy group for the points of $\mathrm{orb}(\sigma)\subset U_{\sigma}$.
Set $M_{\sigma}=M/(M\cap \sigma^{\perp})$.
We claim that $I_{\sigma}\subset T$ is a subtorus and 
\begin{equation}\label{eqn:20220725}
\mathrm{Hom}(I_{\sigma},\mathbb G_m)=M_{\sigma}.
\end{equation}
Indeed, since $M\cap \sigma^{\perp}\subset M$ is saturated, $M_{\sigma}$ is a free $\mathbb Z$-module.
Hence we get an exact sequence of free $\mathbb Z$-modules:
$$
0\to M\cap \sigma^{\perp}\to M\to M_{\sigma}\to 0.
$$
This yields the exact sequence of tori (cf. \cite[Thm 12.9]{Milne})
$$
1\to \mathrm{Spec}\mathbb C[M_{\sigma}]\to \mathrm{Spec}\mathbb C[M]\to \mathrm{Spec}\mathbb C[M\cap \sigma^{\perp}]\to 1.
$$
Hence we have $I_{\sigma}=\mathrm{Spec}\mathbb C[M_{\sigma}]$.
This shows that $I_{\sigma}$ is a torus such that \eqref{eqn:20220725}.

\subsection{Smooth modifications}

\begin{lem}\label{lem:20211030111}
Given an equivariant compactification $\overline{A}$, we may take a smooth equivariant compactification $\hat{A}$ and an equivariant morphism $\hat{A}\to \overline{A}$.
\end{lem}

{\it Proof.}\
Let $\overline{T}$ be the Zariski closure of $T\cdot 0_A\subset A$ in $\overline{A}$.
Then $T\subset \overline{T}$ is an equivariant compactification.
We have a canonical map $\psi:(\overline{T}\times A)/T\to \overline{A}$.
We have a map $T\times \overline{T}\to \overline{T}$ which extends the group morphism $T\times T\to T$.
Let $\tilde{T}\to \overline{T}$ be the normalization.
This is an isomorphism over $T$.
Hence $T\subset \tilde{T}$.
Since $T\times \tilde{T}$ is normal, the composite $T\times\tilde{T}\to T\times \overline{T}\to \overline{T}$ factors $\tilde{T}\to \overline{T}$.
Hence we get $T\times\tilde{T}\to\tilde{T}$, which extends the group morphism $T\times T\to T$.
Hence $\tilde{T}$ is an equivariant compactification of $T$.
Then by Sumihiro's theorem, $T\subset \tilde{T}$ is a torus embedding.
Hence by \cite[p. 23]{Oda}, there exists a smooth equivariant compactification $\hat{T}$ and an equivariant morphism $\hat{T}\to \tilde{T}$.
Set $\hat{A}=(\hat{T}\times A)/T$.
Then we have $A$-equivariant morphisms $\hat{A}\to (\overline{T}\times A)/T\to \overline{A}$.
Note that $\hat{A}$ is a smooth equivariant compactification.
Compare with \cite[Thm 1.1]{knop}.
\hspace{\fill} $\square$

\begin{lem}\label{lem:202206031}
Given a smooth equivariant compactification $\overline{A}$, we may take a smooth projective equivariant compactification $\hat{A}$ and an equivariant morphism $\hat{A}\to \overline{A}$.
\end{lem}

{\it Proof.}\
By Lemma \ref{lem:20200906}, we have $\overline{A}=(\overline{T}\times A)/T$.
By \cite[Prop. 2.17]{Oda}, there exists an equivariant modification $\hat{T}\to \overline{T}$ such that $\hat{T}$ is projective.
By \cite[p. 23]{Oda}, we may assume that $\hat{T}$ is smooth.
Then $\hat{A}=(A\times \hat{T})/T$ is smooth.
We claim that this is also projective.
Indeed, let $D=\sum_jn_jD_j$ be a $T$-invariant, very ample divisor on $\hat{T}$.
We set $E=\sum_jn_j(D_j\times A)/T$.
Then $E$ is a divisor on $\hat{A}$.
Let $p:\hat{A}\to A_0$ be the canonical projection.
Let $\{ U_i\}$ be a Zariski open covering of $A_0$ such that $p^{-1}(U_i)=\hat{T}\times U_i$.
Then $p^{-1}(U_i)\cap E=D\times U_i$, which shows that $\mathcal{O}_{\hat{A}}(E)$ is $p$-very ample relative to $A_0$ (cf. \cite[Def. 13.52]{Gortz}).
Since $A_0$ is projective, there exists an ample line bundle $L$ on $A_0$ such that $p^*L(E)$ is ample on $\hat{A}$ (cf. \cite[Prop. 13.65]{Gortz}). 
Hence $\hat{A}$ is projective.
\hspace{\fill} $\square$

\begin{lem}\label{lem:202109231}
Let $\bar{A}$ and $\tilde{A}$ be equivariant compactifications.
Then there exists a smooth equivariant compactification $\hat{A}$ such that equivariant morphisms $\hat{A}\to\bar{A}$ and $\hat{A}\to \tilde{A}$ exist.
\end{lem}

{\it Proof.}\
Let $A\subset A\times A$ be the diagonal, which induces $A\subset \bar{A}\times \tilde{A}$.
Let $V\subset \bar{A}\times \tilde{A}$ be the Zariski closure of this.
Then $V$ is an $A$-equivariant compactification with equivariant maps $V\to \bar{A}$ and $V\to\tilde{A}$.
By Lemma \ref{lem:20211030111}, there exists a smooth equivariant compactification $\hat{A}$ with an equivariant morphism $\hat{A}\to V$.
The composites of this with $V\to \bar{A}$ and $V\to\tilde{A}$ yield $\hat{A}\to\bar{A}$ and $\hat{A}\to \tilde{A}$.
\hspace{\fill} $\square$

\begin{lem}\label{lem:20220728}
Every semi-abelian variety admits only countably many smooth equivariant compactifications.
\end{lem}

{\it Proof.}\
A smooth equivariant compactification $\overline{A}$ is written as $(\overline{T}\times A)/T$, where $\overline{T}$ is a smooth equivariant compactification (cf. Lemma \ref{lem:20200906}).
By Sumihiro's theorem, every smooth, equivariant compactification $\overline{T}$ of $T$ is a torus embedding associated to some complete fan (\cite[Thm 1.5]{Oda}).
There are only countably many complete fans.
Hence there are only countably many smooth equivariant compactifications of $T$, hence of $A$.
\hspace{\fill} $\square$

\subsection{Orbit closure and isotropy group}

Let $\overline{A}=(\overline{T}\times A)/T$, where $\overline{T}=T\mathrm{emb}\Delta$.
For each $\sigma\in \Delta$, we set $O_{\sigma}=(\mathrm{orb}(\sigma)\times A)/T$.
Then $O_{\sigma}\subset \overline{A}$ is an $A$-orbit.
By \eqref{eqn:202207255}, we have $\overline{A}=\coprod_{\sigma\in\Delta}O_{\sigma}$.
Let $I\subset A$ be the isotropy group for the points of $O_{\sigma}$.
Since $I$ acts trivially on the abelian variety $A_0$, we get $I\subset T$.
Then $I$ is the isotropy group for the points of $\mathrm{orb}(\sigma)\subset \overline{T}$.
Hence $I\subset T$ is a subtorus and \eqref{eqn:20220725} yields
\begin{equation}\label{eqn:202207252}
\mathrm{Hom}(I,\mathbb G_m)=M/(M\cap \sigma^{\perp}).
\end{equation}
Let $V\subset \overline{A}$ be an irreducible Zariski closed set which is $A$-invariant.
Then there exists a unique $A$-orbit $O_{\sigma}\subset \overline{A}$ such that $V$ is the Zariski closure of $O_{\sigma}$.
Suppose $\overline{A}$ is smooth.
Then by \cite[Cor. 1.7]{Oda}, the closure of $\mathrm{orb}(\sigma)\subset \overline{T}$ is smooth.
Hence $V$ is smooth.

\begin{lem}\label{lem:vect}
Let $\overline{A}$ be a smooth equivariant compactification.
Let $V\subset \overline{A}$ be an irreducible Zariski closed set which is $A$-invariant.
Then there exist an $A$-invariant Zariski open set $W\subset \overline{A}$ with $V\subset W$ and an equivariant morphism $p:W\to V$ such that $p:W\to V$ is a total space of a vector bundle over $V$ whose zero-section is the inclusion map $V\hookrightarrow W$.
Moreover $V$ is covered by Zariski open subsets $V_0$ such that the following properties hold:
\begin{enumerate}
\item
$p:W\to V$ is a trivial vector bundle $\mathbb C^l\times V_0\to V_0$ on $V_0$.
\item
Let $O\subset \overline{A}$ be the $A$-orbit such that $V$ is the closure of $O$.
Let $I$ be the isotropy group for the points of $O\subset \overline{A}$.
Then there exists a decomposition $I=\mathbb G_m^l$ such that $(t_1,\ldots,t_l)\in \mathbb G_m^l$ acts $(z_1,\ldots,z_l,y)\in\mathbb C^l\times V_0$ as $(z_1,\ldots,z_l,y)\mapsto (t_1z_1,\ldots,t_lz_l,y)$.
\end{enumerate}
\end{lem}

{\it Proof.}\
We first consider the case $A=T$.
Set $M=\mathrm{Hom}(T,\mathbb G_m)$ and $N=\mathrm{Hom}(\mathbb G_m,T)$.
Let $\Delta$ be the fan defining $\overline{T}$, i.e., $\overline{T}=T\mathrm{emb}(\Delta)$.
By \cite[Prop. 1.6]{Oda}, we may take $\tau\in \Delta$ such that $V$ is the closure of $\mathrm{orb}(\tau)$.
We may take $v_1,,\ldots,v_l\in N$ such that $\tau=\mathbb R_{\geq 0}v_1+\cdots+\mathbb R_{\geq 0}v_l$ and $\tau\cap N=\mathbb Z_{\geq 0}v_1+\cdots+\mathbb Z_{\geq 0}v_l$.
Let $\overline{N}=N/N_{\tau}$, where $N_{\tau}=\mathbb Z(\tau\cap N)$.
Then $\overline{N}$ is the dual of $M\cap \tau^{\perp}$.
We have
$$
\mathrm{orb}(\tau)=\mathrm{Spec}\mathbb C[M\cap \tau^{\perp}].
$$
For each $\sigma\in \Delta$ with $\tau\prec\sigma$, we set $\bar{\sigma}=(\sigma+\mathbb R\tau)/\mathbb R\tau\subset \overline{N}_{\mathbb R}$.
Then $\bar{\sigma}\subset \overline{N}_{\mathbb R}$ is a strongly convex rational polyhedral cone (cf. \cite[Prop. A.8]{Oda}).
Set $\bar{U}_{\bar{\sigma}}=\mathrm{Spec}\mathbb C[(M\cap \tau^{\perp})\cap \bar{\sigma}^{\vee}]$.
The inclusion $\mathbb C[(M\cap \tau^{\perp})\cap \bar{\sigma}^{\vee}]\subset \mathbb C[M\cap\sigma^{\vee}]$ yields $U_{\sigma}\to \bar{U}_{\bar{\sigma}}$.
We have a closed immersion $\bar{U}_{\bar{\sigma}}\hookrightarrow U_{\sigma}$, which is obtained by the surjection $\mathbb C[M\cap\sigma^{\vee}]\to \mathbb C[(M\cap \tau^{\perp})\cap \bar{\sigma}^{\vee}]$ defined by $\chi^m\mapsto 0$ for all $m\in (M\cap \sigma^{\vee})\backslash (M\cap \tau^{\perp}\cap \bar{\sigma}^{\vee})$.
This inclusion $\bar{U}_{\bar{\sigma}}\hookrightarrow U_{\sigma}$ is the section of the map $U_{\sigma}\to \bar{U}_{\bar{\sigma}}$.

We remark that the map $U_{\sigma}\to \bar{U}_{\bar{\sigma}}$ is a trivial vector bundle whose zero-section is the inclusion map $\bar{U}_{\bar{\sigma}}\hookrightarrow U_{\sigma}$.
Indeed, the sequence
$$
0\to N_{\tau}\to N\to \overline{N}\to 0
$$
has a (non-canonical) splitting $N=N_{\tau}\oplus \overline{N}$ such that $\sigma=\tau\oplus\bar{\sigma}$, where we view $\tau\subset N_{\mathbb R}$ as $\tau\subset (N_{\tau})_{\mathbb R}$.
Considering the dual
$$
0\to M\cap \tau^{\perp}\to M\to M_{\tau}\to 0,
$$
we get the associated splitting $M=M_{\tau}\oplus (M\cap \tau^{\perp})$.
Hence we get 
$$
\mathbb C[M\cap \sigma^{\vee}]=\mathbb C[M_{\tau}\cap\tau^{\vee}]\otimes_{\mathbb C}\mathbb C[(M\cap \tau^{\perp})\cap \bar{\sigma}^{\vee}].
$$
We have $\mathrm{Spec}\mathbb C[M_{\tau}\cap\tau^{\vee}]=\mathbb C^l$, hence $U_{\sigma}\simeq\mathbb C^l\times \bar{U}_{\bar{\sigma}}$.
The inclusion $\bar{U}_{\bar{\sigma}}\subset U_{\sigma}$ coincides with $\{ 0\}\times \bar{U}_{\bar{\sigma}}\subset \mathbb C^l\times \bar{U}_{\bar{\sigma}}$.
Hence the map $U_{\sigma}\to \bar{U}_{\bar{\sigma}}$ is a trivial vector bundle whose zero-section is the inclusion map $\bar{U}_{\bar{\sigma}}\hookrightarrow U_{\sigma}$.
Moreover the stabilizer $I$ of $\mathrm{orb}(\tau)$ is identified with $(\mathbb G_m)^l\subset\mathbb C^l$.
Note that the ambiguity of the splitting $N=N_{\tau}\oplus \overline{N}$ yields an isomorphism $\mathbb C^l\times \bar{U}_{\bar{\sigma}}\simeq \mathbb C^l\times \bar{U}_{\bar{\sigma}}$ which is defined by the multiplication given by $\bar{U}_{\bar{\sigma}}\to I$.

Now set $\bar{\Delta}=\{ \bar{\sigma};\ \sigma\in \Delta, \tau\prec\sigma\}$.
Then, by \cite[Cor. 1.7]{Oda}, $\bar{\Delta}$ is a fan in $\overline{N}_{\mathbb R}$ and $V=\mathrm{orb}(\tau)\mathrm{emb}(\bar{\Delta})$.
Hence $V=\cup_{\bar{\sigma}\in \bar{\Delta}}\bar{U}_{\bar{\sigma}}$.
We set
$$
W=\bigcup_{\sigma\in \Delta,\tau\prec\sigma}U_{\sigma}
$$
Then the maps $U_{\sigma}\to \bar{U}_{\bar{\sigma}}$ glue together to get a vector bundle $W\to V$.
Indeed the coordinate transformations of this vector bundle are given by the multiplications by $I$.
Since $U_{\sigma}$ are $T$-invariant, $W$ is $T$-invariant.
This completes the proof when $A=T$.

Next we consider the general case.
We write as $V=(V'\times A)/T$, where $V'\subset \overline{T}$ is a $T$-invariant Zariski closed set.
Then we may take a Zariski open set $W'\subset \overline{T}$ and $p:W'\to V'$ by the previous consideration.
We set $W=(W'\times A)/T$, where $W'$ is $T$-invariant.
Then $W\subset \overline{A}$ is an $A$-invariant Zariski open set and $V\subset W$.
This yields $W\to V$ as desired.
\hspace{\fill} $\square$

\begin{rem}\label{rem:20220728}
Let $\overline{A}$ be a smooth equivariant compactification.
Then the boundary $\partial A$ is a simple normal crossing divisor. 
Indeed for each $x\in \partial A$, we take a unique $A$-orbit $O\subset \overline{A}$ such that $x\in O$ and set $V=\overline{O}$.
Then $V$ is smooth.
We take $W\subset \overline{A}$ as in Lemma \ref{lem:vect}.
Let $V_0\subset V$ be described in Lemma \ref{lem:vect} so that $x\in V_0$.
We may assume $V_0\subset O$.
Then $(\partial A)\cap (\mathbb C^l\times V_0)=\cup_{i=1}^lH_i\times V_0$, where $H_i\subset \mathbb C^l$ is define by $H_i=\{ z_i=0\}$.
Hence $(\partial A)\cap (\mathbb C^l\times V_0)$ is a simple normal crossing divisor on $\mathbb C^l\times V_0$.
Since $x\in\partial A$ is arbitrary, $\partial A$ is simple normal crossing.
\end{rem}

\subsection{Modification of equivariant maps}

The purpose of this subsection is to prove Lemma \ref{lem:20210513}.

\begin{lem}\label{lem:20210923}
Let $\psi:N\to N'$ be surjective.
Let $\Delta$ be a fan in $N$ and let $\Delta'$ be a fan in $N'$.
Set $\Delta''=\{ \sigma\cap\psi^{-1}(\sigma');\sigma\in \Delta, \sigma'\in\Delta'\}$, where we continue to write $\psi:N_{\mathbb R}\to N'_{\mathbb R}$.
Then $\Delta''$ is a fan in $N$.
\end{lem}

Before going to prove this lemma, we remark the followings.
Given a convex rational polyhedral cone $C\subset N_{\mathbb R}$, we have $(C^{\vee})^{\vee}=C$ (cf. \cite[Thm A.1]{Oda}).
For another convex rational polyhedral cone $C'\subset N_{\mathbb R}$, we have $(C\cap C')^{\vee}=C^{\vee}+(C')^{\vee}$ (cf. \cite[Thm A.1]{Oda}).
In particular, $C\cap C'$ is a convex rational polyhedral cone.
Moreover if $C$ is strongly convex, then $C\cap C'$ is strongly convex.
By \cite[Prop. A.9]{Oda}, the following two statements are equivalent:
\begin{enumerate}
\item
The subset $C\cap C'\subset C$ is a face of $C$.
\item
For every $x,x'\in C$, if $x+x'\in C\cap C'$, then both $x\in C\cap C'$ and $x'\in C\cap C'$.
\end{enumerate}
The non-trivial part is the implication (2)$\Longrightarrow$(1).
The converse is trivial.
Indeed suppose $C\cap C'$ is a face of $C$.
Then there exists $u\in C^{\vee}$ such that $C\cap C'=C\cap\{ u=0\}$.
If $x,x'\in C$ satisfy $x+x'\in C\cap C'$, then $u(x)\geq 0$, $u(x')\geq 0$ and $u(x+x')=0$.
Hence $u(x)=u(x')=0$, thus $x,x'\in C\cap C'$.

\medskip

{\it Proof of Lemma \ref{lem:20210923}.}\
(Cf. \cite[Lem. 3.2]{KS}.)
Note that $\psi^{-1}(\sigma')$ is a convex rational polyhedral cone.
Since $\sigma$ is strongly convex, $\sigma\cap\psi^{-1}(\sigma')$ is a strongly convex rational polyhedral cone.
Hence $\Delta''$ is a non-empty collection of strongly convex rational polyhedral cones in $N_{\mathbb R}$.

We first show that a face of $\sigma\cap\psi^{-1}(\sigma')$ is contained in $\Delta''$.
We claim
\begin{equation}\label{eqn:20220727}
(\sigma\cap\psi^{-1}(\sigma'))^{\vee}=\sigma^{\vee}+\psi^*((\sigma')^{\vee}),
\end{equation}
where $\psi^*:M'_{\mathbb R}\to M_{\mathbb R}$ is the dual map.
Indeed, by $(\sigma\cap\psi^{-1}(\sigma'))^{\vee}=\sigma^{\vee}+(\psi^{-1}(\sigma'))^{\vee}$, it is enough to show
\begin{equation*}\label{eqn:20220726}
(\psi^{-1}(\sigma'))^{\vee}=\psi^*((\sigma')^{\vee}).
\end{equation*}
Note that $(\psi^{-1}(\sigma'))^{\vee}\supset \psi^*((\sigma')^{\vee})$ is trivial.
To show the converse, let $u\in (\psi^{-1}(\sigma'))^{\vee}$.
Let $K\subset N_{\mathbb R}$ be the kernel of $\psi$.
Then $K\subset \psi^{-1}(\sigma')$.
Hence $K\subset \{ u\geq 0\}$.
By $-K=K$, we get $K\subset \{ u=0\}$.
Hence there exists $v\in M'_{\mathbb R}$ such that $u=v\circ\psi$.
By $u\in (\psi^{-1}(\sigma'))^{\vee}$, we have $\sigma'\subset \{ v\geq 0\}$.
Hence $v\in (\sigma')^{\vee}$.
This shows $u\in \psi^*((\sigma')^{\vee})$, hence $(\psi^{-1}(\sigma'))^{\vee}\subset \psi^*((\sigma')^{\vee})$.
Thus we get \eqref{eqn:20220727}.

Now a face of $\sigma\cap\psi^{-1}(\sigma')$ is defined by $\sigma\cap\psi^{-1}(\sigma')\cap\{u=0\}$ for some $u\in (\sigma\cap\psi^{-1}(\sigma'))^{\vee}$.
By \eqref{eqn:20220727}, we have $u=v+v'\circ \psi$ for some $v\in \sigma^{\vee}$ and $v'\in(\sigma')^{\vee}$.
We claim that 
\begin{equation}\label{eqn:202207271}
\sigma\cap\psi^{-1}(\sigma')\cap\{u=0\}=(\sigma\cap\{v=0\})\cap\psi^{-1}(\sigma'\cap\{ v'=0\}).
\end{equation}
Indeed the relation `$\supset$' is trivial.
To show the converse, let $x\in \sigma\cap\psi^{-1}(\sigma')\cap\{u=0\}$.
Then we have $v(x)\geq 0$ and $v'\circ\psi(x)\geq 0$.
By $u(x)=0$, we have $v(x)=0$ and $v'\circ\psi(x)=0$.
Hence $x\in (\sigma\cap\{v=0\})\cap\psi^{-1}(\sigma'\cap\{ v'=0\})$.
This shows \eqref{eqn:202207271}.
Set $\tau=\sigma\cap\{v=0\}$ and $\tau'=\sigma'\cap\{ v'=0\}$.
Then $\tau$ and $\tau'$ are faces of $\sigma$ and $\sigma'$, respectively.
Hence $\tau\in \Delta$ and $\tau'\in\Delta'$.
By \eqref{eqn:202207271}, we have
$$\sigma\cap\psi^{-1}(\sigma')\cap\{u=0\}=\tau\cap\psi^{-1}(\tau')\in\Delta''.$$

Next we show that $(\sigma_1\cap\psi^{-1}(\sigma_1'))\cap(\sigma_2\cap\psi^{-1}(\sigma_2'))$ is a face of $\sigma_1\cap\psi^{-1}(\sigma_1')$.
We take $x,y\in \sigma_1\cap\psi^{-1}(\sigma_1')$ such that $x+y\in (\sigma_1\cap\psi^{-1}(\sigma_1'))\cap(\sigma_2\cap\psi^{-1}(\sigma_2'))$.
Then we have $x+y\in \sigma_1\cap\sigma_2$.
Since $\sigma_1\cap\sigma_2$ is a face of $\sigma_1$, and $x,y\in\sigma_1$, we get $x,y\in\sigma_1\cap\sigma_2$.
Also we have $\psi(x+y)=\psi(x)+\psi(y)\in \sigma_1'\cap\sigma_2'$.
Since $\sigma_1'\cap\sigma_2'$ is a face of $\sigma_1'$, and $\psi(x),\psi(y)\in\sigma_1'$, we get $\psi(x),\psi(y)\in\sigma_1'\cap\sigma_2'$.
Hence $x,y\in (\sigma_1\cap\psi^{-1}(\sigma_1'))\cap(\sigma_2\cap\psi^{-1}(\sigma_2'))$.
By \cite[Prop. A.9]{Oda}, we have 
$$(\sigma_1\cap\psi^{-1}(\sigma_1'))\cap(\sigma_2\cap\psi^{-1}(\sigma_2'))\prec (\sigma_1\cap\psi^{-1}(\sigma_1')).$$
Hence $\Delta''$ is a fan.
\hspace{\fill} $\square$

\begin{lem}\label{lem:20210920}
Let $\psi:N\to N'$.
Let $\Delta$ be a complete fan in $N$.
Then there exist a non-singular, complete fan $\Delta'$ in $N'$ and a finite subdivision $\tilde{\Delta}$ of $\Delta$ such that for all $\sigma\in \tilde{\Delta}$, there exists $\sigma'\in\Delta'$ such that $\psi(\sigma)=\sigma'$.
\end{lem}

This is contained in the proof of \cite[Prop. 4.4]{AbKaru}.
We give a proof for the sake of completeness.
Before going to prove this lemma, we remark the followings.
If $\Delta_1,\ldots,\Delta_l$ are complete fans in $N$, then $\Delta=\{\sigma_1\cap\cdots\cap\sigma_l;\ \sigma_1\in\Delta_1, \ldots,\sigma_l\in\Delta_l\}$ is a complete fan in $N$.
Indeed for the case $l=2$, Lemma \ref{lem:20210923} applied to the identity map $N\to N$ yields that $\Delta$ is a fan.
The same holds for all $l$ by induction.
The completeness is obvious.
Each $\sigma_1\in \Delta_1$ is a union of cones in $\Delta$.
Indeed, we have 
$$\sigma_1=\bigcup_{\sigma_2\in \Delta_2,\ldots,\sigma_l\in\Delta_2}\sigma_1\cap\sigma_2\cap\cdots\cap\sigma_l,$$ 
for $\Delta_2,\ldots,\Delta_l$ are complete.

\medskip

{\it Proof of Lemma \ref{lem:20210920}.}\
For each $\sigma\in\Delta$, we may take a complete fun $\Delta'_{\sigma}$ in $N'$ such that $\psi(\sigma)$ is a union of cones in $\Delta'_{\sigma}$.
Indeed, by \cite[Thm A.3]{Oda}, there exist strongly convex rational polyhedral cones $\tau_1,\ldots,\tau_l$ such that $\psi(\sigma)=\tau_1\cup\cdots\cup\tau_l$.
By \cite[Prop. 1.12]{Oda}, there exist complete fans $\Delta_1',\ldots,\Delta_l'$ such that $\tau_i\in\Delta_i'$.
We set $\Delta'_{\sigma}=\{\sigma'_1\cap\cdots\cap\sigma_l';\sigma_i'\in\Delta_i'\}$.
Then by Lemma \ref{lem:20210923}, $\Delta'_{\sigma}$ is a complete fan.
Each $\tau_i\in \Delta_i'$ is a union of cones in $\Delta'_{\sigma}$.
Hence $\psi(\sigma)$ is a union of cones in $\Delta'_{\sigma}$, as desired.

We set $\Delta'=\{ \cap_{\sigma\in\Delta}\kappa_{\sigma};\ \kappa_{\sigma}\in \Delta'_{\sigma}\}$.
Then by Lemma \ref{lem:20210923}, $\Delta'$ is a complete fan in $N'$.
For all $\sigma\in \Delta$, $\psi(\sigma)$ is a union of cones in $\Delta'$.
By replacing $\Delta'$ by its finite subdivision, we may assume that $\Delta'$ is non-singular.

We set $\tilde{\Delta}=\{ \psi^{-1}(\tau)\cap \sigma;\ \tau\in\Delta', \sigma\in \Delta\}$.
Then by Lemma \ref{lem:20210923}, $\tilde{\Delta}$ is a complete fan in $N$.
We shall show that these $\tilde{\Delta}$ and $\Delta'$ satisfy our requirement.
We have $\psi(\psi^{-1}(\tau)\cap \sigma)=\tau\cap\psi(\sigma)$.
Hence it is enough to show 
$$\tau\cap\psi(\sigma)\in \Delta'.$$
We may take $\tau_1,\ldots,\tau_k\in \Delta'$ such that $\psi(\sigma)=\tau_1\cup\cdots\cup\tau_k$.
Hence we have 
$$\tau\cap \psi(\sigma)=(\tau\cap\tau_1)\cup\cdots\cup(\tau\cap\tau_k).$$
Each $\tau\cap\tau_i$ is a face of $\tau$.
We take $x,y\in\tau$ such that $x+y\in \tau\cap \psi(\sigma)$.
Then we have $x+y\in \tau\cap\tau_i$ for some $\tau_i$.
Since $\tau\cap\tau_i$ is a face of $\tau$, we have $x,y\in \tau\cap\tau_i$.
Hence $x,y\in \tau\cap \psi(\sigma)$.
By \cite[Prop. A.9]{Oda}, $\tau\cap \psi(\sigma)$ is a face of $\tau$.
Hence $\tau\cap \psi(\sigma)\in \Delta'$.
This completes the proof of our lemma.
\hspace{\fill} $\square$

\begin{lem}\label{lem:20210513}
Let $\overline{A}$ be a smooth equivariant compactification and let $B\subset A$ be a semi-abelian subvariety.
Then there exists an equivariant compactification $\hat{A}=(\hat{T}\times A)/T$ with an equivariant morphism $\hat{A}\to \overline{A}$ such that the following properties hold:
\begin{enumerate}
\item There exists a smooth equivariant compactification $\widehat{A/B}$ such that an equivariant morphism $p:\hat{A}\to \widehat{A/B}$ exists.
\item 
For every $x\in \hat{A}$, if a semi-abelian subvariety $C\subset A$ satisfies $C\cdot p(x)=p(x)$, then $C\cdot x\subset B\cdot x$ as subsets of $\hat{A}$.
\end{enumerate}
\end{lem}

{\it Proof.}\
We have canonical expressions $0\to T_A\to A\to A_0\to 0$ and $0\to T_B\to B\to B_0\to 0$.
Then we have the canonical expression
$$
0\to T_A/T_B\to A/B\to A_0/B_0\to 0.$$
Set $N=\mathrm{Hom}(\mathbb G_m,T_A)$.
Let $\Delta$ be the fan in $N$ such that $\overline{T}_A=T_A\mathrm{emb}(\Delta)$, where $\overline{A}=(\overline{T}_A\times A)/T_A$.
Set $N'=\mathrm{Hom}(\mathbb G_m,T_A/T_B)$.
Then we have $\psi:N\to N'$.
By Lemma \ref{lem:20210920}, we may find a non-singular, complete fan $\Delta'$ in $N'$ and a finite subdivision $\hat{\Delta}$ of $\Delta$ such that for all $\sigma\in \hat{\Delta}$, there exists $\sigma'\in \Delta'$ such that $\psi(\sigma)= \sigma'$.
Set $\hat{A}=(\hat{T}_A\times A)/T_A$ and $\widehat{A/B}=(\widehat{T_A/T_B}\times A/B)/(T_A/T_B)$, where $\hat{T}_A=T_A\mathrm{emb}(\hat{\Delta})$ and $\widehat{T_A/T_B}=T_A/T_B\mathrm{emb}(\Delta')$.
Then $\widehat{T_A/T_B}$ is smooth.
Hence $\widehat{A/B}$ is smooth.
By \cite[p. 19]{Oda}, we get equivariant morphisms $p:\hat{A}\to \widehat{A/B}$ and $\hat{A}\to \overline{A}$.

We shall show that $p:\hat{A}\to \widehat{A/B}$ satisfies the assertion (2).
Let $x\in\hat{A}$.
Set $I_x=\{ a\in A; a\cdot x=x\}$.
We take $\sigma\in\hat{\Delta}$ such that $x\in ( \mathrm{orb}(\sigma)\times A)/T\subset \hat{A}$. 
Set $M_{\sigma}=M\cap\sigma^{\perp}$, where $M=\mathrm{Hom}(T_A,\mathbb G_m)$.
By \eqref{eqn:202207252}, we have 
$$\mathrm{Hom}(I_{x},\mathbb G_m)=M/M_{\sigma},$$
where $I_x\subset T_A$ is a subtorus.
Similarly, set $J_{p(x)}=\{ a\in A/B; a\cdot p(x)=p(x)\}$.
We take $\sigma'\in\Delta'$ such that 
\begin{equation}\label{eqn:202207274}
\psi(\sigma)=\sigma'.
\end{equation}
Let $\mathrm{orb}(\sigma')\subset \widehat{T_A/T_B}$ be the corresponding orbit.
Denoting by $q:\hat{T}_A\to \widehat{T_A/T_B}$ the canonical map, we have $q(\mathrm{orb}(\sigma))=\mathrm{orb}(\sigma')$ (cf. \cite[p. 56]{Fulton}, \cite[Lem 3.3.21]{Cox}).
In particular, we have $p(x)\in(\mathrm{orb}(\sigma')\times (A/B))/(T_A/T_B)\subset \widehat{A/B}$.
Hence by \eqref{eqn:202207252}, we have 
$$\mathrm{Hom}(J_{p(x)},\mathbb G_m)=M'/M'_{\sigma'},$$
where $M'=\mathrm{Hom}(T_A/T_B,\mathbb G_m)$ and $M'_{\sigma'}=M'\cap(\sigma')^{\perp}$.
The quotient map $T_A\to T_A/T_B$ induces a morphism $I_x\to J_{p(x)}$ of tori.
This induces $\mathrm{Hom}(J_{p(x)},\mathbb G_m)\to\mathrm{Hom}(I_x,\mathbb G_m)$, namely $M'/M'_{\sigma'}\to M/M_{\sigma}$.
This fits into the following commutative diagram:
\begin{equation*}
\begin{CD}
0@>>> M'_{\sigma'}@>>> M'@>>> M'/M'_{\sigma'} @>>> 0\\
@. @VVV @VVV @VVV  \\
0@>>> M_{\sigma} @>>> M@>>> M/M_{\sigma}@>>> 0
\end{CD}
\end{equation*}

We claim that the map $I_x\to J_{p(x)}$ is surjective.
Indeed, by \eqref{eqn:202207274}, we have $M'_{\sigma'}=M'\cap M_{\sigma}$.
Hence the map $M'/M'_{\sigma'}\to M/M_{\sigma}$ is injective.
Thus we get the injection $\mathrm{Hom}(J_{p(x)},\mathbb G_m)\hookrightarrow \mathrm{Hom}(I_x,\mathbb G_m)$.
Since the image of the map $I_x\to J_{p(x)}$ is a subtorus of $J_{p(x)}$ (cf. \cite[Prop B, p. 54]{Hump}), this image should coincides with $J_{p(x)}$.
Thus the map $I_x\to J_{p(x)}$ is surjective.

Now let $C\subset A$ satisfies $C\cdot p(x)$.
Then $(C+B)/B\subset J_{p(x)}$.
Hence $(C+B)/B\subset (I_x+B)/B$, so $C\subset I_x+B$.
Thus $C\cdot x\subset (I_x+B)\cdot x=B\cdot x$.
This completes the proof.
\hspace{\fill} $\square$

\subsection{Logarithmic tangent space}

We recall the isomorphism $A\times\mathrm{Lie}A\simeq TA$ from \eqref{eqn:20211108}.
Let $\overline{A}$ be a smooth equivariant compactification of $A$.
Then $\partial A$ is a simple normal crossing divisor on $\overline{A}$ (cf. Remark \ref{rem:20220728}).

\begin{lem}\label{lem:202111111}
The isomorphism \eqref{eqn:20211108} extends to an isomorphism
$$
\bar{\psi}:\overline{A}\times\mathrm{Lie}A\to
T\overline{A}(-\log (\partial A)).
$$
\end{lem}

{\it Proof.}\
(Cf. \cite[Prop. 5.4.3]{NW})
We denote by
$$\bar{m}:\overline{A}\times A\to \overline{A}$$
the natural action defined by $(x,a)\mapsto x+a$.
By $\bar{m}^*(\partial A)=(\partial A)\times A$, we have
\begin{equation*}
T(\overline{A}\times A)(-\log ((\partial A)\times A)) \to T\overline{A}(-\log (\partial A)).
\end{equation*}
We have subbundle
$$
\overline{A}\times TA\subset T(\overline{A}\times A)(-\log ((\partial A)\times A)).
$$
Thus we get
\begin{equation}\label{eqn:20211111}
\overline{A}\times TA\to T\overline{A}(-\log (\partial A)).
\end{equation}
We restrict this to $\overline{A}\times\{0_A\}\subset \overline{A}\times A$ to get 
$$
\bar{\psi}:\overline{A}\times\mathrm{Lie}A\to
T\overline{A}(-\log (\partial A)),
$$
which extends \eqref{eqn:20211108}.
We shall show that $\bar{\psi}$ is an isomorphism.

Let $x\in \overline{A}$ be an arbitrary point.
It is enough to show that the induced map 
$$
\bar{\psi}_{x}:\mathrm{Lie}A\to T_{x}\overline{A}(-\log (\partial A))
$$
is injective.
Let $O\subset \overline{A}$ be a unique $A$-orbit such that $x\in O$.
Let $I\subset A$ be the isotropy group for $x\in \overline{A}$.
Then we have an isomorphism $O\simeq A/I$.
Set $V=\overline{O}$.
By Lemma \ref{lem:vect}, there exist an $A$-invariant Zariski open $W\subset \overline{A}$ with $V\subset W$ and an $A$-equivariant morphism $p:W\to V$.
We take a Zariski open neighbourhood $V_0\subset V$ of $x$ such that $V_0\subset O$ and $p^{-1}(V_0)=\mathbb A^l\times V_0$.
Set $W_0=p^{-1}(V_0)$.
Then $W_0$ is $I$-invariant and $(t_1,\ldots,t_l)\in I$ acts $(y_1,\ldots,y_l,v)\in \mathbb A^l\times V_0$ as $(t_1y_1,\ldots,t_ly_l,v)$.
We denote this action by
$$\mu:W_0\times I\to W_0.$$ 
Note that $\mu$ is the restriction of $\bar{m}$.
The restriction of \eqref{eqn:20211111} on $W_0\times TI\subset A\times TA$ induces 
$$W_0\times TI\to TW_0(-\log (\partial A)).$$
By the dual, we get
\begin{equation}\label{eqn:202207281}
\mu^*\Omega^1_{W_0}(\log (\partial A))\to \Omega^1_{(W_0\times I)/W_0}.
\end{equation}
The image of the sections $dy_1/y_1,\ldots,dy_l/y_l\in \Omega^1_{W_0}(\log (\partial A))$ on $W_0$ under the map \eqref{eqn:202207281} are $dt_1/t_1,\ldots,dt_l/t_l$.
This shows that the map \eqref{eqn:202207281} is surjective over $W_0\times I$.
Considering the dual on $(x,e_I)\in W_0\times I$, we observe that the restriction
$$
\bar{\psi}_{x}|_{\mathrm{Lie}I}:\mathrm{Lie}I\to T_{x}\overline{A}(-\log (\partial A))
$$
is injective.

Next we have the following commutative diagram:
\begin{equation*}
\begin{CD}
W\times A@>\bar{m}>> W\\
@VqVV    @VVpV
\\
V\times (A/I) @>>\bar{m}_V> V
\end{CD}
\end{equation*}
This induces the following  commutative diagram:
\begin{equation*}
\begin{CD}
\mathrm{Lie}A@>\bar{m}_*>> T_{x}W\\
@Vq_*VV    @VVp_*V
\\
\mathrm{Lie}(A/I) @>>(\bar{m}_V)_*> T_{x}V
\end{CD}
\end{equation*}
By $x\in A/I$, we note that $(\bar{m}_V)_*$ is an isomorphism.
Hence $\ker \bar{m}_{*}\subset \mathrm{Lie}I$.
Since $\bar{m}_{*}$ factors as
$$
\mathrm{Lie}A\overset{\bar{\psi}_{x}}{\longrightarrow} T_{x}\overline{A}(-\log (\partial A))\longrightarrow T_{x}W,
$$
we get $\ker \bar{\psi}_{x}\subset \mathrm{Lie}I$.
This shows that $\bar{\psi}_{x}$ is injective.
\hspace{\fill} $\square$

\subsection{An application of faithfully flat descent}

\begin{lem}\label{lem:20211030}
Let $V$ and $\Sigma$ be varieties.
Let $Z\subset V$ be a closed subscheme.
Let $Y\subset V\times \Sigma$ be a closed subscheme such that $Y_s=Z$ for all $s\in \Sigma$.
Then $Z\times \Sigma\subset Y$ as closed subschemes of $V\times \Sigma$.
\end{lem}

{\it Proof.}\
It is enough to prove the case that $V$ and $\Sigma$ are affine.
Thus we assume $V=\Spec R$ and $\Sigma=\Spec S$.
Let $I_Z\subset R$ and $I_Y\subset R\otimes_{\mathbb C}S$ be the ideals associated to $Z\subset V$ and $Y\subset V\times \Sigma$, respectively.
To show $I_Y\subset I_Z\otimes_{\mathbb C}S$, we take $\varphi\in I_Y$.
Then we may write as 
$$\varphi=f_1\otimes g_1+\cdots +f_l\otimes g_l,$$
where $f_1,\ldots,f_l\in R$ and $g_1,\ldots,g_l\in S$.
We take this  expression so that $l$ is minimum.
We define a subset $L\subset \mathbb C^l$ such that $(a_1,\ldots,a_l)\in L$ iff $a_1f_1+\cdots+a_lf_l\in I_Z$.
Then $L\subset \mathbb C^l$ is a linear subspace.
To prove $L=\mathbb C^l$, we assume contrary $L\subsetneqq \mathbb C^l$.
Then we may take a non-zero $(\lambda_1,\ldots,\lambda_l)\in\mathbb C^l$ such that $\lambda_1a_1+\cdots+\lambda_la_l=0$ for all $(a_1,\ldots,a_l)\in L$.
We may assume without loss of generality that $\lambda_1=1$.
Now for all $s\in \Sigma$, we have $Y_s=Z$.
Hence $g_1(s)f_1+\cdots+g_l(s)f_l\in I_Z$.
Hence $(g_1(s),\ldots,g_l(s))\in L$, so $\lambda_1g_1(s)+\cdots+\lambda_lg_l(s)=0$ for all $s\in \Sigma$.
Since $S$ is an integral domain, we have $\lambda_1g_1+\cdots+\lambda_lg_l=0$ in $S$.
By $\lambda_1=1$, we have $g_1=-(\lambda_2g_2+\cdots+\lambda_lg_l)$.
Hence $\varphi=(f_2-\lambda_2 f_1)\otimes g_2+\cdots+ (f_l-\lambda_lf_1)\otimes g_l$.
This contradicts to the minimality of $l$.
Hence $L=\mathbb C^l$.
Now we have $(1,0,\ldots,0)\in L$, hence $f_1\in I_Z$.
Similarly, we have $f_i\in I_Z$ for all $i$.
Hence $\varphi\in  I_Z\otimes_{\mathbb C}S$.
Thus $I_Y\subset I_Z\otimes_{\mathbb C}S$.
This shows $Z\times \Sigma\subset Y$.
\hspace{\fill} $\square$

\medskip

Let $A$ be a semi-abelian variety and let $B\subset A$ be a semi-abelian subvariety.
Let $\Sigma$ be a variety.
Then we get a projection $p:A\times \Sigma\to (A/B)\times \Sigma$, which yields the following fiber product:
\begin{equation}\label{eqn:20211030}
\begin{CD}
B\times A\times \Sigma@>q_2>> A\times \Sigma \\
@Vq_1VV @VVpV   \\
A\times \Sigma @>>p> (A/B)\times \Sigma
\end{CD}
\end{equation}
Here for $(b,a,s)\in B\times A\times \Sigma$, we have $q_1(b,a,s)=(a,s)$ and $q_2(b,a,s)=(b+a,s)$.

\begin{lem}\label{lem:202110301}
Let $Z\subset A\times \Sigma$ be a closed subscheme which is $B$-invariant in the sense that $b+Z=Z$ for all $b\in B$.
Then $q_1^*(Z)=q_2^*(Z)$ as closed subschemes of $B\times A\times \Sigma$.
\end{lem}

{\it Proof.}\
Note that $q_1^*(Z)=B\times Z$.
Let $\pi:B\times A\times \Sigma\to B$ be the first projection.
Then for $b\in B$, we have $q_2^*(Z)\cap \pi^{-1}(b)=(-b)+Z=Z$.
Hence by Lemma \ref{lem:20211030}, we have $q_1^*(Z)\subset q_2^*(Z)$.
We have an isomorphism $\tau:B\times A\times \Sigma\to B\times A\times \Sigma$ defined by $\tau(b,a,s)=(b,a-b,s)$.
By $q_1=q_2\circ \tau$, we have $\tau^*q_2^*(Z)=B\times Z$.
We have $q_1\circ \tau(b,a,s)=(a-b,s)$.
Hence for $b\in B$, we have $\tau^*q_1^*(Z)\cap \pi^{-1}(b)=b+Z=Z$.
Hence by Lemma \ref{lem:20211030}, we have $\tau^*q_2^*(Z)\subset  \tau^*q_1^*(Z)$.
Hence $q_2^*(Z)\subset q_1^*(Z)$.
Thus we get $q_1^*(Z)=q_2^*(Z)$.
\hspace{\fill} $\square$

\begin{lem}\label{lem:quotient}
Let $A$ be a semi-abelian variety and let $B\subset A$ be a semi-abelian subvariety.
Let $\Sigma$ be a variety and let $Z\subset A\times \Sigma$ be a closed subscheme which is $B$-invariant in the sense that $b+Z=Z$ for all $b\in B$.
Then there exists a closed subscheme $Z_0\subset (A/B)\times \Sigma$ such that $Z=p^*(Z_0)$, where $p:A\times \Sigma\to (A/B)\times \Sigma$ is the projection.
Moreover, the natural map $Z\to Z_0$ is a categorical quotient, i.e., all $B$-invariant morphism $Z\to Y$ factor uniquely through $Z\to Z_0$.
\end{lem}

{\it Proof.}\
To show the existence of $Z_0$, we apply the theory of faithfully flat descent.
By Lemma \ref{lem:202110301}, we have $q_1^*(Z)=q_2^*(Z)$ as closed subschemes of $B\times A\times \Sigma$.
The cocycle condition on $B\times B\times A\times \Sigma$ (cf. \cite[(14.21.1)]{Gortz}) reduces to the obvious identity $\mathrm{id}_{B\times B\times Z}\circ \mathrm{id}_{B\times B\times Z}=\mathrm{id}_{B\times B\times Z}$.
Since the closed immersion $Z\hookrightarrow A\times \Sigma$ is affine, this descent data yields a closed subscheme $Z_0\subset (A/B)\times \Sigma$ such that $Z=p^*(Z_0)$ (cf. \cite[Cor 14.86]{Gortz}).
See also \cite[Remark 2.6.2. (iv)]{Brion}.
The restriction $p|_Z:Z\to Z_0$ is the base change of $p$, hence a $B$-torsor.
Hence $Z\to Z_0$ is a categorical quotient (cf. \cite[Proposition 2.6.4]{Brion}).
\hspace{\fill} $\square$

\section{Flattening via blowing-ups}\label{sec:b}

\begin{lem}\label{lem:20210417}
Let $S$ be a scheme of finite type over $\mathbb C$.
Let $X\subset \mathbb P^N\times S$ be a closed subscheme such that $X\to S$ is flat.
Assume that $S$ is reduced.
Let $U\subset S$ be a dense Zariski open set.
Then the scheme theoretic closure of $X|_U$ is $X$.
In particular $\supp X|_U\subset \supp X$ is dense.
\end{lem}

{\it Proof.}\
Let $Y\subset \mathbb P^N\times S$ be the scheme theoretic closure of $X|_U$.
Then $Y\subset X$ and $Y|_U=X|_U$.
Let $P$ be the Hilbert polynomial for each $X_s$ where $s\in S$.
Take $s\in S-U$.
Let $P'$ be the Hilbert polynomial of $Y_s$.
Then by $Y_s\subset X_s$, we have $P'\leq P$.
On the other hand, we have $P\leq P'$ by upper semi-continuity of the Hilbert polynomial.
Hence $P'=P$.
Since $S$ is reduced, we observe that $Y\to S$ is flat.
Moreover we have $Y_s=X_s$ for all $s\in S$.
Hence $Y=X$.
See \cite[Prop. 14.26]{Gortz}.
\hspace{\fill} $\square$

\begin{lem}\label{lem:202110273}
Let $X\subset \mathbb P^N\times S$ be a closed subscheme where $S$ is a variety.
Let $Z\subset S$ be an irreducible Zariski closed set.
Then there exist a projective birational modification $\varphi:S'\to S$ and a non-empty Zariski open set $U\subset S$ with the following properties.
\begin{enumerate}
\item 
$U\cap Z\not=\emptyset$.
\item
$\varphi^{-1}(U)\to U$ is an isomorphism.
\item
Let $p:X\to S$ be the projection and let $X'\subset \mathbb P^N\times S'$ be the scheme theoretic closure of $p^{-1}(U)\subset X\times_SS'$.
Then $X'|_{Z'}\to Z'$ is flat, where $Z'\subset S'$ is the Zariski closure of $\varphi^{-1}(Z\cap U)\subset S'$ and $X'|_{Z'}=X'\times_{S'}Z'$.
\end{enumerate}
\end{lem}

{\it Proof.}\
This follows from \cite[Lemma 3.1]{Y} as follows.
For each $s\in S$, let $P_{X_s}$ be the Hilbert polynomial of $X_s\subset \mathbb P^N$.
Let $\mathcal{P}$ be the set of all Hilbert polynomials of the closed subschemes of $\mathbb P^N$.
Then $\mathcal{P}$ is a totally ordered set by $P_1\leq P_2$ iff $P_1(m)\leq P_2(m)$ for all large $m$.
Moreover $\mathcal{P}$ is well-ordered (cf. \cite[Lemma 8.2]{Y}). 
Note that $\{ P_{X_s}\}_{s\in S}$ is a finite set (cf. \cite[p. 201, Step 2]{Ser}).
Let $P_{\max}$ be the maximal element in $\{ P_{X_s}\}_{s\in S}$.
We set
\begin{equation*}
T=
\begin{cases}
\{ s\in S;\ P_{X_s}=P_{\max}\} \quad & \text{if $X\not= \emptyset$,} \\
\emptyset  \quad & \text{if $X=\emptyset$.}
\end{cases}
\end{equation*}
Then $T\subset S$ is a Zariski closed subset (cf. \cite[Lemma 3.1]{Y}).

The proof of our lemma is by the transfinite induction on $P_{\max}\in \mathcal{P}$.
So we assume that our lemma is true when $P_{\max}<P$ and consider the case $P_{\max}=P$.
If $X=\emptyset$, then our lemma is obvious.
Hence in the following, we assume $X\not=\emptyset$.
Suppose $Z\subset T$.
Then for all $s\in Z$, we have $P_{X_s}=P$.
Since $Z$ is integral, $X|_Z\to Z$ is flat (cf. \cite[III, Thm 9.9]{H}).
Hence our lemma is true for $S'=S$ and $U=S$. 
So we consider the case $Z\not\subset T$.
Then $T\not= S$.
By \cite[Lemma 3.1]{Y}, there exists a closed subscheme $\mathcal{T}$ with $\supp \mathcal{T}=T$ such that if $\hat{P}=\max\{ P_{\hat{X}_s}\}_{s\in \hat{S}}$, then $\hat{P}<P$, where $\hat{S}=\mathrm{Bl}_{\mathcal{T}}S$ and $\hat{X}\subset \mathbb P^N\times \hat{S}$ is the strict transform.
This strict transform is defined under the isomorphism $\mathbb P^N\times \hat{S}=\mathrm{Bl}_{\mathbb P^N\times \mathcal{T}}(\mathbb P^N\times S)$.
Let $\hat{Z}\subset \hat{S}$ be the strict transform.

Now by the induction hypothesis, there exist a projective birational modification $\psi:S'\to \hat{S}$ and a non-empty Zariski open set $U_0\subset \hat{S}$ with the following properties.
\begin{itemize}
\item 
$U_0\cap \hat{Z}\not=\emptyset$.
\item
$\psi^{-1}(U_0)\to U_0$ is an isomorphism.
\item
Let $\hat{p}:\hat{X}\to \hat{S}$ be the projection and let $\tilde{X}\subset \mathbb P^N\times S'$ be the scheme theoretic closure of $\hat{p}^{-1}(U_0)\subset \hat{X}\times_{\hat{S}}S'$.
Then $\tilde{X}|_{Z'}\to Z'$ is flat, where $Z'\subset S'$ is the Zariski closure of $\psi^{-1}(\hat{Z}\cap U_0)\subset S'$.
\end{itemize}
We denote by $\varphi :S'\to S$ the composite of $S'\to\hat{S}\to S$.
Note that $\hat{S}\to S$ is an isomorphism over $S-T$.
Hence we may consider $S-T\subset S$ as an open subset of $\hat{S}$.
Set $U=(S-T)\cap U_0$.
Then we may consider $U$ as an open subset of $S$.
Then we have $U\cap Z\not=\emptyset$ and $\varphi^{-1}(U)\to U$ is an isomorphism.
Since $Z$ is irreducible, the Zariski closure of $\varphi^{-1}(U\cap Z)$ in $S'$ is equal to $Z'$.
Let $X'\subset \mathbb P^N\times S'$ be the scheme theoretic closure of $p^{-1}(U)\subset X\times_SS'$.
Then we have $X'\subset \tilde{X}$, hence $X'|_{Z'}\subset \tilde{X}|_{Z'}$.
On the other hand, we have $X'|_{\varphi^{-1}(U)}=\tilde{X}|_{\varphi^{-1}(U)}$.
Hence by Lemma \ref{lem:20210417}, the scheme theoretic closure of $X'|_{\varphi^{-1}(U)\cap Z'}\subset \tilde{X}|_{Z'}$ is $\tilde{X}|_{Z'}$.
Hence $\tilde{X}|_{Z'}\subset X'|_{Z'}$, hence $\tilde{X}|_{Z'}= X'|_{Z'}$.
Hence $X'|_{Z'}\to Z'$ is flat.
Hence by the transfinite induction, the proof is completed.
\hspace{\fill} $\square$

\begin{lem}
In Lemma \ref{lem:202110273}, if $S$ is smooth, then $S'$ is taken as smooth.
\end{lem}

{\it Proof.}\
We apply Lemma \ref{lem:202110273} to get $S'$ and $U$.
Suppose $S'$ is not smooth.
Then we take a smooth modification $S''\to S'$.
Since $S$ is smooth, we may assume that $S''\to S'$ is an isomorphism over $\varphi^{-1}(U)$.
Let $\psi:S''\to S$ be the induced map.
Then $\psi:S''\to S$ and $U\subset S$ satisfy the properties (1) and (2) in the statement of Lemma \ref{lem:202110273}.
We prove (3).
We consider $U\subset S''$.
Let $X''\subset \mathbb P^N\times S''$ be the scheme theoretic closure of $p^{-1}(U)\subset X\times_SS''$.
Then we have $X''\subset X'\times_{S'}S''\subset X\times_SS''$.
This implies the closed immersion $X''|_{Z''}\subset (X'\times_{S'}S'')|_{Z''}$, where $Z''\subset S''$ is the Zariski closure of $\psi^{-1}(Z\cap U)\subset S''$.
Since $Z''\to S'$ factors $Z'\to S'$, we have $(X'\times_{S'}S'')|_{Z''}=(X'|_{Z'})\times Z''$.
Hence we get $X''|_{Z''}\subset (X'|_{Z'})\times Z''$.
Note that $X''|_{Z''\cap U}=(X'|_{Z'})\times (Z''\cap U)=p^{-1}(U)|_{Z''\cap U}$.
Since $(X'|_{Z'})\times Z''\to Z''$ is flat, Lemma \ref{lem:20210417} yields $X''|_{Z''}= (X'|_{Z'})\times Z''$.
Hence $X''|_{Z''}\to Z''$ is flat.
We replace $S'$ by $S''$ to completes the proof.
\hspace{\fill} $\square$

\end{document}